\documentclass[final,3p]{elsarticle}
\journal{Computers and Mathematics with Applications (CAMWA)}

\usepackage{tikz}
\usetikzlibrary{math}
\usetikzlibrary{spy}
\usetikzlibrary{calc}
\usetikzlibrary{patterns}
\usetikzlibrary{external}
\usepackage{graphicx}
\graphicspath{{Images/}{../Images/}}
\usepackage{subcaption}
\usepackage{amssymb}
\usepackage{mathrsfs}
\usepackage{mathtools}
\newcommand{\Hilbert}{H} 
\newcommand{\level}{{\varrho}} 
\newcommand{\maxlevel}{{\varrho_{\rm max}}} 
\newcommand{\minlevel}{{\varrho_{\rm min}}} 
\newcommand{\sv}{\varsigma} 
\newcommand{\partition}{\mathcal{P}} 

\usepackage{amsthm}
\newtheorem*{example*}{Example}
\newtheorem{remark*}{Remark}
\newtheorem{remark}{Remark}
\usepackage[lined,boxed,commentsnumbered,linesnumbered]{algorithm2e}

\SetCommentSty{mycommfont}
\usepackage{color}
\usepackage{subfiles}
\usepackage[hidelinks]{hyperref}
\usepackage{booktabs}
\usepackage{multirow}
\biboptions{sort&compress}
\begin{document}
	\begin{frontmatter}
		\title{Error-estimate-based Adaptive Integration For Immersed Isogeometric Analysis}
		%
		\author[eindhoven,pavia]{Sai~C~Divi \corref{mycorrespondingauthor}}
		\cortext[mycorrespondingauthor]{Corresponding author}
		\ead{s.c.divi@tue.nl}
		\author[eindhoven]{Clemens~V~Verhoosel}
		\ead{C.V.Verhoosel@tue.nl}
		\author[pavia]{Ferdinando~Auricchio}
		\ead{auricchio@unipv.it}
		\author[pavia]{Alessandro~Reali}
		\ead{alessandro.reali@unipv.it}
		\author[eindhoven]{E~Harald~van~Brummelen}
		\ead{E.H.v.Brummelen@tue.nl}
		\address[eindhoven]{Department of Mechanical Engineering, Eindhoven University of Technology, 5600MB Eindhoven, The Netherlands}
		\address[pavia]{Department of Civil Engineering and Architecture, University of Pavia, 27100 Pavia, Italy}

\begin{abstract}
	The Finite Cell Method (FCM) together with Isogeometric analysis (IGA) has been applied successfully in various problems in solid mechanics, in image-based analysis, fluid-structure interaction and in many other applications. A challenging aspect of the isogeometric finite cell method is the integration of cut cells. In particular in three-dimensional simulations the computational effort associated with integration can be the critical component of a simulation. A myriad of integration strategies has been proposed over the past years to ameliorate the difficulties associated with integration, but a general optimal integration framework that suits a broad class of engineering problems is not yet available. In this contribution we provide a thorough investigation of the accuracy and computational effort of the octree integration scheme. We quantify the contribution of the integration error using the theoretical basis provided by Strang's first lemma. Based on this study we propose an error-estimate-based adaptive integration procedure for immersed isogeometric analysis. Additionally, we present a detailed numerical investigation of the proposed optimal integration algorithm and its application to immersed isogeometric analysis using two- and three-dimensional linear elasticity problems.
\end{abstract}

	\end{frontmatter}
	%
	\newpage
	\tableofcontents
	%
	\newpage

\section{Introduction}
Immersed finite element methods -- such as, \emph{e.g.}, the finite cell method (FCM) \cite{rank2008}, CutFEM \cite{hansbo2002} and immersogeometric analysis \cite{schillinger2011, rank2012, schillinger2012} -- have been demonstrated to be suitable for computational problems for which the performance of mesh-fitting finite element methods is impeded by complications in the meshing procedure. In recent years, immersed finite element methods have been successfully combined with isogeometric analysis (IGA) \cite{hughes2005}, a spline-based higher-order finite element framework targeting the integration of finite element analysis (FEA) and computer aided design (CAD). On the one hand, immersed methods provide the possibility to conduct IGA on volumetric domains based on a CAD boundary surface representations \cite{schillinger2012} or voxelized geometry data \cite{verhoosel2015}. On the other hand, immersed methods provide a natural framework for the consideration of trimming curves in IGA \cite{rank2012,schmidt2012,ruess2013,ruess2014,guo2015,guo2017,marussig2018,bauer2017,beer2015}. Immersed IGA has been applied successfully to various problems in, amongst others, structural and solid mechanics \cite{schillinger2015, duster2017}, fluid-structure interactions \cite{kamensky2015, hsu2015, xu2018} and scan-based analysis \cite{verhoosel2015, ruess2012}.

In comparison to mesh-fitting finite element methods, immersed methods require special treatment of various computational aspects. A prominent computational challenge that is inherent to immersed finite element methods is the integration over cut-cells, a problem that closely relates to the special treatment of discontinuous integrands in enriched finite element methods such as XFEM and GFEM. Since the geometry of the computational problem is captured by the integration procedure rather than by the ambient mesh in which the domain is immersed, cut-cell integration techniques must be capable of adequately capturing a wide range of cut-cell configurations. A myriad of dedicated integration procedures with this capability has been developed over the years in the context of immersed FEM (see \cite{schillinger2015} for a review) and enriched FEM (see \cite{belytschko2009}), which can be categorized as:
\begin{itemize}
	\item \textit{Octree subdivision:} The general idea of octree (or quadtree in 2D) integration is to capture the geometry of a cut-cell by recursively bisecting sub-cells that intersect with the boundary of the domain. At every level of this recursion, sub-cells that are completely inside the domain are preserved, while sub-cells that are completely outside of the domain are discarded. This cut-cell subdivision strategy was initially proposed in the context of the finite cell method in Ref.~\cite{duster2008} and is generally appraised for its simplicity and robustness with respect to cut-cell configurations. Octree integration has been widely adopted in immersed FEM, see, \emph{e.g.}, Refs.~\cite{schillinger2015,duster2017,kamensky2015,verhoosel2015}. Various generalizations and improvements to the original octree procedure have been proposed, of which the consideration of tetrahedral cells  \cite{varduhn2016,stavrev2016}, the reconstruction of the immersed boundary by tessellation of the lowest level of bisectioning \cite{verhoosel2015}, and the consideration of variable integration schemes for the sub-cells \cite{alireza2013}, are particularly noteworthy. Despite the various improvements to the original octree strategy, a downside of the technique remains the number of integration sub-cells that results from the procedure, especially in three-dimensional cases.
	
	\item \textit{Cut-cell reparametrization:} Accurate cut-cell integration schemes can be obtained by modifying the geometry parametrization of cut-cells in such a way that the immersed boundary is fitted. This strategy was original developed in the context of XFEM by decomposing cut elements into various sub-cells containing cut-cells with only one curved side and then to alter the geometry mapping related to the curved sub-cell to obtain a higher-order accurate integration scheme \cite{fries2010}. This concept has also been considered in the context of implicitly defined geometries (level sets) \cite{fries2016}, the NURBS-enhanced finite element method (NEFEM) \cite{sevilla2008,sevilla2011} and the Cartesian grid finite element method (cgFEM) \cite{nadal2013}. In the context of the finite cell method the idea of cut element reparametrization has been adopted as part of the \emph{smart octree integration} strategy \cite{kudela2016,hubrich2017}, where a boundary fitting procedure is employed at the lowest level of octree bisectioning in order to obtain higher-order integration schemes for cut-cells with curved boundaries. Reparametrization procedures have the potential to yield accurate integration schemes at a significantly lower computational cost than octree procedures, but generally compromises in terms of robustness with respect to cut-cell configurations.

  \item \textit{Polynomial integration}: Provided that one can accurately evaluate integrals over cut-cells (for example using octree integration), it is possible to construct computationally efficient integration rules for specific classes of integrands. In the context of immersed finite element methods it is of particular interest to derive efficient cut-cell integration rules for polynomial functions. The two most prominent methods to integrate polynomial functions over cut-cells are \emph{moment fitting techniques} \cite{mousavi2011,joulaian2016,hubrich2017,hubrich2019}, in which integration point weights and (possibly) positions are determined in order to yield exact quadrature rules, and \emph{equivalent polynomial methods} \cite{ventura2006,alireza2019}, in which a non-polynomial (\emph{e.g.}, discontinuous) integrand is represented by an equivalent polynomial which can then be treated using standard integration procedures. Such methods have been demonstrated to yield efficient quadrature rules for a range of scenarios. A downside of such techniques is the need for the evaluation of the exact integrals (using an adequate cut-cell integration procedure) in order to determine the optimized integration rules. This can make the construction of such quadrature rules computationally expensive, which makes that they are of particular interest mainly in the context of time-dependent and non-linear problems, for which the construction of the integration rule is only considered as a pre-processing operation (for each cut-cell) and the optimized integration rule can then be used throughout the simulation.

  \item \textit{Boundary integral representation}: Depending on the problem under consideration, it can be possible to reformulate volumetric integrals over cut-cells by equivalent boundary integrals. This reformulation, which has been proposed in the context of XFEM in Ref.~\cite{ventura2009} and in the immersed FEM setting in Ref.~\cite{jonsson2017}, is advantageous from a computational effort point of view, as the equivalent boundary integrals are generally less costly to evaluate. Moreover, provided that an adequate description of the boundary surface is available, higher-order accurate integration evaluations are obtained. A downside of integration techniques of this kind is that they require reformulation of the volume integrals, which makes them less general than standard quadrature rules.
	
\end{itemize}
In the selection of an appropriate cut-cell integration scheme one balances integration between robustness, accuracy and expense with respect to cut-cell configurations. If one requires a method that automatically treats a wide range of cut-cell configurations and is willing to pay the price in terms of (higher-order) accuracy and computational effort, octree integration is the compelling option. If constraints are imposed from an accuracy and computational expense point of view and one has some control over the range of configurations to be considered, alternative techniques such as cut-cell reparametrization are attractive. The need to balance between robustness, accuracy and expense has driven the development of hybrid integration schemes, such as \emph{smart octree integration} \cite{kudela2016} and \emph{adaptive moment-fitting} \cite{hubrich2019}, which allow one to attain an integration procedure with the desired properties.

Balancing accuracy with computational effort does not necessarily require the consideration of hybrid integration procedures, but can also be achieved by controlling the parameters of the integration procedures (listed above). This is particularly the case for octree integration, the accuracy of which can be controlled by the bisectioning depth and the integration orders used on the sub-cells. Optimization of these parameters to reduce the computational expense of octree integration without compromising its robustness with respect to configurations was proposed in Ref.~\cite{alireza2013}, where an algorithm is proposed to select the integration order on the different levels of sub-cells. It is demonstrated that reducing the integration order with increasing bisection depth significantly reduces the number of integration points, without unacceptably compromising the accuracy.

The idea of Ref.~\cite{alireza2013} to reduce integration orders on certain levels of the octree subdivision is in agreement with the theory of finite elements. Strang's first lemma \cite{strang1973,strang1972,ern2013} provides a framework to incorporate integration effects in the error analysis of finite element methods. This lemma indicates that integration does not need to be exact in order to attain (optimal) convergence \cite{strang1973}\footnote{As formulated by Strang and Fix in 1973 \cite{strang1973}: \textit{``What degree of accuracy in the integration formula is required for convergence? It is not required that every polynomial which appears be integrated exactly."}}. This lemma has been considered for the analysis of CutFEM, see, \emph{e.g.,} Refs.~\cite{burman2016,burman2017}. In fact, the idea of reduced integration in finite element methods has been studied extensively over the last decades, with applications in the analysis of plates and shells \cite{zienkiewicz1971} and mixed finite element methods \cite{hughes1978reduced,hughes1978mixed} being prominent examples.

In this manuscript we propose an error-estimation-based adaptive algorithm to obtain optimal octree integration rules for cut-cells. A sub-cell evaluation of the integration error in accordance with Strang's first lemma \cite{strang1973,strang1972,ern2013} is combined with a computational-cost estimate based on the number of integration cells to obtain a refinement indicator that optimally balances computational effort and accuracy. The proposed algorithm deviates from the one proposed in Ref.~\cite{alireza2013} in that it is directly based on integration error evaluations. The effectiveness of the adaptive integration procedure is demonstrated in the context of the isogeometric finite-cell framework. Besides the ability to obtain optimal integration rules for cut-cells, the developed adaptive integration procedure provides a tool to postulate rules of thumb for the selection of integration orders at different levels of bisectioning, and to assess the quality of integration rules derived using alternative algorithms such as that proposed in Ref.~\cite{alireza2013, taghipour2018}.

The paper outline is as follows. In Section~\ref{sec:fcm} we introduce the finite cell method, with a focus on the application of Strang's lemma to estimate integration errors. In Section~\ref{sec:integration_error} we study the influence of the octree integration scheme on the computational cost, and we define and evaluate the integration error of a sub-cell in a cut element. Based on the evaluated integration error, we present an optimization algorithm. The optimal integration procedure is investigated for two- and three-dimensional numerical test cases with various cut element configurations in Section~\ref{sec:integrationerroranalysis}. The developed algorithm is then applied to a two- and three-dimensional immersed isogeometric analysis problem in Section \ref{sec:numerical_examples}. Finally, concluding remarks are presented in Section \ref{sec:conclusion}.

\section{The finite cell method} \label{sec:fcm}
Because the integration procedure proposed in this work applies to a wide range of finite-cell simulations, in Section~\ref{sec:formulation} we first introduce the finite cell method in abstract form. Based on this abstract problem setting, Section~\ref{sec:inconsistencyerror} presents an error analysis that incorporates integration errors. In Section~\ref{sec:poissonproblem} the Poisson problem is considered to exemplify the abstract derivations.

\subsection{Formulation}\label{sec:formulation}
\begin{figure}
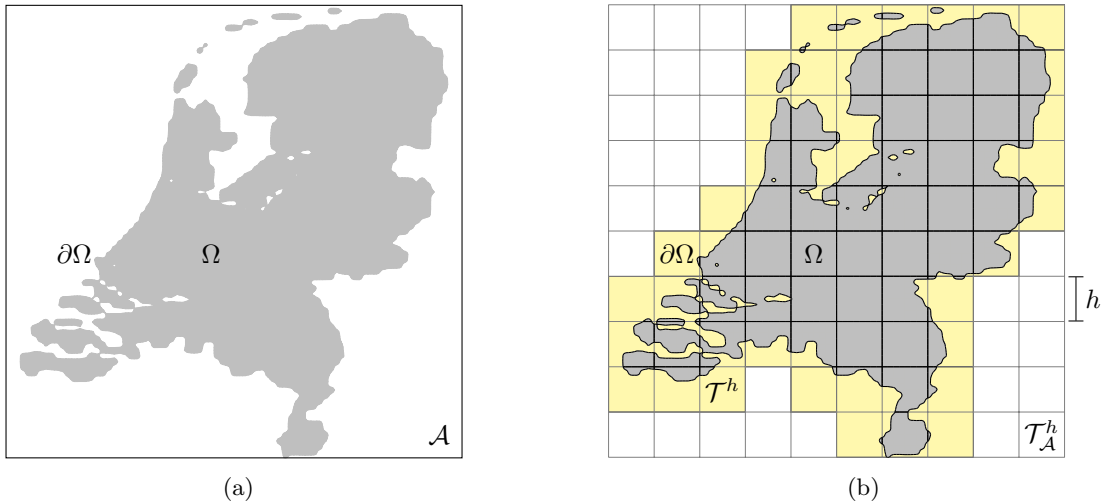

	\centering
	\begin{subfigure}[b]{0.5\textwidth}
		\centering
		\begin{tikzpicture}
  \tikzmath{\wdomain=2;
            \wambient=(3/2)*\wdomain;
            \h=0.3*\wdomain;
            }

    \begin{scope}[scale=3*\wdomain,line join=bevel,ultra thin,fill=black!25,draw=black!25]
      \input{tikz/nl_integration_subcells.dat}
    \end{scope}
    \draw[-](0,0) rectangle (2*\wambient,2*\wambient);
	\node at (4.5*\h,4.5*\h) {$\Omega$};
	\node at (1.5*\h,4.5*\h) {$\partial \Omega$};
	\node at (9.5*\h,0.5*\h) {$\mathcal{A}$};

\end{tikzpicture}
		\caption{}
	\end{subfigure}%
	\begin{subfigure}[b]{0.5\textwidth}
		\centering
		\begin{tikzpicture}
  \tikzmath{\wdomain=2;
            \wambient=(3/2)*\wdomain;
            \h=0.3*\wdomain;
            }

    \draw[ultra thin,fill=yellow!50!white,opacity=0.8]  (0*\h,1*\h) rectangle (3*\h,4*\h);
	\draw[ultra thin,fill=yellow!50!white,opacity=0.8]  (1*\h,4*\h) rectangle (3*\h,5*\h);
	\draw[ultra thin,fill=yellow!50!white,opacity=0.8]  (2*\h,5*\h) rectangle (3*\h,6*\h);
	\draw[ultra thin,fill=yellow!50!white,opacity=0.8]  (3*\h,5*\h) rectangle (7*\h,7*\h);
	\draw[ultra thin,fill=yellow!50!white,opacity=0.8]  (3*\h,7*\h) rectangle (6*\h,9*\h);
	\draw[ultra thin,fill=yellow!50!white,opacity=0.8]  (4*\h,9*\h) rectangle (10*\h,10*\h);
	\draw[ultra thin,fill=yellow!50!white,opacity=0.8]  (10*\h,9*\h) rectangle (9*\h,7*\h);
	\draw[ultra thin,fill=yellow!50!white,opacity=0.8]  (10*\h,7*\h) rectangle (8*\h,6*\h);
	\draw[ultra thin,fill=yellow!50!white,opacity=0.8]  (10*\h,6*\h) rectangle (8*\h,5*\h);
	\draw[ultra thin,fill=yellow!50!white,opacity=0.8]  (9*\h,5*\h) rectangle (7*\h,4*\h);
	\draw[ultra thin,fill=yellow!50!white,opacity=0.8]  (8*\h,4*\h) rectangle (6*\h,3*\h);
	\draw[ultra thin,fill=yellow!50!white,opacity=0.8]  (8*\h,3*\h) rectangle (7*\h,2*\h);
	\draw[ultra thin,fill=yellow!50!white,opacity=0.8]  (8*\h,2*\h) rectangle (5*\h,0*\h);
	\draw[ultra thin,fill=yellow!50!white,opacity=0.8]  (4*\h,1*\h) rectangle (5*\h,2*\h);
	\draw[ultra thin,fill=yellow!50!white,opacity=0.8]  (3*\h,2*\h) rectangle (5*\h,3*\h);
	\draw[ultra thin,fill=yellow!50!white,opacity=0.8]  (3*\h,3*\h) rectangle (4*\h,4*\h);
    \begin{scope}
      \draw[step=\h,gray,ultra thin] (0,0) grid (2*\wambient,2*\wambient);
      \draw[|-|] (2*\wambient+0.25*\h,3*\h) -- node[right]{$h$} ++(0,\h);
    \end{scope}


    \begin{scope}[scale=3*\wdomain,line join=bevel,ultra thin,fill=black!25,draw=black!25]
      \input{tikz/nl_integration_subcells.dat}
    \end{scope}
	\node at (4.5*\h,4.5*\h) {$\Omega$};
	\node at (1.5*\h,4.5*\h) {$\partial \Omega$};
    \node at (2.5*\h,1.5*\h) {$\mathcal{T}^h$};
    \node at (9.5*\h,0.5*\h) {$\mathcal{T}^h_\mathcal{A}$};
    
    \begin{scope}[scale=3*\wdomain,line join=bevel]
      \input{tikz/nl_integration_elements.dat}
    \end{scope}

\end{tikzpicture}
		\caption{}
	\end{subfigure}
	\caption{Schematic representation of (a) the physical domain $\Omega$ (gray) with boundary $\partial \Omega$ which is embedded in the ambient domain $\mathcal{A}$, and (b) the ambient domain mesh $\mathcal{T}^h_\mathcal{A}$ and background mesh $\mathcal{T}^h$ (yellow), with mesh size parameter $h$.}
	\label{fig:generaldomain}
\end{figure}

We consider a domain $\Omega \subset \mathbb{R}^d$ ($d \in \{ 2,3\}$) with boundary $\partial \Omega$ as shown in Figure~\ref{fig:generaldomain}. The boundary is split into a part on which Dirichlet conditions are imposed, denoted by $\partial \Omega_D$, and a complementary part to which Neumann conditions apply, denoted by $\partial \Omega_N$. 
We consider a problem described by a field variable $u$.
Let $W\ni{}u$ denote a suitable ambient space for the solution, equipped with a Hilbert structure corresponding to the inner product $(\cdot,\cdot)_W$ and the association norm
$\|\cdot\|_W$. Similarly, $V$ is a Hilbert space with inner product $(\cdot,\cdot)_V$ and norm $\|\cdot\|_V$, which encompasses the test space for the weak formulation of the problem under consideration. To accommodate the Dirichlet boundary conditions, let the spaces $W_0\subset{}W$ and $V_0\subset{}V$ be composed of functions that vanish at the Dirichlet boundary in a suitable manner. Moreover, let $\ell_g\in{}W$ denote a lift of the Dirichlet data. We consider a weak formulation of the generic form:
\begin{equation}
\left\{ \begin{aligned} 
&\text{Find} \: u \in \ell_g+W_0 \: \text{such that:} \\
&a(u,v) = b(v) \quad \forall v \in V_0,
\end{aligned}\right. 
\label{eq:weakform}
\end{equation}
We assume that the bilinear form $a:W \times V \to \mathbb{R}$ is continuous 
on $W\times{}V$ and weakly coercive on $W_0\times{}V_0$, and the linear form $b:V \to \mathbb{R}$ is continuous. The weak formulation~\eqref{eq:weakform} is then well-posed; see, e.g., \cite{ern2013}.

The finite cell method provides a general framework for constructing finite dimensional subspaces $W^h \subset W$ and $V^h \subset V$, where the superscript $h$ refers 
to a mesh parameter. The subspaces $W^h,V^h$ are subordinate to a regular 
mesh $\mathcal{T}^h_{\mathcal{A}}$, covering an ambient (embedding) domain $\mathcal{A} \supset \Omega$. The collection of elements $K\in\mathcal{T}^h_{\mathcal{A}}$ that intersect with the physical domain comprise the background mesh:
\begin{equation}
\mathcal{T}^h := \{ K \in \mathcal{T}^h_\mathcal{A}: K \cap \Omega \neq \emptyset \}
\label{eq:backgroundmesh}
\end{equation}
By trimming the elements in the background mesh, a mesh for the interior of the 
domain $\Omega$ is obtained: 
\begin{equation}
\mathcal{T}^h_{\Omega} := \{ K \cap \Omega :  K \in \mathcal{T}^h \}
\label{eq:interiormesh}
\end{equation}
Similarly, meshes for the Dirichlet and Neumann boundaries are defined as:
\begin{align}
\mathcal{T}^h_{\partial \Omega_D} := \{ K \cap \partial \Omega_D :  K \in \mathcal{T}^h \} \quad \mbox{and} \quad \mathcal{T}^h_{\partial \Omega_N} := \{ K \cap \partial \Omega_N :  K \in \mathcal{T}^h \}
\label{eq:boundarymeshes}
\end{align}
Note that elements in~$\mathcal{T}^h_{\partial \Omega_D}$ and~$\mathcal{T}^h_{\partial \Omega_N}$ are manifolds of co-dimension 1.

We consider a B-spline basis of degree $k$ and regularity $\alpha$ constructed over the ambient domain using the Cox-De Boor recursion formula \cite{piegl2012}. The span of this B-spline basis is denoted as 
\begin{equation}
\mathcal{S}^{k}_{\alpha}(\mathcal{A}) 
= 
\{ N\in{}C^{\alpha}(\mathcal{A}):N|_{K}\in{}P^k(K),\,\forall{}K\in\mathcal{T}^h_{\mathcal{A}}\}
\end{equation}
with $P^k(K)$ the collection of $d$-variate polynomials on~$K$. The approximation spaces in the finite cell method are obtained by restricting the B-splines 
in~$\mathcal{S}^{k}_{\alpha}(\mathcal{A})$ to the domain~$\Omega$:
\begin{equation}
\label{eq:WhVh}
W^h=V^h=\{N|_{\Omega}:N\in{}\mathcal{S}^{k}_{\alpha}(\mathcal{A}) \}
\end{equation}
A basis for~$W^h,V^h$ follows immediately from the restriction of the B-spline basis 
for~$\mathcal{S}^{k}_{\alpha}(\mathcal{A})$.

It is generally infeasible to impose restrictions on $W^h,V^h$ according to~\eqref{eq:WhVh} to form subspaces of~$W_0,V_0$ that retain suitable approximation properties. Hence, the finite cell method generally relies on weak imposition of the
Dirichlet boundary conditions, typically by means of Nitsche's method~\cite{nitsche1971}.
Accordingly, the bilinear and linear forms for the approximation problem are adapted to 
weakly incorporate the Dirichlet boundary conditions, giving rise to a Galerkin formulation of the form:
\begin{equation}
\left\{ \begin{aligned} 
&\text{Find} \: u^h \in W^h \: \text{such that:} \\
&a^h(u^h, v^h) = b^h(v^h) \quad \forall v^h \in V^h, 
\label{eq:galerkin}
\end{aligned}\right.
\end{equation}
where the superscript $h$ on the bilinear form $a^h : W^h \times V^h \to \mathbb{R}$ and linear form~$b^h : V^h \to \mathbb{R}$ indicate an explicit dependence on the mesh parameter. In general, the bilinear and linear forms in the finite cell method comprise contributions from both the interior and the boundaries. 
Corresponding to this typical setting, in the remainder of this work we assume the operators to be of the form 
\begin{subequations}
	\begin{align}
	a^h(u^h,v^h) &= \int_{\Omega} {A}^h_\Omega (u^h,v^h)(\boldsymbol{x})\,{\rm d}V + \int_{\partial \Omega_D} {A}^h_{\partial \Omega_D}( u^h,v^h )(\boldsymbol{x})\,{\rm d}S, \\  
	b^h(v^h) &= \int_{\Omega} {B}^h_\Omega (v^h )(\boldsymbol{x})\,{\rm d}V  
	+ \int_{\partial \Omega_D} {B}^h_{\partial \Omega_D}( v^h )(\boldsymbol{x})\,{\rm d}S
	+ \int_{\partial \Omega_N} {B}^h_{\partial \Omega_N}( v^h )(\boldsymbol{x})\,{\rm d}S,
	\end{align}
	\label{eq:discreteintegrals}%
\end{subequations}
where the integrands ${A}^h_{\Omega}$ and ${A}^h_{\partial \Omega_D}$ map $W^h\times{}V^h$ into integrable functions on the domain $\Omega$ and the Dirichlet boundary section, $\partial\Omega_D$, respectively. Similarly, ${B}^h_{\Omega}$, ${B}^h_{\partial \Omega_D}$ and ${B}^h_{\partial \Omega_N}$ map  $V^h$ into integrable functions on the domain $\Omega$ and its boundary sections, $\partial \Omega_D$ and~$\partial \Omega_N$. In the setting of~\eqref{eq:galerkin}, the spaces $W^h,V^h$ are equipped with (generally mesh-dependent) norms $\|\cdot\|_{W^h},\|\cdot\|_{V^h}$. Assuming that the bilinear form $a^h$ is continuous and weakly coercive on $W^h\times{}V^h$ and the linear form $b^h$ is continuous on~$V^h$, the approximation problem~\eqref{eq:galerkin} is well posed \cite[ch.~2]{ern2013}.

In practice, the integrals in equation \eqref{eq:discreteintegrals} are approximated by means of quadrature rules, leading to the definition of the quadrature-dependent bilinear and linear forms
\begin{subequations}
	\begin{align}
	{a}^h_\mathcal{Q}(u^h,v^h) &= 
	\sum \limits_{K \in \mathcal{T}^h_\Omega} 
	\sum \limits_{l=1}^{l_K} \omega_K^l {A}^h_\Omega (u^h,v^h)(\boldsymbol{x}_K^l)
	+ 
	\sum \limits_{K \in \mathcal{T}^h_{\partial \Omega_D}} 
	\sum \limits_{l=1}^{l_K} \omega_K^l {A}^h_{\partial \Omega_D} 
	(u^h,v^h)(\boldsymbol{x}_K^l),\\  
	{b}^h_{\mathcal{Q}}(v^h) &= 
	\sum \limits_{K \in \mathcal{T}^h_\Omega} 
	\sum \limits_{l=1}^{l_K} \omega_K^l {B}^h_\Omega (v^h)(\boldsymbol{x}_K^l)
    + 
    \sum \limits_{K \in \mathcal{T}^h_{\partial \Omega_D}} \sum \limits_{l=1}^{l_{K}} 	
    \omega_K^l {B}^h_{\partial \Omega_D} (v^h)(\boldsymbol{x}_K^l)
    \nonumber\\ 
    &\phantom{=\quad}
    + \sum \limits_{K \in \mathcal{T}^h_{\partial \Omega_N}} \sum \limits_{l=1}^{l_{K}} 	\omega_K^l {B}^h_{\partial \Omega_N} (v^h)(\boldsymbol{x}_K^l),
	\end{align} \label{eq:quadratureforms}%
\end{subequations}
where for each (interior or boundary) element $K$, the set 
$\{(\boldsymbol{x}_K^l,\omega_K^l)\}_{l=1}^{l_K}$ represents a quadrature rule, \emph{i.e.}, a suitable collection of pairs of integration points in~$K$ and corresponding weights. We denote by $\mathcal{Q}$ the complete integration scheme, \emph{i.e.}, the collection of all the (interior and boundary) element-wise quadrature rules.
The Galerkin problem corresponding to the approximate (bi-)linear forms $a^h_\mathcal{Q}$ and $b^h_{\mathcal{Q}}$ then writes:
\begin{equation}
\left\{ \begin{aligned} 
&\text{Find} \: {u}^{h}_{\mathcal{Q}} \in W^h \: \text{such that} \\
&a^{h}_{\mathcal{Q}}({u}^h_{\mathcal{Q}},v^h) = b^{h}_{\mathcal{Q}}(v^h) \quad 
\forall v^h \in V^h.
\end{aligned}\right. \label{eq:approxproblem}
\end{equation}
We assume that the integration scheme $\mathcal{Q}$ transfers the coercivity and boundedness properties of $a^h,b^h$ to $a^h_{\mathcal{Q}},b^h_{\mathcal{Q}}$,
so that~\eqref{eq:approxproblem} is well posed. We note that 
the solution ${u}^h_{\mathcal{Q}}$ generally deviates from that 
of~\eqref{eq:galerkin} on account of the inexactness of the integration rules.
%

\subsection{Finite cell error analysis} \label{sec:inconsistencyerror}
The error of the approximated finite cell solution, ${u}^h_{\mathcal{Q}}$, computed using the Galerkin problem \eqref{eq:approxproblem} is composed of two parts, \emph{viz.}: \emph{i)} the discretization error, defined as the difference between the exact solution to~\eqref{eq:weakform}, $u$, and the approximate solution to~\eqref{eq:galerkin} in the absence of integration errors, $u^h$; and \emph{ii)} the inconsistency error related to the integration procedure, which is defined as the difference between the approximate solution in the absence of integration errors, $u^h$, and the approximate solution to \eqref{eq:approxproblem} with integration errors, ${u}^{h}_{\mathcal{Q}}$. In this section we present the error analysis in an abstract setting. A concrete example is provided in Section~\ref{sec:poissonproblem}.
To provide a setting for the error analysis, we denote by $W(h)=\operatorname{span}{}\{u\}\oplus{}W^h$ the
linear space containing the actual solution to~\eqref{eq:weakform} and the approximation
space $W^h$. We equip $W(h)$ with a norm $\|\cdot\|_{W(h)}$ such that $\|w^h\|_{W(h)}=\|w^h\|_{W^h}$ for all $w^h\in{}W^h$ and $\|u\|_{W(h)}\leq{}c\,\|u\|_W$, \emph{i.e.}, the solution of~\eqref{eq:weakform} is continuously embedded in~$W(h)$; cf.~\cite[Ch.~2]{ern2013}. We assume that the solution~$u$ of~\eqref{eq:weakform} is suitably regular,
so that the bilinear form $a^h$ admits a continuous extension to $W(h)$, \emph{i.e.}, there exists a constant $\mathscr{C}>0$ such that:
\begin{equation}
\|a^h\|_{W(h),V^h}:=\sup_{(w,v^h)\in{}W(h)\times{}V^h}\frac{|a^h(w,v^h)|}{\|w\|_{W(h)}\|v\|_{V^h}}\leq{}\mathscr{C}.
\label{eq:ahnorm}
\end{equation}
In this setting, an upper bound for the error $u-u^h_{\mathcal{Q}}$ is provided by Strang's first lemma~\cite{strang1972}:
\begin{equation}
\begin{aligned}
 \left\| u - {u}^h_{\mathcal{Q}} \right\|_{W(h)} 
 \leq & 
 \left( 1 + \frac{\left\| a^h \right\|_{W(h),V^h}}{\alpha^h}  \right) \left\| u - \mathcal{I}^h u \right\|_{W(h)} +  & & \\ 
 &  \frac{1}{\alpha^h} \left(\sup_{v^h \in {V}^h}{\frac{\left|  b^h (v^h) - b^{h}_{\mathcal{Q}}(v^h) \right| }{ \left\| v^h \right\|_{V^h}}} 
 + \sup_{v^h \in {V}^h}{\frac{\left|  a^h(\mathcal{I}^h u,v^h) - a^{h}_{\mathcal{Q}}(\mathcal{I}^h u,v^h) \right| }{ \left\| v^h \right\|_{V^h}}} \right),
\end{aligned}
\label{eq:strang1}
\end{equation}
where $\mathcal{I}^h u$ represents the best approximation of $u$ in $W^h$, \emph{i.e.},
\begin{equation}
\mathcal{I}^h u = \underset{w^h \in{} W^h}{\rm arg\,min} \| u - w^h \|
\end{equation}
and $\alpha^h$ denotes the inf-sup constant of the bilinear form $a^h:W^h\times{}V^h\to\mathbb{R}$:
\begin{equation}
\alpha^h = \inf \limits_{w^h \in W^h\rule{0pt}{8pt}} \sup \limits_{v^h \in V^h}  
\frac{ a^h(w^h,v^h) }{\| w^h \|_{W^h}  \| v^h \|_{V^h} }.
\label{eq:alphah}
\end{equation}
Proof of~\eqref{eq:strang1} is standard; see, for instance, \cite[Lemma 2.27]{ern2013}. The first term in~\eqref{eq:strang1} represents an upper bound to the discretization error. The second term in~\eqref{eq:strang1} bounds the consistency error, \emph{i.e.}, the  integration error.

Equation~\eqref{eq:discreteintegrals} conveys that, in principle, integration errors emerge for both the volume integrals and the surface integrals. Since the focus of this work is on the optimization of the integration orders over octree sub-cells, in the remainder of this work we will restrict our considerations to the integration errors associated with the volume integrals, thereby implicitly assuming that the boundary integrals are evaluated exactly. When neglecting the integration errors associated with the boundary integrals, the error term in \eqref{eq:strang1} associated with the inexact integration of the linear form is bounded by
\begin{align}
\sup_{v^h \in {V}^h}{\frac{\left|  b^h(v^h) - b^{h}_{\mathcal{Q}}(v^h) \right| }{ \left\| v^h \right\|_{V^h}}} &= 
\sup_{v^h \in {V}^h}{\frac{\left| \sum \limits_{K \in \mathcal{T}^h_\Omega} \left( \int_{K} B^h_\Omega( v^h)(\boldsymbol{x})\,{\rm d}V 
-  
\sum \limits_{l=1}^{l_K} \omega^l_K B^{h}_\Omega(v^h)( \boldsymbol{x}_K^l)  \right)\right| }{ \left\| v^h \right\|_{V^h}}} \nonumber \\ 
& \leq {\sum \limits_{K \in \mathcal{T}^h_\Omega} 
\sup_{v^h \in {V}^h}
\frac{  \left|\int_{K} B^h_\Omega(v^h)(\boldsymbol{x})\,{\rm d}V -  \sum \limits_{l=1}^{l_K} \omega^l_K B^{h}_\Omega (v^h)( \boldsymbol{x}_K^l)\right| }{ \left\| v^h \right\|_{V^h}}}.
 \label{eq:upperboundb1}
\end{align}
Note that the supremum is always considered over a function space with the zero function excluded, \emph{e.g.}, $v^h \in V^h \setminus \{0\}$. For notational brevity we omit the zero exclusion in the remainder of this manuscript. Under the (non-restrictive) assumption that $B^h_{\Omega}$ is a local operator in the sense that
\begin{equation}
\big(B^h_{\Omega}(v^h|_K)\big)\Big|_K\text{ is well defined and } 
\big(B^h_{\Omega}(v^h)\big)\Big|_K=\big(B^h_{\Omega}(v^h|_K)\big)\Big|_K,
\label{eq:locality}
\end{equation}
\emph{e.g.}, if $B^h_{\Omega}$ corresponds to a differential operator, the summands in the ultimate expression in~\eqref{eq:upperboundb1} are bounded by the 
element-integration-error indicators
\begin{equation}
e_K^b = \sup_{v^h_K \in {V}^h_K}{\frac{  \left|\int_{K} B^h_\Omega( v^h_K)(\boldsymbol{x})\,{\rm d}V -  \sum \limits_{l=1}^{l_K} \omega^l_K B^{h}_\Omega(v^h_K)( \boldsymbol{x}_K^l) \right| }{ \left\| v^h_K \right\|_{V^h_K}}}, 
\label{eq:elem_interrorb}
\end{equation}
where $V^h_K=\{v^h|_K,\,v^h\in{}V^h\}$ is the restriction of $V^h$ to the element~$K$, equipped with a norm $\| \cdot \|_{V^h_K}$ such that $\|v^h|_K\|_{V^h_K}\leq\|v^h\|_{V^h}$ for all $v^h\in{}V^h$. The element-integration-error indicators $e_K^b$ can be computed element-wise and provide a bound for the integration error in the linear form:
\begin{equation}
\sup_{v^h \in {V}^h}{\frac{\left|  b^h(v^h) - b^{h}_{\mathcal{Q}}(v^h) \right| }{ \left\| v^h \right\|_{V^h}}}
\leq
\sum \limits_{K \in \mathcal{T}^h_\Omega{}}e^b_K
\label{eq:upperboundb}
\end{equation}
Similarly, the integration error for the bilinear form in~\eqref{eq:strang1} is bounded as
\begin{align}
\sup_{v^h \in {V}^h}{\frac{\left|  a^h(\mathcal{I}^h u,v^h) - a^{h}_{\mathcal{Q}}(\mathcal{I}^h u,v^h) \right| }{ \left\| v^h \right\|_{V^h}}}  
&\leq  
\sum \limits_{K \in \mathcal{T}^h_\Omega} e_K^{a}  \label{eq:upperbounda}
\end{align}
with the element-integration-error indicators defined as
\begin{align}
e_K^a = 
\sup_{v^h_K \in {V}^h_K}{\frac{ \left|\int_{K} A^h_\Omega(\mathcal{I}^h u,v^h)(\boldsymbol{x}_K)\,{\rm d}V -  \sum \limits_{l=1}^{l_K} \omega^l_K A^{h}_\Omega(\mathcal{I}^h u,v^h )( \boldsymbol{x}_K^l)\right| }{ \left\| v^h_K \right\|_{V^h_K}}}. \label{eq:elem_interrora}
\end{align}
Substitution of the bounds \eqref{eq:upperboundb} and  \eqref{eq:upperbounda} into the error estimate \eqref{eq:strang1} then yields
\begin{equation}
\begin{aligned}
 \left\| u - {u}^h_{\mathcal{Q}} \right\|_{W(h)} 
 \leq  
 \left( 1 + \frac{\left\| a^h \right\|_{W(h),V^h}}{\alpha^h}  \right) \left\| u - \mathcal{I}^h u \right\|_{W(h)} +  
   \frac{1}{\alpha^h} 
\sum \limits_{K \in \mathcal{T}^h_\Omega} \left(  e_K^b +  e_K^a \right) ,
\end{aligned}
\label{eq:strang1estimatefinal}
\end{equation}
which conveys that the element-integration-error indicators \eqref{eq:elem_interrorb} and \eqref{eq:elem_interrora} yield control over the inconsistency error. This motivates the development of an algorithm to minimize the element-integration-error indicators; see Section~\ref{sec:integrationerror}. 

\subsection{Application to the Poisson problem}\label{sec:poissonproblem}
We apply the general theory presented in Sections~\ref{sec:formulation}--\ref{sec:inconsistencyerror} to the particular
case of the Dirichlet--Poisson problem:
\begin{equation}
\left\{ \begin{aligned} 
-\Delta u &= f \quad \mbox{in} ~ \Omega, \\
u &= g  \quad \mbox{on} ~ \partial \Omega,
\end{aligned}\right. 
\label{eq:stronglaplace}
\end{equation}
with $f:\Omega \to \mathbb{R}$ and $g: \partial \Omega \to \mathbb{R}$ 
exogenous data. A suitable ambient space for the weak formulation of~\eqref{eq:stronglaplace} is provided by $W=H^1(\Omega)$, equipped with
the usual $\Hilbert^1$-norm and inner product.
 Denoting by $\ell_g\in{}\Hilbert^1(\Omega)$ a lift of the Dirichlet data such that $\ell_g|_{\partial\Omega}=g$ and by $W_0=\Hilbert^1_0(\Omega)$ the subspace of functions in~$\Hilbert^1(\Omega)$ that vanish at the boundary in the trace sense, the weak formulation of~\eqref{eq:stronglaplace} assumes the form~\eqref{eq:weakform} with
\begin{equation}
a(u,v)=\int_{\Omega}\nabla{}u\cdot\nabla{}v\,{\rm d}V,\qquad
b(v)=\int_{\Omega}f\,v\,{\rm d}V.\qquad
\end{equation}
We denote by $V^h$ a finite-dimensional subspace of $\Hilbert^1(\Omega)$ according to the construction in Section~\ref{sec:formulation}. We equip $V^h$ with the mesh-dependent norm
\begin{equation}
\|v^h\|_{V^h}^2=\int_{\Omega}|\nabla{}v^h|^2\,{\rm d}V+\frac{c_0}{h}\int_{\partial\Omega}(v^h)^2\,{\rm d}S.
\label{eq:Vhnorm}
\end{equation}
The symmetric Galerkin formulation in~$V^h$ with weak enforcement of the Dirichlet conditions by means of Nitsche's method conforms to~\eqref{eq:galerkin} with 
$a^h:V^h\times{}V^h\to\mathbb{R}$ and $b^h:V^h\to\mathbb{R}$:
\begin{subequations}
\begin{align}
a^h(u^h,v^h)&=\int_{\Omega}\nabla{}u^h\cdot\nabla{}v^h\,{\rm d}V
-\int_{\partial\Omega}v^h\partial_nu^h\,{\rm d}S
-\int_{\partial\Omega}u^h\partial_nv^h\,{\rm d}S
+\frac{c_1}{h}\int_{\partial\Omega}u^hv^h\,{\rm d}S,
\\
b^h(v^h)&=\int_{\Omega}f{}\,v^h\,{\rm d}V
-\int_{\partial\Omega}g\,\partial_nv^h\,{\rm d}S
+\frac{c_1}{h}\int_{\partial\Omega}g\,v^h\,{\rm d}S.
\end{align}
\end{subequations}
For suitably chosen constants $c_0$ and $c_1$, the bilinear form is bounded and 
coercive on $V^h\times{}V^h$ and the linear form $b^h$ is bounded on~$V^h$.
The bilinear form $a^h$ does not admit a continuous extension to all $\Hilbert^1(\Omega)$, owing to the fact that the normal derivatives $\partial_n(\cdot)$ that appear in $a^h$ are not properly defined for all functions in~$\Hilbert^1(\Omega)$. We therefore assume that the solution of the Dirichlet--Poisson problem~\eqref{eq:stronglaplace} is suitably regular, \emph{e.g.}, such that $\partial_nu\in{}L^2(\partial\Omega)$. We introduce the composite space $W(h)=\operatorname{span}\{u\}\oplus{}V^h$ and equip $W(h)$ with the norm $\|\cdot\|_{W(h)}$, corresponding to the extension of the norm in~\eqref{eq:Vhnorm} to~$W(h)$. Trace theory conveys that $\|u\|_{W(h)}\leq{}c_h\|u\|_{\Hilbert^1(\Omega)}$ for some ($h$\nobreakdash-dependent) constant $c_h>0$ and, hence, the embedding of $\operatorname{span}\{u\}$ into~$W(h)$ is continuous. Under the aforementioned regularity assumption on~$u$, the bilinear form $a^h$ admits a continuous extension to $W(h)\times{}V^h$, and \eqref{eq:ahnorm}\nobreakdash--\eqref{eq:alphah} apply.
To elucidate the notation for the quadrature-dependent bilinear and linear forms for the Dirichlet--Poisson problem, we note that in this case:
\begin{equation}
A^h_{\Omega}(u^h,v^h)(\boldsymbol{x})=\nabla{}u^h(\boldsymbol{x})\cdot\nabla{}
v^h(\boldsymbol{x}),
\qquad
B^h_{\Omega}(v^h)(\boldsymbol{x})=f(\boldsymbol{x})\,v^h(\boldsymbol{x}).
\end{equation}
One easily verifies that both $A^h_{\Omega}$ and $B^h_{\Omega}$ satisfy the locality condition~\eqref{eq:locality}.

	\section{Optimized octree cut-cell integration} \label{sec:integration_error}
	
	Based on the error analysis discussed in Section~\ref{sec:inconsistencyerror}, in this section we develop an algorithm to optimize the distribution of integration points over octree-subdivided cut-cells. In Section~\ref{sec:octree_partitioning} we first introduce the considered octree procedure and study its complexity with respect to the most relevant parameters. In Section~\ref{sec:integrationerror} we then express the integration-error estimate in an operator-independent form, making it applicable to a class of operators. The per-cut-element evaluation of the operator-independent integration error is also discussed in this section. Finally, the developed optimization procedure in the form of Algorithm~\ref{alg:adaptive_integration} is presented in Section~\ref{sec:thealgorithm}.
	
	\begin{figure}
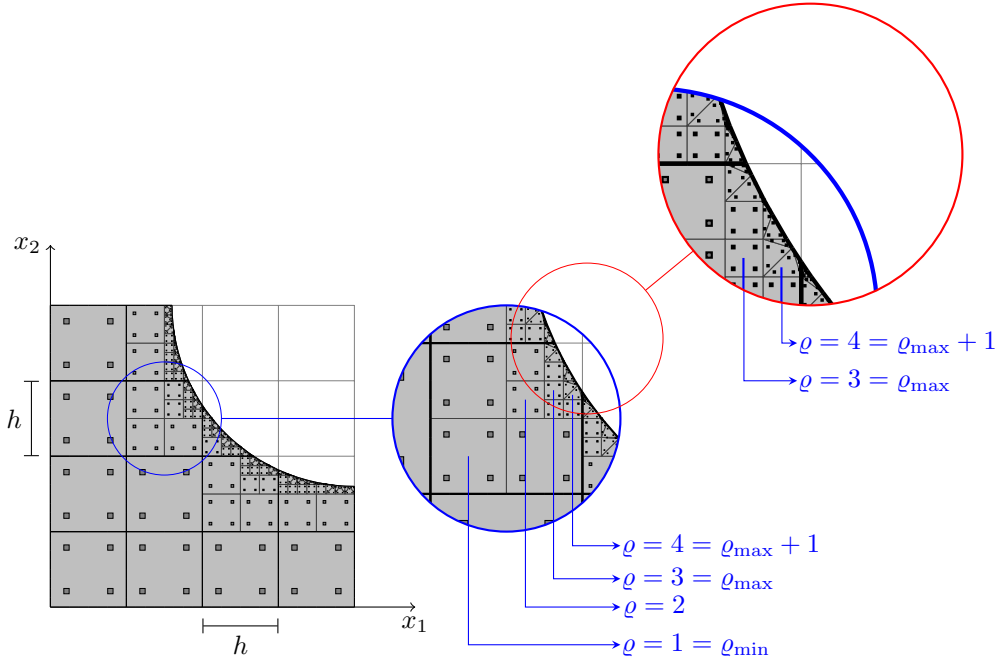

		\centering
		\begin{tikzpicture}
[spy using outlines={black, circle, magnification=2, size=4cm, connect spies,red}]
  \tikzmath{\wdomain=2;
            \wambient=2*\wdomain;
            \h=0.5*\wdomain;
            \Rc=0.4*\wdomain;
            \Rsw=1.273/2*\wdomain;
            \Rse=0.5*\wdomain;
            \Rne=1.333/2*\wdomain;
            \Rnw=0.4*\wdomain;
            }

  \begin{scope}[spy using outlines={black, circle, magnification=2, size=3*\h cm, connect spies,blue}]

    \begin{scope}
      \draw[step=\h,gray,ultra thin] (-\wambient,-\wambient) grid (0,0);
      \draw[|-|] (-\wambient-0.25*\h,-2*\h) -- node[left]{$h$} ++(0,\h);
      \draw[|-|] (-2*\h,-\wambient-0.25*\h) -- node[below]{$h$} ++(\h,0);
    \end{scope}

    \draw[->] (-\wambient,-\wambient) -- (0.2*\wambient,-\wambient) node[below] {$x_1$} coordinate(x axis);
    \draw[->] (-\wambient,-\wambient) -- (-\wambient,0.2*\wambient) node[left] {$x_2$} coordinate(y axis);

    \begin{scope}[scale=2*\wdomain,line join=bevel,ultra thin,fill=black!25,draw=black!75]
      \input{tikz/octree_integration_subcells.dat}
    \end{scope}

    \begin{scope}[scale=2*\wdomain,line join=bevel]
      \input{tikz/octree_integration_elements.dat}
    \end{scope}

    \begin{scope}[scale=2*\wdomain,line width=0.2pt,inner sep=0pt,minimum size=1pt,rectangle,draw=black!100,fill=black!50]
      \input{tikz/octree_integration_points.dat}
    \end{scope}

    \spy on (-2.5*\h, -1.5*\h) in node(zoom1) at (2*\h,-1.5*\h);

  \end{scope}

  \spy on (zoom1.north east) in node(zoom2) at (6, 2);

  \draw[->,>=stealth] (1.5*\h,-2*\h) [blue,circle,draw,inner sep=0pt,minimum size=2pt]{} -- ++(-90:2.5*\h) -- ++(0:2*\h) node[right]{$\varrho = 1 = \varrho_{\rm min}$};
  \draw[->,>=stealth] (2.25*\h,-1.25*\h) [blue,circle,draw,inner sep=0pt,minimum size=2pt]{} -- ++(-90:2.75*\h) -- ++(0:1.25*\h) node[right]{$\varrho = 2$};
  \draw[->,>=stealth] (2.62*\h,-1.125*\h) [blue,circle,draw,inner sep=0pt,minimum size=2pt]{} -- ++(-90:2.5*\h) -- ++(0:0.88*\h) node[right]{$\varrho = 3 = \varrho_{\rm max}$};
  \draw[->,>=stealth] (2.87*\h,-1.187*\h) [blue,circle,draw,inner sep=0pt,minimum size=2pt]{} -- ++(-90:2*\h) -- ++(0:0.63*\h) node[right]{$\varrho = 4 = \varrho_{\rm max} + 1$};

  \draw[->,>=stealth] (5.125*\h,\h) [blue,circle,draw,inner sep=0pt,minimum size=2pt]{} -- ++(-90:2*\h) -- ++(0:0.7*\h) node[right]{$\varrho = 3 = \varrho_{\rm max}$};
  \draw[->,>=stealth] (5.625*\h,\h) [blue,circle,draw,inner sep=0pt,minimum size=2pt]{} -- ++(-90:1.5*\h) -- ++(0:0.2*\h) node[right]{$\varrho = 4 = \varrho_{\rm max} + 1$};

\end{tikzpicture}
		\caption{Schematic representation of the octree integration procedure with tessellation at the lowest level of bisectioning (see \ref{sec:tessellation}).}
		\label{fig:octreeintegration}
	\end{figure}
	
	\subsection{Octree partitioning}\label{sec:octree_partitioning}
	To construct an explicit parametrization of the geometry (including its boundary) we herein consider the octree procedure proposed in Ref.~\cite{verhoosel2015}, which is illustrated in Figure~\ref{fig:octreeintegration} for a two-dimensional cut-cell. In this procedure an element in the background mesh that intersects the boundary of the domain is bisected into $2^d$ sub-cells. If a sub-cell is completely inside the domain, it is preserved in the partitioning of the cut-cell, whereas a sub-cell is discarded if it is completely outside of the domain. This bisectioning procedure is recursively applied to the sub-cells that intersect the boundary, until $\maxlevel$-times bisected sub-cells are obtained. At the lowest level of bisectioning, \emph{i.e.}, for the $\maxlevel$-times bisected sub-cells, a tessellation procedure is applied to construct a partitioning of the sub-cells that intersect with the domain boundary, resulting in an additional level, \emph{i.e.}, $\maxlevel + 1$. On the one hand, this tessellation procedure provides an $\mathcal{O}( h^2 / 2^{2 \maxlevel} )$ accurate parametrization of the interior volume \cite{verhoosel2015}, while, on the other hand, it provides a parametrization for the trimmed surface. See \ref{sec:tessellation} for details regarding the employed tessellation procedure.
	
	Since the octree procedure provides a parametrization of the cut-cell and integration rules are generally available for all sub-cells \cite{alireza2013,taghipour2018}, cut-cell integration can be performed by agglomeration of all sub-cell quadrature points. The accuracy of the cut-cell integration scheme can then be controlled through the selection of the quadrature rules on the sub-cells. In particular, polynomials can be integrated exactly over the cut-cell when Gauss quadrature of the appropriate order is used on all sub-cells. An example of a cut-cell integration scheme based on Gauss quadrature for third order polynomials is illustrated in Figure~\ref{fig:octreeintegration}, where the gray squares represent the Gauss points, and the relative size of the squares is representative for the integration weights.
	
	Selection of Gauss points of optimal order on all sub-cells is evidently very attractive, in the sense that an adequate cut-cell integration scheme can be constructed by specification of the bisectioning depth $\maxlevel$ and the order of polynomials to be integrated exactly. However, as observed from Figure~\ref{fig:octreeintegration}, the obtained cut-cell integration rule is generally not computationally efficient, in the sense that the majority of integration points is formed on small sub-cells, as a consequence of which their relative contribution to the overall integral (observed from the size of the squares) is generally limited. As studied in Ref.~\cite{alireza2013}, the computational efficiency of the cut-cell integration scheme can be improved by reducing the number of integration points on lower bisectioning levels. Although this reduction decreases the accuracy of the integration scheme, the obtained improvement in computational effort associated with the lower number of integration points outweighs this disadvantage. The algorithm proposed in Ref.~\cite{alireza2013} targets the optimization of this balance between accuracy and computational effort.
	
	\begin{figure}
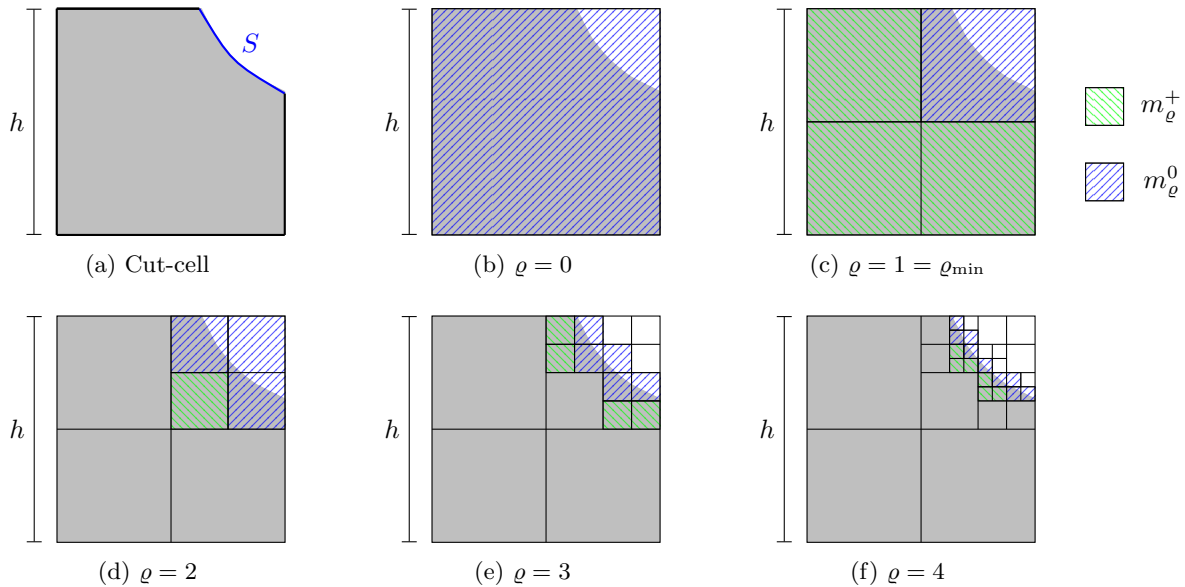

		\begin{subfigure}{0.3\textwidth}
			\centering
			\begin{tikzpicture}
  \tikzmath{\wdomain=1.5;
            \wambient=(3/2)*\wdomain;
            \h=1.0*\wdomain;
            }
        
  \begin{scope}[spy using outlines={black, circle, magnification=4, size=8*\h cm, connect spies,blue}]
    \coordinate (A) at (0.0*\wdomain,0.0*\wdomain);
    \coordinate (B) at (2.0*\wdomain,0.0*\wdomain);
    \coordinate (C) at (2.0*\wdomain,2.0*\wdomain);
    \coordinate (D) at (0.0*\wdomain,2.0*\wdomain);

    \begin{scope}[scale=2*\wdomain,line join=bevel,ultra thin,fill=black!25,draw=black!25]
      \input{tikz/octree_level3_subcells.dat}
    \end{scope}
    \draw[thick,color=blue] (1.25*\wdomain,2.0*\wdomain) .. controls (1.53*\wdomain,1.53*\wdomain) .. (2.0*\wdomain,1.25*\wdomain);
    \node at (1.7*\wdomain,1.7*\wdomain) {\color{blue} $S$};
    \draw[thick] (1.25*\wdomain,2.0*\wdomain) -- (D);
    \draw[thick] (A) -- (B);
    \draw[thick] (B) -- (2.0*\wdomain,1.25*\wdomain);
    \draw[thick] (1.25*\wdomain,2.0*\wdomain) -- (D);
    \draw[thick] (D) -- (A);
    \draw[|-|] (-0.2*\wdomain,0.0*\wdomain) -- (-0.2*\wdomain,2.0*\wdomain) node[midway, left] {$h$};

  \end{scope}
    
\end{tikzpicture}
			\caption{Cut-cell}
			\label{fig:subcells}
		\end{subfigure}%
		\begin{subfigure}{0.3\textwidth}
			\centering
			\begin{tikzpicture}
  \tikzmath{\wdomain=1.5;
            \wambient=(3/2)*\wdomain;
            \h=1.0*\wdomain;
            }
        
  \begin{scope}[spy using outlines={black, circle, magnification=4, size=8*\h cm, connect spies,blue}]
    \coordinate (A) at (0.0*\wdomain,0.0*\wdomain);
    \coordinate (B) at (2.0*\wdomain,0.0*\wdomain);
    \coordinate (C) at (2.0*\wdomain,2.0*\wdomain);
    \coordinate (D) at (0.0*\wdomain,2.0*\wdomain);

    \begin{scope}[scale=2*\wdomain,line join=bevel,ultra thin,fill=black!25,draw=black!25]
      \input{tikz/octree_level3_subcells.dat}
    \end{scope}
    
	\draw[ultra thin, pattern=north east lines,pattern color = blue,opacity=0.8]  (A)rectangle (C);

    \draw[-](0,0) rectangle (2*\wdomain,2*\wdomain);
    
    \draw[|-|] (-0.2*\wdomain,0.0*\wdomain) -- (-0.2*\wdomain,2.0*\wdomain) node[midway, left] {$h$};

  \end{scope}
    
\end{tikzpicture}
			\caption{$\level = 0$}
			\label{fig:subcells_level0}
		\end{subfigure}%
		\begin{subfigure}{0.3\textwidth}
			\centering
			\begin{tikzpicture}
  \tikzmath{\wdomain=1.5;
            \wambient=(3/2)*\wdomain;
            \h=1.0*\wdomain;
            }
        
  \begin{scope}[spy using outlines={black, circle, magnification=4, size=8*\h cm, connect spies,blue}]
    \coordinate (A) at (0.0*\wdomain,0.0*\wdomain);
    \coordinate (B) at (2.0*\wdomain,0.0*\wdomain);
    \coordinate (C) at (2.0*\wdomain,2.0*\wdomain);
    \coordinate (D) at (0.0*\wdomain,2.0*\wdomain);
    \coordinate (A0) at (1.0*\wdomain,0.0*\wdomain);
	\coordinate (B0) at (1.0*\wdomain,2.0*\wdomain);
	\coordinate (AB0) at (1.0*\wdomain,1.0*\wdomain);
	\coordinate (C0) at (0.0*\wdomain,1.0*\wdomain);
	\coordinate (D0) at (2.0*\wdomain,1.0*\wdomain);

    \begin{scope}[scale=2*\wdomain,line join=bevel,ultra thin,fill=black!25,draw=black!25]
      \input{tikz/octree_level3_subcells.dat}
    \end{scope}
    
	\draw[ultra thin, pattern=north west lines,pattern color = green,opacity=0.8]  (A)rectangle (D0);
	\draw[ultra thin, pattern=north west lines,pattern color = green,opacity=0.8]  (C0)rectangle (B0);
	\draw[ultra thin, pattern=north east lines,pattern color = blue,opacity=0.8]  (AB0)rectangle (C);

	\draw[-] (A0) --(B0);
	\draw[-] (C0) --(D0);
    \draw[-](0,0) rectangle (2*\wdomain,2*\wdomain);
    
    \draw[|-|] (-0.2*\wdomain,0.0*\wdomain) -- (-0.2*\wdomain,2.0*\wdomain) node[midway, left] {$h$};

  \end{scope}
    
\end{tikzpicture}
			\caption{$\level = 1 = \minlevel$}
			\label{fig:subcells_level1}
		\end{subfigure}%
		\begin{subfigure}{0.1\textwidth}
			\begin{tikzpicture}
			\draw[ultra thin, pattern=north west lines,pattern color = green,opacity=0.8]  (0,1) rectangle (0.5,1.5);
			\draw[ultra thin, pattern=north east lines,pattern color = blue,opacity=0.8]  (0,0) rectangle (0.5,0.5);
			\draw[-](0,0) rectangle (0.5,0.5);
			\draw[-](0,1) rectangle (0.5,1.5);
			\node at (1,0.25) {$m^0_{\level}$};
			\node at (1,1.25) {$m^{+}_{\level}$};
			\end{tikzpicture}
		\end{subfigure}\\[12pt]
		\begin{subfigure}{0.3\textwidth}
			\centering
			\begin{tikzpicture}
  \tikzmath{\wdomain=1.5;
            \wambient=(3/2)*\wdomain;
            \h=1.0*\wdomain;
            }
        
  \begin{scope}[spy using outlines={black, circle, magnification=4, size=8*\h cm, connect spies,blue}]
    \coordinate (A) at (0.0*\wdomain,0.0*\wdomain);
    \coordinate (B) at (2.0*\wdomain,0.0*\wdomain);
    \coordinate (C) at (2.0*\wdomain,2.0*\wdomain);
    \coordinate (D) at (0.0*\wdomain,2.0*\wdomain);
    \coordinate (A0) at (1.0*\wdomain,0.0*\wdomain);
	\coordinate (B0) at (1.0*\wdomain,2.0*\wdomain);
	\coordinate (AB0) at (1.0*\wdomain,1.0*\wdomain);
	\coordinate (C0) at (0.0*\wdomain,1.0*\wdomain);
	\coordinate (D0) at (2.0*\wdomain,1.0*\wdomain);
	\coordinate (A1) at (1.5*\wdomain,1.0*\wdomain);
	\coordinate (B1) at (1.5*\wdomain,2.0*\wdomain);
	\coordinate (AB1) at (1.5*\wdomain,1.5*\wdomain);
	\coordinate (AB2) at (1.5*\wdomain,1.75*\wdomain);
	\coordinate (C1) at (1.0*\wdomain,1.5*\wdomain);
	\coordinate (D1) at (2.0*\wdomain,1.5*\wdomain);
    \begin{scope}[scale=2*\wdomain,line join=bevel,ultra thin,fill=black!25,draw=black!25]
      \input{tikz/octree_level3_subcells.dat}
    \end{scope}
    
	\draw[ultra thin, pattern=north west lines,pattern color = green,opacity=0.8]  (AB0)rectangle (AB1);
	\draw[ultra thin, pattern=north east lines,pattern color = blue,opacity=0.8]  (A1)rectangle (C);
	\draw[ultra thin, pattern=north east lines,pattern color = blue,opacity=0.8]  (C1)rectangle (B1);
	
	\draw[-] (A1) --(B1);
	\draw[-] (C1) --(D1);
	\draw[-] (A0) --(B0);
	\draw[-] (C0) --(D0);
    \draw[-](0,0) rectangle (2*\wdomain,2*\wdomain);
	\draw[|-|] (-0.2*\wdomain,0.0*\wdomain) -- (-0.2*\wdomain,2.0*\wdomain) node[midway, left] {$h$};
  \end{scope}
    
\end{tikzpicture}
			\caption{$\level = 2$}
			\label{fig:subcells_level2}
		\end{subfigure}%
		\begin{subfigure}{0.3\textwidth}
			\centering
			\begin{tikzpicture}
  \tikzmath{\wdomain=1.5;
            \wambient=(3/2)*\wdomain;
            \h=1.0*\wdomain;
            }
        
  \begin{scope}[spy using outlines={black, circle, magnification=4, size=8*\h cm, connect spies,blue}]
    \coordinate (A) at (0.0*\wdomain,0.0*\wdomain);
    \coordinate (B) at (2.0*\wdomain,0.0*\wdomain);
    \coordinate (C) at (2.0*\wdomain,2.0*\wdomain);
    \coordinate (D) at (0.0*\wdomain,2.0*\wdomain);
    \coordinate (A0) at (1.0*\wdomain,0.0*\wdomain);
	\coordinate (B0) at (1.0*\wdomain,2.0*\wdomain);
	\coordinate (C0) at (0.0*\wdomain,1.0*\wdomain);
	\coordinate (D0) at (2.0*\wdomain,1.0*\wdomain);
	\coordinate (A1) at (1.5*\wdomain,1.0*\wdomain);
	\coordinate (B1) at (1.5*\wdomain,2.0*\wdomain);
	\coordinate (AB1) at (1.5*\wdomain,1.5*\wdomain);
	\coordinate (AB2) at (1.5*\wdomain,1.75*\wdomain);
	\coordinate (C1) at (1.0*\wdomain,1.5*\wdomain);
	\coordinate (D1) at (2.0*\wdomain,1.5*\wdomain);
	\coordinate (A2) at (1.75*\wdomain,1.0*\wdomain);
	\coordinate (B2) at (1.75*\wdomain,2.0*\wdomain);
	\coordinate (AB21) at (1.75*\wdomain,1.25*\wdomain);
	\coordinate (AB22) at (1.75*\wdomain,1.5*\wdomain);
	\coordinate (AB23) at (1.75*\wdomain,1.75*\wdomain);
	\coordinate (C2) at (1.25*\wdomain,1.5*\wdomain);
	\coordinate (CD1) at (1.25*\wdomain,1.75*\wdomain);
	\coordinate (D2) at (1.25*\wdomain,2.0*\wdomain);
	\coordinate (E2) at (1.5*\wdomain,1.25*\wdomain);
	\coordinate (F2) at (2.0*\wdomain,1.25*\wdomain);
	\coordinate (G2) at (1.0*\wdomain,1.75*\wdomain);
	\coordinate (H2) at (2.0*\wdomain,1.75*\wdomain);
    \begin{scope}[scale=2*\wdomain,line join=bevel,ultra thin,fill=black!25,draw=black!25]
      \input{tikz/octree_level3_subcells.dat}
    \end{scope}
	\draw[ultra thin, pattern=north west lines,pattern color = green,opacity=0.8]  (A2)rectangle (F2);
	\draw[ultra thin, pattern=north west lines,pattern color = green,opacity=0.8]  (A1)rectangle (AB21);
	\draw[ultra thin, pattern=north west lines,pattern color = green,opacity=0.8]  (C1)rectangle (CD1);
	\draw[ultra thin, pattern=north west lines,pattern color = green,opacity=0.8]  (G2)rectangle (D2);
	\draw[ultra thin, pattern=north east lines,pattern color = blue,opacity=0.8]  (E2)rectangle (D1);
	\draw[ultra thin, pattern=north east lines,pattern color = blue,opacity=0.8]  (C2)rectangle (AB23);
	\draw[ultra thin, pattern=north east lines,pattern color = blue,opacity=0.8]  (CD1)rectangle (B1);
	
	\draw[-] (A2) -- (B2);
	\draw[-] (C2) -- (D2);
	\draw[-] (E2) -- (F2);
	\draw[-] (G2) -- (H2);
	\draw[-] (A1) --(B1);
	\draw[-] (C1) --(D1);
	\draw[-] (A0) --(B0);
	\draw[-] (C0) --(D0);
    \draw[-](0,0) rectangle (2*\wdomain,2*\wdomain);
	\draw[|-|] (-0.2*\wdomain,0.0*\wdomain) -- (-0.2*\wdomain,2.0*\wdomain) node[midway, left] {$h$};
  \end{scope}
    
\end{tikzpicture}
			\caption{$\level = 3$}
			\label{fig:subcells_level3}
		\end{subfigure}%
		\begin{subfigure}{0.3\textwidth}
			\centering
			\begin{tikzpicture}
  \tikzmath{\wdomain=1.5;
            \wambient=(3/2)*\wdomain;
            \h=1.0*\wdomain;
            }
        
  \begin{scope}[spy using outlines={black, circle, magnification=4, size=8*\h cm, connect spies,blue}]
    \coordinate (A) at (0.0*\wdomain,0.0*\wdomain);
    \coordinate (B) at (2.0*\wdomain,0.0*\wdomain);
    \coordinate (C) at (2.0*\wdomain,2.0*\wdomain);
    \coordinate (D) at (0.0*\wdomain,2.0*\wdomain);
    \coordinate (A0) at (1.0*\wdomain,0.0*\wdomain);
	\coordinate (B0) at (1.0*\wdomain,2.0*\wdomain);
	\coordinate (C0) at (0.0*\wdomain,1.0*\wdomain);
	\coordinate (D0) at (2.0*\wdomain,1.0*\wdomain);
	\coordinate (A1) at (1.5*\wdomain,1.0*\wdomain);
	\coordinate (B1) at (1.5*\wdomain,2.0*\wdomain);
	\coordinate (AB1) at (1.5*\wdomain,1.5*\wdomain);
	\coordinate (AB2) at (1.5*\wdomain,1.75*\wdomain);
	\coordinate (C1) at (1.0*\wdomain,1.5*\wdomain);
	\coordinate (D1) at (2.0*\wdomain,1.5*\wdomain);
	\coordinate (A2) at (1.75*\wdomain,1.0*\wdomain);
	\coordinate (B2) at (1.75*\wdomain,2.0*\wdomain);
	\coordinate (AB21) at (1.75*\wdomain,1.25*\wdomain);
	\coordinate (AB22) at (1.75*\wdomain,1.5*\wdomain);
	\coordinate (AB23) at (1.75*\wdomain,1.75*\wdomain);
	\coordinate (C2) at (1.25*\wdomain,1.5*\wdomain);
	\coordinate (CD1) at (1.25*\wdomain,1.75*\wdomain);
	\coordinate (D2) at (1.25*\wdomain,2.0*\wdomain);
	\coordinate (E2) at (1.5*\wdomain,1.25*\wdomain);
	\coordinate (F2) at (2.0*\wdomain,1.25*\wdomain);
	\coordinate (G2) at (1.0*\wdomain,1.75*\wdomain);
	\coordinate (H2) at (2.0*\wdomain,1.75*\wdomain);
	\coordinate (A3) at (1.625*\wdomain,1.25*\wdomain);
	\coordinate (B3) at (1.625*\wdomain,1.75*\wdomain);
	\coordinate (AB31) at (1.625*\wdomain,1.375*\wdomain);
	\coordinate (AB32) at (1.625*\wdomain,1.5*\wdomain);
	\coordinate (AB33) at (1.625*\wdomain,1.625*\wdomain);
	\coordinate (AB24) at (1.75*\wdomain,1.375*\wdomain);
	\coordinate (AB3) at (1.5*\wdomain,1.625*\wdomain);
	\coordinate (C3) at (1.875*\wdomain,1.25*\wdomain);
	\coordinate (D3) at (1.875*\wdomain,1.5*\wdomain);
	\coordinate (E3) at (1.375*\wdomain,1.5*\wdomain);
	\coordinate (F3) at (1.375*\wdomain,2.0*\wdomain);
	\coordinate (EF31) at (1.375*\wdomain,1.625*\wdomain);
	\coordinate (EF32) at (1.375*\wdomain,1.75*\wdomain);
	\coordinate (G3) at (1.5*\wdomain,1.375*\wdomain);
	\coordinate (H3) at (2.0*\wdomain,1.375*\wdomain);
	\coordinate (I3) at (1.25*\wdomain,1.625*\wdomain);
	\coordinate (J3) at (1.75*\wdomain,1.625*\wdomain);
	\coordinate (K3) at (1.25*\wdomain,1.875*\wdomain);
	\coordinate (L3) at (1.5*\wdomain,1.875*\wdomain);
	%
    \begin{scope}[scale=2*\wdomain,line join=bevel,ultra thin,fill=black!25,draw=black!25]
      \input{tikz/octree_level3_subcells.dat}
    \end{scope}
    \draw[ultra thin, pattern=north west lines,pattern color = green,opacity=0.8]  (E2) rectangle (AB24);
    \draw[ultra thin, pattern=north west lines,pattern color = green,opacity=0.8]  (G3) rectangle (AB32);
    \draw[ultra thin, pattern=north west lines,pattern color = green,opacity=0.8]  (C2) rectangle (AB3);
    \draw[ultra thin, pattern=north west lines,pattern color = green,opacity=0.8]  (I3) rectangle (EF32);
    \draw[ultra thin, pattern=north east lines,pattern color = blue,opacity=0.8] (AB21) rectangle (H3);
    \draw[ultra thin, pattern=north east lines,pattern color = blue,opacity=0.8] (AB31) rectangle (D3);
    \draw[ultra thin, pattern=north east lines,pattern color = blue,opacity=0.8] (AB1) rectangle (AB33);
    \draw[ultra thin, pattern=north east lines,pattern color = blue,opacity=0.8] (EF31) rectangle (AB2);
    \draw[ultra thin, pattern=north east lines,pattern color = blue,opacity=0.8] (CD1) rectangle (L3);
    \draw[ultra thin, pattern=north east lines,pattern color = blue,opacity=0.8] (K3) rectangle (F3);
    
    \draw[-] (A3) -- (B3);
    \draw[-] (C3) -- (D3);
    \draw[-] (E3) -- (F3);
    \draw[-] (G3) -- (H3);
    \draw[-] (I3) -- (J3);
    \draw[-] (K3) -- (L3);
	\draw[-] (A2) -- (B2);
	\draw[-] (C2) -- (D2);
	\draw[-] (E2) -- (F2);
	\draw[-] (G2) -- (H2);
	\draw[-] (A1) --(B1);
	\draw[-] (C1) --(D1);
	\draw[-] (A0) --(B0);
	\draw[-] (C0) --(D0);
    \draw[-](0,0) rectangle (2*\wdomain,2*\wdomain);
	\draw[|-|] (-0.2*\wdomain,0.0*\wdomain) -- (-0.2*\wdomain,2.0*\wdomain) node[midway, left] {$h$};
  \end{scope}
    
\end{tikzpicture}
			\caption{$\level = 4$}
			\label{fig:subcells_level4}
		\end{subfigure}
	\caption{Schematic illustration of the sub-cell formation using the octree partitioning for a single cut-element. The green color represents the sub-cells at level $\level$ that are completely inside the domain, $m^{+}_{\level}$, and the blue cells represent the intersected sub-cells of the level $\level$, $m^{0}_{\level}$.}
	\label{fig:subcell_count}
	\end{figure}

	Since the balance between accuracy and computational effort also forms the basis of the error-estimation-based optimization procedure in this work it is important to understand the distribution of the number of octree sub-cells over the levels of bisectioning. To establish a scaling relation for the number of sub-cells we consider a random cut-cell of size $h$ with trimmed surface area $S$, as illustrated in Figure~\ref{fig:subcells}. The octree subdivision approximates the surface area at the maximum octree depth as
	\begin{align}
	S &\approx m^0_{\maxlevel} \bar{s} \left( \frac{h}{2^{\maxlevel}} \right)^{d-1}   &  & \maxlevel \geq {\minlevel} \geq 1
	\label{eq:surfacerelation}
	\end{align}
	where $m^0_{\maxlevel}$ denotes the average number of sub-cells at octree depth $\maxlevel$ that is intersected by the trimmed boundary, and where $\bar{s}$ is the average surface area of a randomly cut unit cube in $d$ dimensions. This surface approximation is valid under the condition that
	\begin{align}
	\maxlevel \geq {\minlevel}  = {\rm ceil}\left(  \frac{\log_2{\left( \sv \right)}}{1-d} \right) \geq 1,  
	\end{align}
	with surface fraction $\sv= \frac{S}{ \bar{s} h^{d-1} }$. For example, for the cut-cell illustrated in Figure~\ref{fig:subcell_count}, $\maxlevel \geq \minlevel = 1$.
	
	From \eqref{eq:surfacerelation} the number of intersected elements at level $\maxlevel$ follows to be:
	\begin{align}
	m^0_{\maxlevel} \approx 
	\begin{cases} 
	1 & 0 \leq \maxlevel <  {\minlevel} \\
	\sv  2^{\maxlevel(d-1)}  & \maxlevel \geq {\minlevel \geq 1}
	\end{cases} 
	\label{eq:mintersection}
	\end{align}
	Denoting the number of sub-cells at level $\level$ that are completely inside of the domain by $m_{\level}^{+}$ (see Figure~\ref{fig:subcell_count}), and noting that the average volume of a randomly cut unit cube is $\bar{v}=\frac{1}{2}$, the cut-cell volume can be approximated as
	\begin{align}
	V &\approx m_{\maxlevel}^0  \frac{1}{2} \left( \frac{h}{2^{\maxlevel}} \right)^d +   \sum_{\level = 1}^{\maxlevel} m^+_{\level} \left( \frac{h}{2^{\level}} \right)^d  & \maxlevel \geq 1, 
	\end{align}
	which is divided by the volume of the background cell to obtain the cut-cell volume fraction:
\begin{subequations}
	\begin{align}
	\eta  \approx \eta_{\maxlevel} &= m_{\maxlevel}^0  2^{-(\maxlevel d+1)} +   \sum_{\level = 1}^{\maxlevel} m^+_{\level} \left( \frac{1}{2^{\level}} \right)^d   &  &\maxlevel \geq 1,\\
	&=  \sv  2^{-(\maxlevel+1)}  +   \sum_{\level = 1}^{\maxlevel} m^+_{\level} \left( \frac{1}{2^{\level}} \right)^d  & &\maxlevel \geq \minlevel \geq 1.
	\end{align}
\end{subequations}
	Assuming that $\maxlevel$ adequately resolves the volume fraction of the cut-cell, it holds that $\eta_{\maxlevel} \approx \eta_{\maxlevel-1}$ with $\maxlevel \geq \minlevel+1$ and
	\begin{align}
	\eta_{\maxlevel} -  \eta_{\maxlevel-1}  &= - \sv  2^{-(\maxlevel+1)}  +   m^+_{\maxlevel} {2^{-\maxlevel d}} \approx 0.
	\end{align}
	For the number of sub-cells at bisectioning level $\maxlevel$, we then obtain
	\begin{align}
	m^+_{\maxlevel} &\approx \sv  2^{\maxlevel(d-1)-1}  &    &\maxlevel > \minlevel.
	\end{align}
	Since it follows from \eqref{eq:mintersection} that the interface localizes in a single sub-cell if $\level \leq \minlevel$ ($m_\level^0 = 1 \: \level < \minlevel, m_{\minlevel}^0 = \sv 2^{\minlevel (d-1)} \approx 1$), the number of preserved sub-cells across all levels is given by
	\begin{align}
	m_{\level}^+ \approx 
	\begin{cases}
	0 &  \level \leq \minlevel ~ \land ~ \eta \leq \frac{1}{2}\\
	2^{d} - 1   & \level \leq \minlevel ~ \land ~ \eta > \frac{1}{2}\\
	\sv  2^{\level(d-1)-1} & \level > \minlevel
	\end{cases}
	\label{eq:mpositive}
	\end{align}
	where the volume fraction $\eta$ dictates which side of the interface corresponds to the interior domain. For example, for the cut element displayed in Figure~\ref{fig:subcell_count} the volume fraction is larger than a half, which implies that all but one of the sub-cells at level $\level=1$ are preserved. If the complement of the element would be considered, none of the level $\level=1$ sub-cells would be included in the partitioning. 
	
	If we denote the average number of integration points per octree sub-cell by $\bar{q}_{\level}$, the total number of integration points per level follows as
	\begin{align}
	n_{\level} \approx 
	\begin{cases}
	0 &  \level \leq \minlevel ~ \land ~ \eta \leq \frac{1}{2}\\
	\bar{q}_{\level}( 2^{d} - 1 )   & \level \leq \minlevel ~ \land ~ \eta > \frac{1}{2}\\
	\bar{q}_{\level}  \sv  2^{\level(d-1)-1} &  \minlevel < \level \leq \maxlevel\\
	\bar{q}_{\level} \bar{t}_d \sv  2^{\maxlevel(d-1)}  & \level = \maxlevel+1
	\end{cases}
	\label{eq:pointsperlevel}
	\end{align}
	where $\bar{t}_{d}$ is the number of sub-cells in which the lowest bisectioning level $\level = \maxlevel + 1$ is tessellated (see \ref{sec:tessellation}). The total number of integration points for a level $\maxlevel$ octree partitioning then follows as
	\begin{align}
	N_{\maxlevel} &\approx 
	\begin{cases}
	\sv \left( \bar{q}_{\maxlevel+1} \bar{t}_d 2^{\maxlevel(d-1)} + \sum \limits_{\level=\minlevel+1}^{\maxlevel} \bar{q}_{\level} 2^{\level(d-1)-1}  \right) & \eta \leq \frac{1}{2},  \\
	\sv \left( \bar{q}_{\maxlevel+1} \bar{t}_d 2^{\maxlevel(d-1)} + \sum \limits_{\level=\minlevel+1}^{\maxlevel} \bar{q}_{\level} 2^{\level(d-1)-1}  \right) + \sum \limits_{\level = 1}^{\minlevel} \bar{q}_{\level} \left( 2^{d} - 1 \right) & \eta > \frac{1}{2}.
	\end{cases}
	\label{eq:numberofpoints}
	\end{align}
This approximation reflects the strong dependence of the number of integration points on the surface fraction $\sv$. When $\sv$ is very small, $\minlevel$ is very large, as a result of which a significant number of integration points for non-intersected sub-cells can be present (depending on the volume fraction). When $\sv$ is large, or, more generally when $\maxlevel \gg \minlevel$, the number of non-intersected sub-cells is negligible. For the case that the Gauss order is selected equal on all levels, \emph{i.e.}, $\bar{q}_\level = q_{\rm line}^d$ for $\level \leq \maxlevel$ (with $q_{\rm line}$ the number of integration points along a one-dimensional line segment) and $\bar{q}_{\maxlevel+1} = q_{{\rm tes},d}$, the number of integration points \eqref{eq:numberofpoints} reduces to:
	\begin{align}
	N_{\maxlevel} &\approx  
	\sv \left( {q}_{{\rm tes},d} \bar{t}_d  +  \frac{{q}_{\rm line}^d}{ 2-2^{2-d}} \right) 2^{\maxlevel(d-1)} & &\maxlevel \gg \minlevel. \label{eq:approx_numberofpoints}
	\end{align}
	This estimate conveys that the number of integration points scales linearly with the surface fraction, and exponentially with the octree depth. This exponential scaling depends on the number of spatial dimensions. In two dimensions the number of integration points doubles while increasing the octree depth by one, while in three dimensions it quadruples. Hence, in comparison to the two-dimensional setting, the number of integration points in three dimensions scales much more dramatically with the octree depth. A numerical study of the number of integration points in relation to the scaling relations presented in this section is provided in Section~\ref{sec:integrationerroranalysis}.

  \subsection{Cut-cell integration errors}\label{sec:integrationerror}
  In Section~\ref{sec:inconsistencyerror} the application of Strang's first lemma to the quantification of integration errors in the finite-cell method was discussed. Equations \eqref{eq:upperboundb} and \eqref{eq:upperbounda} convey that the integration error can be estimated by the summation of the operator-dependent element-integration-error indicators $e_K^{{b}}$ and $e_K^{{a}}$. This section discusses the evaluation of these errors and their localization to the integration sub-cells. In order to evaluate the cell-wise integration error, in Section~\ref{sec:integrationdefinition} the integration errors are expressed as the product of an operator-independent integration error and a scaling term depending on the operator. Subsequently, in Section~\ref{sec:integrationevaluation} the numerical evaluation of the operator-independent cut-cell integration error is discussed.

\subsubsection{Integration error definition}\label{sec:integrationdefinition}
The element-integration-error indicators $e_K^{b}$ and $e_K^{a}$, expressed by equations~\eqref{eq:elem_interrorb} and \eqref{eq:elem_interrora} respectively, depend on the operations $A^h_\Omega$ and $B^h_\Omega$. It is generally desirable to apply a single integration scheme for all terms and, hence, to have uniform control over $e_K^{a}$ and $e_K^{b}$. To enable a uniform treatment of both $e_K^{a}$ and $e_K^{b}$, we first note that the integrals in the numerators of \eqref{eq:elem_interrorb} and \eqref{eq:elem_interrora} constitute linear functionals on $V^h_{K}$. By the Riesz-representation theorem, there exist elements $T^{a}, T^{b} \in V^h_{K}$ such that
\begin{align}
\int \limits_{K} T^{a} v^h_K \,{\rm d}V &= \int \limits_{K} A^h_\Omega\left( \mathcal{I}^h u, v^h_K \right)(\boldsymbol{x}_K) \,{\rm d}V, \label{eq:riesz_representation_a} \\
\int \limits_{K} T^{b} v^h_K \,{\rm d}V &= \int \limits_{K} B^h_\Omega\left( v^h_K \right)(\boldsymbol{x}_K) \,{\rm d}V, \label{eq:riesz_representation_b}
\end{align}
for all $v^h_K \in V^h_K$. We now proceed under the assumption that the error introduced in the numerical integration is equivalent for the original functionals and their Riesz-representation, \emph{i.e.}, the difference in applying the integral rate to the left- and right-hand-side members of \eqref{eq:riesz_representation_a} and \eqref{eq:riesz_representation_b} is negligible. It then holds that
\begin{align}
e_K^{a} & \leq \sup_{v^h_K \in V^h_K}{ \frac{\left|\int_{K} T^{a} v^h_K \,{\rm d}V -  \sum \limits_{l=1}^{l_K} \omega^l_K \left(T^{a} v^h_K \right) (\boldsymbol{x}_K^l)\right| }{ \| v^h_K \|_{V^h_K}}} \nonumber \\
& \leq \| T^a \|_{L^2(K)} \sup_{T^a,v^h_K \in V^h_K}{ \frac{\left|\int_{K} T^{a} v^h_K \,{\rm d}V -  \sum \limits_{l=1}^{l_K} \omega^l_K \left(T^{a} v^h_K \right) (\boldsymbol{x}_K^l)\right| }{ \| T^a \|_{L^2(K)} \| v^h_K \|_{V^h_K}}}, \label{eq:interror_riesz_a}
\end{align}
for \eqref{eq:elem_interrora} and the following similar expression applies to \eqref{eq:elem_interrorb}
\begin{align}
e_K^{b} & \leq \| T^b \|_{L^2(K)} \sup_{T^b,v^h_K \in V^h_K}{ \frac{\left|\int_{K} T^{b} v^h_K \,{\rm d}V -  \sum \limits_{l=1}^{l_K} \omega^l_K \left(T^{b} v^h_K \right) (\boldsymbol{x}_K^l)\right| }{ \| T^b \|_{L^2(K)} \| v^h_K \|_{V^h_K}}}. \label{eq:interror_riesz_b}
\end{align}
For $T^a,v^h_K$ in the polynomial space $V^h_K$, the product $T^a v^h_K$ resides in the double-degree polynomial space $P_K$. By virtue of the equivalence of norms in finite dimensional spaces, $P_K$ can be equipped with a norm $\| \cdot \|_{P_K}$ such that
\begin{align}
\| T^a \|_{L^2(K)} \| v^h_K \|_{V^h_K} \geq \frac{1}{\mathcal{C}_K}  \| T^a v^h_K \|_{P_K},
\end{align}
for a constant $\mathcal{C}_K \geq 1$. A similar inequality holds for $T^b$. Defining $\mathcal{C}_K^a = \mathcal{C}_K \| T^a \|_{L^2(K)}$ and $\mathcal{C}_K^b = \mathcal{C}_K \| T^b \|_{L^2(K)}$, it follows from \eqref{eq:interror_riesz_a} and \eqref{eq:interror_riesz_b} that $e_K^a \leq \mathcal{C}_K^a e_K^p$ and $e_K^b \leq \mathcal{C}_K^b e_K^p$ with
\begin{align}
e_K^{p} = \sup_{p_K \in P_K}{\frac{  \left|\int_{K} p_K(\boldsymbol{x}_K)\,{\rm d}V -  \sum \limits_{l=1}^{l_K} \omega^l_K {p_K}(\boldsymbol{x}_K^l)\right| }{ \| p_K \|_{P_K} } } \label{eq:polyerror}.
\end{align}
Substitution of the operation-independent element integration errors \eqref{eq:polyerror} in the global error bound \eqref{eq:strang1estimatefinal} yields
\begin{equation}
 \left\| u - {u}^h_{\mathcal{Q}} \right\|_{W(h)} \leq \left( 1 + \frac{\left\| a^h \right\|_{W(h),V^h}}{\alpha^h}  \right) \left\| u - \mathcal{I}^h u \right\|_{W(h)} +  
\frac{1}{\alpha^h} 
\sum \limits_{K \in \mathcal{T}^h_\Omega} \left( \mathcal{C}_K^b +  \mathcal{C}_K^a
\right)   e_K^p. \label{eq:errorboundglobal}
\end{equation}
The error bound \eqref{eq:errorboundglobal} conveys that the integration errors can be controlled through the operator-independent elemental integration errors $e_K^p$.

\subsubsection{Integration error evaluation}\label{sec:integrationevaluation}
The element integration error \eqref{eq:polyerror} can be approximated by discretization of the space $P_K$. We select the finite dimensional space $P_K$ as the tensor-product space of univariate polynomials of order $k$ over the element $K$, \emph{i.e.},
\begin{align}
p_K(\boldsymbol{x}_K) &= \prod \limits_{i=1}^d \sum \limits_{j=0}^{k} \alpha_{i,j} x_{K,i}^j = \boldsymbol{\Phi}(\boldsymbol{x}_K)^T \mathbf{a},
\label{eq:polyspace}
\end{align}
with $\boldsymbol{\Phi}(\boldsymbol{x}_K) \in \mathbb{R}^{n_p}$ the basis of monomials of size $n_p = (k+1)^d$ and $\mathbf{a}\in \mathbb{R}^{n_p}$ the corresponding coefficients. The space $P_K$ is equipped with the Sobolev norm $\| \cdot \|_{\Hilbert^{r}}$ with $r \geq 0$. Using the polynomial basis \eqref{eq:polyspace}, the element integration error \eqref{eq:polyerror} can be expressed as
\begin{align}
e_K^{p} &= \sup_{p_K \in \mathbb{P}_K }{\frac{  \left|\int_{K} p_K(\boldsymbol{x}_K)\,{\rm d}V -  \sum \limits_{l=1}^{l_K} \omega^l_K {p_K}(\boldsymbol{x}_K^l)\right| }{ \| p_K \|_{\Hilbert^{r}(K)} } } \nonumber \\
&= {\sqrt{ \sup_{p_K \in \mathbb{P}_K } \frac{ \left|\int_{K} p_K(\boldsymbol{x}_K)\,{\rm d}V -  \sum \limits_{l=1}^{l_K} \omega^l_K {p_K}(\boldsymbol{x}_K^l)\right|^2 }{ \| p_K \|_{\Hilbert^{r}(K)}^2 }} } \nonumber \\
&= {\sqrt{ \sup_{\mathbf{a} \in \mathbb{R}^{n_p} } \frac{ \mathbf{a}^T ( \boldsymbol{\xi}- \bar{\boldsymbol{\xi}}  ) ( \boldsymbol{\xi}- \bar{\boldsymbol{\xi}}  )^T \mathbf{a} }{  \mathbf{a}^T \mathbf{G} \mathbf{a} }} } \nonumber \\
&= {\sqrt{ \lambda_{\rm max}  } }. \label{eq:polyerrorapprox}
\end{align}
where $\lambda_{\rm max}$ is the largest eigenvalue of the generalized eigenvalue problem
\begin{align}
&\left[   \left( \boldsymbol{\xi} - \bar{\boldsymbol{\xi}} \right) \left( \boldsymbol{\xi} - \bar{\boldsymbol{\xi}} \right)^T \right] \mathbf{v}_i = \lambda_i \mathbf{G} \mathbf{v}_i &  &i=1,\ldots,n_p,
\label{eq:eigenvalueproblem}
\end{align}
with eigenvalues $\lambda_i$ and eigenvectors $\mathbf{v}_i$, and with the vectors $\boldsymbol{\xi}$ and $\bar{\boldsymbol{\xi}}$ defined as
\begin{align}
\boldsymbol{\xi} &= \int_{K} \boldsymbol{\Phi}\,{\rm d}V  &  \bar{\boldsymbol{\xi}}  &= \sum \limits_{l=1}^{l_K} \omega^l_K {\boldsymbol{\Phi}}(\boldsymbol{x}_K^l). \label{eq:xivectors}
\end{align}
The matrix $\mathbf{G}$ in \eqref{eq:eigenvalueproblem} is defined as the Gramian matrix for the basis functions $\boldsymbol{\Phi}$ associated with the $\Hilbert^r(K)$ Sobolev space.

Since the left-hand-side matrix in the eigenvalue problem \eqref{eq:eigenvalueproblem} is the dyadic product of two vectors its rank is equal to one. As a consequence, the eigenvalue problem \eqref{eq:eigenvalueproblem} only has one non-zero eigenvalue, with the corresponding eigenvector being in the direction $\mathbf{G}^{-1} ( \boldsymbol{\xi} - \bar{\boldsymbol{\xi}} )$. Since the Gramian matrix is positive definite, this single non-zero eigenvalue is positive and therefore equal to the maximum eigenvalue. This maximum eigenvalue and its corresponding eigenvector are equal to
\begin{align}
 \lambda_{\rm max} &= \| \boldsymbol{\xi} - \bar{\boldsymbol{\xi}} \|^2_{\mathbf{G}^{-1}}, &   
 \mathbf{v}_{\rm max} &= \frac{\mathbf{G}^{-1} ( \boldsymbol{\xi} - \bar{\boldsymbol{\xi}}  )}{  \| \boldsymbol{\xi} - \bar{\boldsymbol{\xi}} \|_{\mathbf{G}^{-1}} } ,
\end{align}
with
\begin{align}
 \| \boldsymbol{\xi} - \bar{\boldsymbol{\xi}} \|_{\mathbf{G}^{-1}}=  \sqrt{  ( \boldsymbol{\xi} - \bar{\boldsymbol{\xi}} )^T  \mathbf{G}^{-1} ( \boldsymbol{\xi} - \bar{\boldsymbol{\xi}} )  } . 
\label{eq:gramianinvnorm}
\end{align}
Substitution of $\lambda_{\rm max}$ in equation \eqref{eq:polyerrorapprox} expresses the element integration error in terms of the vector $\boldsymbol{\xi} - \bar{\boldsymbol{\xi}}$ as:
\begin{align}
  e_K^{p} = \| \boldsymbol{\xi} - \bar{\boldsymbol{\xi}} \|_{\mathbf{G}^{-1}}.
\end{align}
The polynomial corresponding to this error, \emph{i.e.}, the function in $P_K$ yielding the largest integration error follows directly from the maximum eigenvector $\mathbf{v}_{\rm max}$ as
\begin{align}
  p_{K,{\rm max}} =  \mathbf{\Phi}^T \mathbf{v}_{\rm max} = \displaystyle \frac{\mathbf{\Phi}^T  \mathbf{G}^{-1} ( \boldsymbol{\xi} - \bar{\boldsymbol{\xi}}  )}{  \| \boldsymbol{\xi} - \bar{\boldsymbol{\xi}} \|_{\mathbf{G}^{-1}} }. \label{eq:evaluatedpolynomial}
\end{align}
Note that the scaling of this function with respect to the norm \eqref{eq:gramianinvnorm} is chosen such that
\begin{align}
  \|  p_{K,{\rm max}}  \|_{\Hilbert^r(K)}^2 = \mathbf{v}_{\rm max}^T \mathbf{G} \mathbf{v}_{\rm max} =  \displaystyle \frac{( \boldsymbol{\xi} - \bar{\boldsymbol{\xi}}  )^T \mathbf{G}^{-1} ( \boldsymbol{\xi} - \bar{\boldsymbol{\xi}}  )}{  \| \boldsymbol{\xi} - \bar{\boldsymbol{\xi}} \|^2_{\mathbf{G}^{-1}} } = 1,
\end{align}
and hence
\begin{align}
e_K^{p} &= \left|\int_{K} p_{K,{\rm max}}(\boldsymbol{x}_K)\,{\rm d}V -  \sum \limits_{l=1}^{l_K} \omega^l_K {p_{K,{\rm max}}}(\boldsymbol{x}_K^l)\right| .
\label{eq:normedintegrationerror}
\end{align}

\subsection{The cut element quadrature optimization algorithm}\label{sec:thealgorithm}
The ability to evaluate the maximal cut-cell integration error and the corresponding integrand that leads to this error serves as the basis for the quadrature-optimization procedure developed in this work. The error \eqref{eq:normedintegrationerror} is, in general, controlled by the octree depth, $\maxlevel$, and the quadrature rules selected on all sub-cells. In the remainder of this work we assume the octree depth to be fixed at a level where the partitioning error is negligible.

For the space of polynomials considered above, \emph{i.e.}, those of equation~\eqref{eq:polyspace}, the vector $\boldsymbol{\xi}$ of basis function integrals can be evaluated exactly on the exact or approximated geometry and by appropriate selection of the Gauss integration order. As discussed in Section~\ref{sec:octree_partitioning}, the usage of exact Gauss integration on all sub-cells is impractical from a computational effort point of view in large scale finite-cell simulations. The conceptual idea behind the quadrature optimization procedure developed here is to determine integration rules with a significantly smaller number of integration points, and to use the evaluable integration error \eqref{eq:normedintegrationerror} to find the optimal distribution of these points over the sub-cells.

The developed optimization procedure is intended as a per-element pre-processing operation, which results in optimized quadrature rules for all cut elements in a finite cell simulation. In order to conduct this pre-processing operation, for each cut element the exact shape function integrals $\boldsymbol{\xi}$ and the Gramian $\mathbf{G}$ must be determined. These computations are relatively expensive, but especially in time-dependent or non-linear problems -- where the optimised quadrature scheme will be re-used many times -- the computational gain in the simulations outweighs the effort in this pre-processing operation. The computational effort will be studied in Section \ref{sec:cellvslevel}.

In Section~\ref{sec:optimizationproblem} we first formalize the integration error minimization problem that serves as the basis for the optimization algorithm. In Section~\ref{sec:algorithmicaspect} the algorithm is then introduced and various algorithmic aspects are discussed.

\subsubsection{The minimization problem}\label{sec:optimizationproblem}
We consider an octree partitioning of the element $K$ as described in Section~\ref{sec:octree_partitioning} and denote this partitioning by $\partition =  \bigcup \limits_{\level=0}^{\maxlevel+1} \mathcal{P}^{\maxlevel}_{K,\level}$, with $\mathcal{P}^{\maxlevel}_{K,\level}$ the sub-cells corresponding to the octree level $\level$. On each sub-cell $\wp \in \partition$ a sequence of quadrature rules is defined:
\begin{align}
  \left\{ \mathcal{Q}_\wp^\imath \mid  \imath = 0 , 1, 2, \ldots  \right\}.
  \label{eq:quadraturesequence}
\end{align}
We assume that the quadrature rules $\mathcal{Q}_\wp^\imath$ are nested, in the sense that the set of polynomials on $\wp$ that are integrated exactly by $\mathcal{Q}_\wp^{\imath+1}$ includes all polynomials that are integrated exactly by $\mathcal{Q}_\wp^\imath$. The quadrature rule for the complete partitioning is defined as 
\begin{align}
  \mathcal{Q}_{\partition}^{\boldsymbol{\imath}} = \bigcup \limits_{\wp \in \partition} \mathcal{Q}_\wp^{\imath_{\wp}},
\end{align}
with $\boldsymbol{\imath} = \{  \imath_{\wp} \mid  \wp \in  \partition\} $ a list of length $m = \# \partition$ with sub-cell quadrature rule indices. For a given type of integration rule, \emph{e.g.}, Gauss integration, the index list $\boldsymbol{\imath} \in \mathbb{N}^m$ fully determines the quadrature rule for the partitioned cut-cell $K$. 

Both the total number of integration points in the cut element, defined by $\#\mathcal{Q}_\wp^\imath$, and the integration error \eqref{eq:normedintegrationerror} can be expressed as functions of the quadrature index list $\boldsymbol{\imath}$:
\begin{subequations}
\begin{align}
e_K^{p,\boldsymbol{\imath}} &= \left|\int_{K} p_{K,{\rm max}}(\boldsymbol{x}_K)\,{\rm d}V -  \sum \limits_{\wp \in \partition} \sum \limits_{ ( \omega_K, \boldsymbol{x}_K ) \in \mathcal{Q}_\wp^{\imath_\wp} } \omega_K {p_{K,{\rm max}}}(\boldsymbol{x}_K)\right|, \label{eq:intergrationerrorsubcells}%
\end{align}
\end{subequations}
with $ p_{K,{\rm max}}$ according to \eqref{eq:evaluatedpolynomial} for the integration rule $\mathcal{Q}_{\partition}^{\boldsymbol{\imath}}$. Using these functions the intended quadrature optimization procedure can be formulated as the constrained minimization problem
\begin{align}
  \underset{ \boldsymbol{\imath} \in  \mathbb{N}^m }{\mbox{minimize}} \quad e_K^{p,\boldsymbol{\imath}} \quad \mbox{subject to} \quad \#\mathcal{Q}_{\partition}^{\boldsymbol{\imath}} = q^\star,
\label{eq:minimizationproblem}
\end{align}
with $q^\star$ the specified number of integration points. Conversely, we can express the optimization problem in terms of the minimization of the number of integration points for a fixed error
\begin{align}
  \underset{ \boldsymbol{\imath} \in  \mathbb{N}^m }{\mbox{minimize}} \quad \#\mathcal{Q}_{\partition}^{\boldsymbol{\imath}} \quad \mbox{subject to} \quad  e_K^{p,\boldsymbol{\imath}} = e^{\star}_{K},
\end{align}
with $e^{\star}_{K}$ the intended error level.

\begin{remark}[Quadrature refinement strategy]
In this work, the integration quadrature per sub-cell is refined by increasing the order of the integration scheme. Alternatively, refinement would be possible by subdivision of a sub-cell and then constructing quadrature rules over the refined sub-cells. Since we consider the integration of smooth functions over the cut-cells, improving the quadrature by increasing its order is more efficient than spatial refinement of the sub-cells.
\end{remark}

\subsubsection{The optimization algorithm}\label{sec:algorithmicaspect}
The developed algorithm to obtain the optimal distribution of integration points over all sub-cells is presented in Algorithm~\ref{alg:adaptive_integration}. This algorithm approximates the minimization problem \eqref{eq:minimizationproblem} through the generation of a sequence of refinement schemes, $\{ \boldsymbol{\imath}^0, \boldsymbol{\imath}^1, \cdots, \boldsymbol{\imath}^\star \}$, where for the initial quadrature rule the lowest order of integration on each sub-cell is considered, \emph{i.e.}, $\boldsymbol{\imath}^0=\mathbf{0}$, and where the optimized quadrature rule (approximately) satisfying the constraint condition is denoted by $\boldsymbol{i}^\star$. 

Given the $r$-th iterate in the optimization procedure, $\boldsymbol{\imath}^r$, the next integration rule in the sequence, $\boldsymbol{\imath}^{r+1}$, is determined in three steps:
\begin{enumerate}
  \item The function corresponding to the maximum integration error, $p_{K,{\rm max}}$, is determined for an integration rule $\boldsymbol{\imath}^r$ using equation \eqref{eq:evaluatedpolynomial}. The corresponding integration error can be found using \eqref{eq:intergrationerrorsubcells}. This step corresponds to the lines \eqref{alg:line:xivector} -- \eqref{alg:line:interror} in Algorithm~\ref{alg:adaptive_integration}.
  \item The integration error \eqref{eq:intergrationerrorsubcells} is localized to the sub-cells in order to form indicator functions representing the sub-cell-wise error reduction per added integration point. The definition of the sub-cell indicators, which is discussed in detail below, ensures that the algorithm approximates the minimization problem \eqref{eq:minimizationproblem}. This step corresponds to the lines \eqref{alg:line:interror} -- \eqref{alg:line:indicator} in Algorithm~\ref{alg:adaptive_integration}.
  \item The integration rule index of the sub-cells with the largest indicators, $\emph{i.e.}$, with the largest reduction in error per added integration point, is increased by one as to reduce the integration error in these sub-cells. The employed marking strategies are discussed below. This step corresponds to the line \eqref{alg:line:marknigstrategy} in Algorithm~\ref{alg:adaptive_integration}.
\end{enumerate}
These steps are repeated until the specified number of integration points is reached, \emph{i.e.}, when $\#\mathcal{Q}_\wp^{\boldsymbol{\imath}} \geq q^\star$. The algorithm is stopped prematurely when the integration rule sequence \eqref{eq:quadraturesequence} is depleted. Specifically, in the considered implementation \cite{nutils} the maximum quadrature orders for the lowest level tessellation is equal to $6$ ($12$ points) for triangles and $7$ ($31$ points) for tetrahedrons. The algorithm termination conditions correspond to line \eqref{alg:line:stopping} in Algorithm~\ref{alg:adaptive_integration}.

\paragraph{Sub-cell indicators}
To form the sub-cell indicators we consider a uniform refinement of the integration indices, \emph{i.e.}, $\imath^{r+1} = \imath^r +1$ for each sub-cell. The error reduction per added integration point in the case of this uniform refinement can be expressed as
\begin{equation}
- \frac{ e_K^{p,\boldsymbol{\imath}^{r+1}} -  e_K^{p,\boldsymbol{\imath}^{r}} }{ \#\mathcal{Q}_{\partition}^{\boldsymbol{\imath}^{r+1}} - \#\mathcal{Q}_{\partition}^{\boldsymbol{\imath}^r}} \leq  \frac{e_K^{p,\boldsymbol{\imath}^{r}}}{ \#\mathcal{Q}_{\partition}^{\boldsymbol{\imath}^{r+1}} - \#\mathcal{Q}_{\partition}^{\boldsymbol{\imath}^{r}} } . \label{eq:decresingerror}
\end{equation}
Using the sub-cell integration error according to
\begin{align}
e_\wp^{p,\imath} &= \left|\int_{\wp} p_{K,{\rm max}}(\boldsymbol{x}_K)\,{\rm d}V - \sum  \limits_{ ( \omega_K, \boldsymbol{x}_K ) \in \mathcal{Q}_\wp^{\imath} } \omega_K {p_{K,{\rm max}}}(\boldsymbol{x}_K) \right|, \label{eq:subcellintegrationerror}
\end{align}
and the reduction in the number of integration points\footnote{\vspace{-15pt}\begin{align*}
\#\mathcal{Q}_{\partition}^{\boldsymbol{\imath}^{r+1}} - \#\mathcal{Q}_{\partition}^{\boldsymbol{\imath}^r} & \geq m \min \left(   \#\mathcal{Q}_\wp^{\imath^{r+1}} - \#\mathcal{Q}_\wp^{\imath^{r}}  \right) \geq m \frac{\min \left(   \#\mathcal{Q}_\wp^{\imath^{r+1}} - \#\mathcal{Q}_\wp^{\imath^{r}}  \right)}{\max \left(   \#\mathcal{Q}_\wp^{\imath^{r+1}} - \#\mathcal{Q}_\wp^{\imath^{r}}  \right)} \max \left(   \#\mathcal{Q}_\wp^{\imath^{r+1}} - \#\mathcal{Q}_\wp^{\imath^{r}}  \right) \\ & \geq m \frac{\min \left(   \#\mathcal{Q}_\wp^{\imath^{r+1}} - \#\mathcal{Q}_\wp^{\imath^{r}}  \right)}{\max \left(   \#\mathcal{Q}_\wp^{\imath^{r+1}} - \#\mathcal{Q}_\wp^{\imath^{r}}  \right)} \left(   \#\mathcal{Q}_\wp^{\imath^{r+1}} - \#\mathcal{Q}_\wp^{\imath^{r}}  \right).\end{align*}}, the error reduction \eqref{eq:decresingerror} is bounded as:
\begin{align}
\frac{e_K^{p,\boldsymbol{\imath}^{r}}}{ \#\mathcal{Q}_{\partition}^{\boldsymbol{\imath}^{r+1}} - \#\mathcal{Q}_{\partition}^{\boldsymbol{\imath}^{r}} }
  &\leq \sum \limits_{ \left( \imath, \wp \right) \in \left( \boldsymbol{\imath}, \partition\right)} \frac{ e_\wp^{p,\imath^{r}}}{ \#\mathcal{Q}_{\partition}^{\boldsymbol{\imath^{r+1}}} - \#\mathcal{Q}_{\partition}^{\boldsymbol{\imath^{r}}}} \nonumber  \\
&\leq \frac{1}{m} \frac{\max \limits_{ \left( \imath, \wp \right) \in \left( \boldsymbol{\imath}, \partition\right)} \left(   \#\mathcal{Q}_\wp^{\imath^{r+1}} - \#\mathcal{Q}_\wp^{\imath^{r}}  \right)  }{\min \limits_{ \left( \imath, \wp \right) \in \left( \boldsymbol{\imath}, \partition\right)} \left(   \#\mathcal{Q}_\wp^{\imath^{r+1}} - \#\mathcal{Q}_\wp^{\imath^{r}}  \right) }  \sum \limits_{ \left( \imath, \wp \right) \in \left( \boldsymbol{\imath}, \partition\right)} \frac{ e_\wp^{p, \imath^{r} }} { \left(   \#\mathcal{Q}_\wp^{\imath^{r+1}} - \#\mathcal{Q}_\wp^{\imath^{r}}  \right) } \nonumber \\
&\leq \frac{c}{m}\sum \limits_{ \left( \imath, \wp \right) \in \left( \boldsymbol{\imath}, \partition\right)} \mathcal{R}_\wp^{\imath^{r}},
\label{eq:estimatorbound}
\end{align}
with the sub-cell indicator function defined as 
\begin{align}
\mathcal{R}_\wp^{\imath} = \frac{e_\wp^{p,\imath }}{\#\mathcal{Q}_\wp^{\imath + 1}-\#\mathcal{Q}_\wp^{\imath}}.
\label{eq:indicator}
\end{align}

The sub-cell indicator \eqref{eq:indicator} weighs the local integration error by the cost of increasing the integration order. Under the assumption that the sub-cell integration error reduces significantly when increasing the order of the integration scheme by one, the numerator in \eqref{eq:indicator} can be interpreted as the error reduction rather than the error itself. This assumption of local error reduction is, however, in general, not satisfied. As a consequence, in the proposed algorithm sub-cells whose integration error does not decrease significantly (or even increases) by raising the order of integration by one would be underrated in the indicator function. Considering the error rather than the error reduction in the indicator provides a more robust measure of the potential to decrease the error in a sub-cell.

\paragraph{Sub-cell marking strategy} \label{sec:refinementmarker}
Based on the upper bound \eqref{eq:estimatorbound} of the error reduction per added integration point, sub-cells are marked for refinement. We herein consider two marking strategies, with $\mathcal{M}$ denoting the set of sub-cells marked for increasing the integration order:
\begin{itemize}
  \item \emph{Sub-cell marking:} In this marking strategy we increase the integration order of the single sub-cell with the largest indicator:
\begin{align}
  \mathcal{M}  = \underset{\left( \imath, \wp \right) \in \left( \boldsymbol{\imath}, \partition\right)}{\rm arg\,max} \left( \mathcal{R}_\wp^{\imath} \right).
\end{align}
  \item \emph{Octree level marking:} Sub-cell marking involves refinement of individual sub-cell in each refinement step. This generally leads to a large number of optimization steps, which is undesirable from a computational effort point of view. In order to expedite the optimization procedure, it is possible to mark multiple sub-cells simultaneously. A natural way of marking multiple sub-cells is to agglomerate all sub-cells in a particular level. This leads to a marking strategy in which we increase the integration order for the level (in the octree partitioning) with the largest indicator:
\begin{align}
  \mathcal{M}  = \underset{ \{ \mathcal{P}^{\level}_{K} \}_{\level=0}^{\maxlevel+1} }{\rm arg\,max} \left( \sum \limits_{\left( \imath, \wp \right) \in \left( \boldsymbol{\imath},  \mathcal{P}^{\level}_{K}\right)} \mathcal{R}_\wp^{\imath} \right).
\end{align}
\end{itemize}

\begin{remark}[Global adaptive integration]
\label{rem:globalmarking}
  The algorithm presented above is framed in an element-by-element setting, \emph{i.e.}, each cut element is considered separately in the pre-processing operation to determine the integration schemes. Since the evaluation of the integration error and its indicators involves the repeated computation of solutions to a linear system with the size equal to the number of supported basis functions, \emph{i.e.}, equation~\eqref{eq:normedintegrationerror}, global application of the algorithm is impractical from a computational effort point of view.

It should be noted that the element-by-element setting considered here does not account for the fact that the operator-dependent constants in the multiplicative decomposition \eqref{eq:polyerror} in principle vary per element. For example, a source term might be negligibly small on an element that is challenging to integrate, and hence improving quadrature on that element would be inefficient from the error approximation point of view. Additional integration points are, however, assigned to such a cut element in the element-by-element strategy employed here, as the operator is not considered in the marking strategy.

Global adaptive integration is feasible, however, by computing the indicators on all elements individually (\emph{i.e.}, solving multiple relatively small local linear systems, rather than a large global system), and then to apply a global marking strategy. Both the sub-cell marking and the level marking strategy can then be applied globally. The potential benefit of such a global marking strategy is that it automatically accounts for the fact that not all cut elements are equally hard to integrate. However, in principle, the operator-dependent constants should then be incorporated in the error indicators. Moreover, one should be aware that parallelization of the pre-processing operation, which is trivial in the per-element setting, is more challenging using the global marking strategy, as communication between the elements is required. The development of a rigorous global integration optimization routine, which would consider the above-mentioned complications, is considered beyond the scope of the current work. 
\end{remark}

\begin{algorithm}[H]
	\DontPrintSemicolon
	\SetNoFillComment
	\SetAlgoLined
	\SetKwComment{Comment}{\color{blue}\#\,\color{blue}}{}
	\SetKwComment{tcc}{\color{blue}\#\,\color{blue}}{}
	
	\SetStartEndCondition{ }{}{}
	\SetKw{KwTo}{to}
	\SetKw{KwFrom}{from}
	\SetKwFor{For}{for}{:}{end}
	\SetKwFor{While}{while}{:}{end}
	\SetKw{Not}{not}
	\SetKw{In}{in}
	\SetKw{True}{True}
	
	\SetKwFunction{basis}{get\_monomial\_basis}
	\SetKwFunction{quadrature}{get\_quadrature}
	\SetKwFunction{marking}{mark\_sub\_cells}
	
	\KwIn{ $\mathcal{P}= \cup_{0}^{\maxlevel+1} \mathcal{P}_K^\level = \cup \wp$, $k$, $marking\_strategy$}
		\KwOut{ $ \mathcal{Q}_{\partition}^{\boldsymbol{\imath}^\star} $  {\Comment*[r]{optimized quadrature rule}}}
		
	\BlankLine
	\BlankLine
		
    \Comment{Initialization}
		$\boldsymbol{\imath} = \boldsymbol{0}$ {\Comment*[r]{quadrature index}}
		
		$\boldsymbol{\Phi} = $\basis{$k$}  {\Comment*[r]{monomial basis}}
		
	\BlankLine
    \Comment{Reference integrals}
		$\xi_i = \int_{\partition} \Phi_i\,{\rm d}V$ \label{alg:line:xivector} {\Comment*[r]{basis function integral vector}}
		
		$G_{ij} = \left(\Phi_i, \Phi_j \right)_{\Hilbert^1}$ \label{alg:line:gramianmatrix} {\Comment*[r]{Gramian matrix}}
		
		\BlankLine

		\While{\Not terminate $\,$}{
			
			\BlankLine
			\Comment{Quadrature rule}
				$\mathcal{Q}_{\partition}^{\boldsymbol{\imath}}=$ \quadrature{$\mathcal{P}$} {\Comment*[r]{quadrature rule}}

				$\bar{\boldsymbol{\xi}}  = \sum \limits_{( \omega_{\partition}, \boldsymbol{x}_{\partition} ) \in \mathcal{Q}_{\partition}^{\boldsymbol{\imath}}} \omega_{\partition} {\boldsymbol{\Phi}}(\boldsymbol{x}_{\partition})$ \label{alg:line:xibarvector} {\Comment*[r]{approximate basis function integral vector}}
			
			\BlankLine
			\Comment{Worst possible function to integrate}
				$\mathbf{v}_{\rm max} = \dfrac{\mathbf{G}^{-1} ( \boldsymbol{\xi} - \bar{\boldsymbol{\xi}}  )}{  \| \boldsymbol{\xi} - \bar{\boldsymbol{\xi}} \|_{\mathbf{G}^{-1}} }$ \label{alg:line:eigenvalue} {\Comment*[r]{eigenvector}}
			
				$  p_{\rm max} = \mathbf{\Phi}^T \mathbf{v}_{\rm max}$ \label{alg:line:functioneval} {\Comment*[r]{eigenfunction}}
			
			\BlankLine
			\Comment{Sub-cell errors and indicators}
				\For{$\wp$ \In $\partition$}{
					$e_\wp^{p,\imath_{\wp}} = \left|\int_{\wp} p_{\rm max}(\boldsymbol{x}_{\partition})\,{\rm d}V - \sum  \limits_{ ( \omega_{\wp}, \boldsymbol{x}_{\wp} ) \in \mathcal{Q}_\wp^{\imath_{\wp}} } \omega_{\wp} {p_{\rm max}}(\boldsymbol{x}_{\wp}) \right|$ {\Comment*[r]{error}} \label{alg:line:interror} 
			
				$ \mathcal{R}_\wp^{\imath_{\wp}} = \frac{e_\wp^{p,\imath_{\wp} }}{\#\mathcal{Q}_\wp^{\imath_{\wp}+1}-\#\mathcal{Q}_\wp^{\imath_{\wp}}} $ \label{alg:line:indicator} {\Comment*[r]{indicator}}}

			\BlankLine
		    \Comment{Sub-cells marked for refinement}
				$\mathcal{M}$ = \marking{$\{ \mathcal{R}_\wp^{\imath_{\wp}} \}_{\wp \in \mathcal{P}}$, marking\_strategy} \label{alg:line:marknigstrategy}

				\For{$\wp$ \In $\mathcal{M} $}{
					${\imath}_{\wp} \rightarrow {\imath}_{\wp} + {1}$
					
					\Comment{stopping criterion}
					\lIf{\Not $\mathcal{Q}_{\mathcal{\wp}}^{\imath}$ \,}{$terminate=$ \True \label{alg:line:stopping}}
				}
		}
		\Return{$\mathcal{Q}_{\partition}^{\boldsymbol{\imath}^\star}$}
		\caption{Adaptive integration algorithm.}
		\label{alg:adaptive_integration}
	\end{algorithm}

\section{Numerical study of the adaptive integration procedure} \label{sec:integrationerroranalysis}
To assess the performance of the developed adaptive integration technique, in this section we study its behavior in terms of integration accuracy versus the number of integration points. We here consider the case of a single cut-cell. The effect of the integration accuracy on actual finite-cell simulations will be studied subsequently in Section~\ref{sec:numerical_examples}.

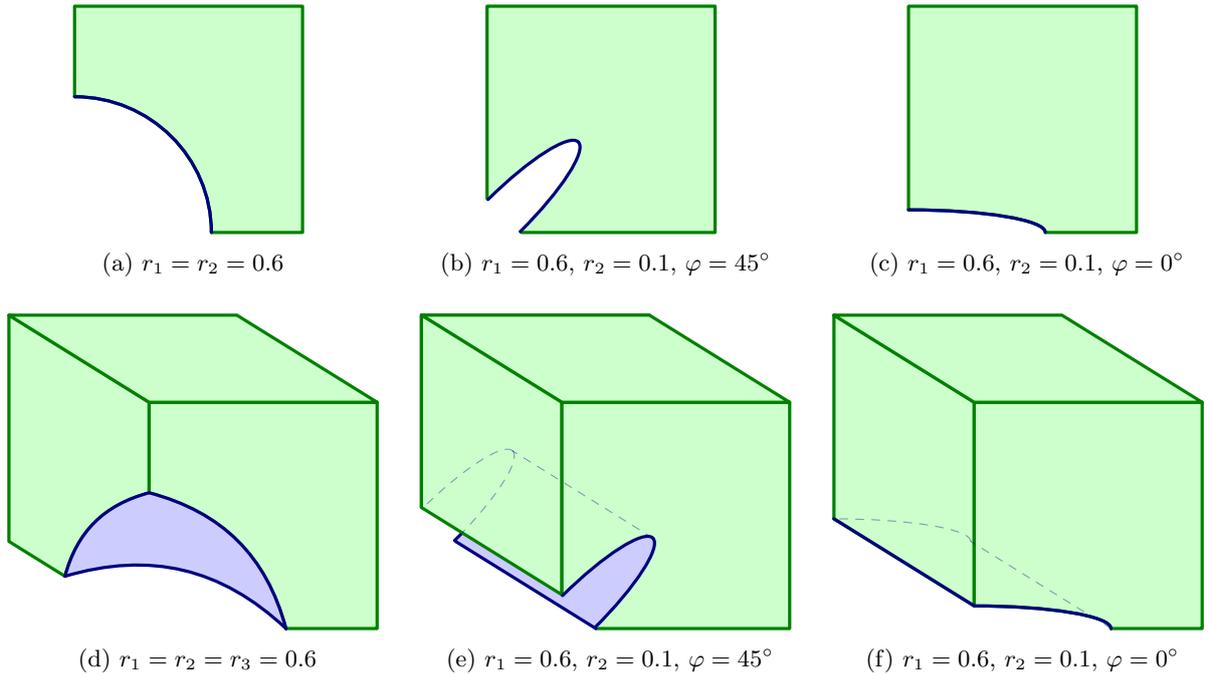
\begin{figure}
	\centering
	\begin{subfigure}[t]{0.33\textwidth}
		\centering
		\begin{tikzpicture}[
scale=3,
line cap=round,
line join=round,
opaque/.style={black,thin,opacity=0.4}
]
	\coordinate (A) at (0,0);
	\coordinate (B) at (1,0);
	\coordinate (C) at (1,1);
	\coordinate (D) at (0,1);

	\coordinate (AB) at (0.6,0);
	\coordinate (AD) at (0,0.6);
	
	\coordinate (origin) at (0.1,0.1);
	\coordinate (x1) at (0.3,0.1);
	\coordinate (x2) at (0.1,0.3);

	\filldraw[draw=green!50!black,very thick,fill=green!20!white] (B) -- (AB) arc[start angle=0, end angle=90, radius=0.6] -- (AD) -- (D) -- (C) -- cycle;
	\draw[draw=blue!50!black,very thick] (AB) arc[start angle=0, end angle=90, radius=0.6];
\end{tikzpicture}
		\caption{$r_{1}=r_{2}=0.6$}
	\end{subfigure}%
	\begin{subfigure}[t]{0.33\textwidth}
		\centering
	    \begin{tikzpicture}[
scale=3,
line cap=round,
line join=round,
opaque/.style={black,thin,opacity=0.4}
]

	\coordinate (A) at (0,0,0);
	\coordinate (B) at (1,0,0);
	\coordinate (C) at (1,1,0);
	\coordinate (D) at (0,1,0);

	\coordinate (AB) at (0.15,0);
	\coordinate (AD) at (0,0.15);
	\coordinate (tip) at (0.4,0.4);

	\filldraw[draw=green!20!white,very thick,fill=green!20!white,rotate=45] (tip) arc(0:90:0.46 and 0.1) -- (AD) -- (D) -- (C) -- cycle;
	\filldraw[draw=green!20!white,very thick,fill=green!20!white,rotate=45] (tip) arc(0:-90:0.46 and 0.1) -- (AB) -- (B) -- (C) -- cycle;
	\draw[draw=green!50!black,very thick] (AD) -- (D) --(C) -- (B) -- (AB);
	\draw[draw=blue!50!black,very thick,rotate=45] (tip) arc(0:90:0.46 and 0.1);
	\draw[draw=blue!50!black,very thick,rotate=45] (tip) arc(0:-90:0.46 and 0.1);
\end{tikzpicture}
		\caption{$r_{1}=0.6$, $r_{2}=0.1$, $\varphi=45^\circ$}
	\end{subfigure}
	\begin{subfigure}[t]{0.33\textwidth}
		\centering
		\begin{tikzpicture}[
scale=3,
line cap=round,
line join=round,
opaque/.style={black,thin,opacity=0.4}
]

	\coordinate (A) at (0,0,0);
	\coordinate (B) at (1,0,0);
	\coordinate (C) at (1,1,0);
	\coordinate (D) at (0,1,0);

	\coordinate (AB) at (0.6,0);
	\coordinate (AD) at (0,0.1);

	\filldraw[draw=green!50!black,very thick,fill=green!20!white,rotate=0] (AB) arc(0:90:0.6 and 0.1) -- (AD) -- (D) -- (C) -- (B) -- cycle;
	\draw[draw=blue!50!black,very thick,rotate=0] (AB) arc(0:90:0.6 and 0.1);
\end{tikzpicture}
		\caption{$r_{1}=0.6$, $r_{2}=0.1$, $\varphi=0^\circ$}
	\end{subfigure}\\[12pt]
	\begin{subfigure}[t]{0.33\textwidth}
		\centering
		\begin{tikzpicture}[
  scale=3,
  line cap=round,
  line join=round,
  opaque/.style={black,thin,opacity=0.4}
]

  \coordinate (A) at (0,0,0);
  \coordinate (B) at (1,0,0);
  \coordinate (C) at (1,1,0);
  \coordinate (D) at (0,1,0);
  \coordinate (E) at (1,0,1);
  \coordinate (F) at (2,0,1);
  \coordinate (G) at (2,1,1);
  \coordinate (H) at (1,1,1);
  
  \coordinate (AE) at (0.4,0,0.4);
  \coordinate (EF) at (1.6,0,1);
  \coordinate (EH) at (1,0.6,1);
  
  \coordinate (origin) at (1.1,0,0.9);
  \coordinate (x1) at (1.28,0,0.9);
  \coordinate (x2) at (1.1,0.18,0.9);
  \coordinate (x3) at (0.9,0,0.7);
    
  \filldraw[draw=green!50!black,very thick,fill=green!20!white] (EH) to[bend right] (AE) -- (A) -- (D) -- (H) -- (EH);
  \filldraw[draw=green!50!black,very thick,fill=green!20!white] (EF) to[bend right] (EH) -- (H) -- (G) -- (F) -- (EF);
  \filldraw[draw=green!50!black,very thick,fill=green!20!white] (H) -- (G) -- (C) -- (D) -- cycle;
  \filldraw[draw=blue!50!black,very thick,fill=blue!20!white] (EF) to[bend right] (EH) to [bend right] (AE) to[bend left] (EF);
%
\end{tikzpicture}
		\caption{$r_{1}=r_{2}=r_{3}=0.6$}
	\end{subfigure}%
	\begin{subfigure}[t]{0.33\textwidth}
		\centering
		\begin{tikzpicture}[
  scale=3,
  line cap=round,
  line join=round,
  opaque/.style={black,thin,opacity=0.4}
]

  \coordinate (A) at (0,0,0);
  \coordinate (B) at (1,0,0);
  \coordinate (C) at (1,1,0);
  \coordinate (D) at (0,1,0);
  \coordinate (E) at (1,0,1);
  \coordinate (F) at (2,0,1);
  \coordinate (G) at (2,1,1);
  \coordinate (H) at (1,1,1);
  
  \coordinate (AB) at (0.15,0,0);
  \coordinate (AD) at (0,0.15,0);
  \coordinate (tip1) at (0.4,0.4,0);
  
  \coordinate (EF) at (1.15,0,1);
  \coordinate (EH) at (1,0.15,1);
  \coordinate (tip2) at (1.4,0.4,1);
  \draw[draw=blue!20!white,fill=blue!20!white,rotate=45] (tip2) arc(0:90:0.46 and 0.1) -- (EH) -- (EF) -- (tip2);
  \draw[draw=blue!20!white,fill=blue!20!white,rotate=45] (tip2) arc(0:-90:0.46 and 0.1) -- (EF) -- (EH) -- (tip2);
  \draw[draw=blue!20!white,fill=blue!20!white,rotate=45] (tip1) arc(0:-90:0.46 and 0.1) -- (AB) -- (EF) -- (tip1);
  \draw[draw=blue!50!black,very thick,rotate=45] (tip1) arc(0:-90:0.46 and 0.1);
  \filldraw[draw=green!20!white,very thick,fill=green!20!white,rotate=45] (tip1) arc(0:90:0.46 and 0.1) -- (AD) -- (D) -- (C) -- cycle;
  \draw[draw=green!50!black,very thick] (AD) -- (D) --(C);
  \filldraw[draw=green!20!white,very thick,fill=green!20!white,rotate=45] (tip2) arc(0:90:0.46 and 0.1) -- (EH) -- (H) -- (G) -- cycle;
  \filldraw[draw=green!20!white,very thick,fill=green!20!white,rotate=45] (tip2) arc(0:-90:0.46 and 0.1) -- (EF) -- (F) -- (G) -- cycle;
  \draw[draw=green!50!black,very thick] (EH) -- (H) --(G) -- (F) -- (EF);
  \filldraw[draw=green!50!black,very thick,fill=green!20!white] (D) --(C) -- (G) -- (H) -- cycle;
  \filldraw[draw=green!50!black,very thick,fill=green!20!white] (AD) --(D) -- (H) -- (EH) -- cycle;
  \draw[draw=blue!50!black,dashed,thin,opacity=0.4,rotate=45] (tip1) arc(0:90:0.46 and 0.1);
  \draw[draw=blue!50!black,dashed,thin,opacity=0.4,rotate=45] (tip1) arc(0:-90:0.46 and 0.1);
  \draw[draw=blue!50!black,very thick,rotate=45] (tip2) arc(0:90:0.46 and 0.1);
  \draw[draw=blue!50!black,very thick,rotate=45] (tip2) arc(0:-90:0.46 and 0.1);
  \draw[draw=blue!50!black,dashed,thin,opacity=0.4] (tip2) -- (tip1);
  \draw[draw=blue!50!black,very thick] (EF) -- (AB);
\end{tikzpicture}
		\caption{$r_{1}=0.6$, $r_{2}=0.1$, $\varphi=45^{\circ}$}
	\end{subfigure}%
	\begin{subfigure}[t]{0.33\textwidth}
		\centering
		\begin{tikzpicture}[
  scale=3,
  line cap=round,
  line join=round,
  opaque/.style={black,thin,opacity=0.4}
]

  \coordinate (A) at (0,0,0);
  \coordinate (B) at (1,0,0);
  \coordinate (C) at (1,1,0);
  \coordinate (D) at (0,1,0);
  \coordinate (E) at (1,0,1);
  \coordinate (F) at (2,0,1);
  \coordinate (G) at (2,1,1);
  \coordinate (H) at (1,1,1);
  
  \coordinate (AB) at (0.6,0,0);
  \coordinate (AD) at (0,0.1,0);
  
  \coordinate (EF) at (1.6,0,1);
  \coordinate (EH) at (1,0.1,1);
  \filldraw[draw=green!50!black,very thick,fill=green!20!white] (AD) -- (D) -- (H) -- (EH) -- (AD);
  \filldraw[draw=green!50!black,very thick,fill=green!20!white,rotate=0] (EF) arc(0:90:0.6 and 0.1) -- (EH) -- (H) -- (G) -- (F) -- (EF);
  \filldraw[draw=green!50!black,very thick,fill=green!20!white] (D) -- (C) -- (G) -- (H) -- (D);
  \draw[draw=blue!50!black,dashed,thin,opacity=0.4,rotate=0] (AB) arc(0:90:0.6 and 0.1);
  \draw[draw=blue!50!black,very thick,rotate=0] (EF) arc(0:90:0.6 and 0.1);
  \draw[draw=blue!50!black,dashed,thin,opacity=0.4] (EF) -- (AB);
  \draw[draw=blue!50!black,very thick] (EH) -- (AD);
\end{tikzpicture}
		\caption{$r_{1}=0.6$, $r_{2}=0.1$, $\varphi=0^{\circ}$}
	\end{subfigure}
	\caption{Schematic representation of the two-dimensional (top row) and three-dimensional (bottom row) cut-element corresponding to the level set function \eqref{eq:levelset}.}
	\label{fig:singlecutcell}
\end{figure}

We consider a $d$-dimensional unit cube $[0,1]^d$ in two and three dimensions as a single element of the background mesh. Throughout this section we exclude an ellipsoid with semi-major axis $r_{1}$    and semi-minor axes $r_{2}$ (with the sphere as the special case $r_{1}=r_{2}$)  that is centered at the origin of the background element and with its major axis residing in the $x_1-   x_2$ plane at an inclination of $\varphi$ with respect to the $x_1$ axis. The resulting cut-element corresponds to the subdomain for which the level set function
\begin{align}
l(\boldsymbol{x}) = \left( \frac{\bar{x}_1}{r_{1}} \right)^2 + \sum_{\delta = 2}^d \left( \frac{\bar{x}_\delta}{r_{2}} \right)^2
\label{eq:levelset}
\end{align}
is positive, where $\bar{x}_1= x_1 \cos{\varphi} - x_2 \sin{\varphi}$, $\bar{x}_2= x_1 \sin{\varphi} +   x_2 \cos{\varphi}$, and, if $d=3$, $\bar{x}_3=x_3$. Schematics of various cut-elements generated using the octree procedure, discussed in Section~\ref{sec:octree_partitioning}, with a maximum octree depth of $\maxlevel=3$ for a spherical exclusion and $\maxlevel=4$ for an elliptical exclusion are displayed in Figure~\ref{fig:singlecutcell}. Here, the maximum octree depth is chosen such that the geometric error is negligible.

\subsection{Equal-order degree integration}\label{sec:equalorder}
\begin{figure}
	\centering
	\begin{tabular}[t]{@{}p{0.45\textwidth}@{}}
		\begin{subfigure}{\hsize}\centering
			\includegraphics[width=\textwidth]{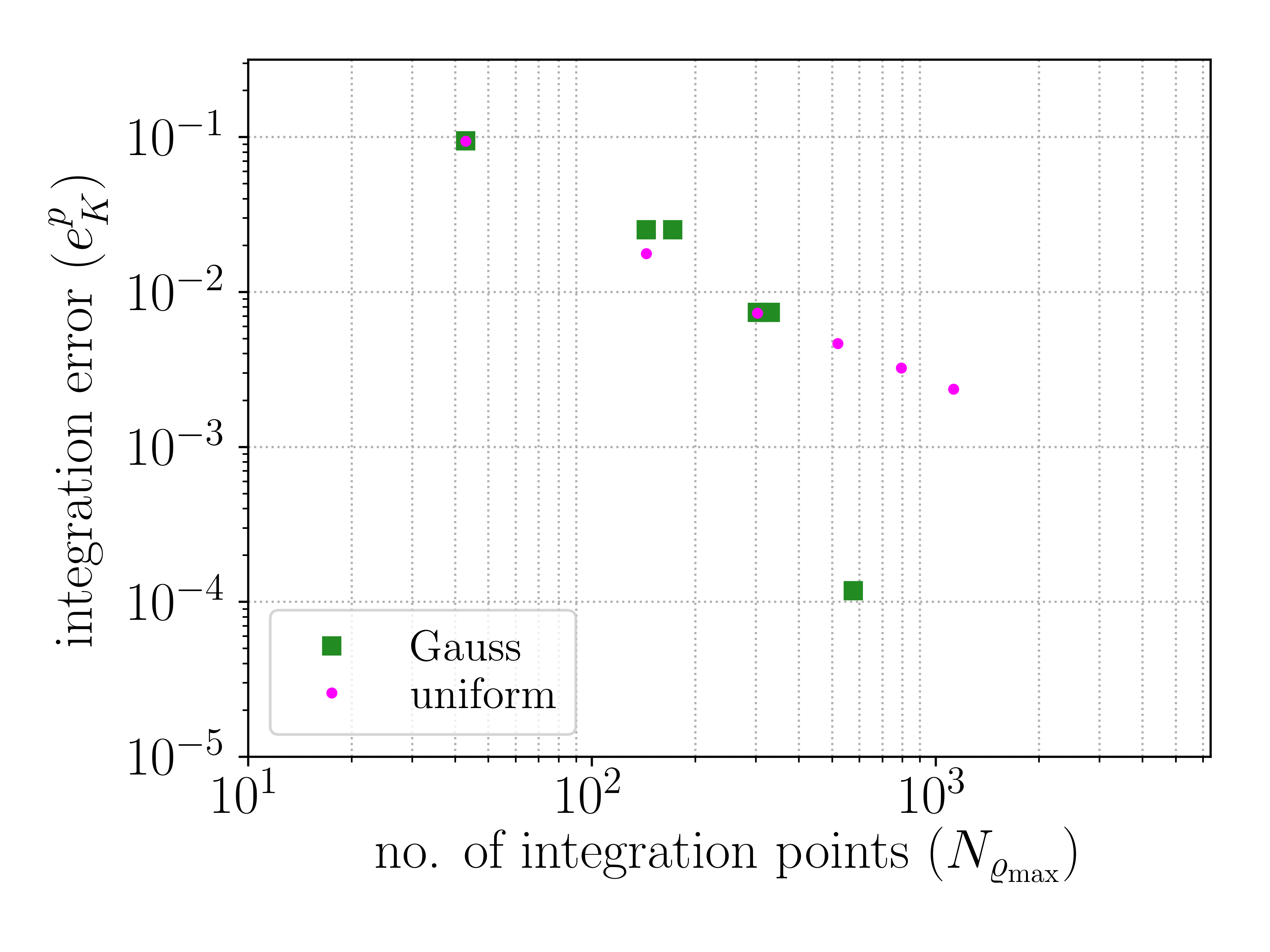}
			\caption{$d=2$, $\maxlevel=3$}
		\end{subfigure}\\[12pt]
		\begin{subfigure}{\hsize}\centering
			\includegraphics[width=\textwidth]{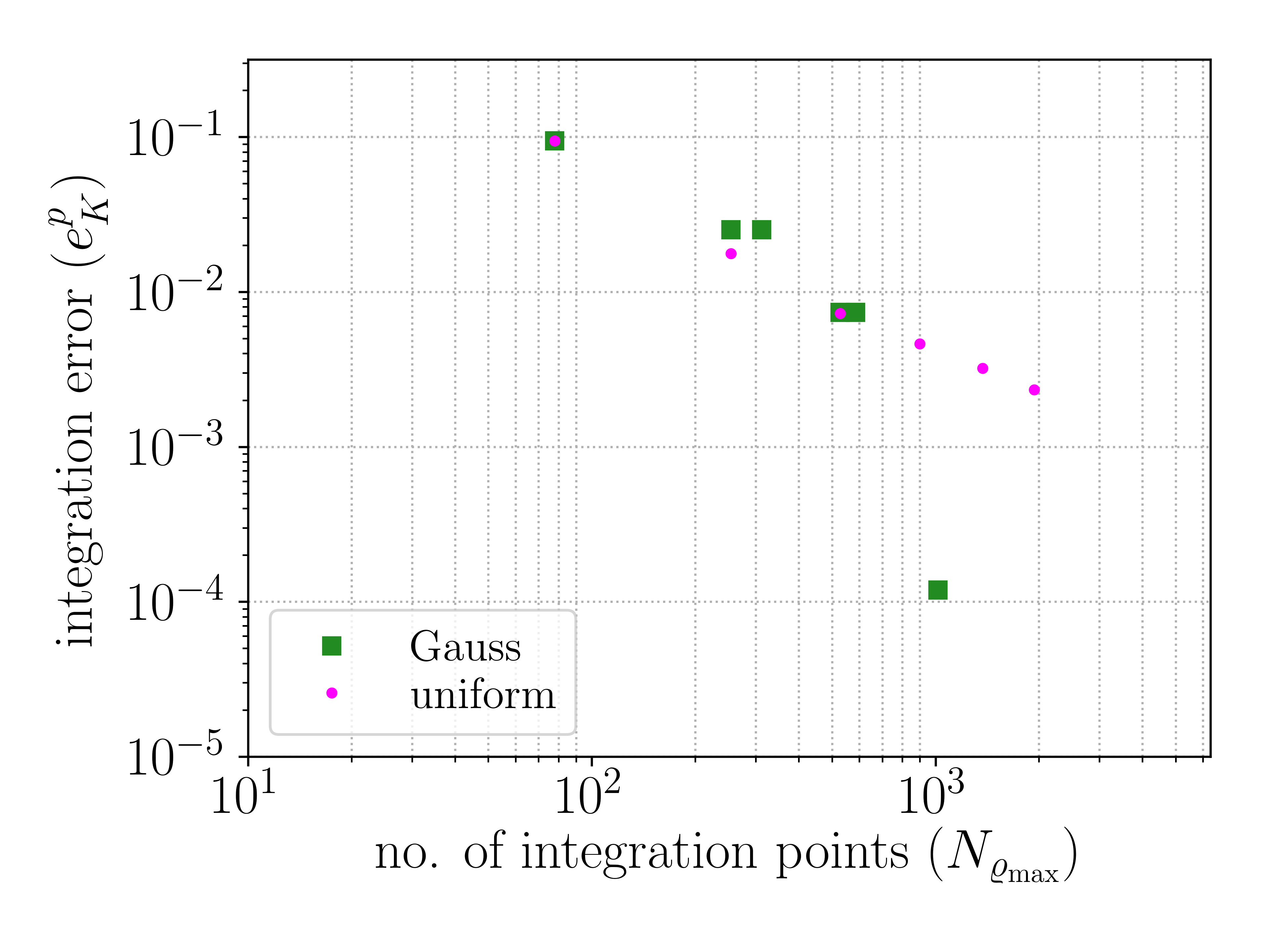}
			\caption{$d=2$, $\maxlevel=4$}
		\end{subfigure}\\[12pt]
		\begin{subfigure}{\hsize}\centering
			\includegraphics[width=\textwidth]{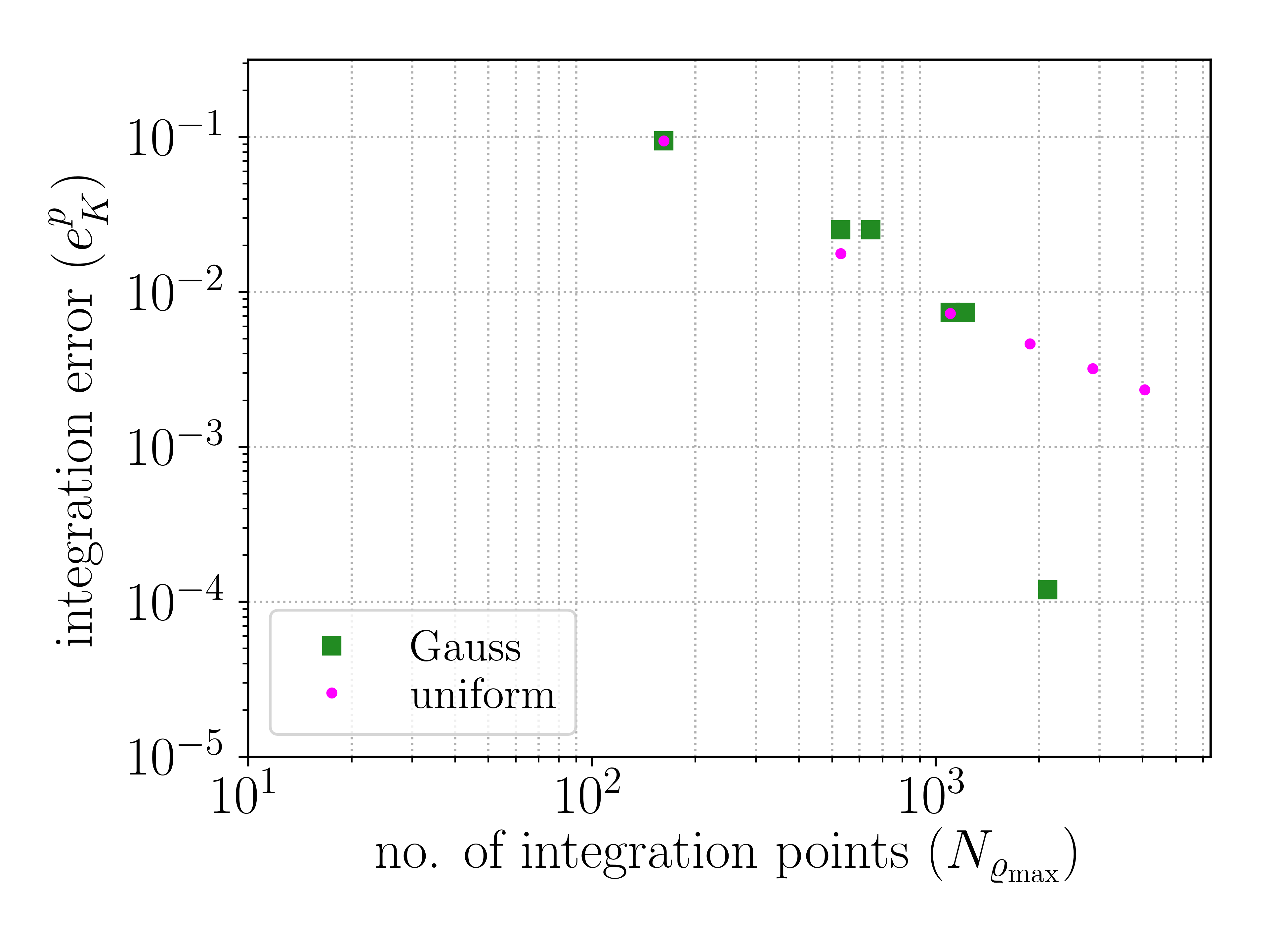}
			\caption{$d=2$, $\maxlevel=5$}
		\end{subfigure}%
	\end{tabular}
	\begin{tabular}[t]{@{}p{0.45\textwidth}@{}}
		\begin{subfigure}{\hsize}\centering
			\includegraphics[width=\textwidth]{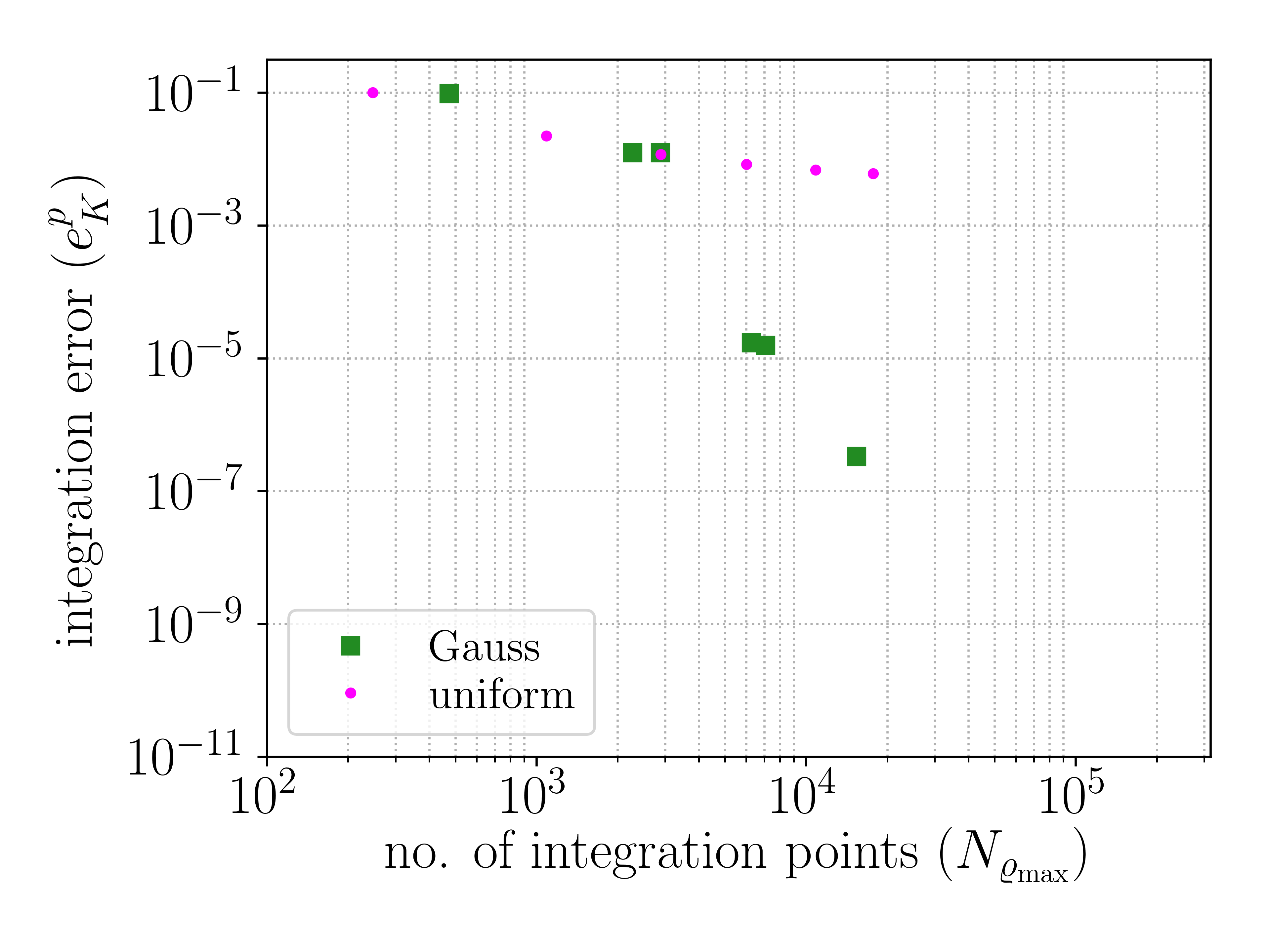}
			\caption{$d=3$, $\maxlevel=2$}
		\end{subfigure}\\[12pt]
		\begin{subfigure}{\hsize}\centering
			\includegraphics[width=\textwidth]{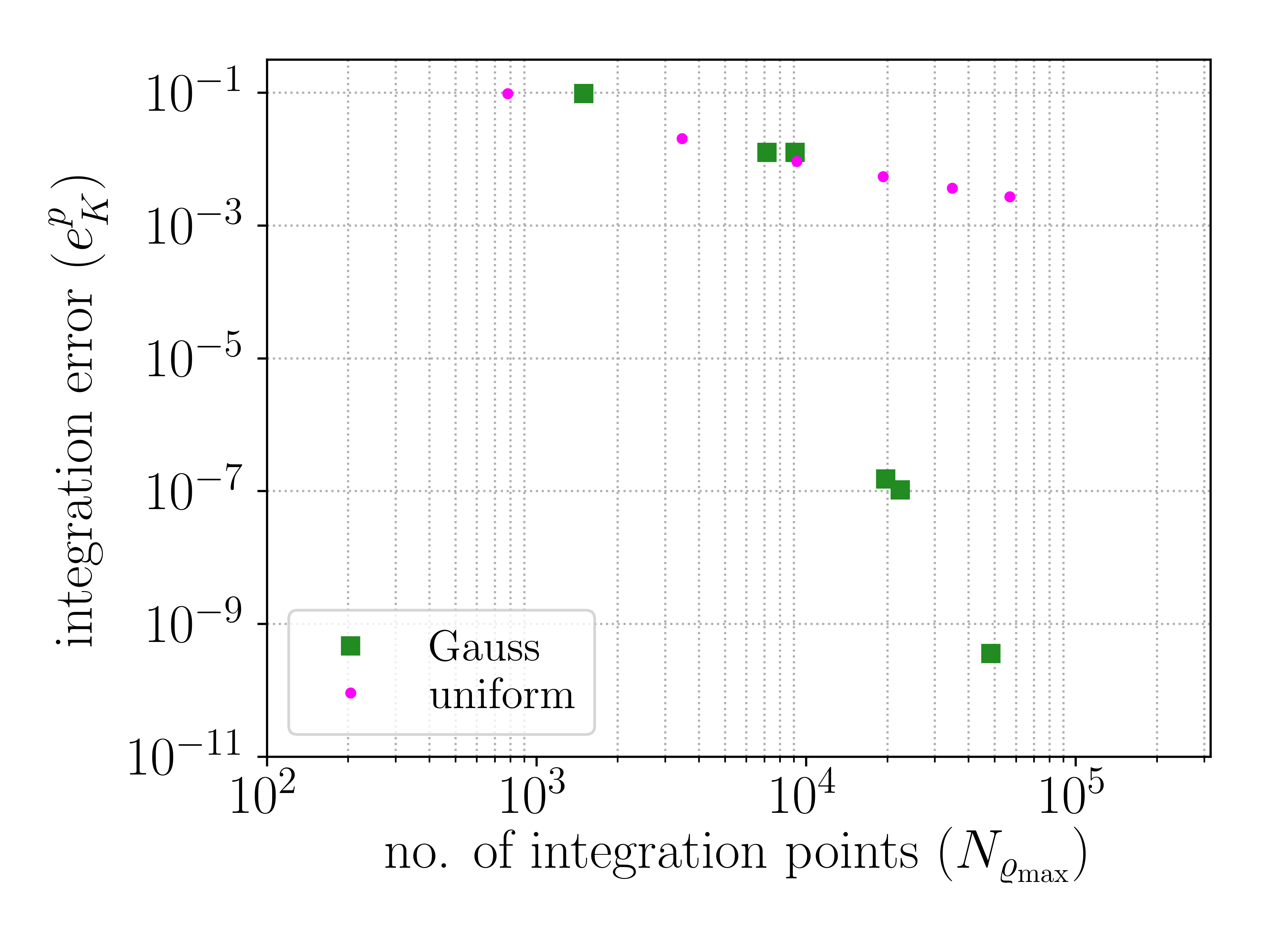}
			\caption{$d=3$, $\maxlevel=3$}
		\end{subfigure}\\[12pt]
		\begin{subfigure}{\hsize}\centering
			\includegraphics[width=\textwidth]{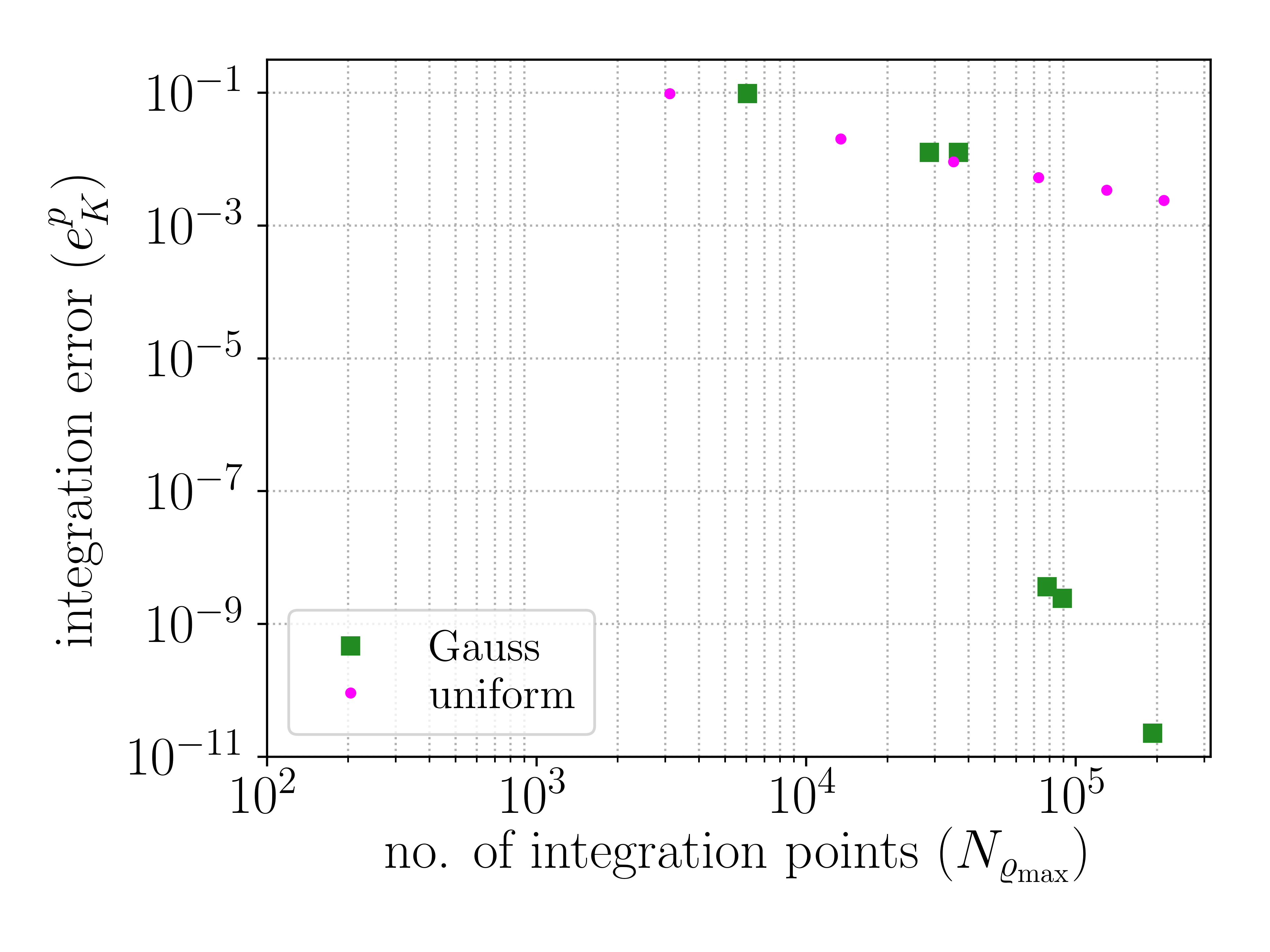}
			\caption{$d=3$, $\maxlevel=4$}
		\end{subfigure}
	\end{tabular}
	\caption{Integration errors for a spherical exclusion with radius $0.6$ using the same integration scheme on all sub-cells.}
	\label{fig:equalorders}
\end{figure}

Before we proceed with the presentation of the results obtained using the adaptive integration procedure, we first consider the case for which the order of the integration scheme is chosen to be same over all integration sub-cells. We consider the case of a spherical exclusion with $r_1=r_2=0.6$. In Figure~\ref{fig:equalorders} the integration error \eqref{eq:normedintegrationerror} is plotted for the cases of $\maxlevel=3,4,5$ in two dimensions and $\maxlevel = 2,3,4$ in three dimensions. In each of the sub plots the integration index $\boldsymbol{\imath}$ is increased uniformly over all sub-cells. For all presented results the polynomial order $k=8$ is considered for two dimensions and $k=5$ is considered for three dimensions. The error \eqref{eq:normedintegrationerror} is defined with respect to the $\Hilbert^1$ norm. The integration error is in all cases computed with respect to the exact integral over the considered partitioning, and hence the geometric error corresponding to the octree partitioning is not represented in the results.

\begin{figure}
	\centering
	\begin{subfigure}[b]{0.49\textwidth}
		\centering
		\includegraphics[width=\textwidth]{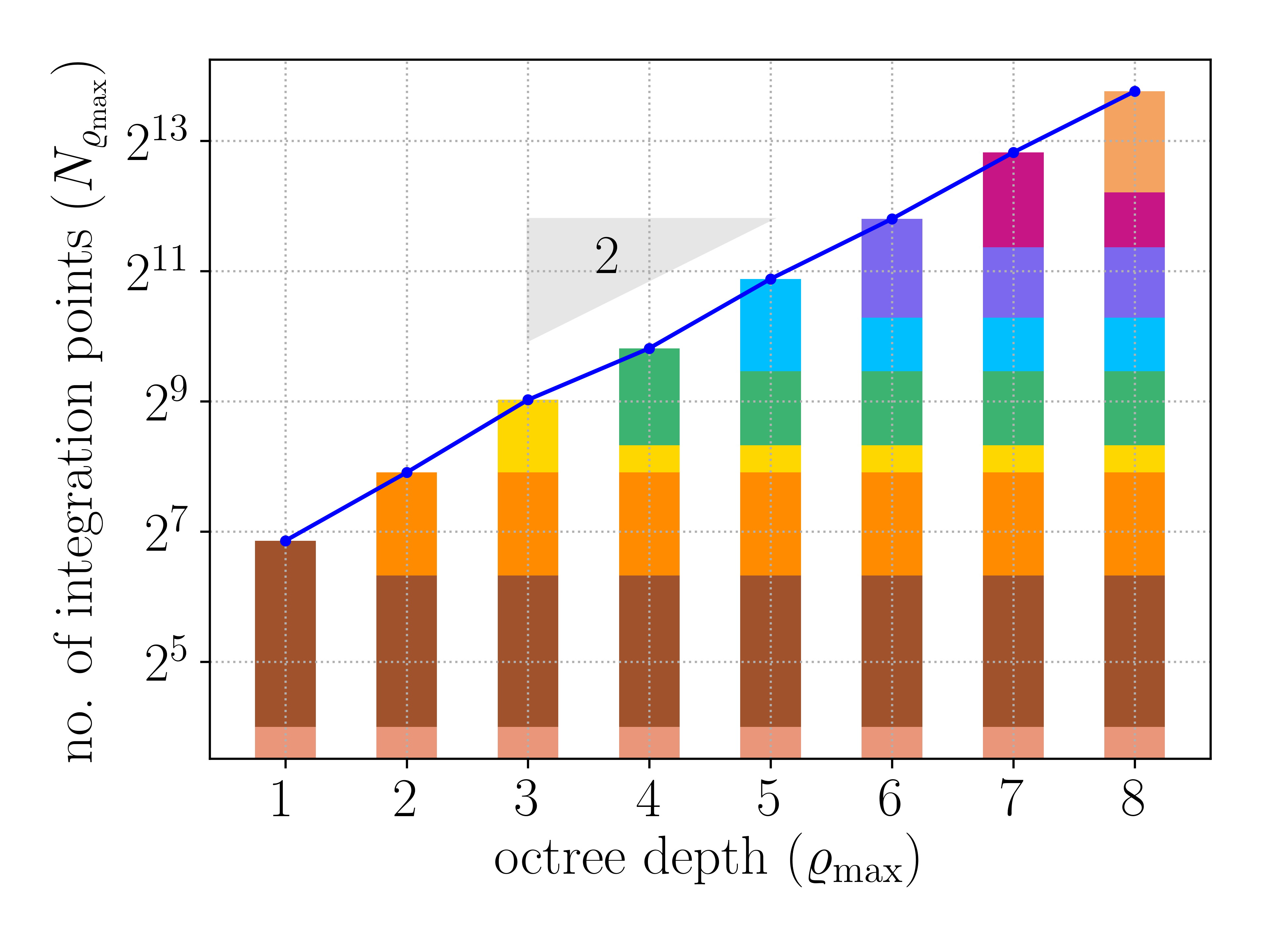}
		\caption{Uniform, $d=2$}
	\end{subfigure}%
	\begin{subfigure}[b]{0.49\textwidth}
		\centering
		\includegraphics[width=\textwidth]{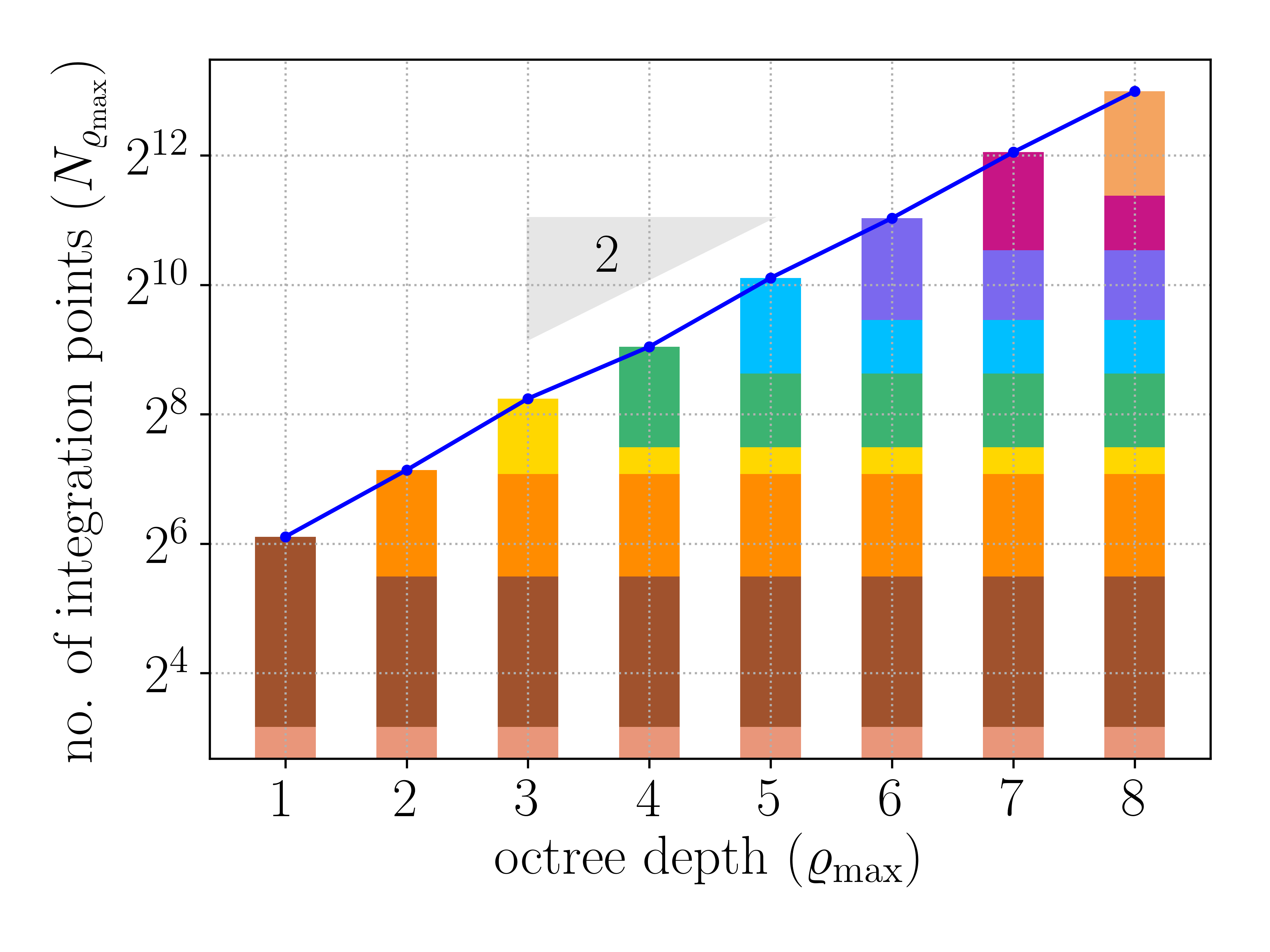}
		\caption{Gauss, $d=2$}
	\end{subfigure}\\[12pt]
	\begin{subfigure}[b]{0.49\textwidth}
		\centering
		\includegraphics[width=\textwidth]{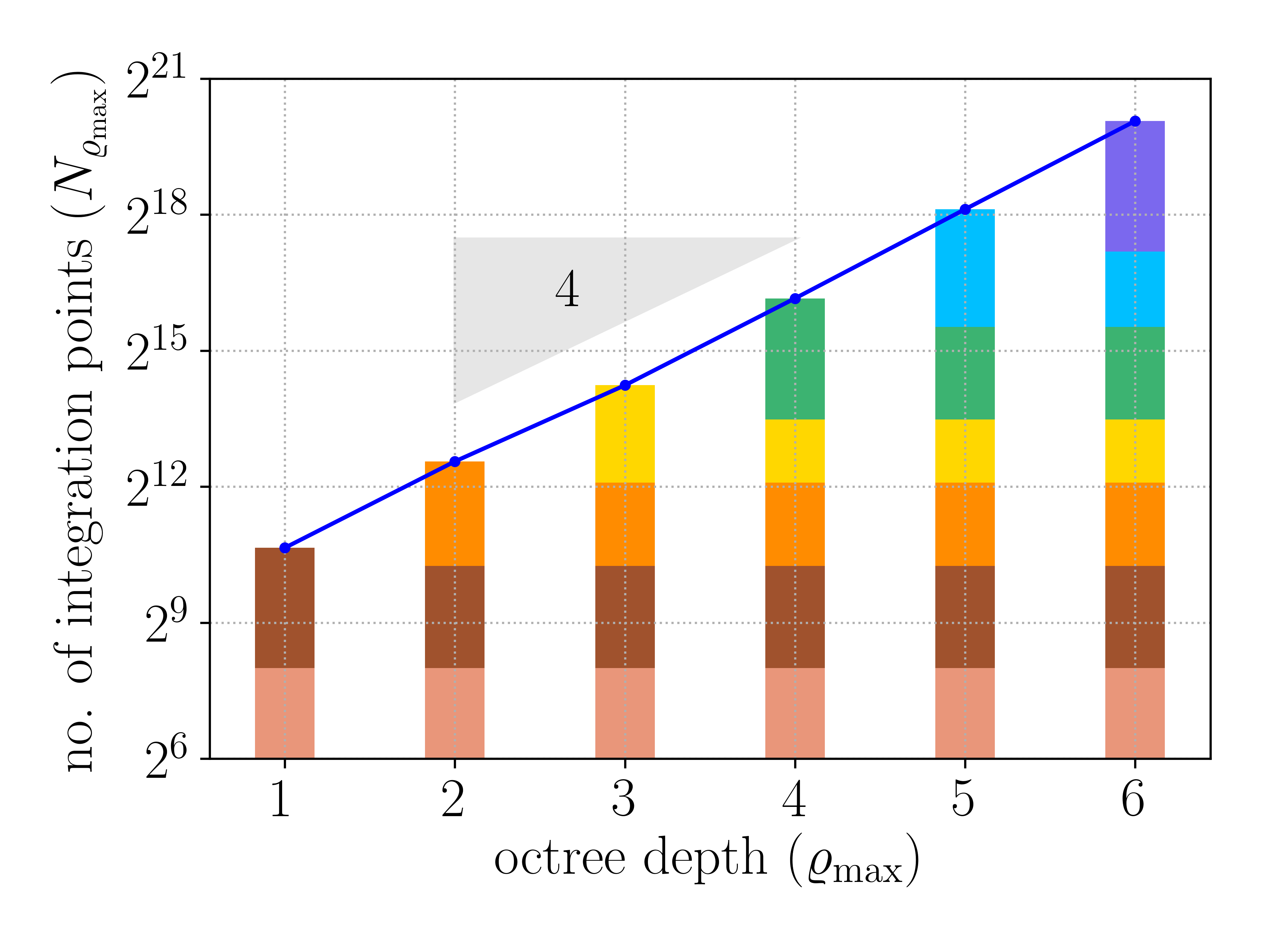}
		\caption{Uniform, $d=3$}
	\end{subfigure}%
	\begin{subfigure}[b]{0.49\textwidth}
		\centering
		\includegraphics[width=\textwidth]{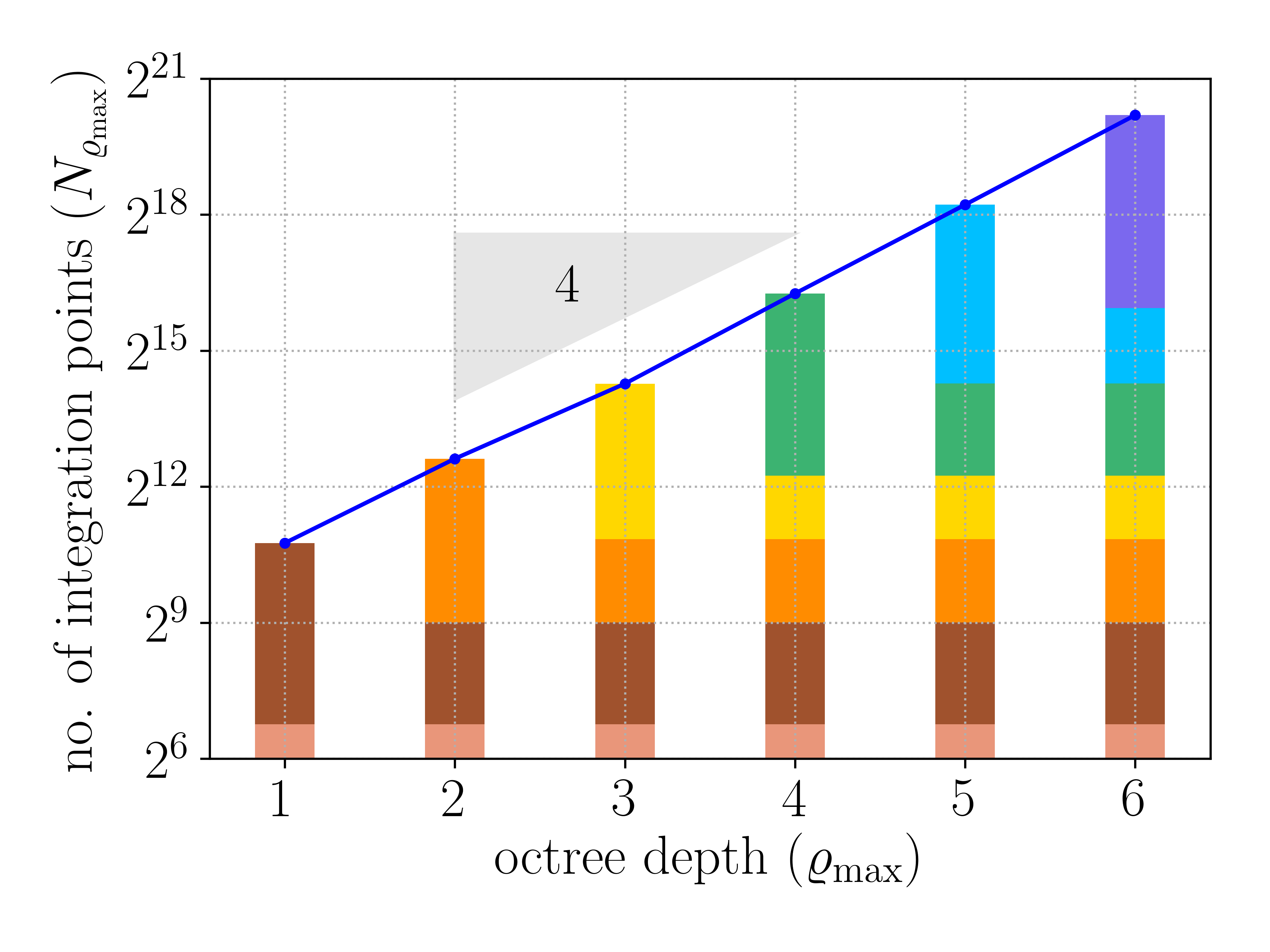}
		\caption{Gauss, $d=3$}
	\end{subfigure}	
	\caption{Dependence of the number of integration points on the octree depth for the case of uniform integration (with 4 point per dimension) on all sub-cells (top row) and Gauss quadrature of order 4 for all sub-cells (bottom row).}
	\label{fig:maxrefine}
\end{figure}

Figure~\ref{fig:equalorders} displays integration errors based on both equal-order Gauss quadrature and equal-order uniform quadrature on all sub-cells. It is observed that for all $\maxlevel$ the observed errors for a particular integration order are similar, which is explained by the fact that the integration error is dominated by contributions from the sub-cells that are already present at the lowest $\maxlevel$ and that the errors are computed with respect to the considered partitioning. In terms of the number of integration points, the octree depth does, however, have a substantial effect. The observed increase in number of integration points is studied in detail in Figure~\ref{fig:maxrefine} for the case of a uniform scheme with $4^2$ points per sub-cell in two dimensions and $4^3$ points per sub-cell in three dimensions and for a Gauss scheme of order $4$ on each sub-cell. This figure conveys that the total number of integration points scales in agreement with the relation \eqref{eq:numberofpoints}. In particular the doubling of the number of points in two dimensions and the quadrupling in three dimensions while increasing the octree by a single level is clearly observed.

From Figure~\ref{fig:equalorders} it is observed that for Gauss orders that are significantly below the order required for exact integration there is not a substantial difference with uniform integration. Once the Gauss order reaches that needed for exact integration, the error observed for Gauss integration is substantially smaller than that using uniform integration. When considering the same integration order on all sub-cells, this advantage is, however, only attained at the expense of introducing a large number of integration points, in particular in cases of high octree depths.
\subsection{Adaptive integration} \label{sec:adaptiveintegration}
\begin{figure}
	\centering
	\begin{subfigure}{\textwidth}
		\includegraphics[width=\textwidth]{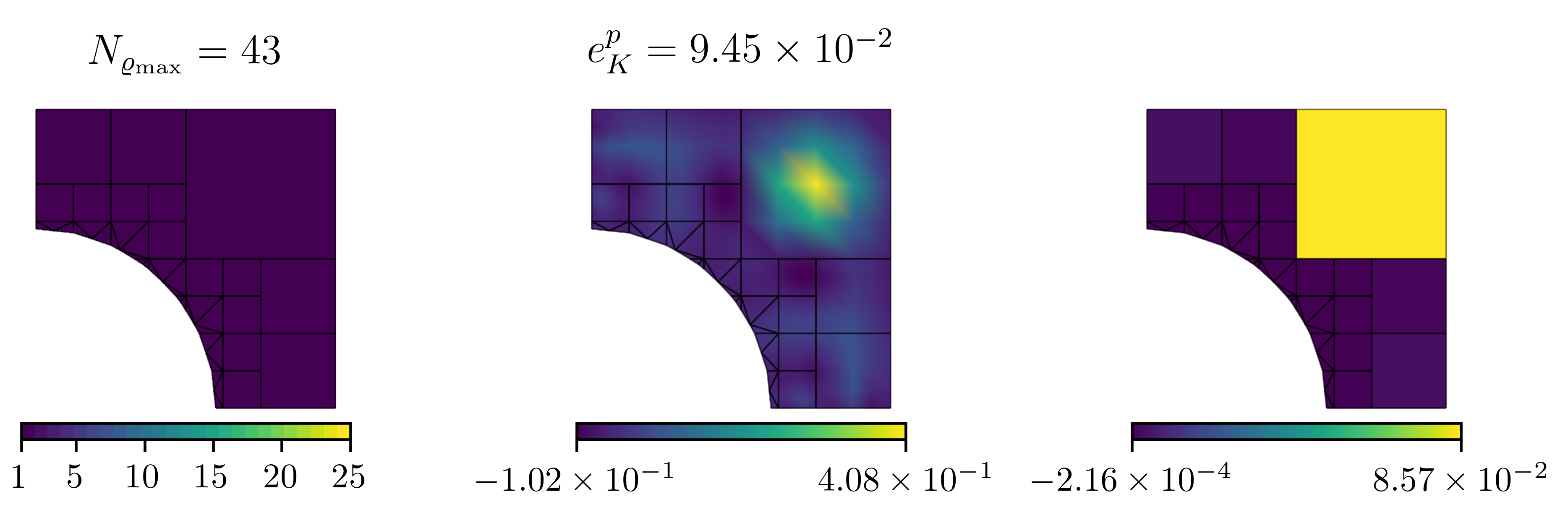}
		\caption{Step 1}
	\end{subfigure}\\[12pt]
	\begin{subfigure}{\textwidth}
		\centering
		\includegraphics[width=\textwidth]{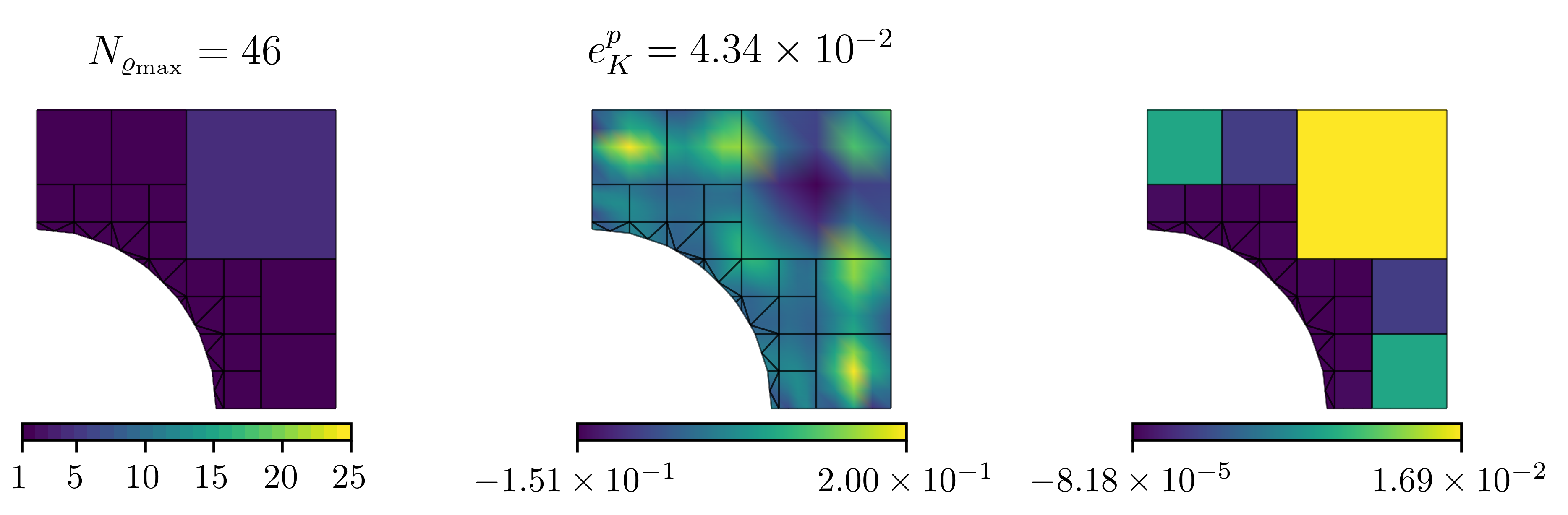}
		\caption{Step 2}
	\end{subfigure}\\[12pt]
	\begin{subfigure}{\textwidth}
		\centering
		\includegraphics[width=\textwidth]{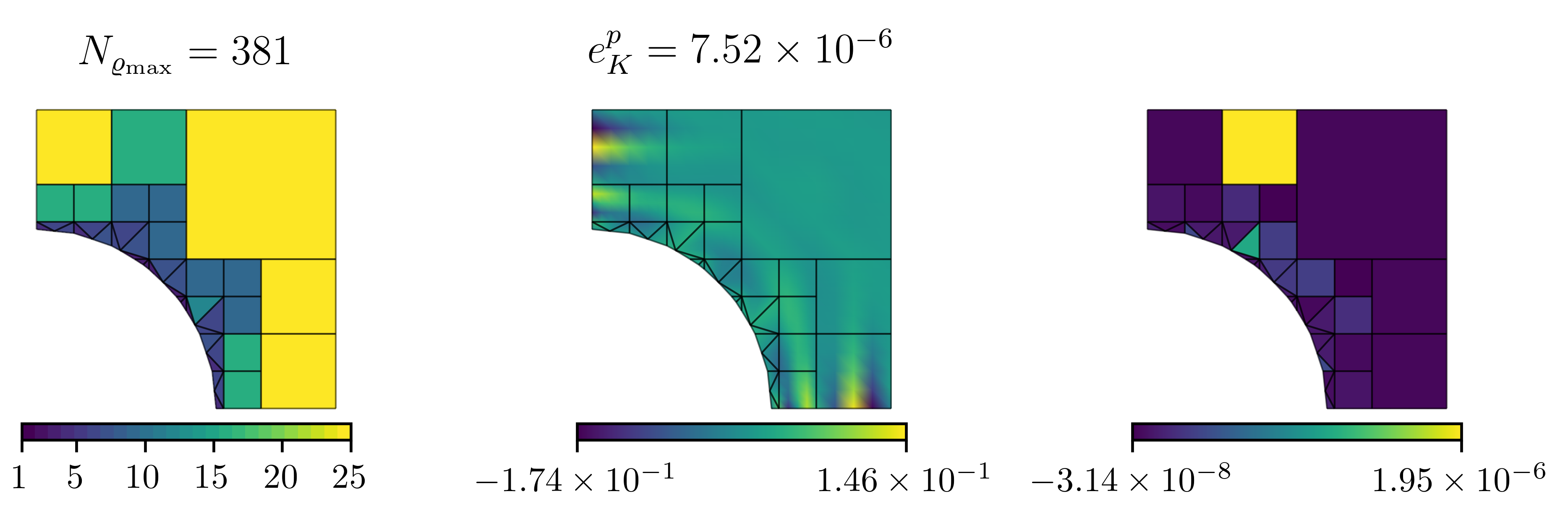}
		\caption{Step 115}
		\label{fig:step17}
	\end{subfigure}
	\caption{Selection of steps in the quadrature optimization procedure for a circular exclusion. (left column)  The number of integration points per sub-cell and the total number of points on the cut-element. (middle column) The integration error and the function leading to that error. (right column) The sub-cell integration error indicators.}
	\label{fig:stagesofalgorithm_level}
\end{figure}

We now consider the adaptive-integration procedure for the circular exclusion setting considered in Section~\ref{sec:equalorder}. In Figure~\ref{fig:stagesofalgorithm_level} various steps in the quadrature-optimization procedure are shows, displaying for each step the number of integration points per sub-cell (left column), error and the function leading to that error (middle column), and the sub-cell error indicators \eqref{eq:normedintegrationerror} (right column). The presented results are generated using the per sub-cell marking strategy. The sequence of steps demonstrates that in each step the worst possible function in terms of integration is determined. In the first step, this corresponds to a function that is large in magnitude at the largest sub-cell. This function leads to an indicator that is highest in that sub-cell, and hence the order of integration on that cell is increased. As a result, in the second step a function that is large in magnitude on the sub-cells surrounding the largest sub-cell is found to be the worst possible in terms of integration error. The larger volume of the $\level=1$ integration cell makes, however, that the largest indicator is still found for that cell. After another increase in integration order on that sub-cell, in the subsequent steps the largest indicators are found on the $\level=2$ sub-cells, which are therefore gradually increased in order.

\begin{figure}
	\centering
	\begin{subfigure}[b]{0.5\textwidth}
		\includegraphics[width=\textwidth]{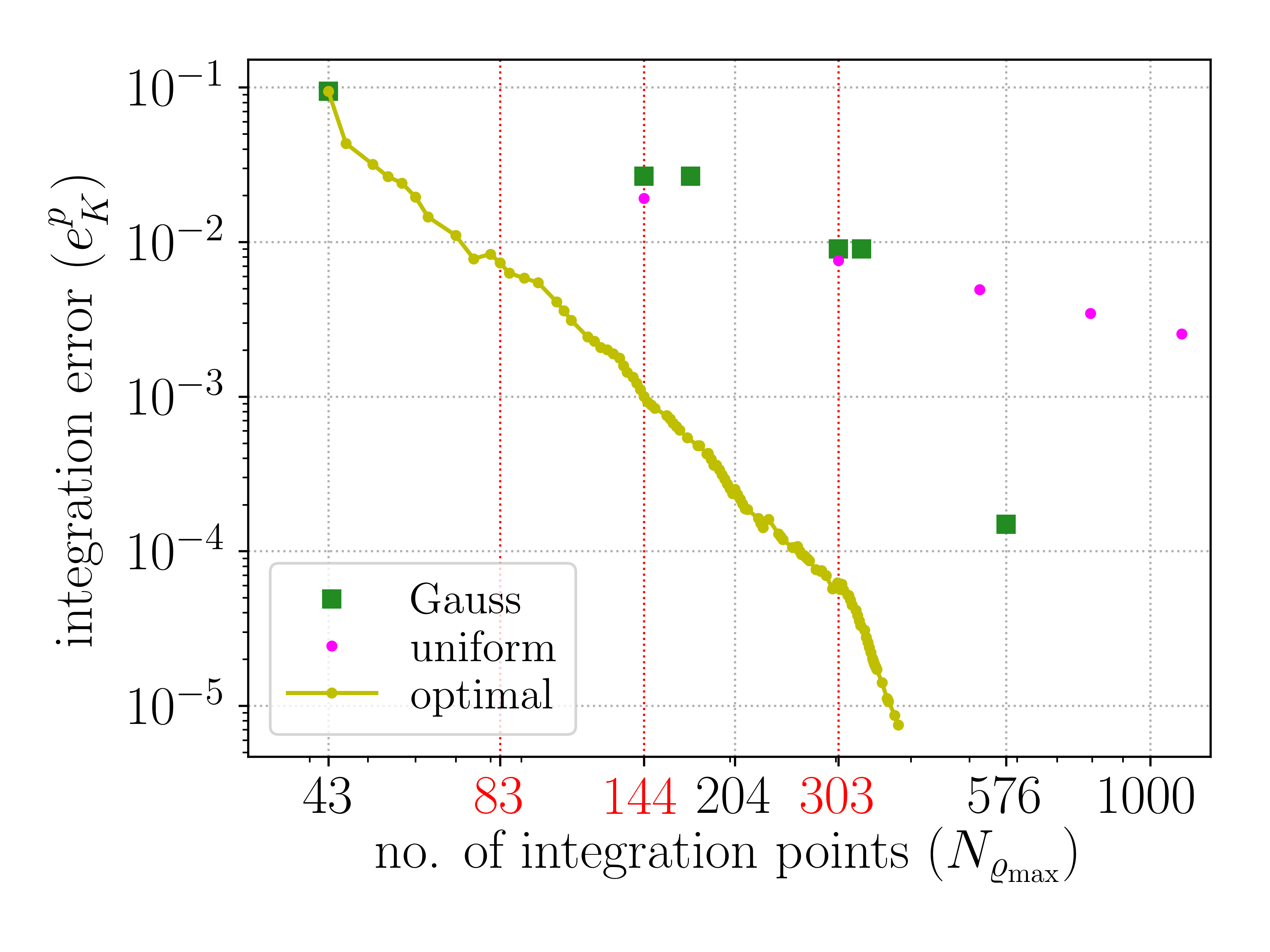}
		\caption{$d=2$}
		\label{fig:optsubcell2d}
	\end{subfigure}%
	\begin{subfigure}[b]{0.5\textwidth}
		\includegraphics[width=\textwidth]{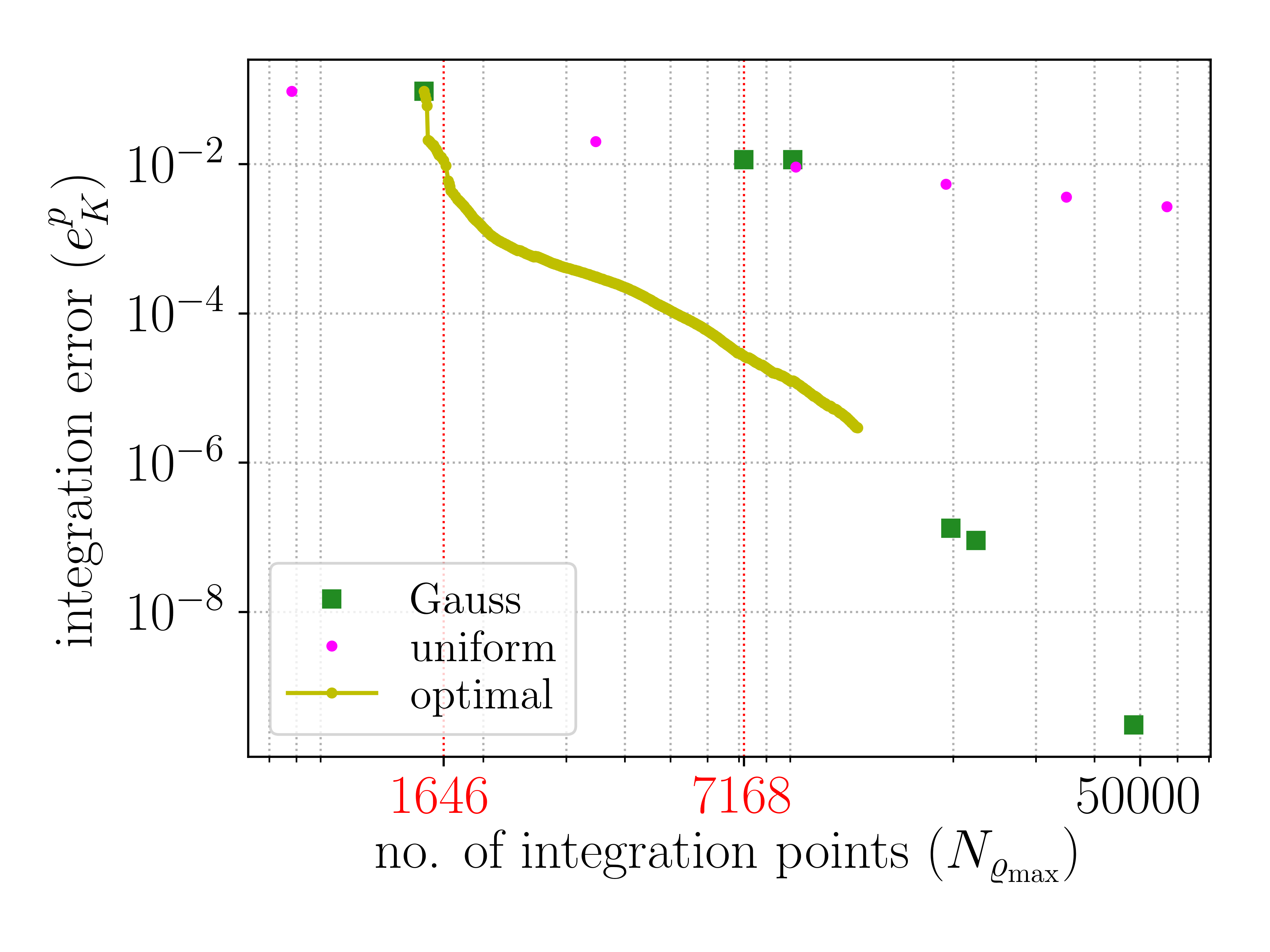}	
		\caption{$d=3$}
		\label{fig:optsubcell3d}
	\end{subfigure}
	\caption{Evolution of the integration error and number of integration points using the per sub-cell marking strategy. The results for equal-order quadrature are shown for reference.}
	\label{fig:percell}
\end{figure}

\begin{figure}
	\centering
	\begin{subfigure}[b]{0.42\textwidth}
		\centering
		\includegraphics[width=0.8\textwidth]{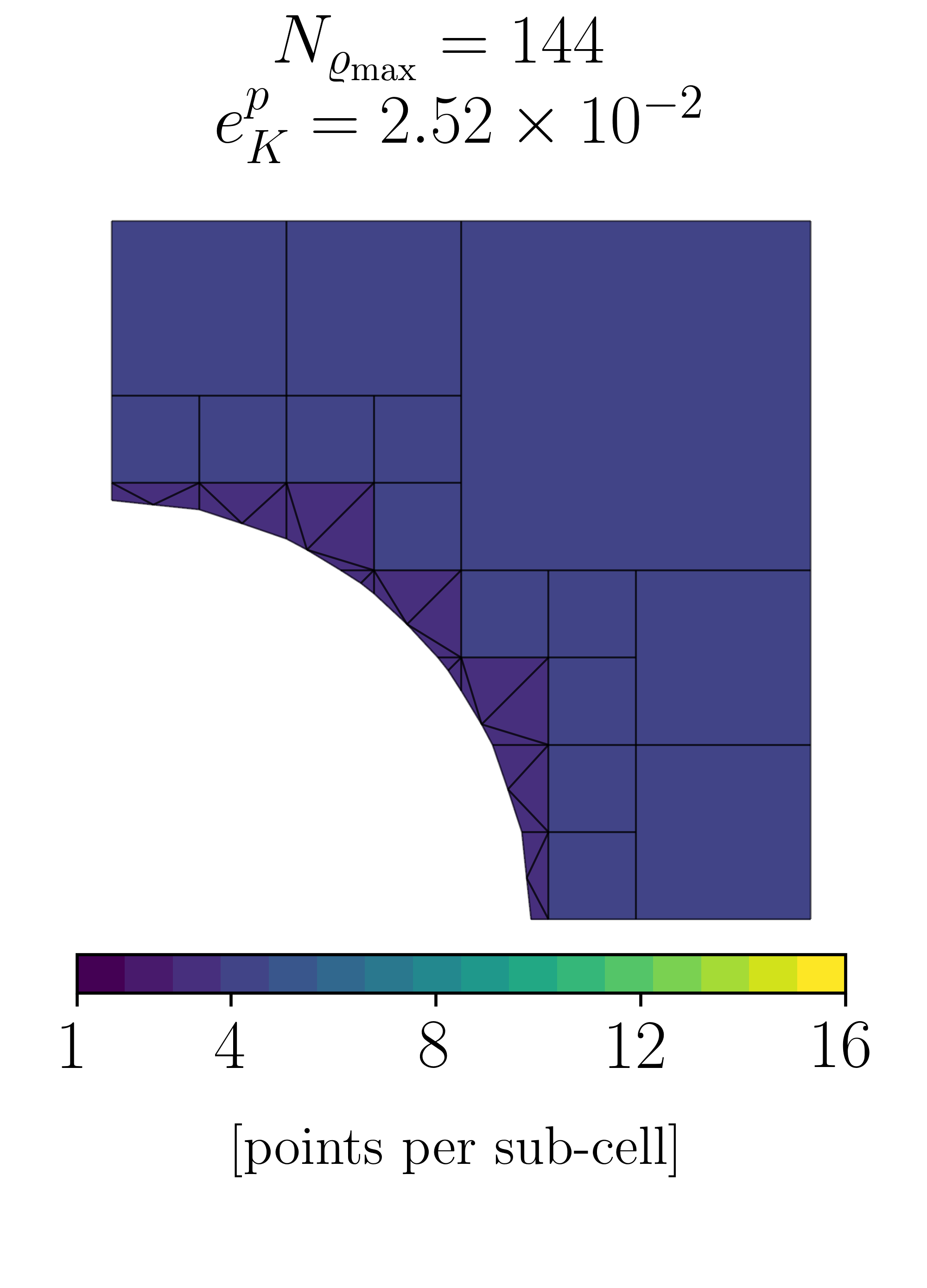}
		\caption{Equal-order Gauss}
		\label{fig:equalorder_gauss}
	\end{subfigure}%
	\begin{subfigure}[b]{0.42\textwidth}
		\centering
		\includegraphics[width=0.8\textwidth]{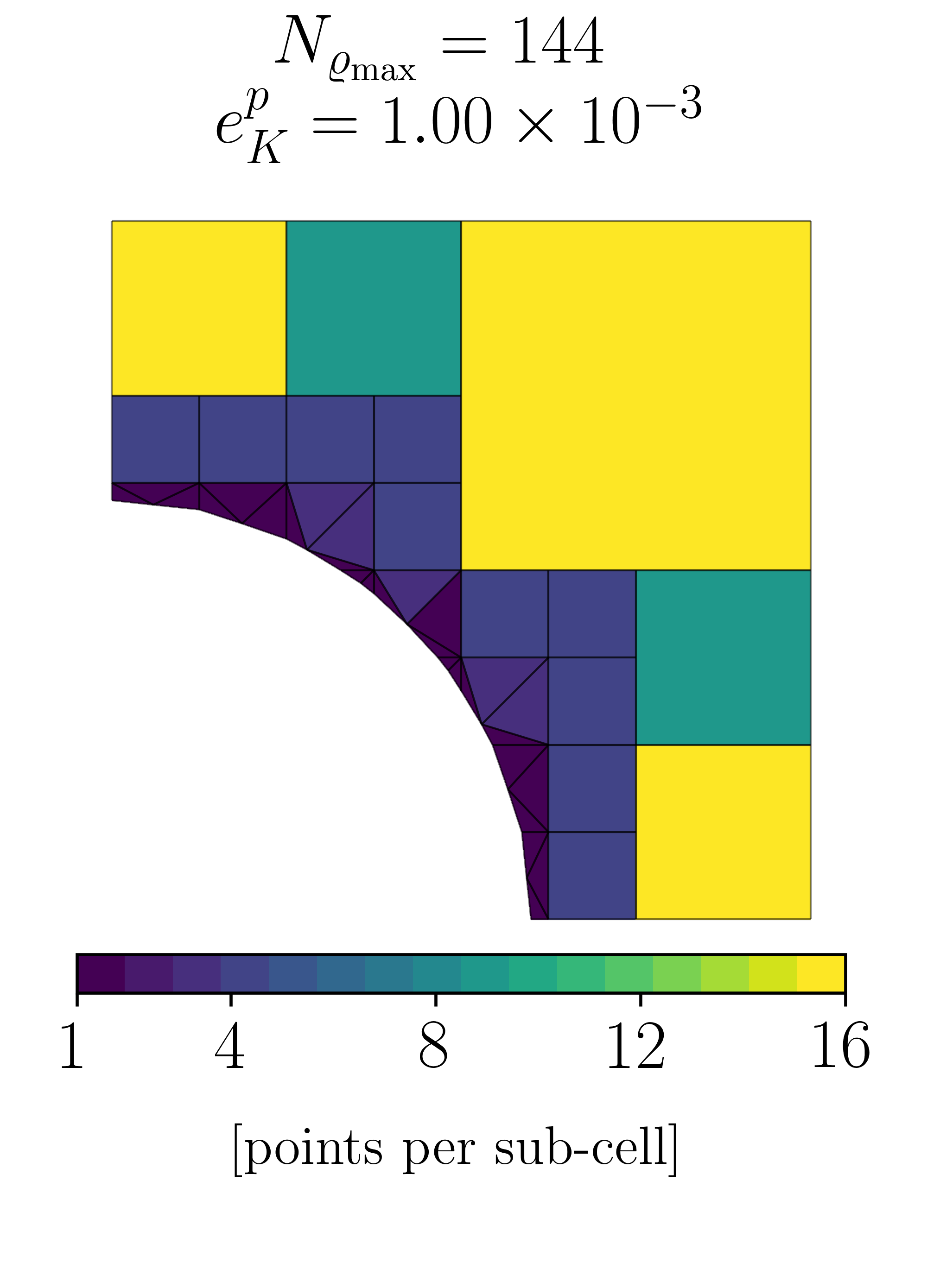}
		\caption{Optimal integration}
		\label{fig:optimal_gauss}
	\end{subfigure}
	\caption{Distribution of integration points over the cut-element for $144$ points in two dimensions, which corresponds to the number of points attained using a second order Gauss scheme on all sub-cells. Note that the error is reduced by a factor of 25 in two dimensions by using the per sub-cell marking strategy.}
	\label{fig:optimaldist2d}
\end{figure}

Figure~\ref{fig:percell} displays the integration error versus the total number of integration points as evolving during the optimization procedure. The displayed results pertain to $\maxlevel=3$ in both two and three dimensions, displaying the equal-order results discussed in Section~\ref{sec:equalorder} for reference. As can be seen, the error per integration point is substantially lower using the adaptive integration procedure. For example, in two dimensions, the error corresponding to the equal second order Gauss scheme is equal to $2.52 \times 10^{-2}$, while that particular integration scheme involves $144$ points. For the same number of points, the error corresponding to the optimized quadrature is equal to $1.00 \times 10^{-3}$, \emph{i.e.}, a factor 25 reduction in error. Figure~\ref{fig:optimaldist2d} displays the distribution of the integration points over the sub-cells for the equal-order Gauss scheme and the optimized case, which clearly demonstrates that the significant reduction in error is achieved by assigning more integration points to the larger sub-cells before introducing additional points in the smaller sub-cells.

\begin{figure}
	\centering
	\begin{subfigure}[b]{0.42\textwidth}
		\centering
		\includegraphics[width=0.8\textwidth]{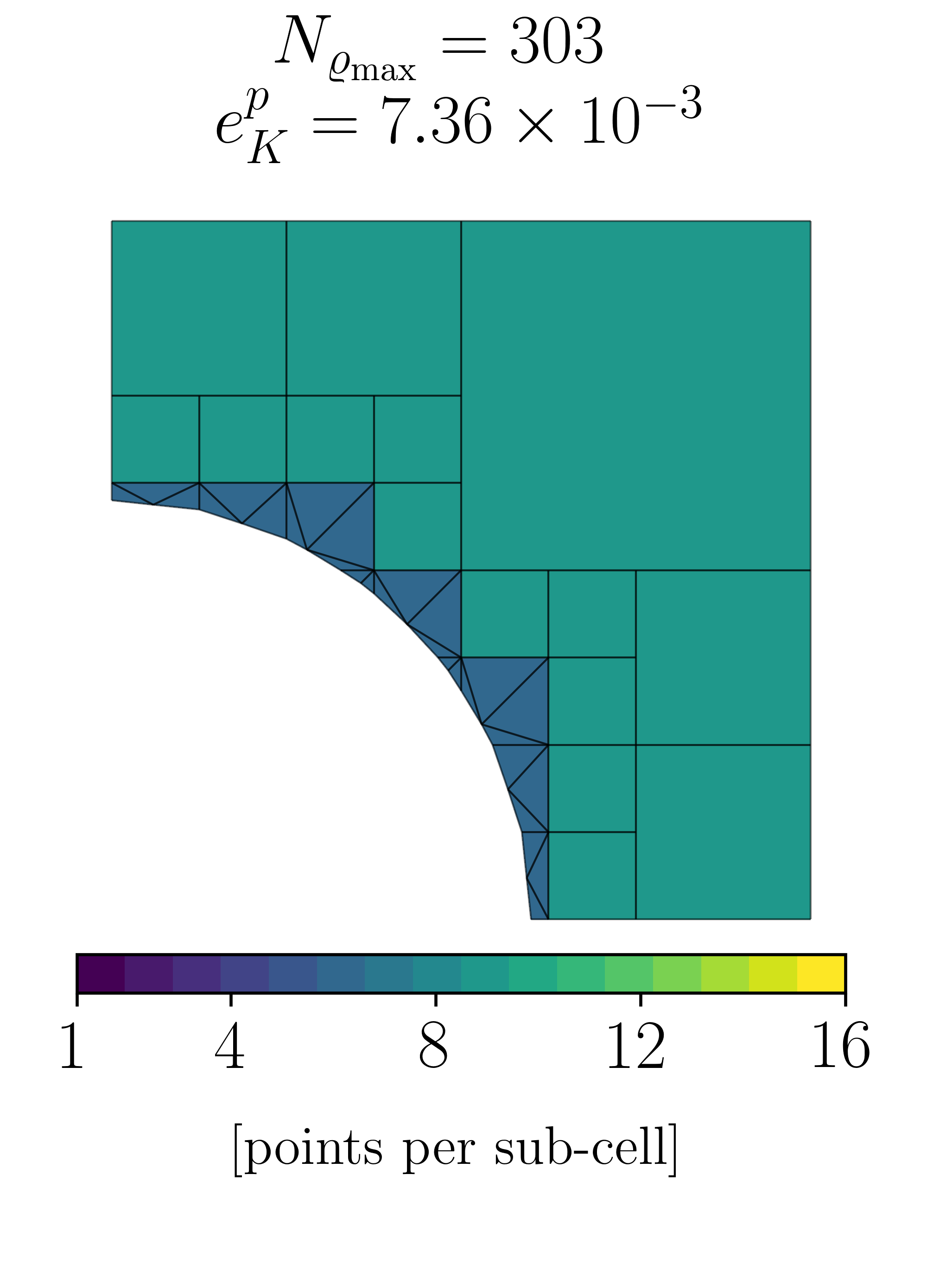}
		\caption{Equal-order Gauss}
		\label{fig:err_equalorder_gauss}
	\end{subfigure}%
	\begin{subfigure}[b]{0.42\textwidth}
		\centering
		\includegraphics[width=0.8\textwidth]{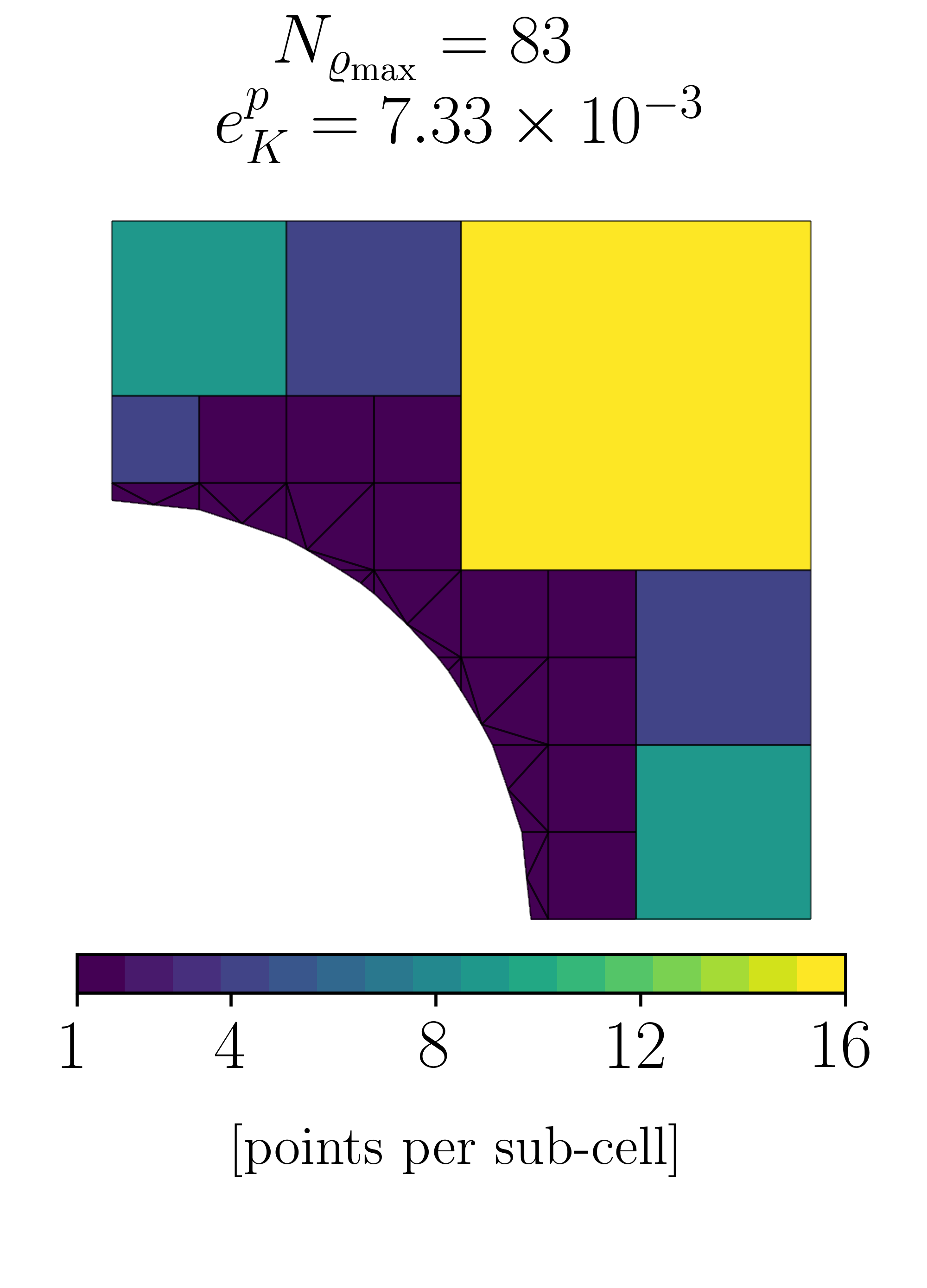}
		\caption{Optimal integration}
		\label{fig:err_optimal_gauss}
	\end{subfigure}
	\caption{Distribution of integration points over the cut-element for an error of approximately $7.3 \times 10^{-3}$. Note that the number of integration points is reduced by a factor of 4 in two dimensions by using the per sub-cell marking strategy.}
	\label{fig:optimalerr2d}
\end{figure}

In general, for a fixed integration error, the number of integration points is substantially lower using the adaptive integration procedure than when using the equal-order schemes. For example, in two dimensions, the error corresponding to the equal fourth order Gauss scheme is equal to $7.35 \times 10^{-3}$, which involves $303$ points. For a similar error, the number of points for the optimized quadrature is equal to $83$, \emph{i.e.}, a reduction by a factor of approximately 4. Figure~\ref{fig:optimalerr2d} displays the distribution of integration points over the sub-cells for the equal-order Gauss scheme and the optimized quadrature, which clearly shows the reduction in number of points for a fixed error. This figure also conveys that the observed reduction factor of 4 is significantly influenced by the initial integration order in the adaptive procedure, which in this case has been set to 1 (for this case the minimum number of points is equal to 43).

Figure~\ref{fig:optsubcell3d} displays the error versus the number of integration points for the three-dimensional spherical exclusion with $\maxlevel=3$. The integration error corresponding to the equal second order Gauss schemes is equal to approximately $1.14 \times 10^{-2}$ and pertains to $7168$ integration points. For the same number of points, the adaptive procedure reduces the error to $2.67 \times 10^{-5}$, which constitutes a reduction factor of approximately $450$. This reduction factor is significantly higher than that observed in the two-dimensional case, despite the similarity in asymptotic scaling rate between the integration error and the number of integration points in the two- and three-dimensional cases. The larger reduction factor in three dimensions is attributed to the significant reduction in error that is achieved during the first integration point optimization steps, in which the integration order on the $\level=1$ sub-cells is increased. The difference between the two-dimensional and three-dimensional setting in this regard is the fact that in three dimensions there is a significant number of such sub-cells, whereas in two dimensions there is only one sub-cell of level $\level=1$. For an integration error $1.139 \times 10^{-2}$, the adaptive procedure requires only $1646$ points, which, as in the two-dimensional case, is a reduction by a factor of approximately $4$. As for the two-dimensional case, this reduction factor is significantly influenced by the integration order with which the adaptive procedure is initiated.

\subsection{Influence of algorithmic settings and cut-cell configurations}
In the previous section the performance of the adaptive integration routine has been studied using representative settings for the cut-cell geometry and algorithmic parameters. In this section we will study the influence of the most prominent parameters on the performance of the adaptive integration technique.

\subsubsection{The marking strategy} \label{sec:cellvslevel}
From the refinement patterns that emerge from the per sub-cell refinement strategies -- such as the ones discussed in the previous section -- it is observed that, as a general trend, integration orders are increased on a per-level basis, \emph{i.e.}, the number of integration points is increased from level $\level$ to $\maxlevel + 1$. This is explained by the fact that the indicators scale with the volume of the sub-cells. Based on this observation it is anticipated that the level-based marking strategy discussed in Section~\ref{sec:refinementmarker} can be very efficient, in the sense that it yields a similar refinement pattern as the per-cell marking, but that it needs fewer iterations by virtue of marking a larger number of sub-cells per step.

\begin{figure}
	\centering
	\begin{subfigure}[b]{0.5\textwidth}
		\includegraphics[width=\textwidth]{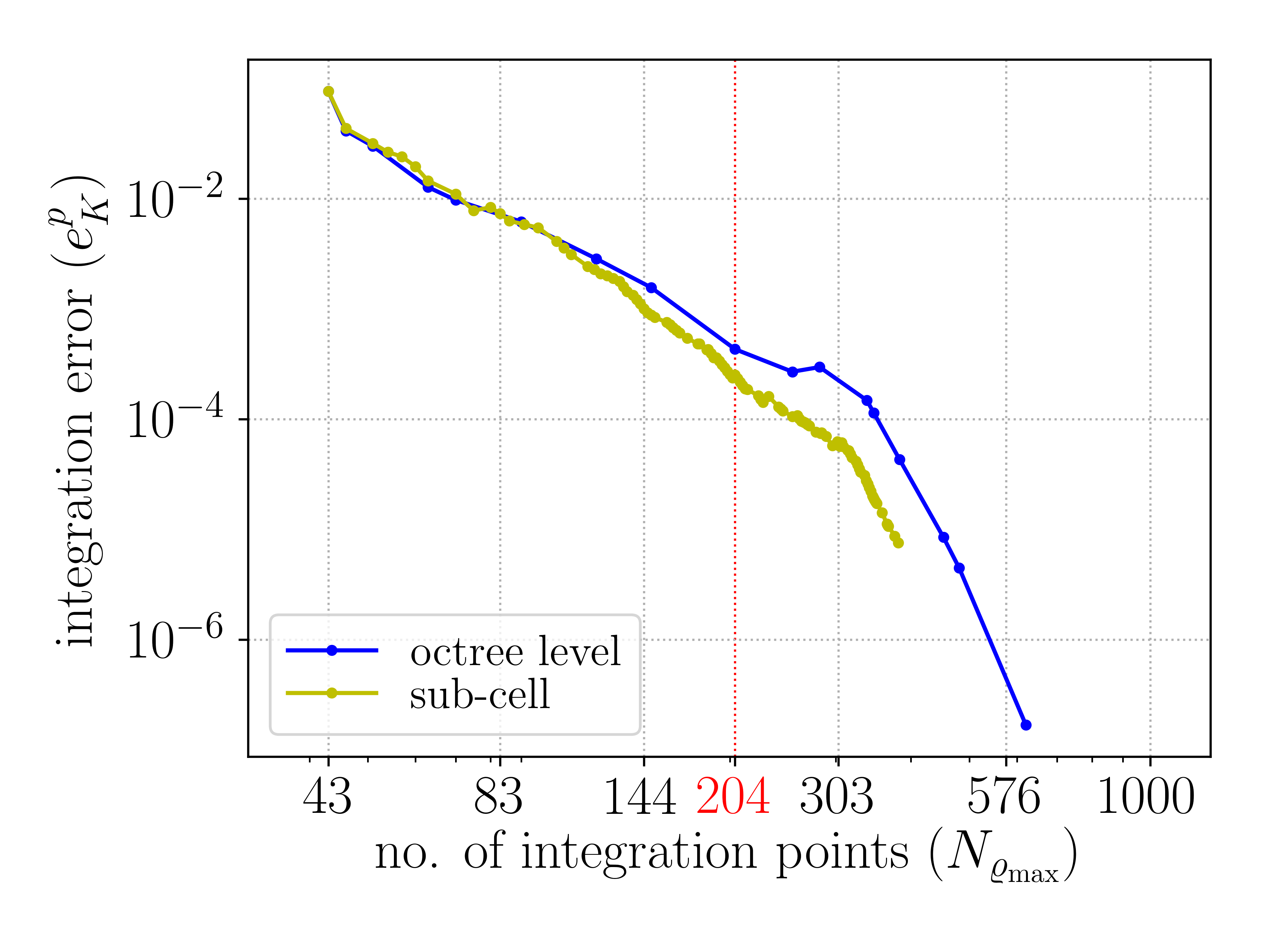}
		\caption{$d=2$}
		\label{fig:cellvslevel2d}
	\end{subfigure}%
	\begin{subfigure}[b]{0.5\textwidth}
		\includegraphics[width=\textwidth]{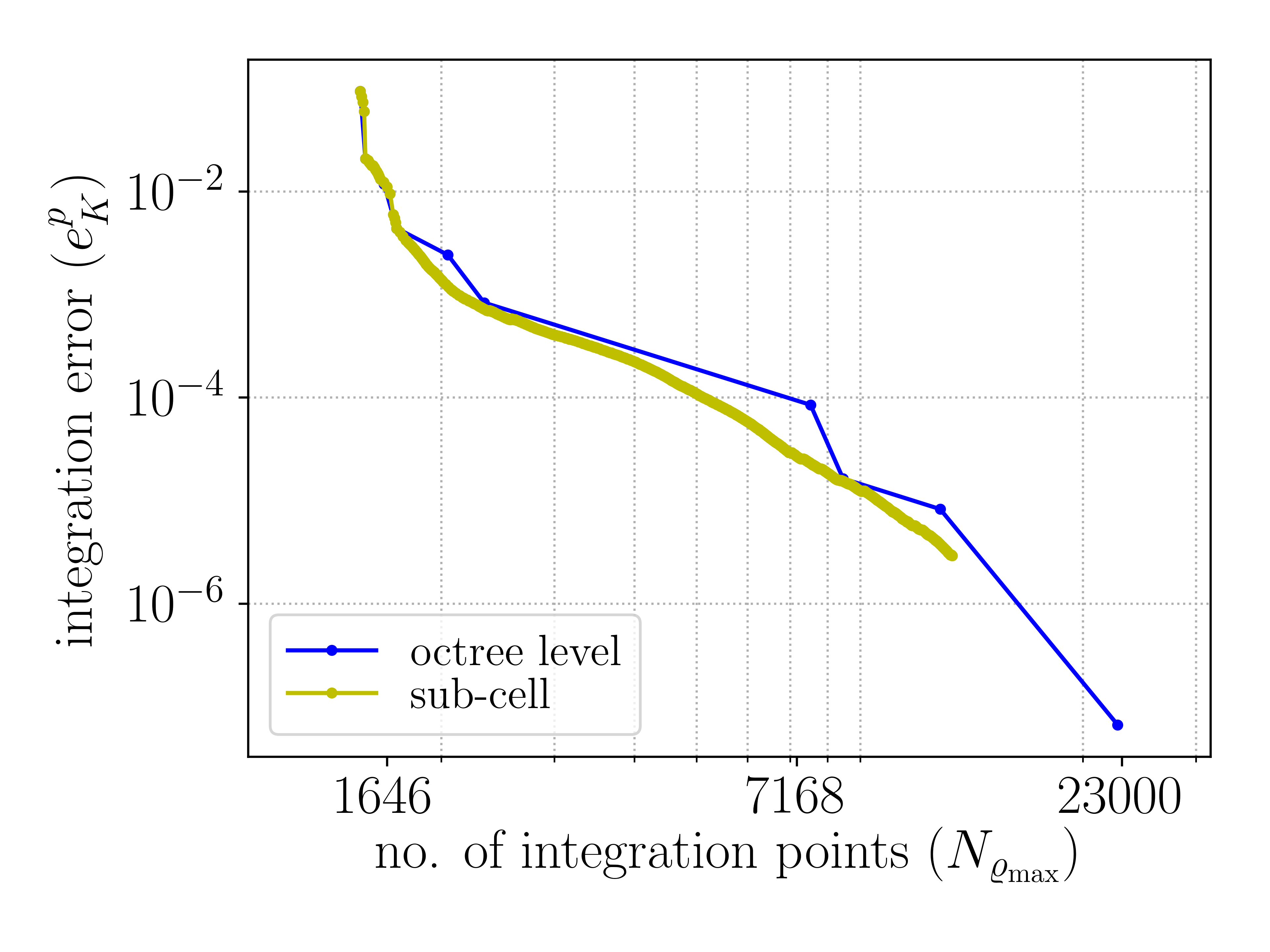}
		\caption{$d=3$}
		\label{fig:cellsvslevel3d}
	\end{subfigure}
	\caption{Comparison of the evolution of the integration error and number of integration points using the per-sub-cell and per-octree-level marking strategy.}
	\label{fig:cellvslevel}
\end{figure}

\begin{figure}
	\centering
	\begin{subfigure}[b]{0.42\textwidth}
	   \centering
	   \includegraphics[width=0.8\textwidth]{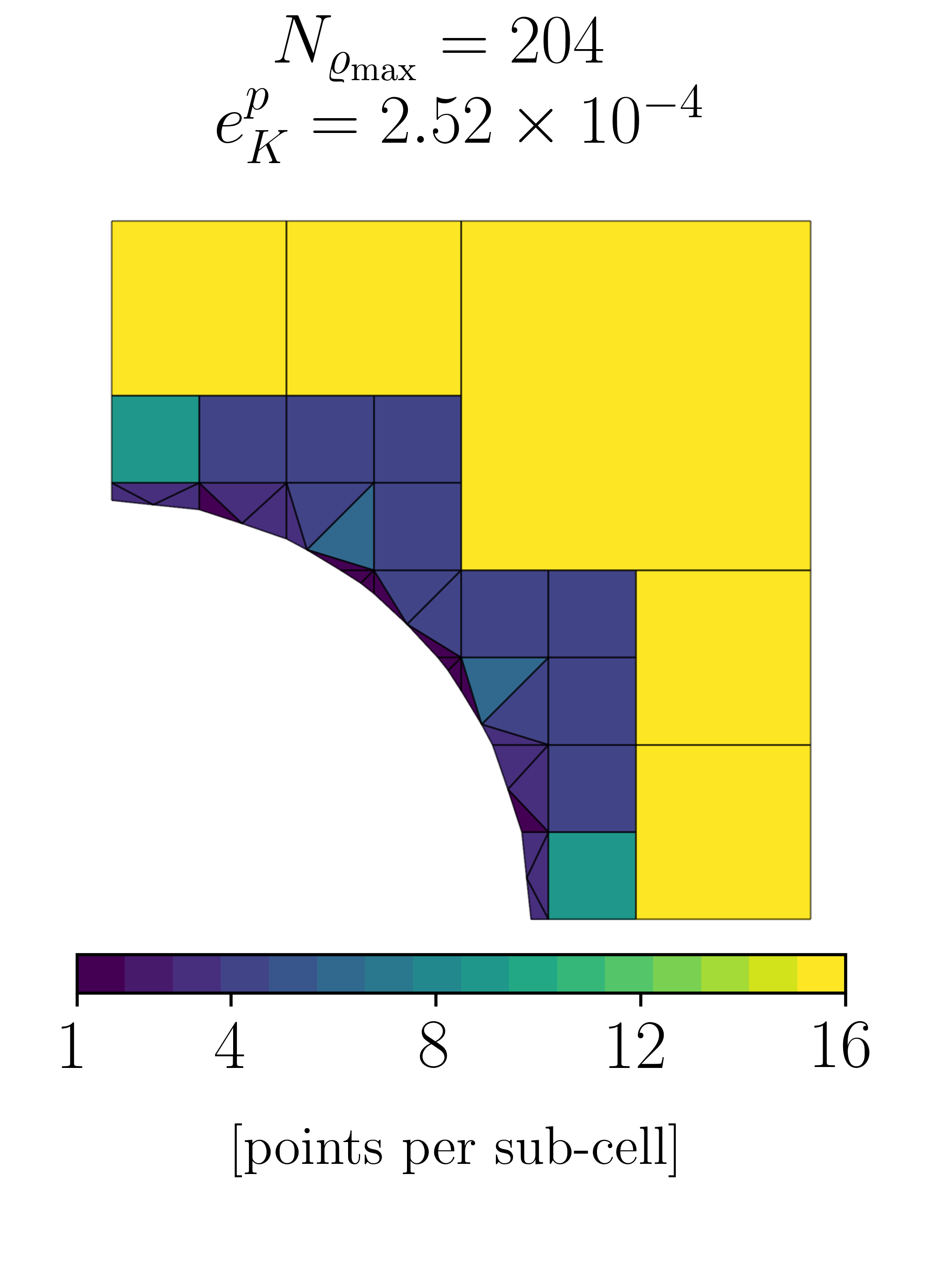}
	   \caption{Sub-cell marking}
	   \label{fig:cell2d}
	\end{subfigure}%
	\begin{subfigure}[b]{0.42\textwidth}
	   \centering
	   \includegraphics[width=0.8\textwidth]{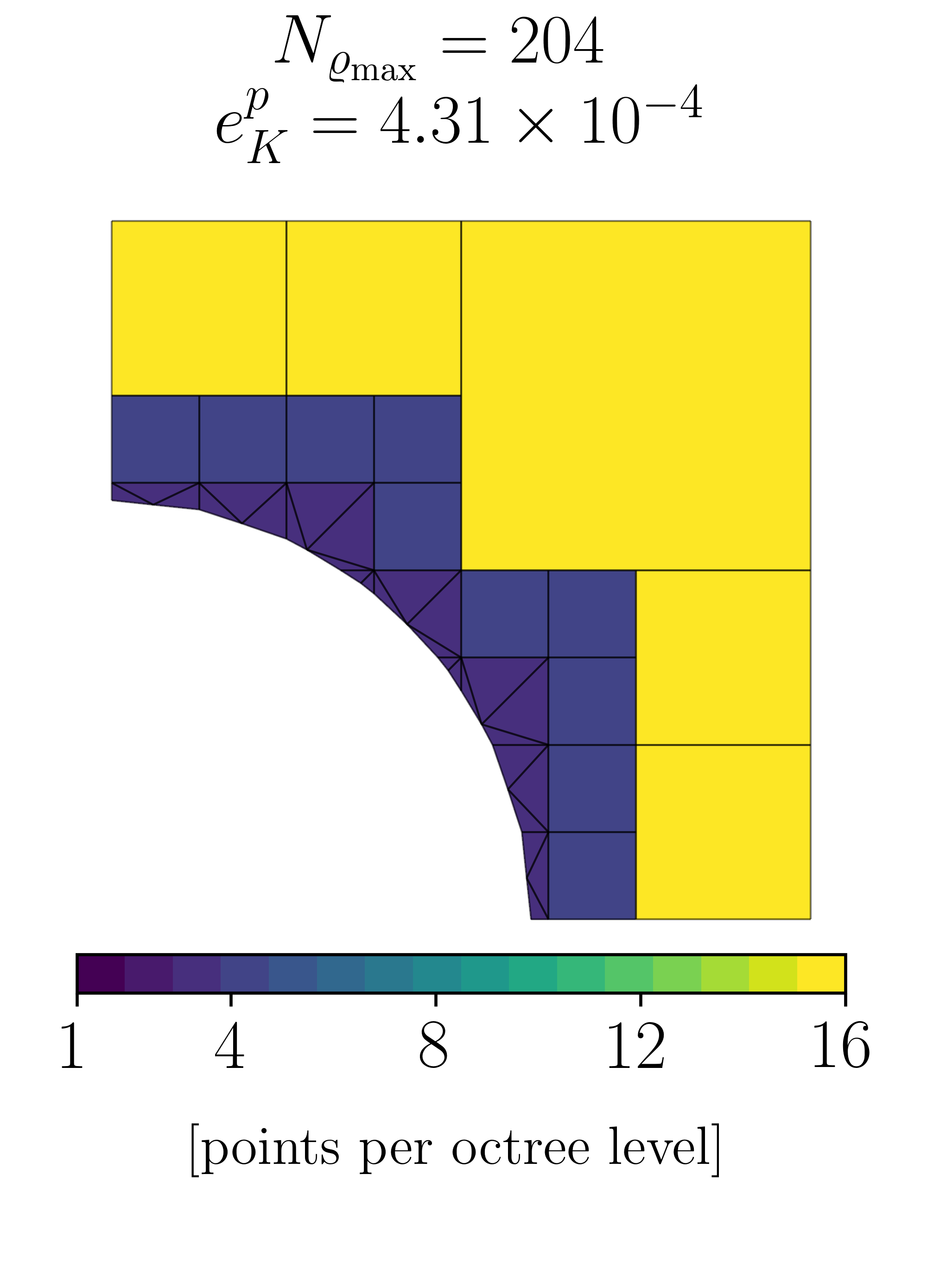}
	   \caption{Level marking}
	   \label{fig:level2d}
	\end{subfigure}
	\caption{Distribution of integration points over a cut-element with $204$ integration points in two dimensions using both sub-cell and level marking.}
	\label{fig:optimaldistcellvslevel}
\end{figure}

In Figure~\ref{fig:cellvslevel} the per-level and per-cell marking strategies are compared for the test case introduced above. In both two and three dimensions it is observed that the per-level marking strategy closely follows the per-cell marking. In particular in the initial steps of the optimization a very close agreement is observed between the marking strategies. In Figure~\ref{fig:optimaldistcellvslevel} we inspect the distribution of points corresponding to the two marking strategies with $204$ points (in the two-dimensional setting). This figure conveys that, indeed, the distribution of integration points between the two markings is very similar.

Although the observed point distributions in Figure~\ref{fig:optimaldistcellvslevel} are similar, the number of iterations required to attain these distributions is significantly different. The per-cell marking strategy requires $54$ iterations, whereas the corresponding result using the per-level marking is achieved in $9$ iterations. Evidently, this reduction in number of iterations translates into a computational-effort advantage for the per-level marking strategy.

In Table~\ref{tab:cputime} we study the computational effort for both marking strategies corresponding to the optimization procedure to attain the results in Figure~\ref{fig:cellvslevel2d}. The reported CPU times are based on a Python implementation using the open source finite element toolbox Nutils \cite{nutils}, which is executed on a four core CPU with a $2.60$ GHz 7th generation Intel Core i5 processor. Although the reported CPU times are highly dependent on a myriad of aspects, it is noted that the computation time of the Gramian matrix (0.4--0.45\,s) is representative for the computational effort involved in the construction of the element contribution to the system matrix. By relating the reported CPU times to this element system matrix construction time, a more meaningful notion of the computational effort is obtained, although, of course, also this relative notion of performance is subject to implementation and architecture considerations.

\begin{table}
	\centering
	\caption{Analysis of the CPU time to attain the results displayed in Figures~\ref{fig:cellvslevel2d}.} \label{tab:cputime}
	\begin{tabular}{ccccc}
		\toprule
		Marking strategy & \multicolumn{2}{c}{sub-cell} & \multicolumn{2}{c}{octree level} \\
		& $n_{\rm call}$ & $t/{\rm call}$[s] & $n_{\rm call}$ & $t/{\rm call}$[s] \\
		\midrule
		Initialization & 1 & 0.5 & 1 & 0.55 \\
		\midrule
		Reference basis function integral ($\boldsymbol{\xi}$) & 1 & 0.1 & 1 & 0.1 \\
		Gramian matrix ($\mathbf{G}$) & 1 & 0.4 & 1 & 0.45 \\
		\midrule
		Iterations & 115 & 0.33 & 17 & 0.35 \\
		\midrule
		Approximate basis function integral ($\boldsymbol{\bar{\xi}}$) & 115 & 0.12 & 17 & 0.1 \\
		Maximum eigenvalue ($\mathbf{v}_{\rm max}$) & 115 & 0.004 & 17 & 0.005 \\
		Localized integration errors ($e_\wp^{p^k}$) & 115 & 0.15 & 17 & 0.15 \\
		Integration error indicator ($\mathcal{R}$) & 115 & 0.01 & 17 & 0.01 \\
		Integration error marker ($\mathcal{M}$) & 115 & 0.0005 & 17 & 0.005  \\
		Post-process & 115 & 0.05 & 17 & 0.05 \\
		\midrule
		\textbf{Total CPU time [s]} & \multicolumn{2}{c}{\textbf{39}} & \multicolumn{2}{c}{\textbf{7}} \\
		\bottomrule
	\end{tabular}
\end{table}

Table~\ref{tab:cputime} conveys that the optimization algorithm has a similar initialization time for both marking strategies. This initialization time pertains to the computation of the reference basis function integrals \eqref{eq:xivectors} and the Gramian matrix ($\mathbf{G}$). Moreover, it is observed that the CPU time per step is also very similar between the two marking strategies. The majority of the computational effort per step resides in the computation of the basis function integrals, $\bar{\boldsymbol{\xi}}$, and in the evaluation of the integration error localized to the sub-cells. In total, the per-level marking strategy is, however, computationally more efficient on account of the significant reduction in number of iterations.

The computation time for the determination of the optimized integration scheme depends primarily on the number of sub-cells and on the dimension of the considered polynomial space. In principle, the computation time of the integral evaluations in the algorithm -- \emph{i.e.}, the computation of the basis function integrals, the Gramian, and the localized errors -- scales linearly with the number of sub-cells. It should be noted that the number of sub-cells is highly dependent on the bisectioning level and on the number of dimensions; see equation \eqref{eq:mpositive}. For example, for the two-dimensional case reported in Figure~\ref{fig:cellvslevel2d} and Table~\ref{tab:cputime} the number of sub-cells is equal to 43, whereas the three-dimensional case in Figure~\ref{fig:cellsvslevel3d} has 341 sub-cells. Since the integral evaluations dominate the overall CPU time of the initialization phase and that of an optimization step, these operations become proportionally more expensive. It is important to note, however, that the number of iterations scales with the number of sub-cells in the case of the per sub-cell marking, whereas it scales with the number of octree levels in the case of the per-level marking. This implies that, with a growing number of sub-cells, the per-level marking strategy becomes more favorable from a computational effort point of view.

In terms of computational effort, the dimension of the polynomial space primarily has an effect on the computation time for the basis function integrals and that for the Gramian. In particular the Gramian computation becomes more expensive with an increase in basis size, as the number of terms to be integrated scales quadratically with this size. The linear system solving step involved in the computation of the worst possible function also increases with an increase in system size, but for all considered systems the computational effort involved in the employed direct solver remains negligible compared to the integral evaluations.

\subsubsection{Influence of the functional setting} \label{sec:comp_k}
From the functional-setting point of view, the distribution of integration points is influenced by both the approximation order of the polynomial space \eqref{eq:polyspace} and by the norm \eqref{eq:evaluatedpolynomial} in which the integrands are considered, encoded by the Gramian $\mathbf{G}$.

\begin{figure}
	\centering
	\begin{subfigure}[b]{0.5\textwidth}
		\includegraphics[width=\textwidth]{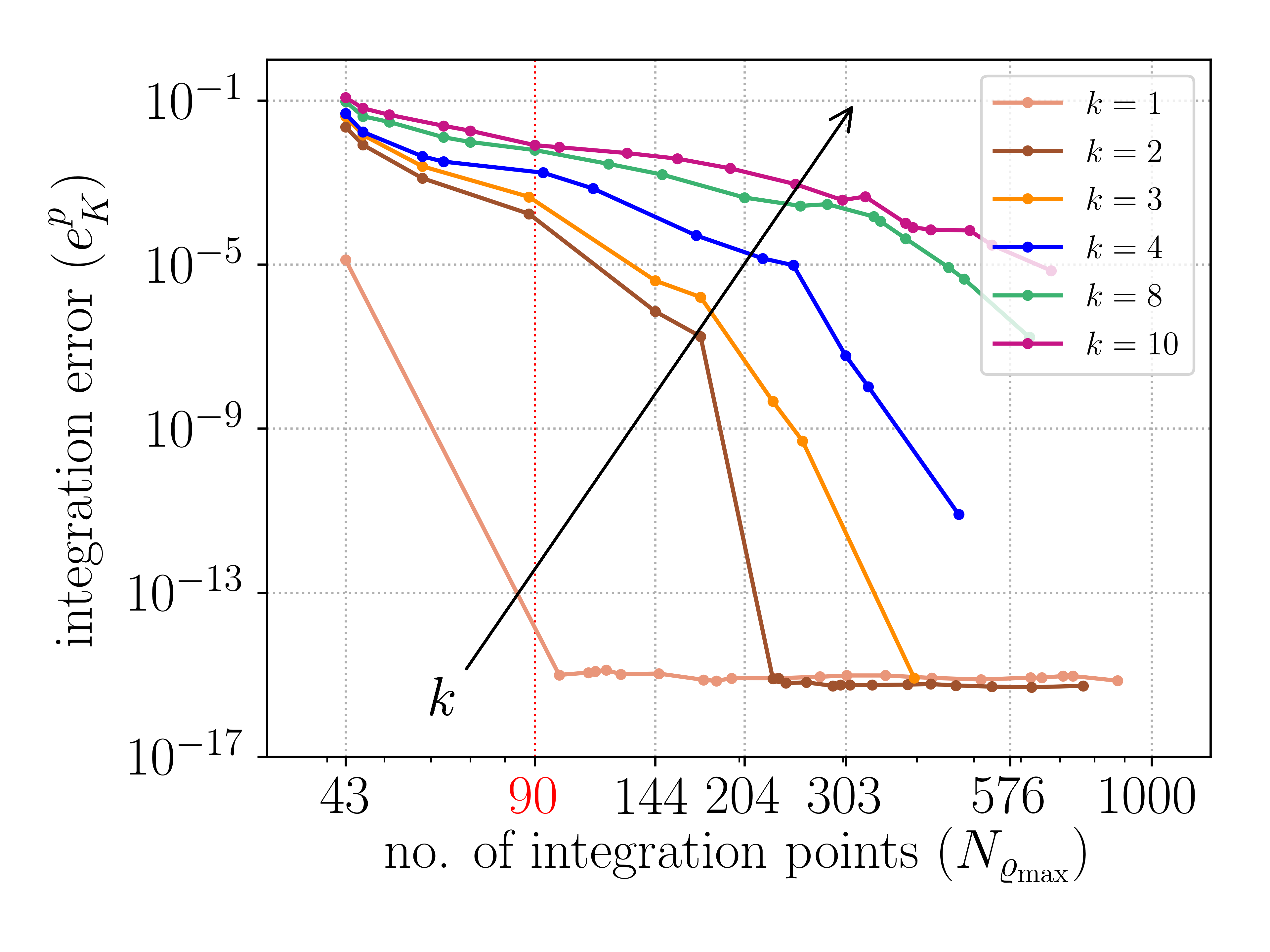}
		\caption{$d=2$}
		\label{fig:conv_degree_2d}
	\end{subfigure}%
	\begin{subfigure}[b]{0.5\textwidth}
		\includegraphics[width=\textwidth]{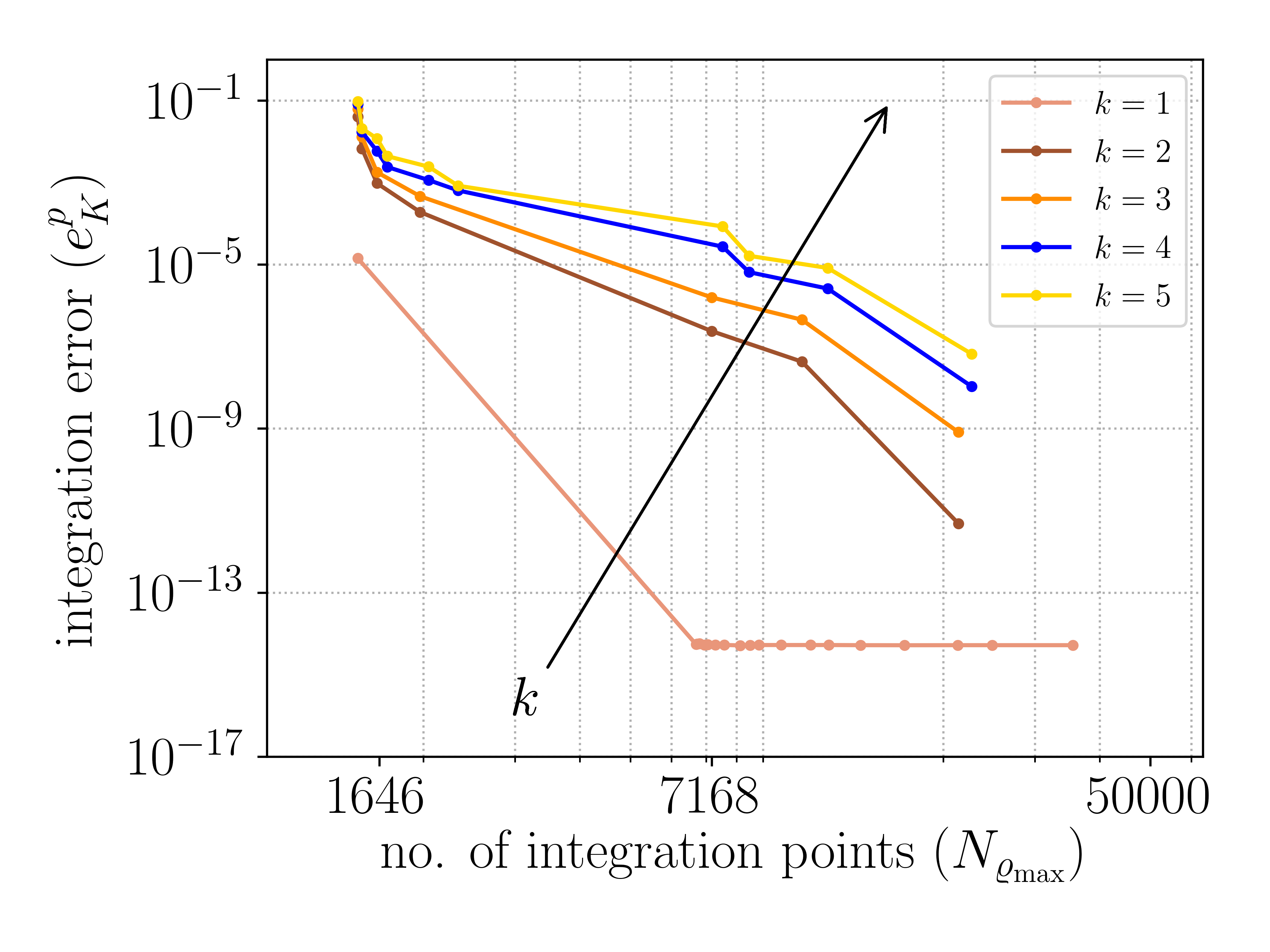}
		\caption{$d=3$}
		\label{fig:conv_degree_3d}
	\end{subfigure}
	\caption{Comparison of the evolution of the integration error versus the number of integration points using the octree level marking strategy for various polynomial functions of order $k$ as an integrand.}
	\label{fig:degreecomp}
\end{figure}

\begin{figure}
	\centering
	\begin{subfigure}[b]{0.33\textwidth}
		\centering
		\includegraphics[width=\textwidth]{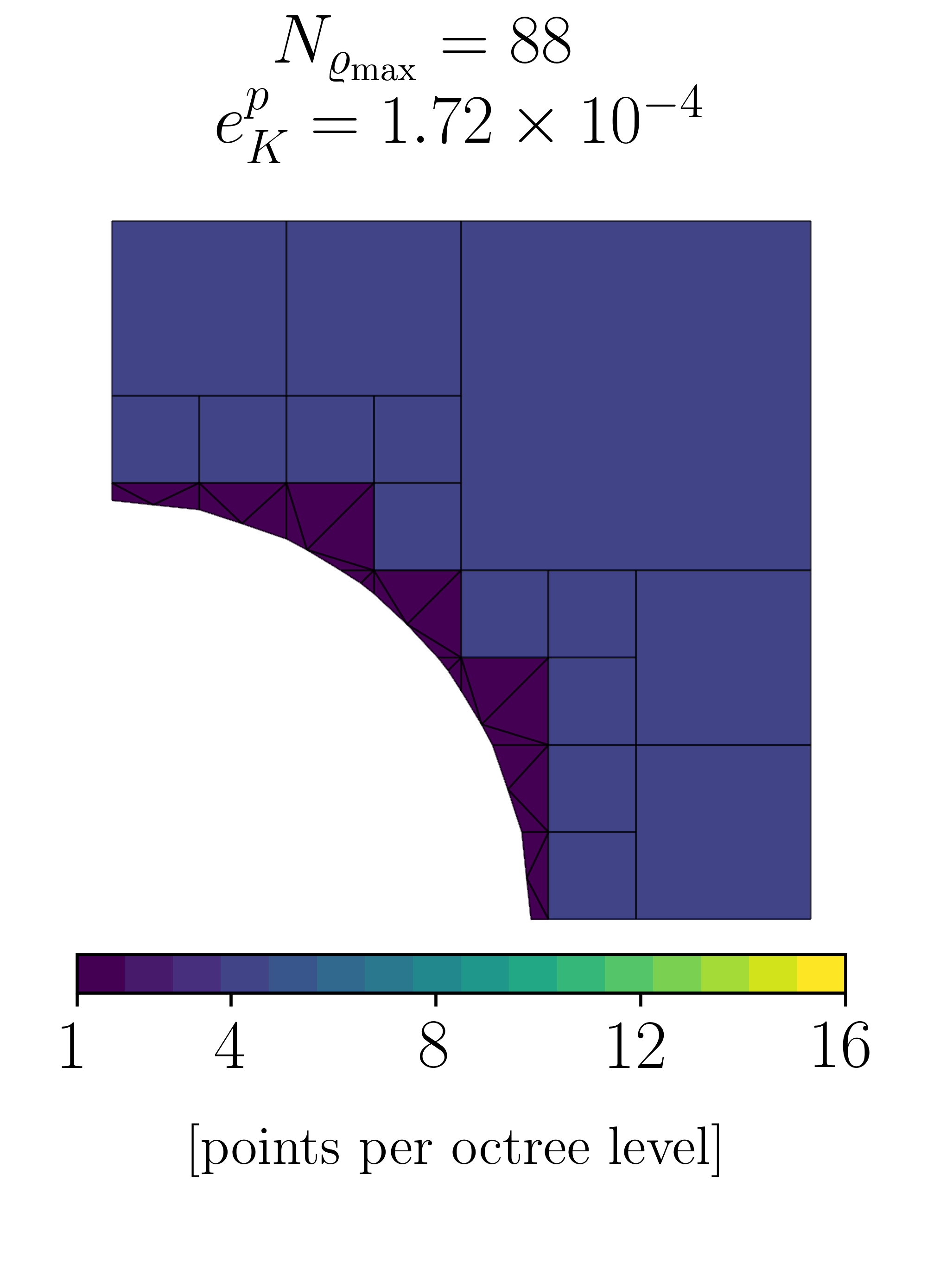}
		\caption{$k=2$}
		\label{fig:npts_k2}
	\end{subfigure}%
	\begin{subfigure}[b]{0.33\textwidth}
		\centering
		\includegraphics[width=\textwidth]{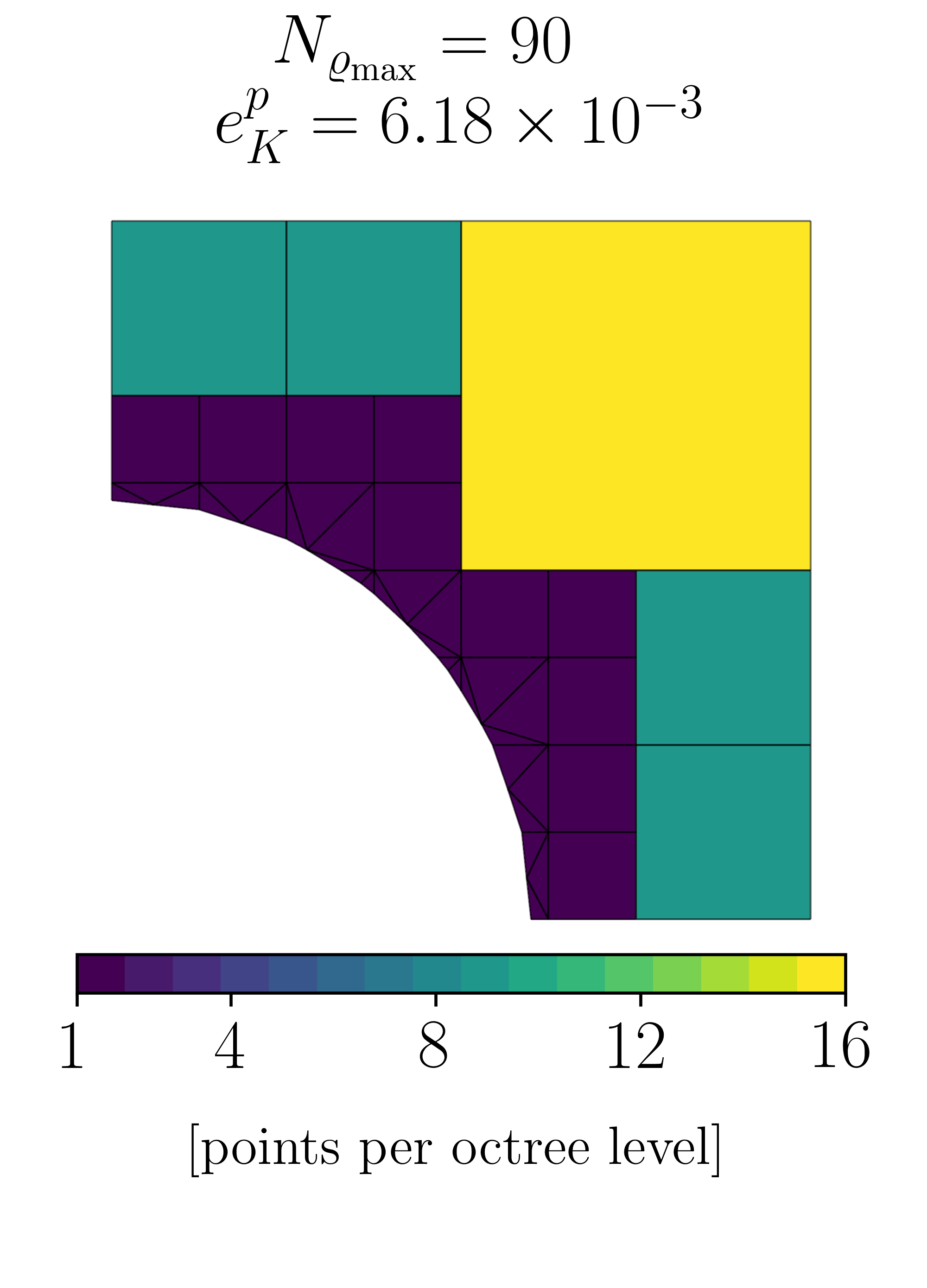}
		\caption{$k=8$}
		\label{fig:npts_k8}
	\end{subfigure}%
	\begin{subfigure}[b]{0.33\textwidth}
		\centering
		\includegraphics[width=\textwidth]{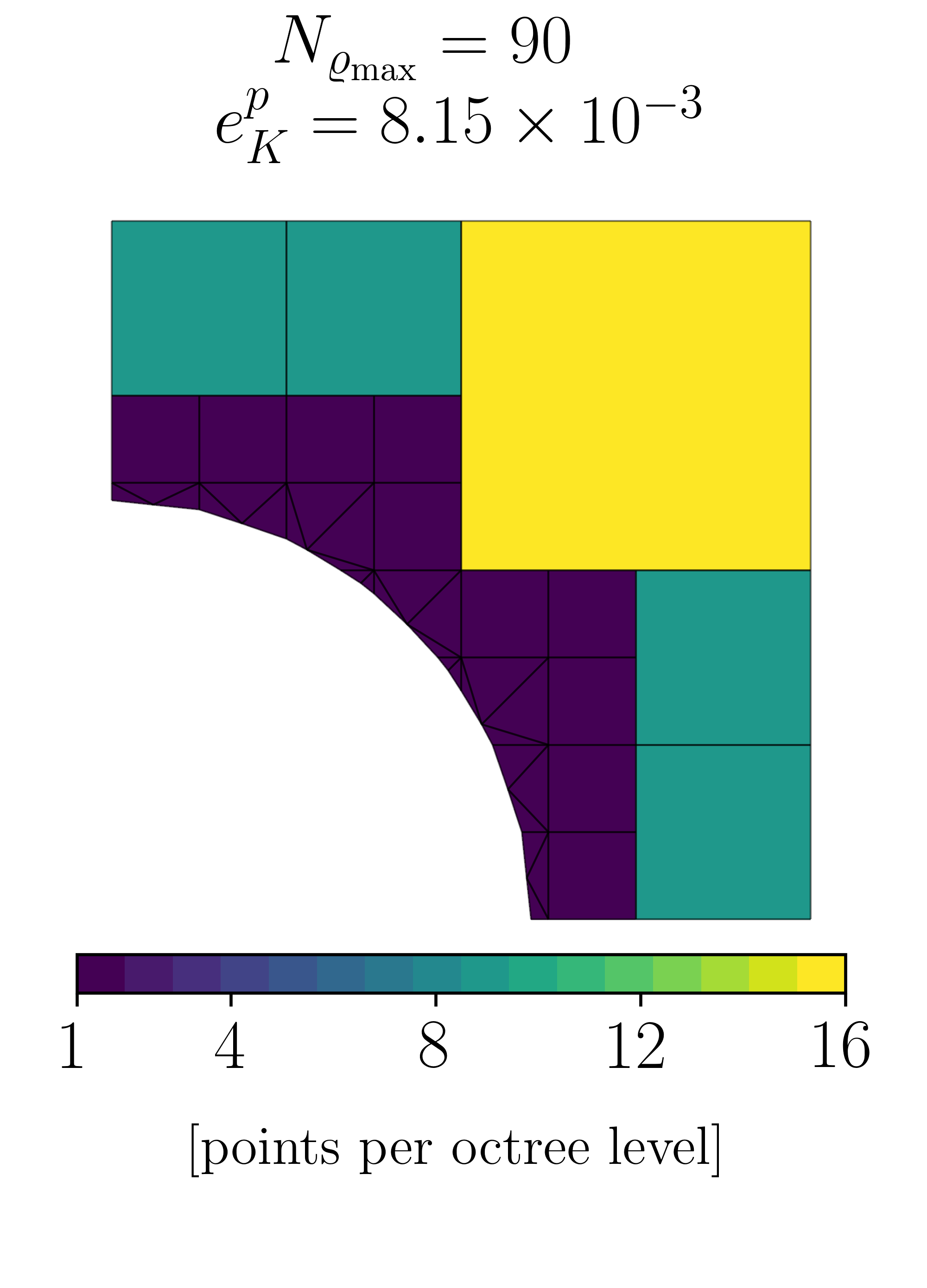}
		\caption{$k=10$}
		\label{fig:npts_k10}
	\end{subfigure}
	\caption{Distribution of approximately 90 integration points over the cut-element in two dimensions using the octree level marking strategy for various polynomial orders $k$.}
	\label{fig:degreecompdist}
\end{figure}

In Figure~\ref{fig:degreecomp} we study the influence of the integration error for various orders of the integrands. Note that the lower bound on the error at approximately $10^{-16}$ corresponds to machine precision errors, and hence these integrands can be interpreted to be exact. As expected, with an increase of the integrand order, a larger number of points is required to attain a certain accuracy. From the relatively low order settings, in particular the case of linear functions, it is observed that the algorithm quickly reaches machine precision, as all sub-cells are marked up until the degree required for exact polynomial integration. Note that in the initial setting, \emph{i.e.}, one point per sub-cell, the $k=1$ integrands are not resolved exactly, because the cross term $x_1 x_2$ (in two dimensions) is not integrated exactly on the triangulated sub-cells. The same applies to the three-dimensional setting.

In Figure~\ref{fig:degreecompdist} the point distributions for the two dimensional case corresponding to approximately 90 points are shown for various orders. It is observed that the distribution of the points over the levels is initially rather uniform for the $k=2$ case, but becomes more dispersed as the order increases. The reason for this is that the higher-order integrands encompass functions that are more localized to the bigger sub-cells than the lower-order integrands. For a fixed number of integration points, the distribution of integration points over the elements is observed to converge with increasing integrand orders.


\begin{figure}
	\centering
	\begin{subfigure}[b]{0.5\textwidth}
		\includegraphics[width=\textwidth]{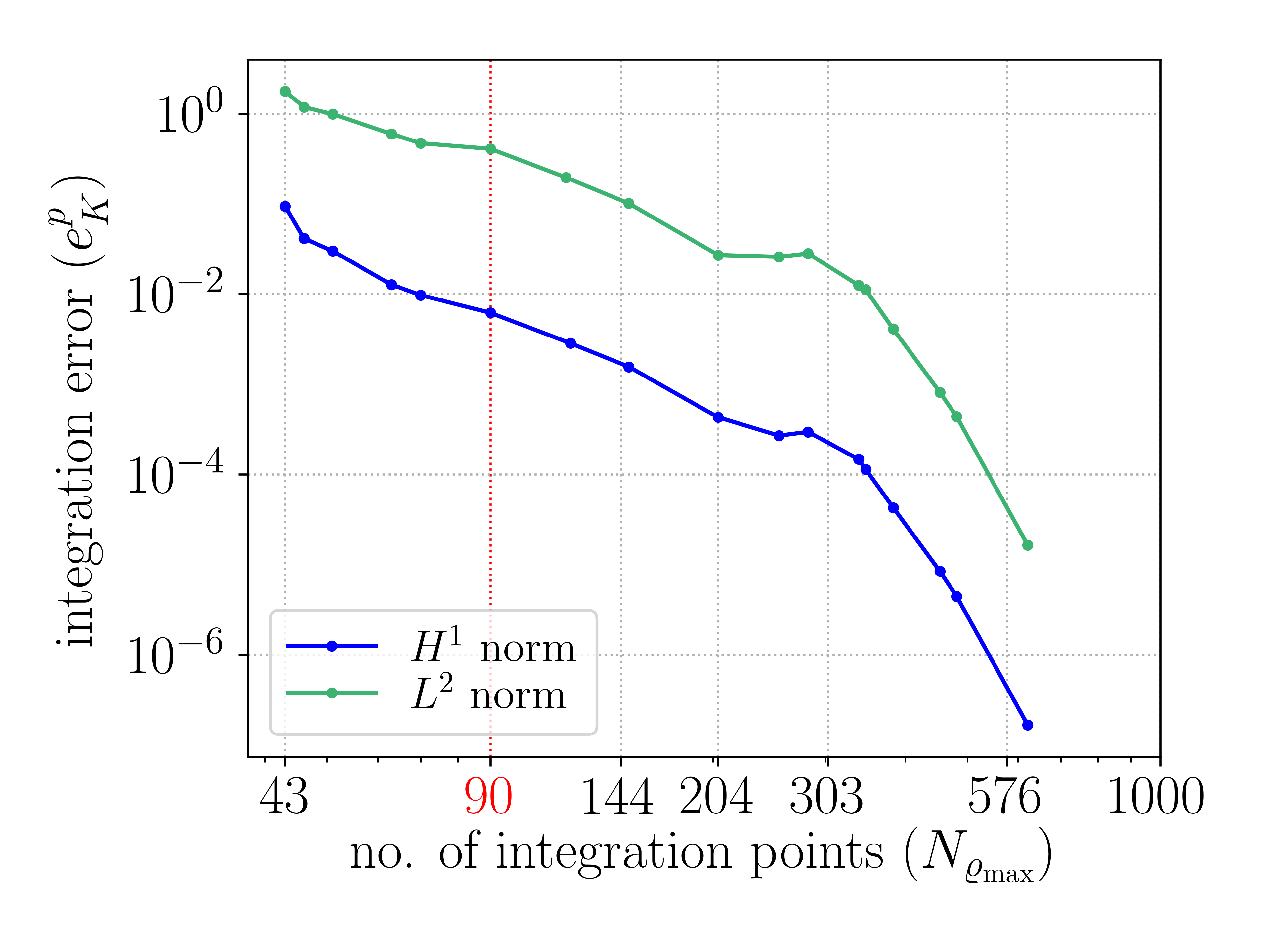}
		\caption{$d=2$}
		\label{fig:h1vsl2_2d}
	\end{subfigure}%
	\begin{subfigure}[b]{0.5\textwidth}
		\includegraphics[width=\textwidth]{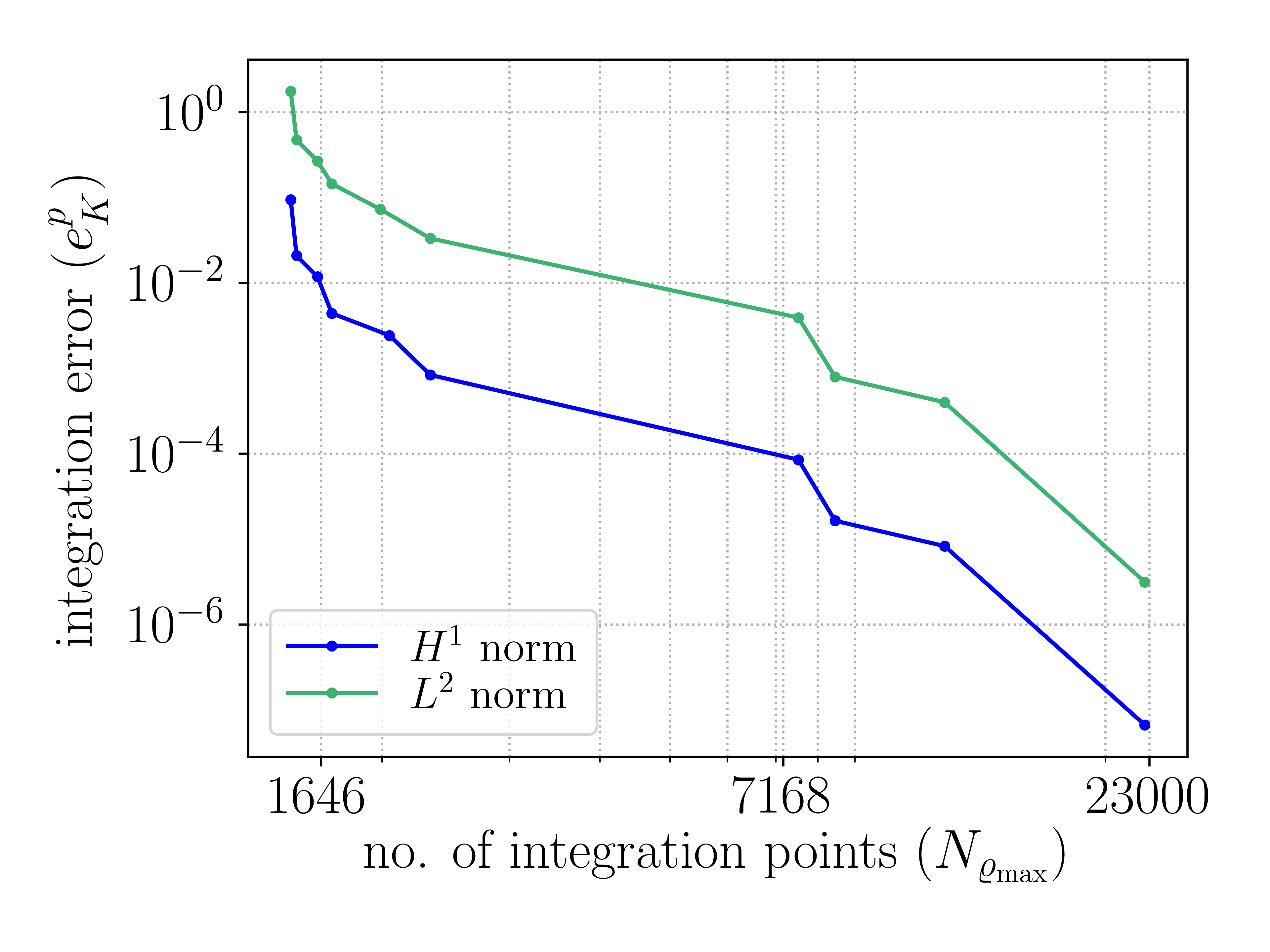}
		\caption{$d=3$}
		\label{fig:h1vsl2_3d}
	\end{subfigure}
	\caption{Comparison of the evolution of the integration error versus the number of integration points using the octree level marking strategy and different settings of the integrand space norms.}
	\label{fig:h1vsl2}
\end{figure}

\begin{figure}
	\centering
	\begin{subfigure}[b]{0.42\textwidth}
		\centering
		\includegraphics[width=0.8\textwidth]{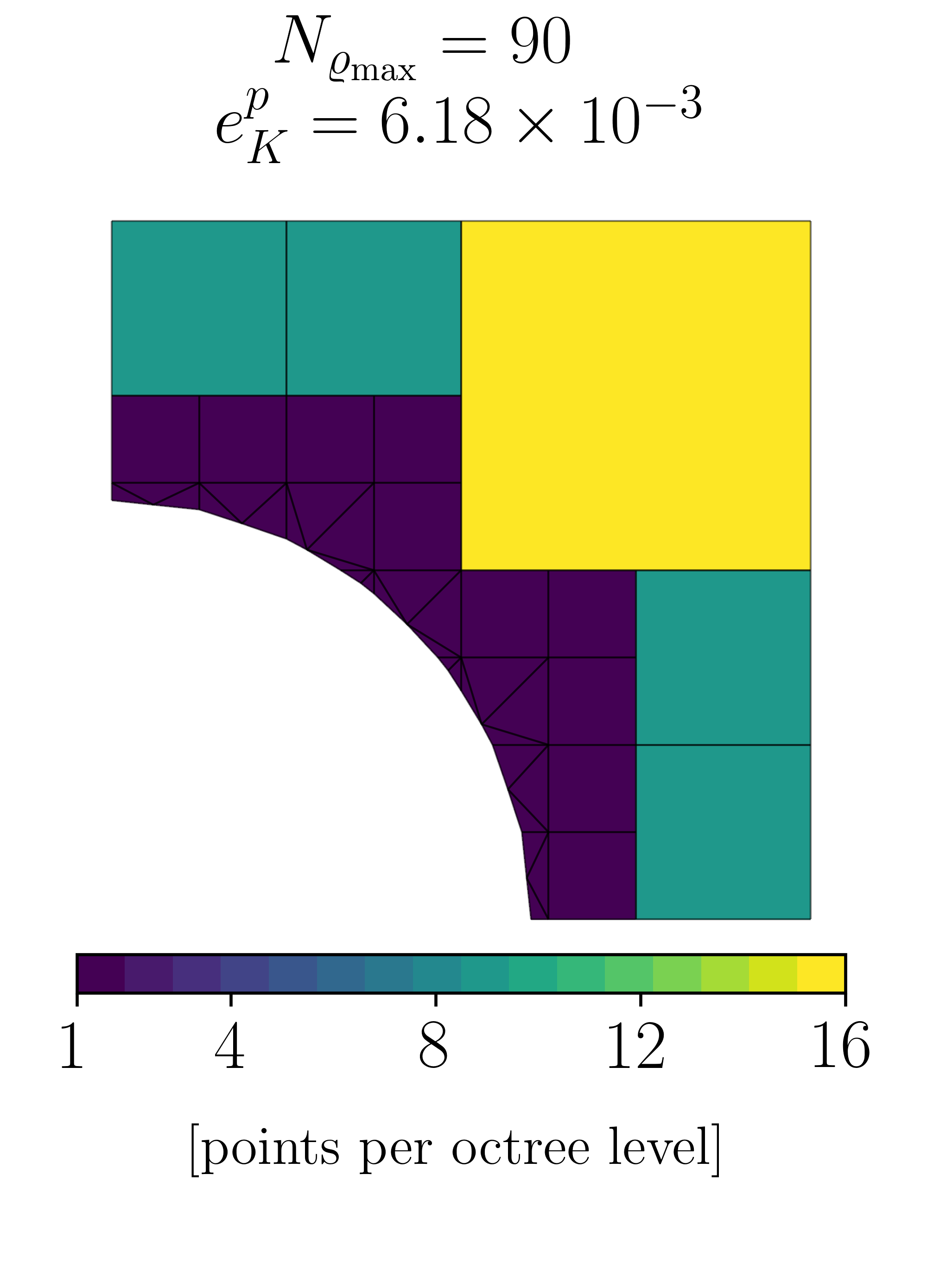}
		\caption{$\Hilbert^1$-normalized}
		\label{fig:npts_h1}
	\end{subfigure}%
	\begin{subfigure}[b]{0.42\textwidth}
		\centering
		\includegraphics[width=0.8\textwidth]{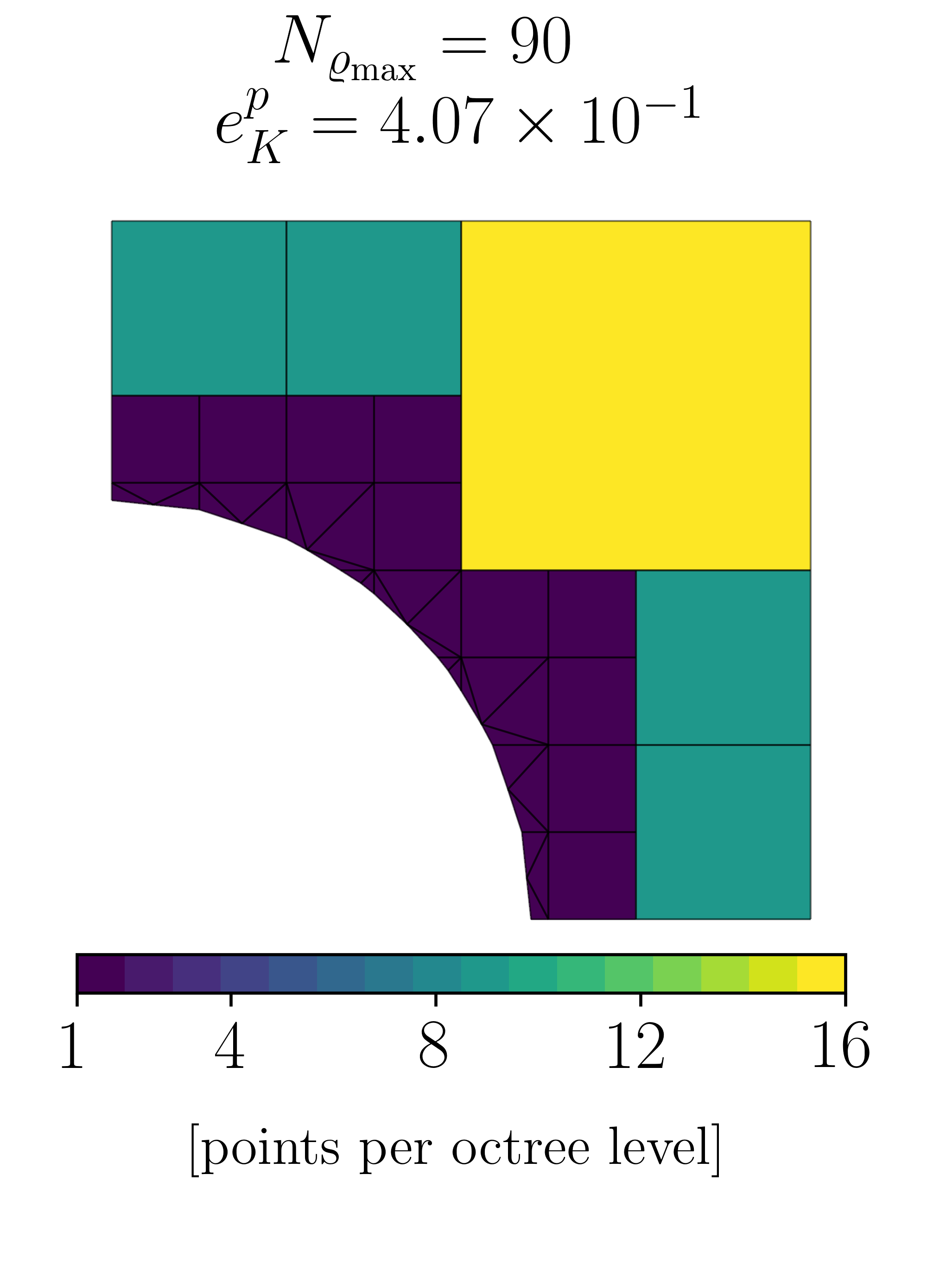}
		\caption{$L^2$-normalized}
		\label{fig:npts_l2}
	\end{subfigure}
	\caption{Distribution of 90 integration points over a two-dimensional cut-element using the octree level marking strategy and integration errors normalized by the $\Hilbert^1$ and $L^2$ norms.}
	\label{fig:normcompdist}
\end{figure}

In Figure~\ref{fig:h1vsl2} the influence of the function space norm \eqref{eq:evaluatedpolynomial} is studied. Although the normalization does affect the overall magnitude of the error, the difference in norm is observed to have a minor effect on the distribution of the integration points. This is confirmed in Figure~\ref{fig:normcompdist} for the two-dimensional setting, where it is observed that for a total number of 90 points, the distribution using either the $L^2$-norm or the $\Hilbert^1$-norm is virtually identical.

\subsubsection{Influence of the cut-element geometry}
The studies presented above were all conducted for cut-elements with spherical exclusions. To study the influence of the geometry of the cut-elements, we here consider the effect of variations in the ellipticity of the exclusions, and in the orientation of the excluded ellipsoids.

\begin{figure}
	\centering
	\begin{subfigure}[b]{0.32\textwidth}
		\centering
		\includegraphics[width=\textwidth]{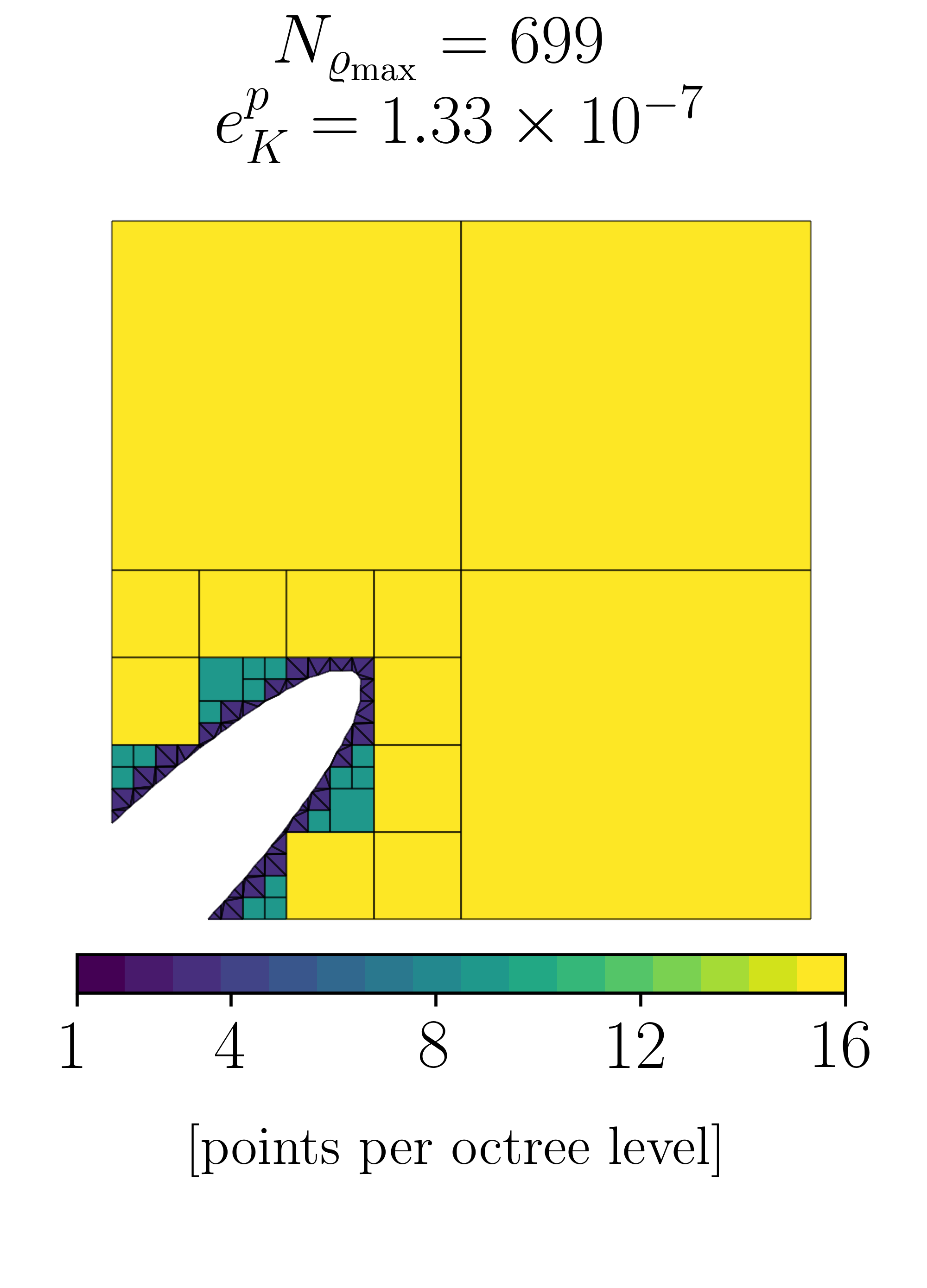}
		\caption{$r_{1} = 0.5$}
		\label{fig:npts_r5}
	\end{subfigure}
	\begin{subfigure}[b]{0.32\textwidth}
		\centering
		\includegraphics[width=\textwidth]{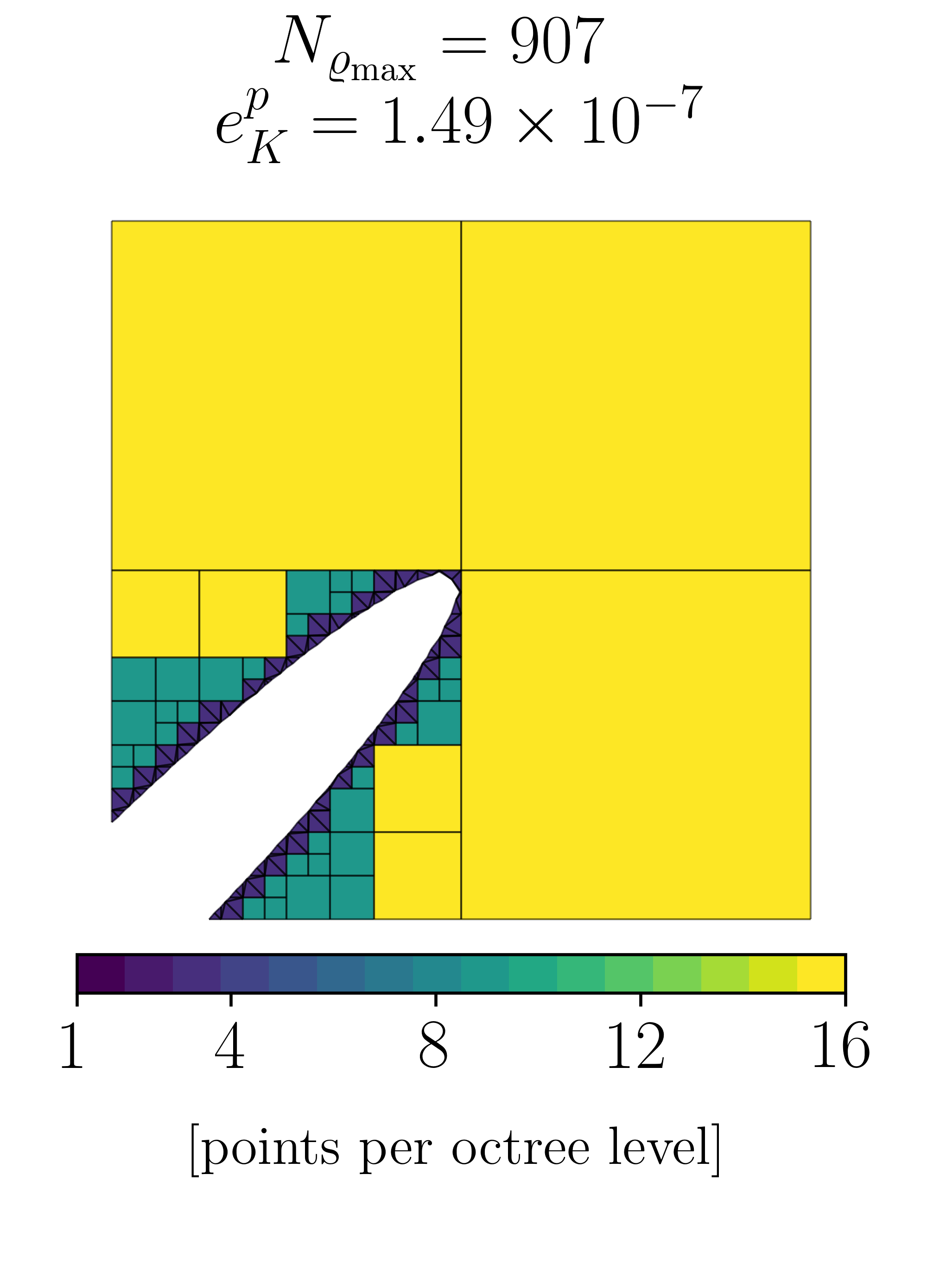}
		\caption{$r_{1} = 0.7$}
		\label{fig:npts_r7}
	\end{subfigure}
	\begin{subfigure}[b]{0.32\textwidth}
		\centering
		\includegraphics[width=\textwidth]{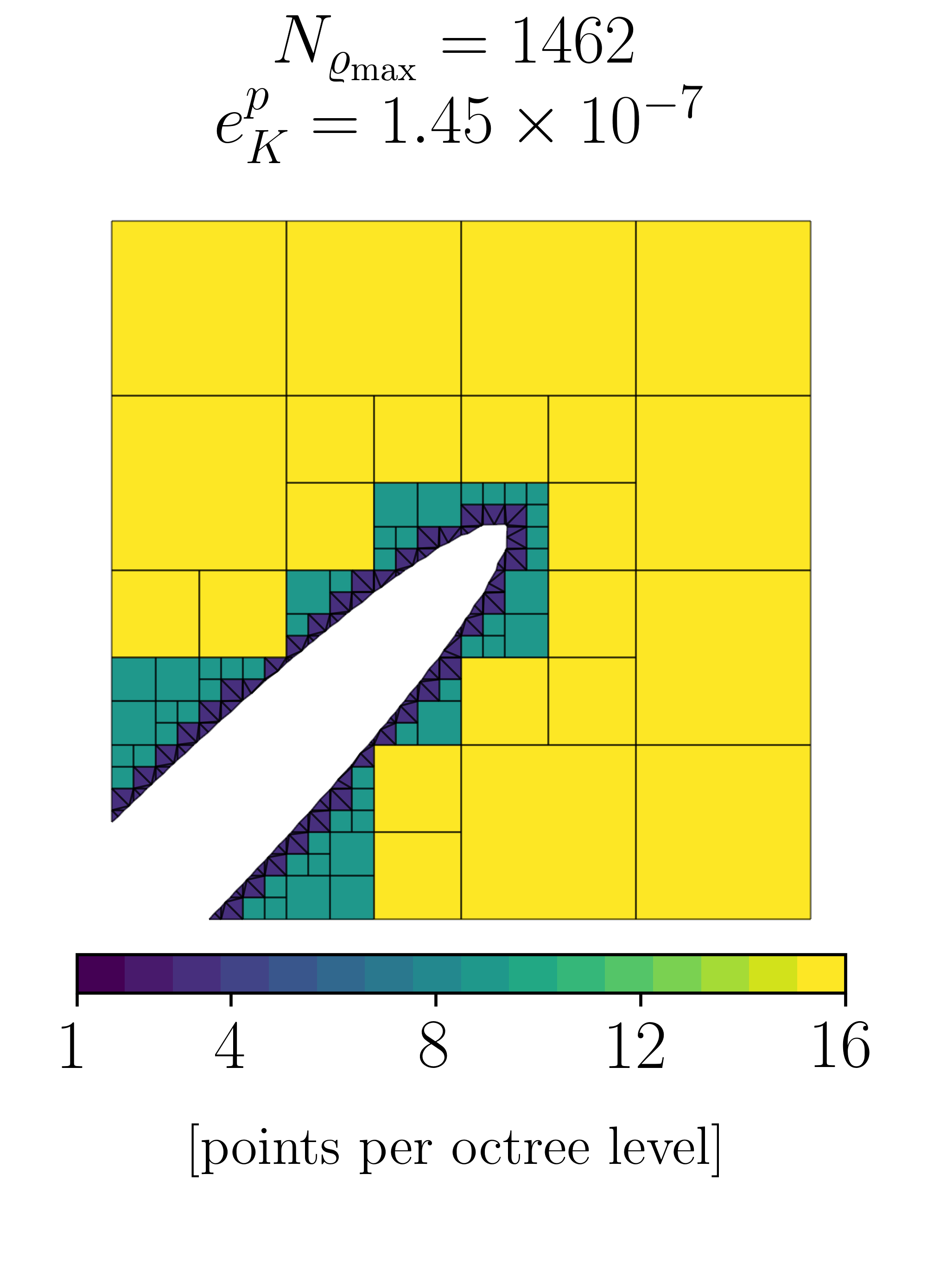}
		\caption{$r_{1} = 0.8$}
		\label{fig:npts_r8}
	\end{subfigure} \\[6pt]
	\caption{Distribution of integration points over cut-elements with various surface to volume ratios $\sv$. The inclination angle for all cases is equal to $\varphi=45^\circ$.}
	\label{fig:optimaldist_convgeom}
\end{figure}

\begin{figure}
	\centering
	\begin{subfigure}[b]{0.5\textwidth}
		\includegraphics[width=\textwidth]{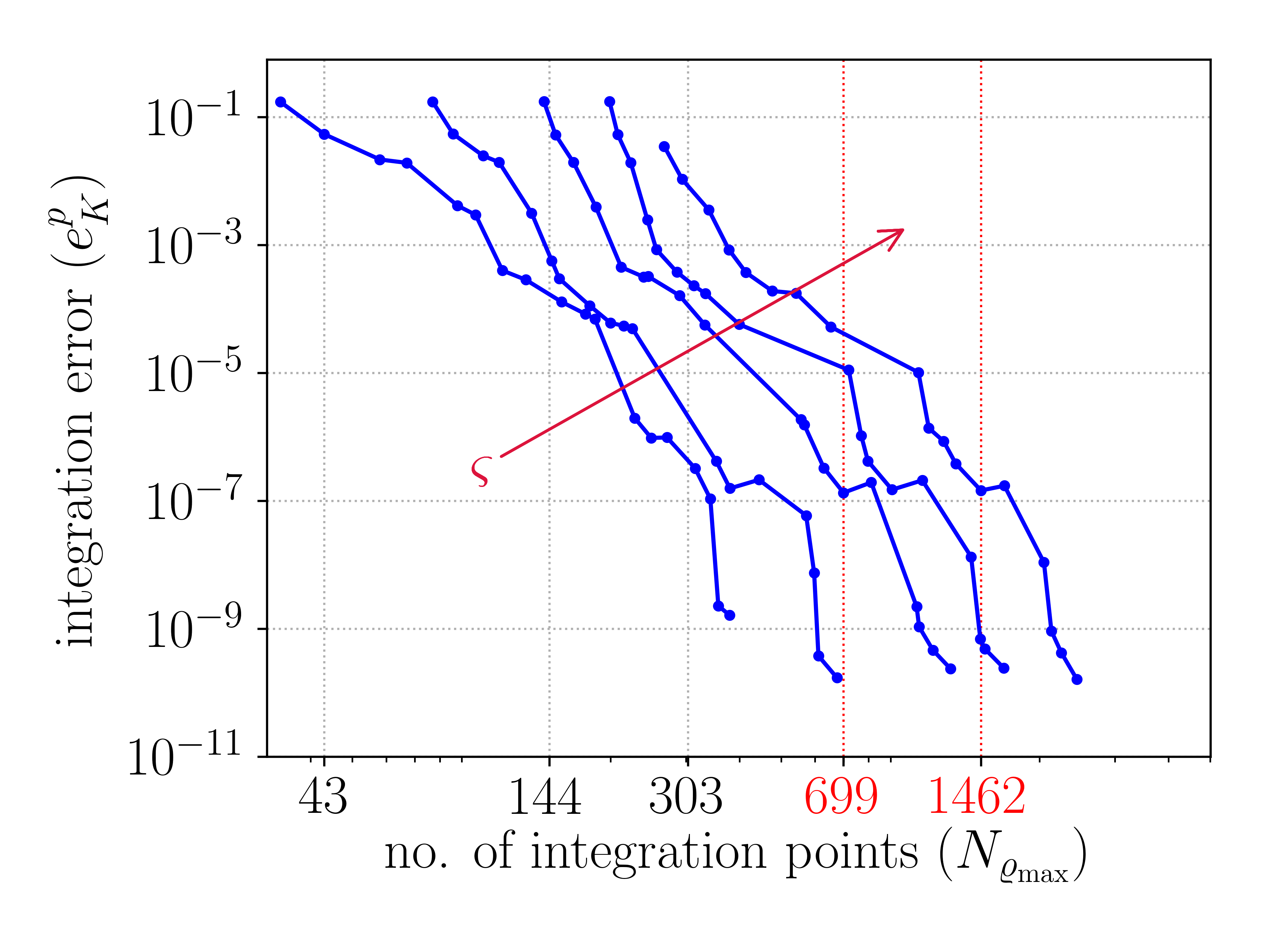}
		\caption{$d=2$}
		\label{fig:vf_2d}
	\end{subfigure}%
	\begin{subfigure}[b]{0.5\textwidth}
		\includegraphics[width=\textwidth]{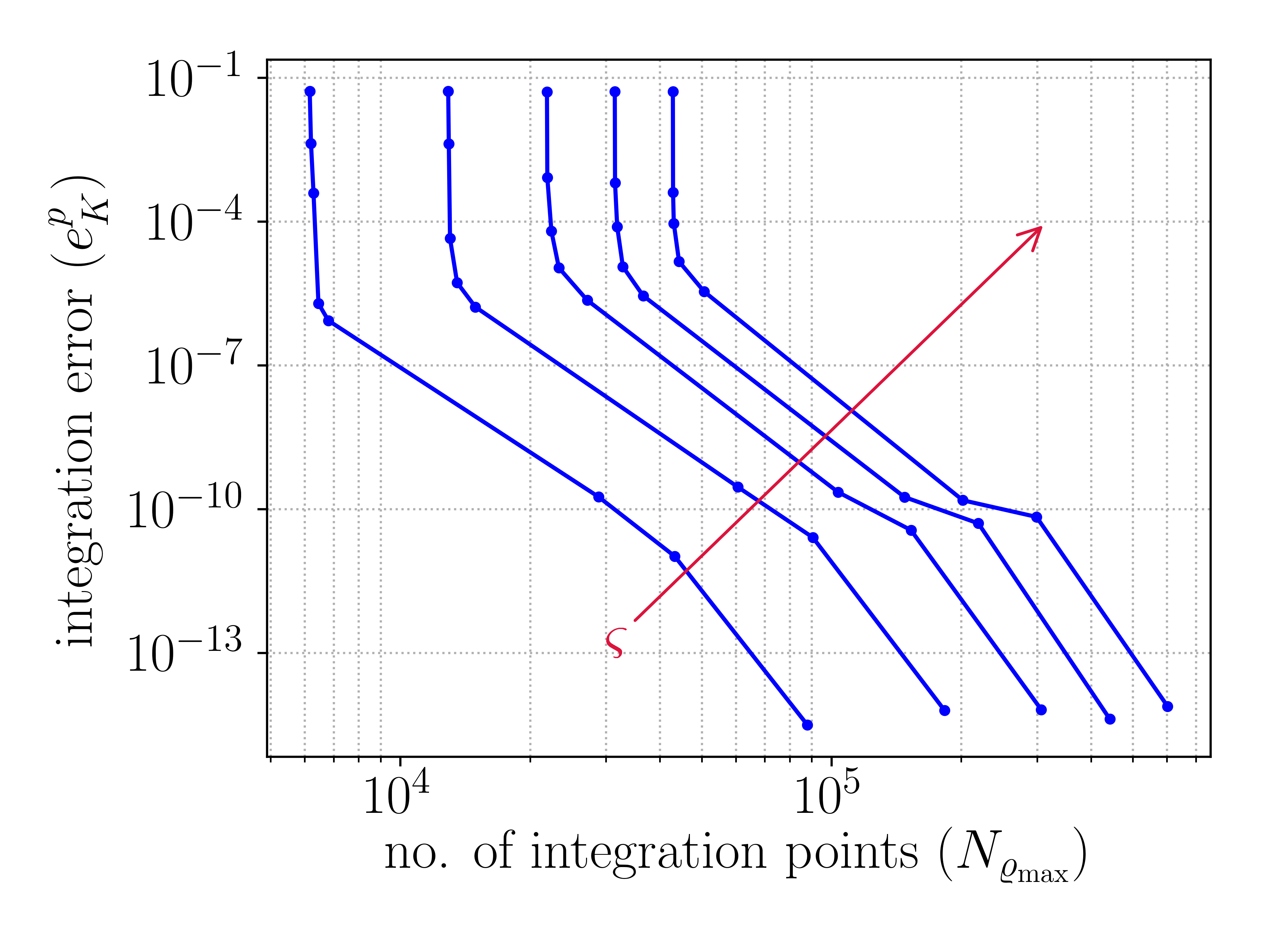}
		\caption{$d=3$}
		\label{fig:vf_3d}
	\end{subfigure}
	\caption{Comparison of the evolution of the integration error versus the number of integration points using the octree level marking strategy for various surface to volume ratios $\sv$.}
	\label{fig:convgeom}
\end{figure}

In Figure~\ref{fig:optimaldist_convgeom} we consider a range of exclusions in a two-dimensional cut-element with differing ellipticity. The corresponding error per integration point using the per-level adaptive integration procedure is displayed in Figure~\ref{fig:convgeom}. It is noted that the rate with which the error decreases is similar for all geometries, but that the number of points to attain a particular error is geometry dependent. This behavior is explained by the fact that although the distribution of the orders over the levels is very similar for each of the geometries, the number of sub-cells in each level is geometry dependent. As elaborated in Section~\ref{sec:octree_partitioning} the number of cells on each level scales with the surface fraction. This scaling is reflected in the results in Figure~\ref{fig:optimaldist_convgeom}. For example, the surface-to-volume ratio of Figure~\ref{fig:npts_r8} is approximately two times that of Figure~\ref{fig:npts_r5}. An increase in number of integration points by a factor of approximately $2$ is also observed from the offset of the corresponding curves in Figure~\ref{fig:convgeom}.

\begin{figure}
	\centering
	\begin{subfigure}[b]{0.5\textwidth}
		\includegraphics[width=\textwidth]{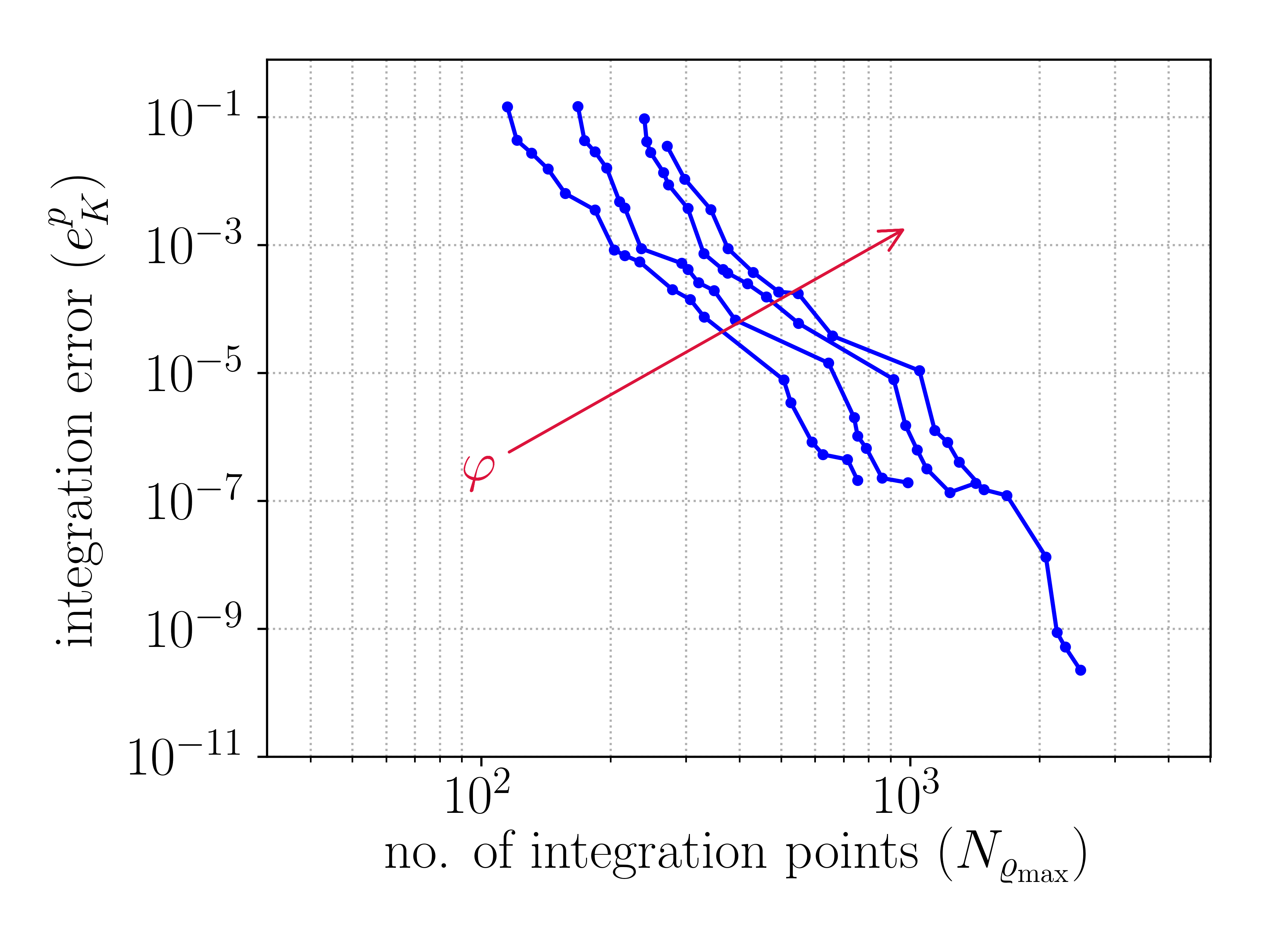}
		\caption{$d=2$}
		\label{fig:phi_2d}
	\end{subfigure}%
	\begin{subfigure}[b]{0.5\textwidth}
		\includegraphics[width=\textwidth]{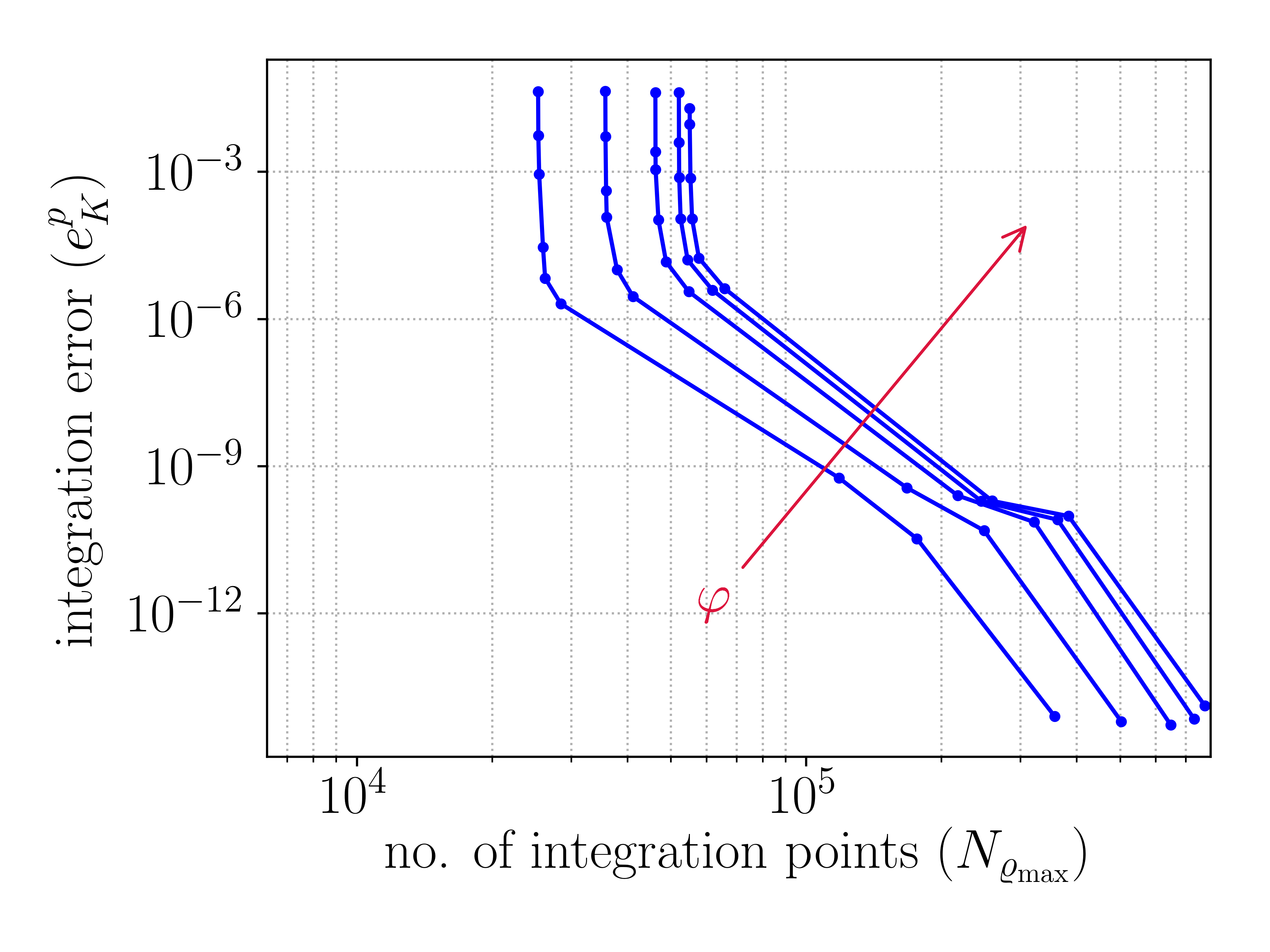}
		\caption{$d=3$}
		\label{fig:phi_3d}
	\end{subfigure}
	\caption{Comparison of the evolution of the integration error versus the number of integration points using the octree level marking strategy for various surface to volume ratios $\sv$ related to an exclusion with $r_1=0.6$, $r_2=0.1$ and varying inclination angle $\varphi$.}
	\label{fig:convgeom_phi}
\end{figure}

Figure~\ref{fig:convgeom_phi} shows the integration error versus the number of integration points for a cut-cell with $r_1=0.6$ and $r_2=0.1$ for various inclination angles $\varphi$. This study confirms the scaling relation with the surface to volume ratio as observed for the ellipticity variations considered above.

\subsection{Manually selected quadrature rules}
Although the computational effort involved in the construction of the optimized quadrature rules is in general acceptable when one wants to re-use a quadrature rule multiple times, some computational effort is involved in this construction. In addition, one has to set up a suitable code to determine the optimal distributions for arbitrary cut-cells. Considering this, one may not be interested in obtaining the optimized distributions of the points, but may instead want a simple rule of thumb to select the quadrature on a cut-cell; see, \emph{e.g.}, Refs.~\cite{alireza2013,alireza2019}. 

\begin{figure}
	\centering
	\begin{subfigure}[b]{0.42\textwidth}
		\centering
		\includegraphics[width=0.8\textwidth]{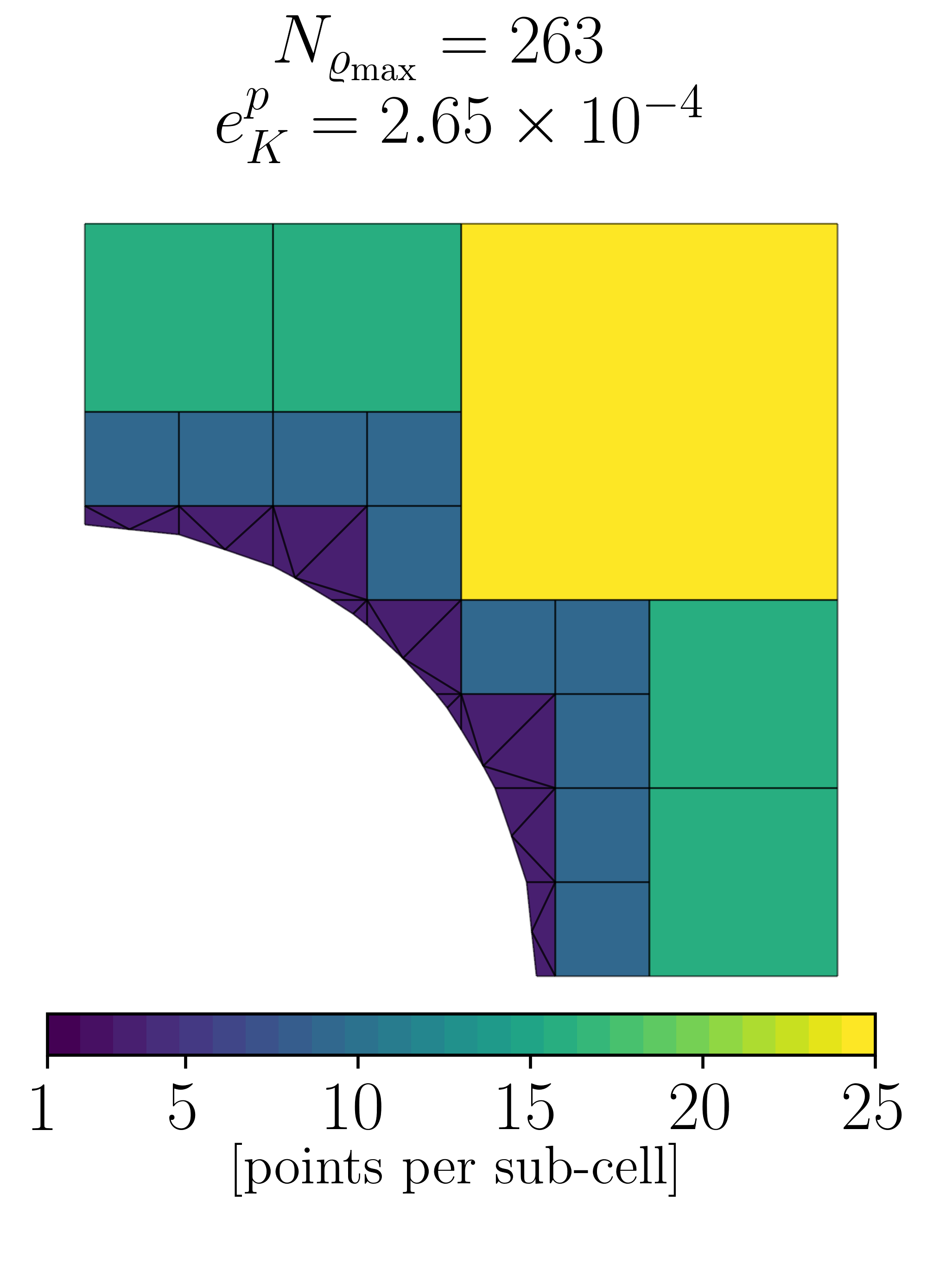}
		\caption{Rule of thumb A.}
		\label{fig:npts_ml263}
	\end{subfigure}%
	\begin{subfigure}[b]{0.42\textwidth}
		\centering
		\includegraphics[width=0.8\textwidth]{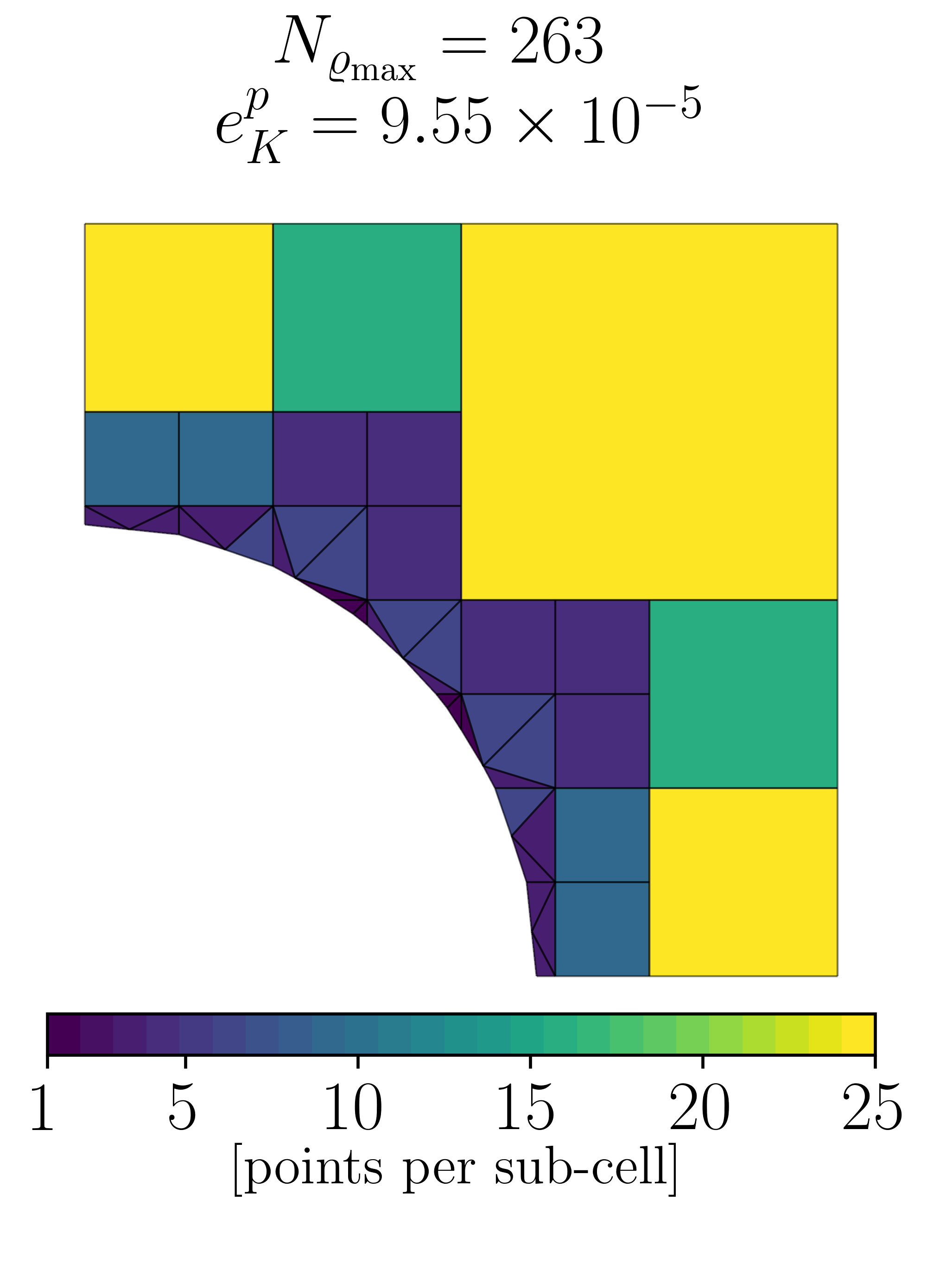}
		\caption{Per-cell optimized}
	\end{subfigure}\\[12pt]
	\begin{subfigure}[b]{0.42\textwidth}
		\centering
		\includegraphics[width=0.8\textwidth]{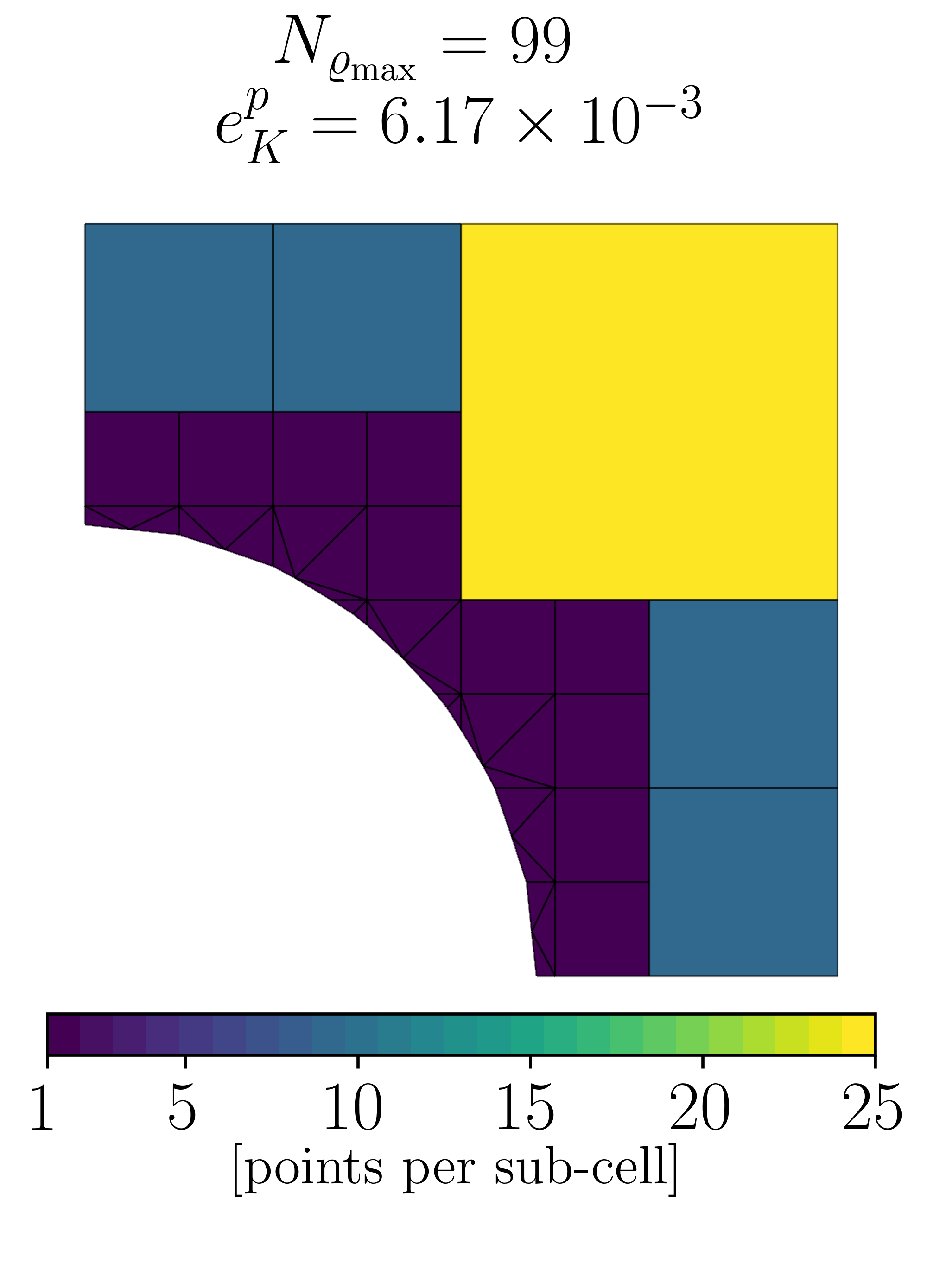}
		\caption{Rule of thumb B.}
		\label{fig:npts_ml99}
	\end{subfigure}%
	\begin{subfigure}[b]{0.42\textwidth}
		\centering
		\includegraphics[width=0.8\textwidth]{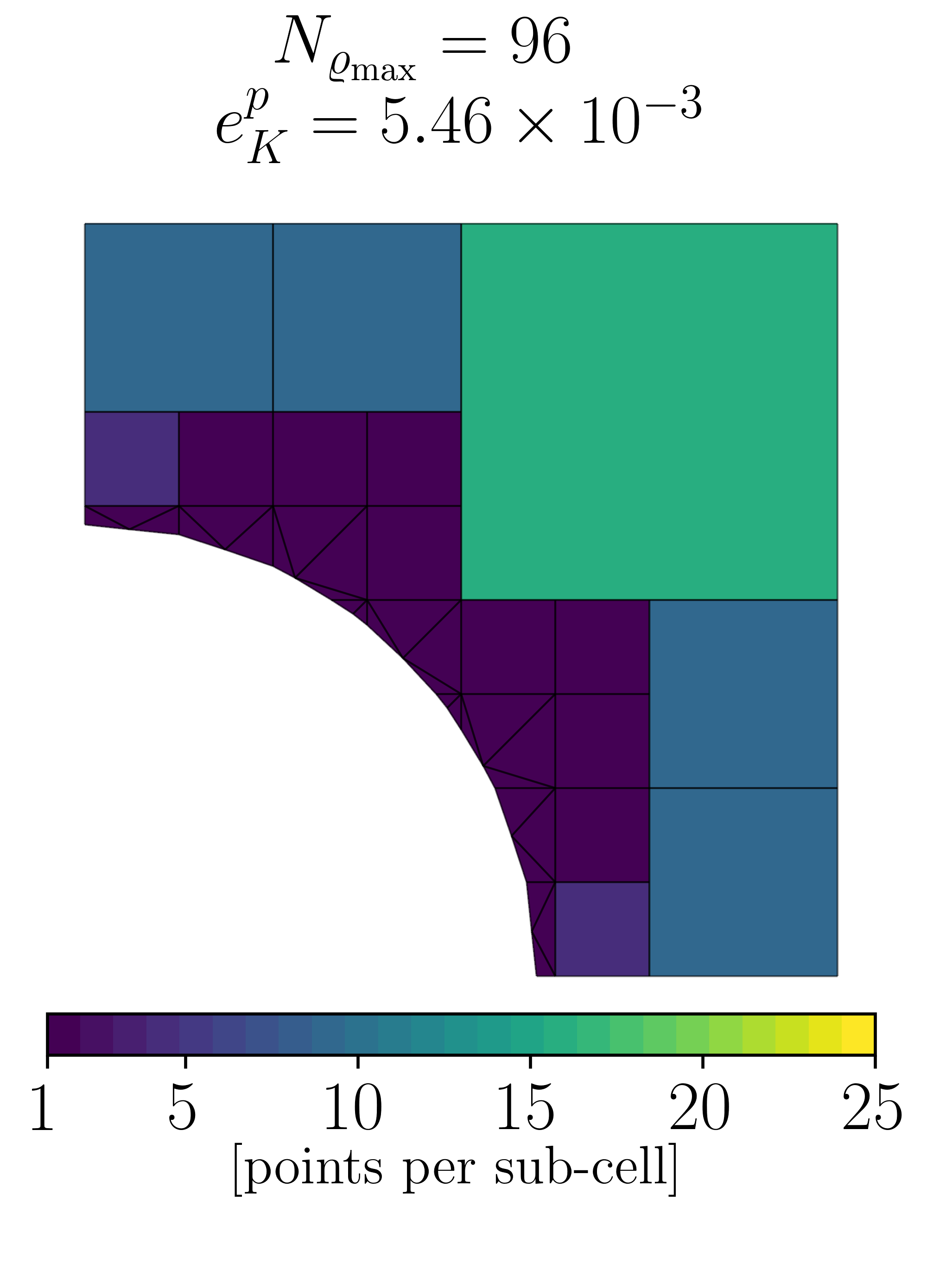}
		\caption{Per-cell optimized}
	\end{subfigure}

	\caption{Distribution of integration points over the cut-element in two dimensions using the manual selection of quadrature rules per-level for a polynomial order $k=8$.}
	\label{fig:mldegreedist}
\end{figure}

The per-level selection of the integration order makes it practical to manually select integration rules that outperform full order integration on all octree levels. We here consider two manual selection strategies based on Gauss integration of the sub-cells:
\begin{itemize}
  \item[A.] \emph{Minimal degree lowering:} In this strategy we set the order of the Gauss scheme on the level $\level=1$ sub-cells to $k_{\rm max}$. We then reduce the number of Gauss points per direction by one for each level, which implies that for $\level \leq \maxlevel$ the Gauss degree is decreased by two between two octree levels. The integration order of the tessellated sub-cells at level $\maxlevel+1$ is set to be two orders lower than that at level $\maxlevel$. Once a Gauss degree of zero is reached this value is maintained for the underlying levels. For example, for the case of $k=8$ and $\maxlevel=3$, the Gauss orders over the levels are set to $8$ for $\level=1$, $6$ for $\level=2$, $4$ for the $\level=\maxlevel$, and $2$ for the $\level=\maxlevel+1$ (see Figure~\ref{fig:npts_ml263}). That is, for an octree-depth of $\maxlevel=3$, the list of levels $[1,2,3,4]$ contains $[8,6,4,2]$ as the corresponding list of the Gauss orders.
  \item[B.] \emph{Uniform degree lowering:} This strategy also starts with selecting the $\level=1$ integration degree as $k_{\rm max}$, but then decreases the degree between two levels in such a way that the single point (degree is zero) is reached at the levels $\maxlevel$, and $\maxlevel+1$. For the case considered above, the orders are set to $8$ for $\level=1$, $4$ for $\level=2$, and $0$ for $\level=\maxlevel$ and  $\level=\maxlevel+1$ (see Figure~\ref{fig:npts_ml99}). That is, for an octree-depth of $\maxlevel=3$, the list of levels $[1,2,3,4]$ contains $[8,4,0,0]$ as the corresponding list of the Gauss orders.
\end{itemize}
Evidently, alternative quadrature rules can be formulated, but a detailed study of such rules is beyond the scope of the current work.

\begin{figure}
	\centering
	\begin{subfigure}{0.5\textwidth}
		\centering
		\includegraphics[width=\textwidth]{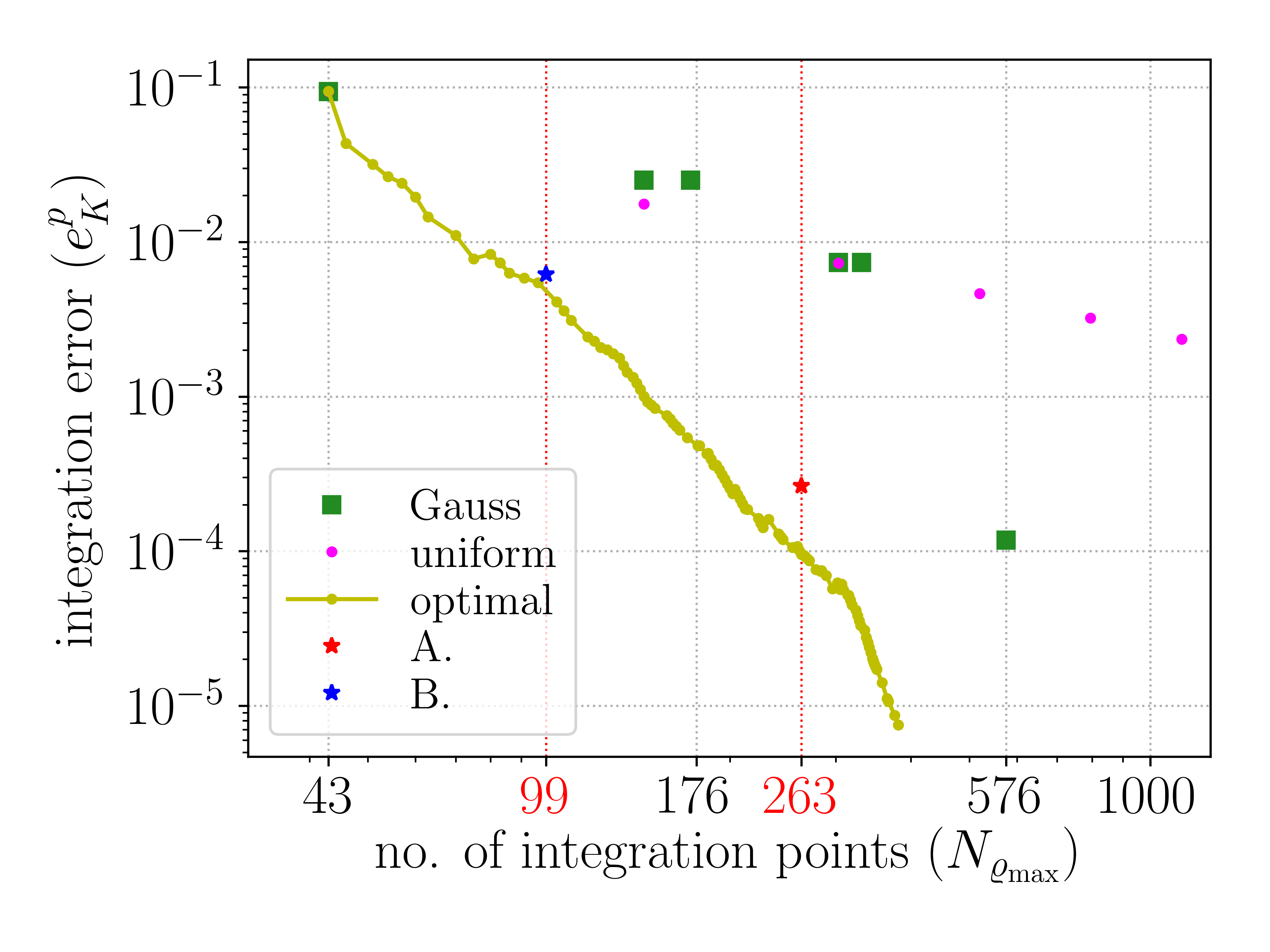}
		\caption{$d=2$}
   	\label{fig:thumbrule2d}
	\end{subfigure}%
	\begin{subfigure}{0.5\textwidth}
	\centering
	\includegraphics[width=\textwidth]{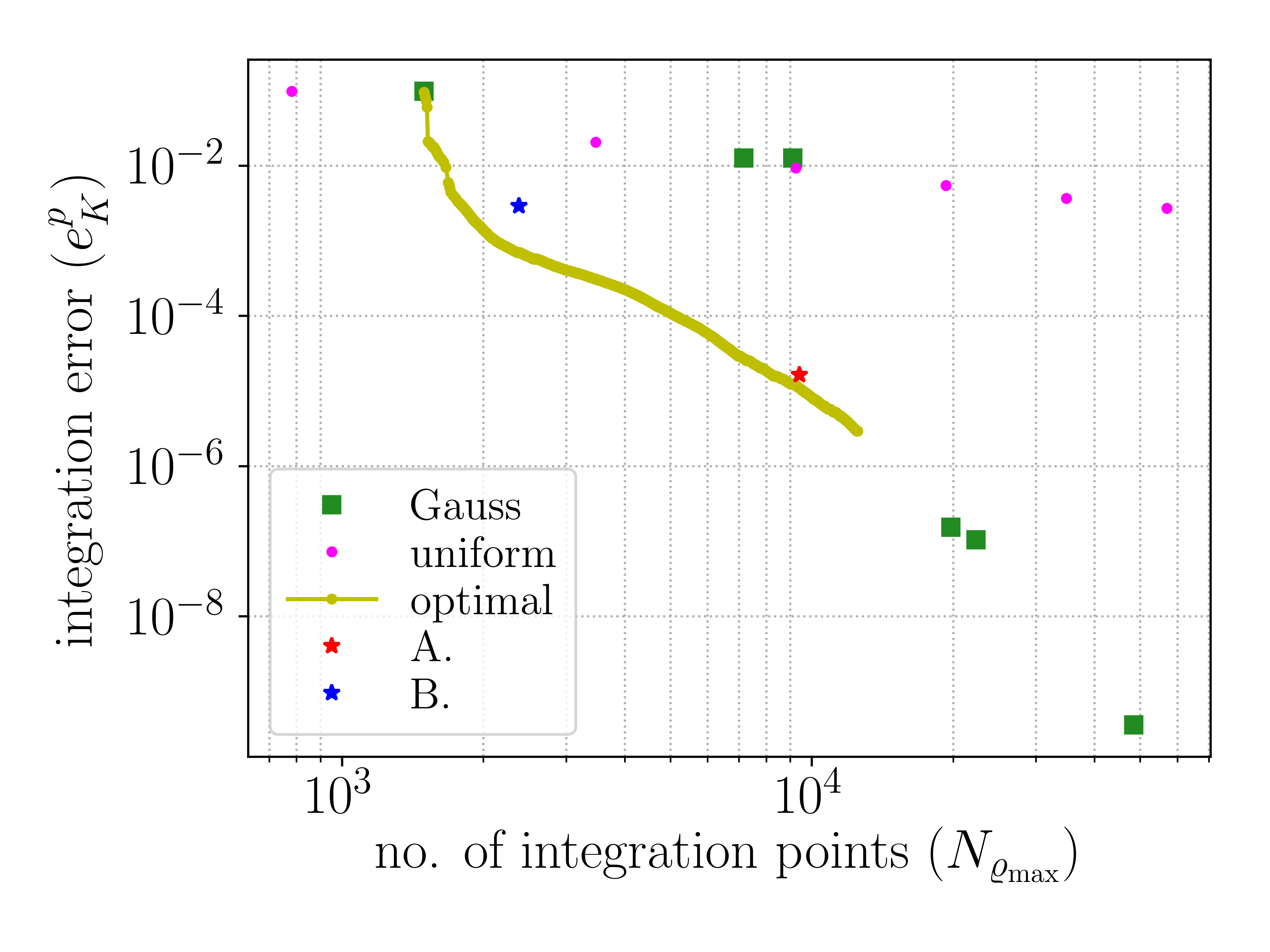}
	\caption{$d=3$}
   	\label{fig:thumbrule3d}
\end{subfigure}%
	\caption{Comparison of the integration error and the number of integration points using various thumb rules.}
	\label{fig:thumbrule}
\end{figure}

The adaptive algorithm developed in this work allows us to asses the suitability of the rules of thumb defined above. As a preliminary study on the effectiveness of these rules of thumb, we again consider the case of a two-dimensional circular exclusion with $\maxlevel = 3$ and a polynomial function of order $k=8$. In Figure~\ref{fig:thumbrule2d} the manually selected schemes are compared to the full order integration schemes and to the per sub-cell optimized integration schemes. This plot conveys that although the manually selected schemes are, as expected, outperformed by the optimized schemes, they generally do provide a dramatic improvement in accuracy per integration point for a fixed number of points. This observed behavior is explained by consideration of the distributions of orders over the sub-cells as shown in Figure~\ref{fig:mldegreedist}, from which it is observed that the per-level decrease of the order of both rules of thumb indeed qualitatively matches the results of the optimization procedure. Evidently, a difference between the two rules of thumb is that the minimal degree lowering strategy (A.) results in a larger number of points and a lower error compared to the uniform degree lowering (B.). Both schemes do, however, reasonably well resemble the optimized distribution patterns. Similar observations for the three-dimensional setting are obtained, as illustrated in Figure~\ref{fig:thumbrule3d}. The displayed three-dimensional results are based on $\maxlevel=3$ and $k=5$.

\section{Application to immersed isogeometric analysis} \label{sec:numerical_examples}
In this section we assess the suitability of the optimized integration schemes in the context of immersogeometric analysis. We consider an elasticity problem in two and three dimensions, which we solve using globally defined values of the integration degrees for each octree level. The integration degrees are obtained by application of the optimization procedure to the entire mesh, starting with a single integration point in all sub-cells. The integration errors computed per octree level are summed over all elements, after which the level with the highest global error is refined in terms of integration degree. It is noted that this global optimization strategy ignores operator-dependent variations between the elements, as discussed in Remark~\ref{rem:globalmarking}.

We consider a linear elasticity problem on a computational domain that is immersed into a background mesh of size $L^d$; see Figure~\ref{fig:le_domain_2d} for $d=2$. Displacements are prescribed on the exterior boundary, while the interior boundary is traction free. In the absence of inertia effects and body forces, the boundary value problem reads as
\begin{equation}
\left\{ \begin{aligned} 
& \text{Find $\boldsymbol{u}$ such that:} \\
&\boldsymbol{\rm div}(\, \boldsymbol{\sigma}(\boldsymbol{u}) \,) = \boldsymbol{0} & &\mbox{in} \: \Omega \\
& \boldsymbol{u} = \boldsymbol{0}  & & \mbox{on} \:  \partial \mathcal{A}_{0}\\
& \boldsymbol{u} = \bar{u} \boldsymbol{n}  & & \mbox{on} \:  \partial \mathcal{A}_{\bar{u}}\\
& \boldsymbol{u} \cdot \boldsymbol{n} = 0  & & \mbox{on} \:  \partial \mathcal{A} \setminus (  \mathcal{A}_{0}  \cup  \mathcal{A}_{\bar{u}} )\\
& \left[ \boldsymbol{I} - \boldsymbol{n} \otimes \boldsymbol{n}  \right]\boldsymbol{\sigma}\boldsymbol{n}= \boldsymbol{0}  & & \mbox{on} \:  \partial \mathcal{A} \setminus (  \mathcal{A}_{0}  \cup  \mathcal{A}_{\bar{u}} )\\
& \boldsymbol{\sigma} \boldsymbol{n}= \boldsymbol{0}  & & \mbox{on} \:  \partial \Omega \setminus \partial \mathcal{A}
\end{aligned}\right. \label{eq:strong_linear_elasticity}
\end{equation}
with stress tensor $\boldsymbol{\sigma} (\boldsymbol{u}) = \lambda \boldsymbol{\rm div} (\boldsymbol{u}) \boldsymbol{I} + 2 \mu \nabla^s \boldsymbol{u}$, where $\nabla^s$ denotes the symmetric gradient operator. Throughout this section, the Lam\'e parameters are set to $\lambda = \frac{1}{2}$ and $\mu = \frac{1}{2}$.

The Dirichlet condition on the external boundary is applied strongly, \emph{i.e.}, by constraining the degrees of freedom related to the boundary displacements. The Galerkin problem corresponding to \eqref{eq:strong_linear_elasticity} then follows as
\begin{equation}
\left\{ \begin{aligned} 
& \text{Find $\boldsymbol{u}^h \in W^h(\Omega)$ such that for all $\boldsymbol{v}^h \in W^h_0(\Omega)$:} \\
& \int_{\Omega} \nabla^s\boldsymbol{v}^h : \boldsymbol{\sigma}(\boldsymbol{u}^h) \,{\rm d}V = 0
\end{aligned}\right. \label{eq:linear_elasticiy}
\end{equation}
with the discrete trial space being a subset of $\Hilbert^1(\Omega)$ and satisfying the Dirichlet boundary conditions. The discrete test space is identical to the trial space, modulo inhomogeneous boundary conditions. All results in this section pertain to full regularity, $C^{p-1}$, B-splines of degree $p=2$. The polynomial order used to compute the optimized integration schemes is taken equal to $k=4$.

As a quantity of interest we consider the constrained modulus
\begin{equation}
 M^\star = \frac{\sigma_{11}}{\varepsilon_{11}} = \frac{L}{A} \frac{F}{\bar{u}},
\end{equation}
where $F$ is the resultant reaction force in normal direction along the displaced boundary, and where $A$ is the area (length in two dimensions) at which this resultant force acts. We normalize the computed constrained modulus by that of a homogeneous square with the same material properties as the considered model, \emph{i.e.}, by $M = \lambda+2\mu$.

For both the two and three dimensional test cases we consider three levels of bisectioning, \emph{i.e.}, $\maxlevel=3$. We represent the selection of the per level integration orders as a list of length $\maxlevel+2$, where the first $\maxlevel+1$ entries refer to the bisectioning levels $0$ to $\maxlevel$. Note that the level $\level=0$ in fact pertains to cells that are not intersected by the boundary. The final entry in the list of orders pertains to the integration order on the tessellated sub-cells at level $\maxlevel+1$.

\subsection{Two-dimensional test case}
\begin{figure}
	\centering
	\usetikzlibrary{decorations.markings}
\begin{tikzpicture}[decoration={markings,mark=between positions 0.25 and 0.75 step 0.125 with {\node [yshift=0.3cm] {$ $};},
	mark=between positions 0.1 and 1.5 step 0.1 with {\draw[very thick] (0pt,0pt) circle (3pt);},}]
  \tikzmath{\wdomain=2;
            \wambient=(3/2)*\wdomain;
            \h=0.3*\wdomain;
            \Rc=0.4*\wdomain;
            \Rsw=1.273/2*\wdomain;
            \Rse=0.5*\wdomain;
            \Rne=1.333/2*\wdomain;
            \Rnw=0.4*\wdomain;
            }

    \begin{scope}[scale=3*\wdomain,line join=bevel,ultra thin,fill=black!25,draw=black!25]
      \input{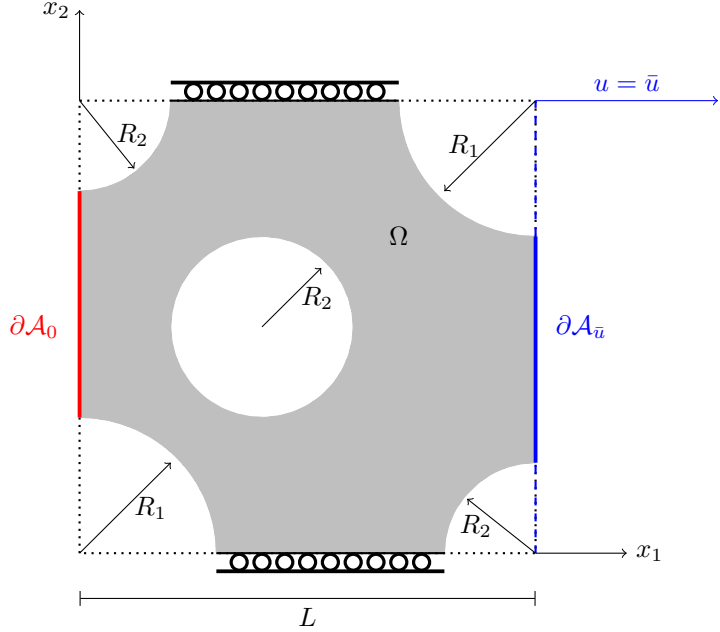}
    \end{scope}
    \draw[dotted,thick](0,0) rectangle (2*\wambient,2*\wambient);
    \draw[dashed,thick,color=blue] (10*\h,0*\h)--(10*\h,10*\h);

    \node at (7*\h,7*\h) {$\Omega$};
    \draw[->,color=blue] (10*\h,10*\h) -- (14*\h,10*\h) node[midway, above] {\color{blue} $u = \bar{u}$};
    \draw[->] (10*\h,0*\h) -- (12*\h,0*\h) node[right] {$x_1$};
    \draw[->] (0*\h,10*\h) -- (0*\h,12*\h) node[left] {$x_2$};
    \draw[->] (0*\h,0*\h) -- (2*\h,2*\h) node[midway, right] {$R_1$};
    \draw[->] (10*\h,10*\h) -- (8*\h,8*\h) node[midway, left] {$R_1$};
    \draw[->] (4*\h,5*\h) -- (5.3*\h,6.3*\h) node[midway, right] {$R_2$};
    \draw[->] (10*\h,0*\h) -- (8.5*\h,1.2*\h) node[midway, left] {$R_2$};    
    \draw[->] (0*\h,10*\h) -- (1.2*\h,8.5*\h) node[midway, right] {$R_2$};
    
    \draw[|-|] (0*\h,-1*\h) -- (10*\h,-1*\h) node[midway, below] {$L$};
    
    \draw[thick] (3*\h,0*\h)--(8*\h,0*\h);
    \path [postaction={decorate}] (3*\h,-0.2*\h) -- (8*\h,-0.2*\h);
    \draw[line width=0.5mm] (3*\h,-0.4*\h)--(8*\h,-0.4*\h);
    
    \draw[thick] (2*\h,10*\h)--(7*\h,10*\h);
	\path [postaction={decorate}] (2*\h,10.2*\h) -- (7*\h,10.2*\h);
	\draw[line width=0.5mm] (2*\h,10.4*\h)--(7*\h,10.4*\h);
	
    \draw[line width=0.5mm,color=red] (0*\h,3*\h)--(0*\h,8*\h);
    \draw[line width=0.5mm, color=blue] (10*\h,2*\h)--(10*\h,7*\h);
    
	\node at (11*\h,5*\h) {\color{blue} $\partial \mathcal{A}_{\bar{u}}$};
	\node at (-1*\h,5*\h) {\color{red} $\partial \mathcal{A}_{0}$};
\end{tikzpicture}
	\caption{Schematic representation of the considered domain and the mechanical loading conditions}
	\label{fig:le_domain_2d}
\end{figure}

\begin{figure}
	\centering
	\begin{subfigure}[b]{0.5\textwidth}
		\centering
		\includegraphics[width=\textwidth]{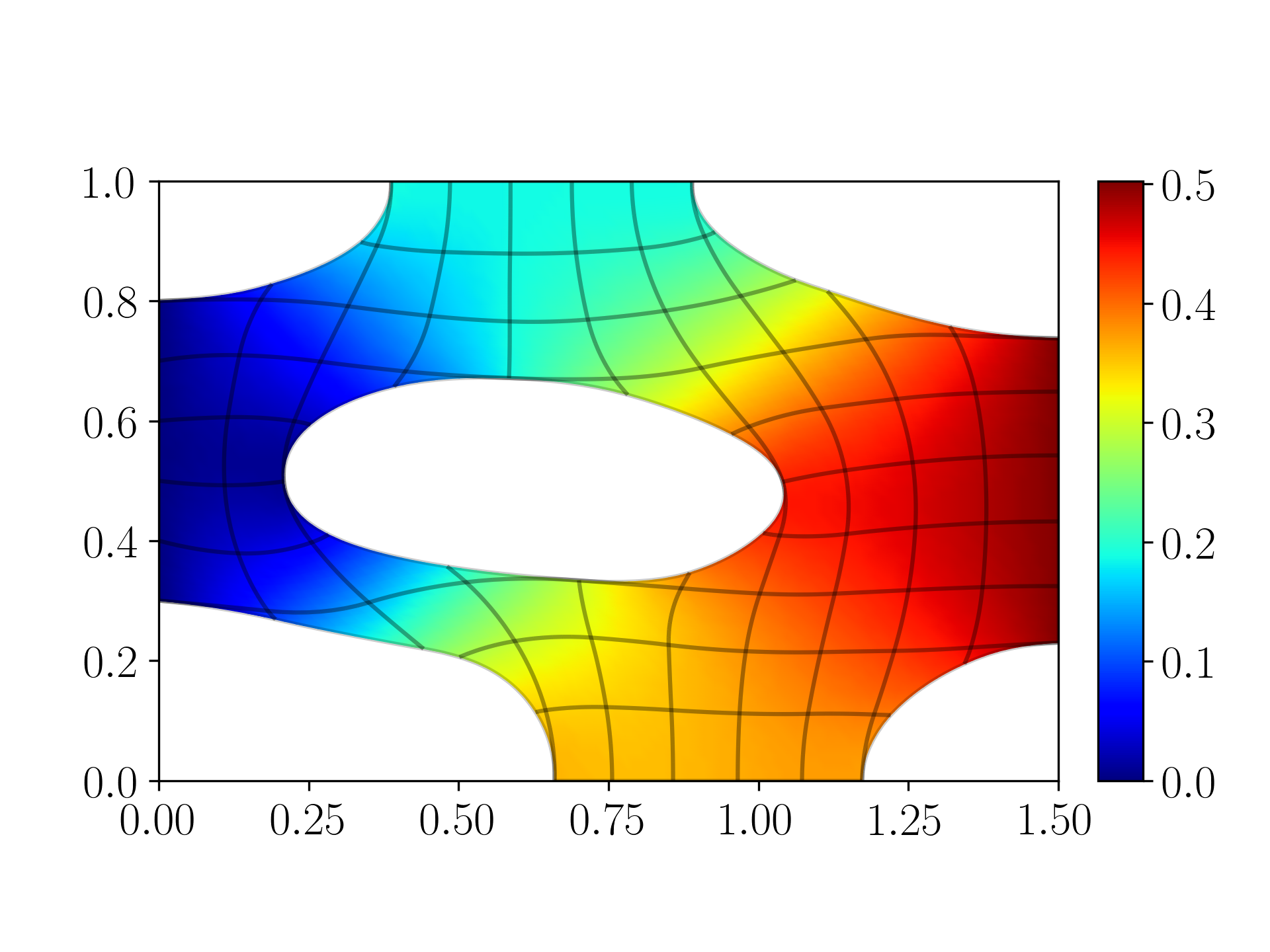}
		\caption{Displacement magnitude}
		\label{fig:displacement_2d}
	\end{subfigure}%
	\begin{subfigure}[b]{0.5\textwidth}
		\centering
		\includegraphics[width=\textwidth]{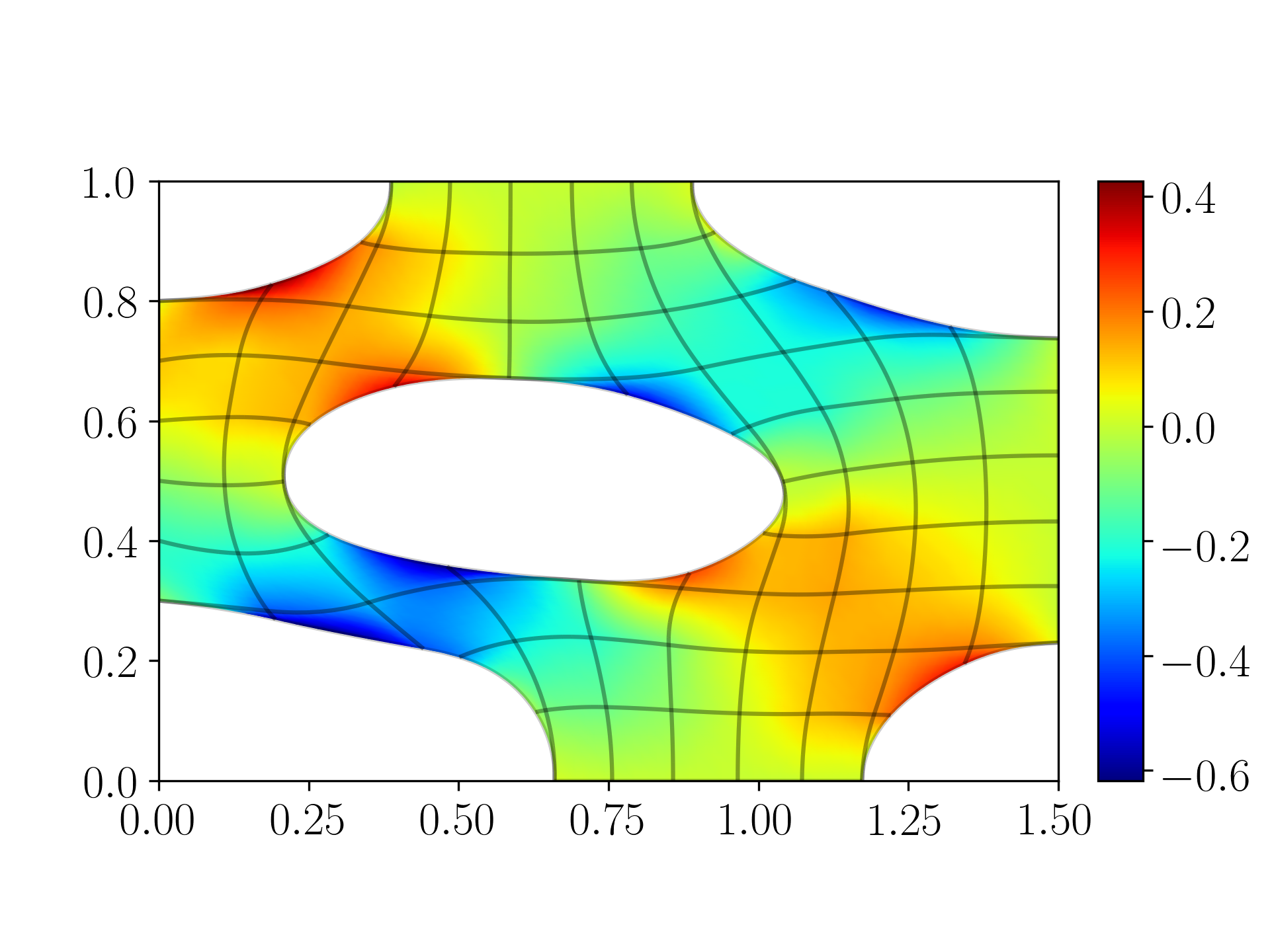}
		\caption{Shear stress $\sigma_{12}$}
		\label{fig:shear_2d}
	\end{subfigure}
	\caption{Solution of the two-dimensional linear elasticity problem, displayed in the deformed configuration.}
	\label{fig:le_results_2d}
\end{figure}

We consider the two-dimensional test case introduced in Ref.~\cite{deprenter2019}, which is illustrated in Figure~\ref{fig:le_domain_2d}. The size of the bounding square is set to $L =1$, and the radii of the exclusions are takes as $R_1 = 0.3$ and $R_2 = 0.2$. The right boundary is displaced by $\bar{u}=\frac{1}{2}$. Figure~\ref{fig:le_results_2d} displays the solution to the elasticity problem \eqref{eq:linear_elasticiy} in terms of the displacement magnitude and the shear stress. The displayed result is based on a $10 \times 10$ ambient domain mesh with uniformly distributed fourth order Gauss integration scheme.

\begin{figure}
	\centering
	\begin{subfigure}[b]{0.5\textwidth}
		\centering
		\includegraphics[width=\textwidth]{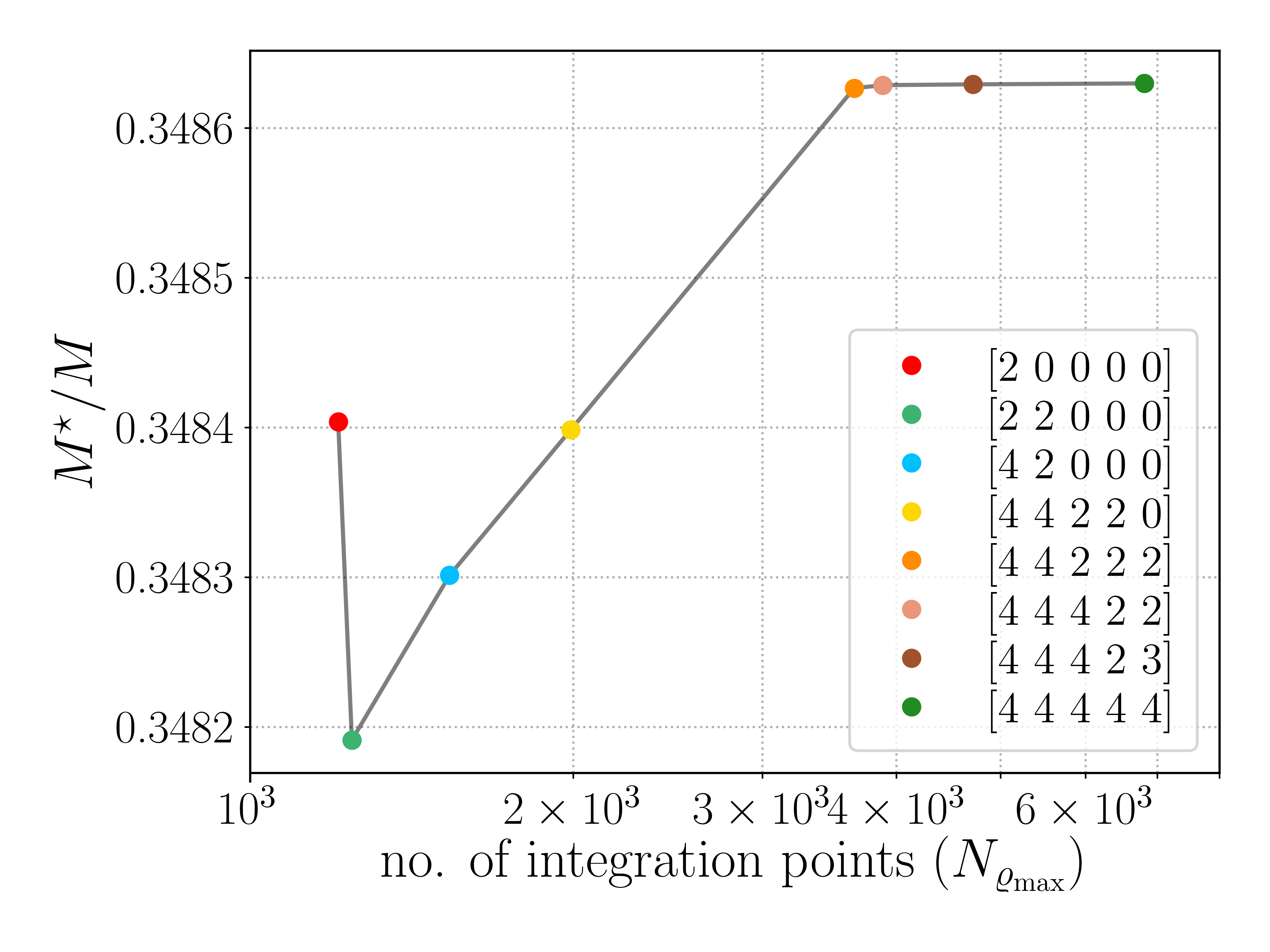}
		\caption{$d=2$}
		\label{fig:linear_elasticity_2d}
	\end{subfigure}%
	\begin{subfigure}[b]{0.5\textwidth}
		\centering
		\includegraphics[width=\textwidth]{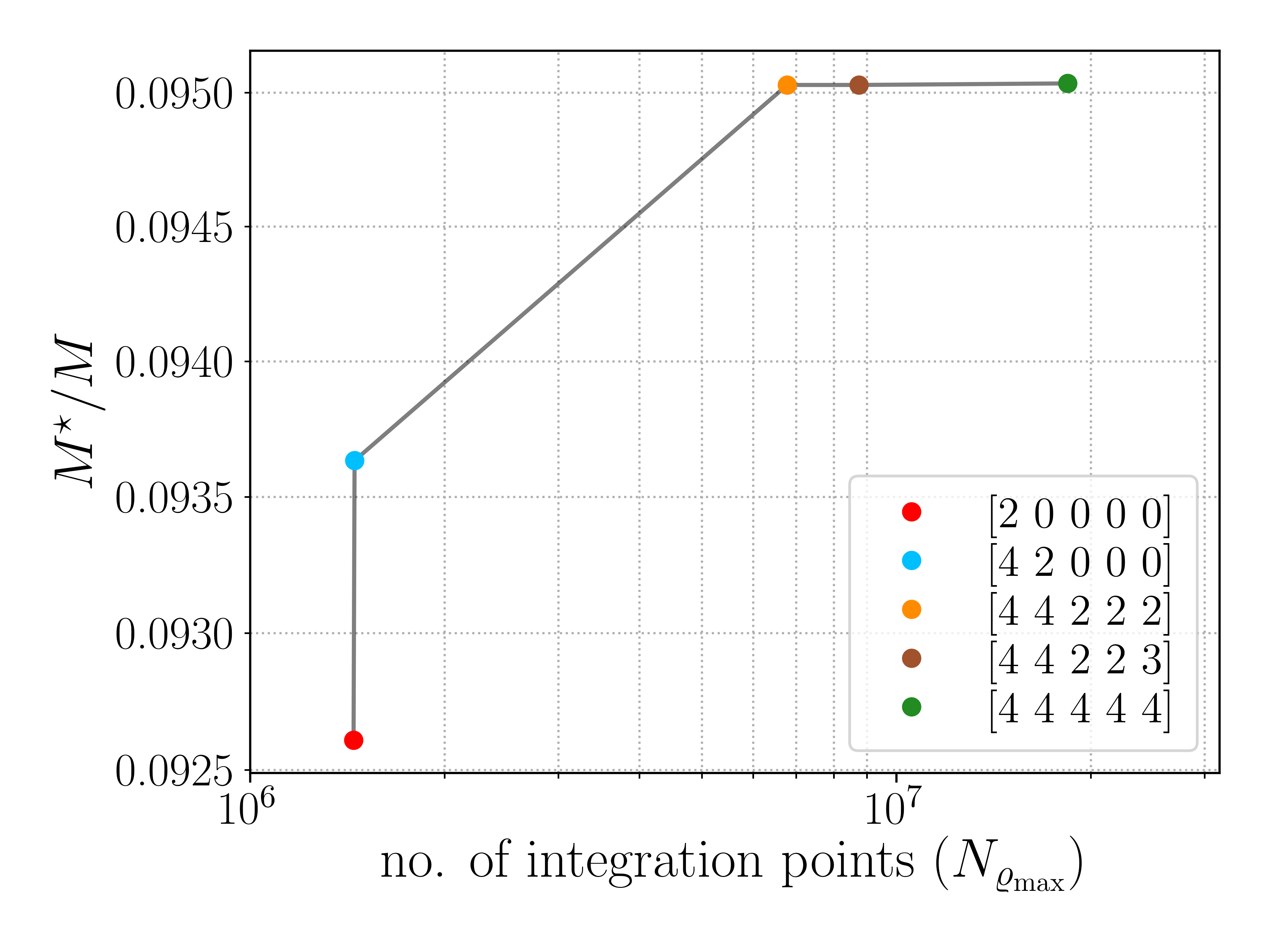}
		\caption{$d=3$}
		\label{fig:linear_elasticity_3d}
	\end{subfigure}
	\caption{Comparison of the constrained modulus for different quadrature rules that have evolved from optimization algorithm for both the two-dimensional and three-dimensional linear elasticity problem.}
	\label{fig:linear_elasticity}
\end{figure}

Figure~\ref{fig:linear_elasticity_2d} displays the computed constrained modulus for different steps in the integration optimization procedure. It is observed that by selecting a second order Gauss scheme on the untrimmed elements, and a single point in each sub-cell of the trimmed elements, yields a modulus of approximately $0.3484 M$, which is already within 0.1\% of the fully resolved integral result. This computation, however, uses a total of only $1208$ points, in comparison to the $6810$ points used for the full integral computation. Upon application of the integration optimization procedure the constrained modulus evidently converges to that using full integration. A virtually identical result as the full integration result is obtained with $3654$ points. For this result fourth order Gauss quadrature is used on the untrimmed cells and on the first bisection level, and second order Gauss quadrature on all lower levels.

\subsection{Three-dimensional test case} \label{sec:immersed3D}
\begin{figure}
	\centering
	\begin{subfigure}[b]{0.5\textwidth}
		\centering
		\includegraphics[width=\textwidth]{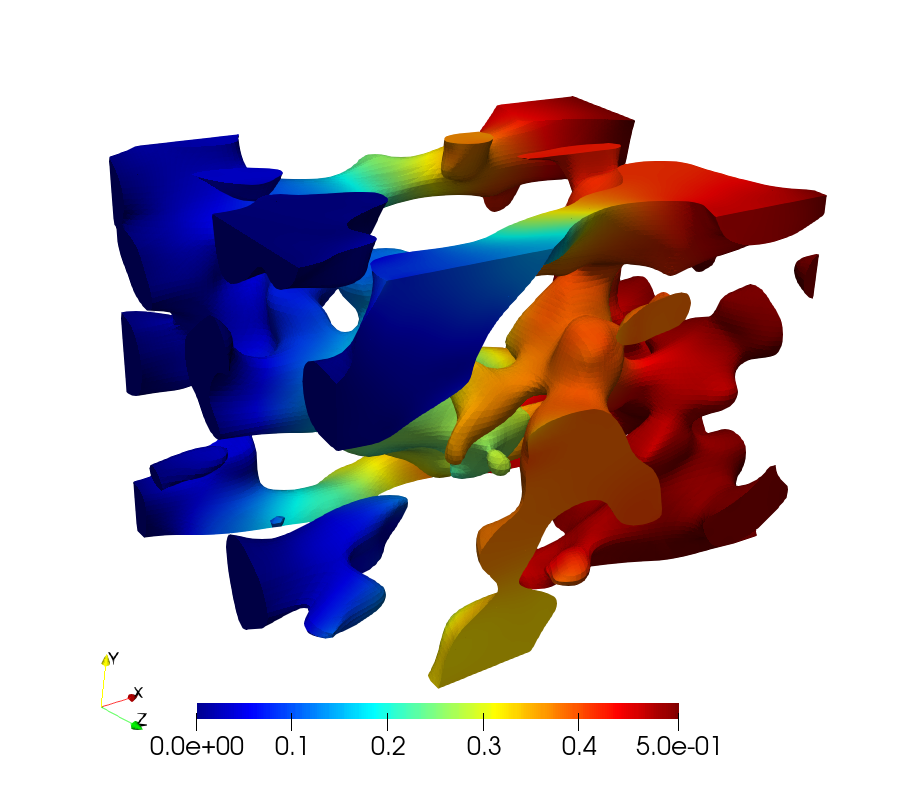}
		\caption{Displacement magnitude}
		\label{fig:displacement_3d}
	\end{subfigure}%
	\begin{subfigure}[b]{0.5\textwidth}
		\centering
		\includegraphics[width=\textwidth]{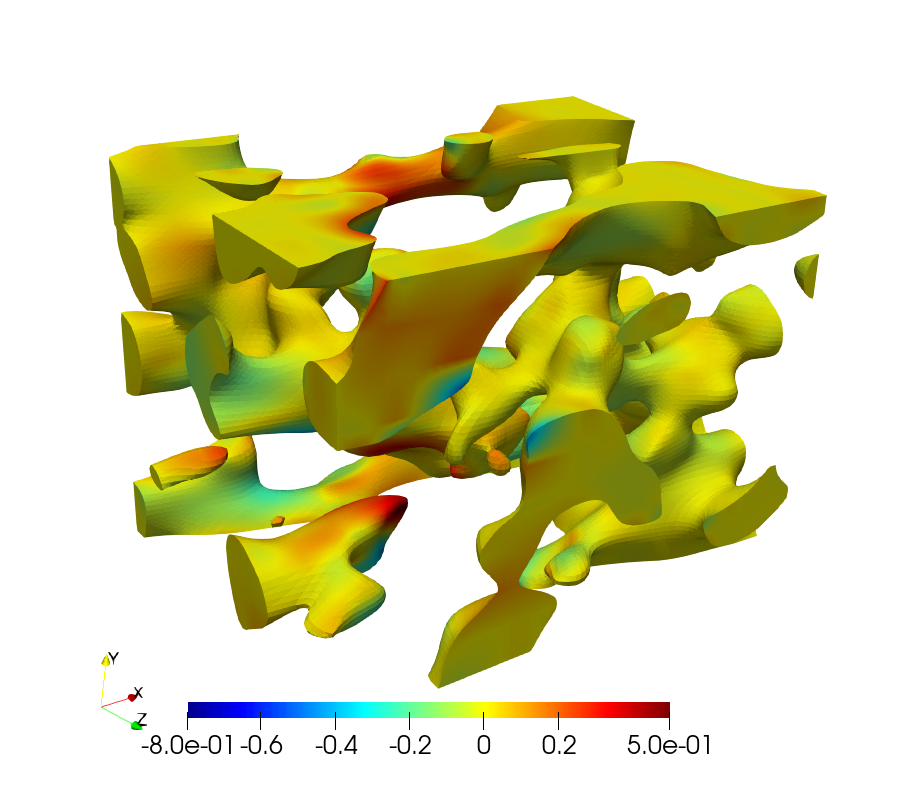}
		\caption{Shear stress $\sigma_{xy}$}
		\label{fig:shear_3d}
	\end{subfigure}
	\caption{Solution of the three-dimensional linear elasticity problem, displayed in the deformed configuration.}
	\label{fig:le_results_3d}
\end{figure}

We consider the elastic deformation of a sintered glass specimen, which is extracted from $\mu$CT-scan data as discussed in Ref.~\cite{hoang2019}\footnote{The scan data used for this simulation can be downloaded from: \newline \href{www.gitlab.tue.nl/20175645/sintered\_scan\_data}{www.gitlab.tue.nl/20175645/sintered\_scan\_data}.}. The size of the bounding box is set to $L=1.5$, and the right boundary is displaced by $\bar{u} = \frac{1}{2}$. Figure~\ref{fig:le_results_3d} displays the solution to the elasticity problem \eqref{eq:linear_elasticiy} in terms of the displacement magnitude and the $\sigma_{xy}$ shear stress. The displayed result is based on a $10 \times 10 \times 10$ ambient domain mesh and uniformly distributed fourth order Gauss integration scheme. As reported in Ref.~\cite{verhoosel2015}, for coarse meshes the basis functions artificially connect disjoint parts of the geometry. This effect is also observed in the shear stress displayed in Figure~\ref{fig:shear_3d}. For the purpose of illustrating the performance of the integration optimization procedure considered in this work, this effect is, however, not prohibitive.

Figure~\ref{fig:linear_elasticity_3d} displays the constrained modulus versus the number of integration points. By using second order quadrature on the untrimmed cells and one point per sub-cell in the trimmed elements, provides a modulus of approximately $0.0925M$, which is within 2.5\% of the fully resolved integral result. However, this result is obtained at the expense of $1.4 \times 10^{6}$ integration points. This is a factor thirteen lower than the $18.4 \times 10^{6}$ points used to compute the fully resolved integral result for a reference. A virtually identical result to the reference result is obtained at $6.8\times 10^6$ points, which is still approximately a factor of three lower than the full integration scheme.

\section{Conclusion}\label{sec:conclusion}
We have developed an algorithm to construct quadrature rules for cut-elements in which the integration points are distributed over the (cut) element in such a way that the integration error is minimized. Strang's first lemma -- which provides an error estimate for the approximate solution including integration errors -- provides a theoretical underpinning of the developed integration procedure. In the setting of this approximation theory, it is found that the integration error is bound by the supremum of integration errors over the considered test space. When the integration error is localized to the (cut) elements, this supremum translates to the consideration of a worst possible function to be integrated over a particular element.  By discretizing the space of functions over an element by means of polynomials, the worst possible function and its associated integration error become computable.

The ability to evaluate the integration error forms the basis of the proposed integration optimization procedure for octree subdivision. The pivotal idea behind the proposed algorithm is to gradually increase the number of integration points by adding integration points to the sub-cells for which the error reduction per added integration point is highest. As the number of sub-cells increases dramatically with an increase in number of dimensions and with the octree depth -- as illustrated by the derived scaling relation for the number of sub-cells -- such an optimization procedure has the potential to significantly reduce the computational effort involved in the integration procedures for immersed methods. The presented numerical simulations demonstrate that, for a given number of integration points, the integration error resulting from the optimization algorithm is, in general, significantly lower than that of the integration scheme considering the same integration orders on all sub-cells. Conversely, when fixing the error, the number of integration points required for the optimized quadrature rule is significantly lower than that of the equal order schemes. The considered simulations demonstrate that the developed optimization algorithm efficiently distributes integration points over cut-elements for a wide range of cut-cell configurations. 

The developed algorithm provides insight into the way in which integration points should be distributed over the cut-elements in order to minimize integration errors. Based on these insights, it is evident that integration rules for which the number of points per sub-cell is decreased with increasing octree depth form a good approximation to the optimized distribution of integration points. The integration point distribution algorithm presented in Ref.~\cite{alireza2013} therefore can also be expected to yield quadrature rules that are a good approximation to the rules that follow from the theoretical framework considered in this contribution.

Although the costs involved in the quadrature optimization algorithm are limited on account of the fact that the employed polynomial basis should only be integrated once using the expensive full order integration scheme, the determined optimized integration rules are particularly attractive when they can be re-used for multiple simulations. This is, for example, the case in time-dependent or non-linear problems. The attained reduction in number of integration points is also highly advantageous when data is to be stored in integration points, such as is the case for example with history data in many non-linear constitutive models.

The developed algorithm is highly generic in the sense that one can use the algorithm to optimize other integration schemes. For example, nested Gauss schemes can be attractive on account of the fact that the points (or locations) are nested in these schemes, which facilitates the re-using of configurations between integration schemes.

When employing the quadrature optimization algorithm one must choose various parameters, most importantly the order of the polynomial integrand to be considered and the norm used for the integrand normalization. In particular the order of polynomials is important with respect to the obtained distribution of integration points, in the sense that selecting this order too low will lead to suboptimal integration results. In general it is possible to determine the order of the integrands based on the problem under consideration and the finite element basis functions being considered.

In this contribution we have restricted ourselves to optimizing the distribution of integration points over the volumes of the cut-elements. The underlying theoretical framework is, however, more general in the sense that other parameters controlling the integration error can be incorporated as well. In the context of the considered octree integration procedure it would be natural to include the octree depth in the optimization procedure. Since both the error reduction associated with the increase in octree depth and the associated computational cost are quantifiable, including this depth as an additional parameter in the optimization algorithm is anticipated to be feasible. Inclusion of this parameter is, however, considered beyond the scope of the current work. The same holds for the incorporation of errors associated with the boundary integrals, in particular the Nitsche terms typically considered in finite cell simulations.

The focus in this work has been on the optimization of the integration error contribution in Strang's first lemma, not taking into account the approximation error. It is noted, however, that the need to optimize the integration error depends on the approximation error. If the approximation error is small compared to the overall error, then there is a need to optimize the integration error. However, if the approximation error is the dominating contribution to the overall error, then there is not a strong need to optimize the integration error. Evaluation of the balance between the approximation error and the integration error is an important topic of further study.

	
\section*{Acknowledgement}
We acknowledge the support from the European Commission EACEA Agency, Framework Partnership Agreement 2013-0043 Erasmus Mundus Action 1b, as a part of the EM Joint Doctorate Simulation in Engineering and Entrepreneurship Development (SEED). All the simulations in this work were performed based on the open source software package Nutils (www.nutils.org) \cite{nutils}. We acknowledge the support of the Nutils team.

	\bibliography{Bibliography/integration}

\begin{thebibliography}{10}
\expandafter\ifx\csname url\endcsname\relax
  \def\url#1{\texttt{#1}}\fi
\expandafter\ifx\csname urlprefix\endcsname\relax\def\urlprefix{URL }\fi
\expandafter\ifx\csname href\endcsname\relax
  \def\href#1#2{#2} \def\path#1{#1}\fi

\bibitem{rank2008}
A.~D{\"u}ster, J.~Parvizian, Z.~Yang, E.~Rank, The finite cell method for
  three-dimensional problems of solid mechanics, Computer Methods in Applied
  Mechanics and Engineering 197~(45-48) (2008) 3768--3782.

\bibitem{hansbo2002}
A.~Hansbo, P.~Hansbo, An unfitted finite element method, based on nitsche’s
  method, for elliptic interface problems, Computer Methods in Applied
  Mechanics and Engineering 191~(47-48) (2002) 5537--5552.

\bibitem{schillinger2011}
D.~Schillinger, E.~Rank, An unfitted hp-adaptive finite element method based on
  hierarchical b-splines for interface problems of complex geometry, Computer
  Methods in Applied Mechanics and Engineering 200~(47-48) (2011) 3358--3380.

\bibitem{rank2012}
E.~Rank, M.~Ruess, S.~Kollmannsberger, D.~Schillinger, A.~D{\"u}ster, Geometric
  modeling, isogeometric analysis and the finite cell method, Computer Methods
  in Applied Mechanics and Engineering 249 (2012) 104--115.

\bibitem{schillinger2012}
D.~Schillinger, L.~Dede, M.~A. Scott, J.~A. Evans, M.~J. Borden, E.~Rank, T.~J.
  Hughes, An isogeometric design-through-analysis methodology based on adaptive
  hierarchical refinement of {NURBS}, immersed boundary methods, and {T}-spline
  {CAD} surfaces, Computer Methods in Applied Mechanics and Engineering 249
  (2012) 116--150.

\bibitem{hughes2005}
T.~J. Hughes, J.~A. Cottrell, Y.~Bazilevs, Isogeometric analysis: {CAD}, finite
  elements, {NURBS}, exact geometry and mesh refinement, Computer Methods in
  Applied Mechanics and Engineering 194~(39-41) (2005) 4135--4195.

\bibitem{verhoosel2015}
C.~V. Verhoosel, G.~Van~Zwieten, B.~Van~Rietbergen, R.~de~Borst, Image-based
  goal-oriented adaptive isogeometric analysis with application to the
  micro-mechanical modeling of trabecular bone, Computer Methods in Applied
  Mechanics and Engineering 284 (2015) 138--164.

\bibitem{schmidt2012}
R.~Schmidt, R.~W{\"u}chner, K.~U. Bletzinger, Isogeometric analysis of trimmed
  {NURBS} geometries, Computer Methods in Applied Mechanics and Engineering 241
  (2012) 93--111.

\bibitem{ruess2013}
M.~Ruess, D.~Schillinger, Y.~Bazilevs, V.~Varduhn, E.~Rank, Weakly enforced
  essential boundary conditions for {NURBS}-embedded and trimmed {NURBS}
  geometries on the basis of the finite cell method, International Journal for
  Numerical Methods in Engineering 95~(10) (2013) 811--846.

\bibitem{ruess2014}
M.~Ruess, D.~Schillinger, A.~I. Oezcan, E.~Rank, Weak coupling for isogeometric
  analysis of non-matching and trimmed multi-patch geometries, Computer Methods
  in Applied Mechanics and Engineering 269 (2014) 46--71.

\bibitem{guo2015}
Y.~Guo, M.~Ruess, Weak dirichlet boundary conditions for trimmed thin
  isogeometric shells, Computers \& Mathematics with Applications 70~(7) (2015)
  1425--1440.

\bibitem{guo2017}
Y.~Guo, M.~Ruess, D.~Schillinger, A parameter-free variational coupling
  approach for trimmed isogeometric thin shells, computational Mechanics 59~(4)
  (2017) 693--715.

\bibitem{marussig2018}
B.~Marussig, T.~J. Hughes, A review of trimming in isogeometric analysis:
  challenges, data exchange and simulation aspects, Archives of computational
  methods in engineering 25~(4) (2018) 1059--1127.

\bibitem{bauer2017}
A.~Bauer, M.~Breitenberger, B.~Philipp, R.~W{\"u}chner, K.~U. Bletzinger,
  Embedded structural entities in {NURBS}-based isogeometric analysis, Computer
  Methods in Applied Mechanics and Engineering 325 (2017) 198--218.

\bibitem{beer2015}
G.~Beer, B.~Marussig, J.~Zechner, A simple approach to the numerical simulation
  with trimmed {CAD} surfaces, Computer Methods in Applied Mechanics and
  Engineering 285 (2015) 776--790.

\bibitem{schillinger2015}
D.~Schillinger, M.~Ruess, The finite cell method: A review in the context of
  higher-order structural analysis of {CAD} and image-based geometric models,
  Archives of Computational Methods in Engineering 22~(3) (2015) 391--455.

\bibitem{duster2017}
A.~D{\"u}ster, E.~Rank, B.~Szab{\'o}, The p-version of the finite element and
  finite cell methods, Encyclopedia of computational Mechanics Second Edition
  (2017) 1--35.

\bibitem{kamensky2015}
D.~Kamensky, M.~C. Hsu, D.~Schillinger, J.~A. Evans, A.~Aggarwal, Y.~Bazilevs,
  M.~S. Sacks, T.~J. Hughes, An immersogeometric variational framework for
  fluid--structure interaction: Application to bioprosthetic heart valves,
  Computer Methods in Applied Mechanics and Engineering 284 (2015) 1005--1053.

\bibitem{hsu2015}
M.~C. Hsu, D.~Kamensky, F.~Xu, J.~Kiendl, C.~Wang, M.~C. Wu, J.~Mineroff,
  A.~Reali, Y.~Bazilevs, M.~S. Sacks, Dynamic and fluid--structure interaction
  simulations of bioprosthetic heart valves using parametric design with
  {T}-splines and fung-type material models, computational Mechanics 55~(6)
  (2015) 1211--1225.

\bibitem{xu2018}
F.~Xu, S.~Morganti, R.~Zakerzadeh, D.~Kamensky, F.~Auricchio, A.~Reali, T.~J.
  Hughes, M.~S. Sacks, M.~C. Hsu, A framework for designing patient-specific
  bioprosthetic heart valves using immersogeometric fluid--structure
  interaction analysis, International journal for numerical methods in
  biomedical engineering 34~(4) (2018) e2938.

\bibitem{ruess2012}
M.~Ruess, D.~Tal, N.~Trabelsi, Z.~Yosibash, E.~Rank, The finite cell method for
  bone simulations: verification and validation, Biomechanics and modeling in
  mechanobiology 11~(3-4) (2012) 425--437.

\bibitem{belytschko2009}
T.~Belytschko, R.~Gracie, G.~Ventura, A review of extended/generalized finite
  element methods for material modeling, Modelling and Simulation in Materials
  Science and Engineering 17~(4) (2009) 043001.

\bibitem{duster2008}
A.~D{\"u}ster, J.~Parvizian, Z.~Yang, E.~Rank, The finite cell method for
  three-dimensional problems of solid mechanics, Computer Methods in Applied
  Mechanics and Engineering 197~(45-48) (2008) 3768--3782.

\bibitem{varduhn2016}
V.~Varduhn, M.~C. Hsu, M.~Ruess, D.~Schillinger, The tetrahedral finite cell
  method: higher-order immersogeometric analysis on adaptive
  non-boundary-fitted meshes, International Journal for Numerical Methods in
  Engineering 107~(12) (2016) 1054--1079.

\bibitem{stavrev2016}
A.~Stavrev, L.~H. Nguyen, R.~Shen, V.~Varduhn, M.~Behr, S.~Elgeti,
  D.~Schillinger, Geometrically accurate, efficient, and flexible quadrature
  techniques for the tetrahedral finite cell method, Computer Methods in
  Applied Mechanics and Engineering 310 (2016) 646--673.

\bibitem{alireza2013}
A.~Abedian, J.~Parvizian, A.~D{\"u}ster, E.~Rank, The finite cell method for
  the {J2} flow theory of plasticity, Finite Elements in Analysis and Design 69
  (2013) 37--47.

\bibitem{fries2010}
K.~W. Cheng, T.~P. Fries, Higher-order {XFEM} for curved strong and weak
  discontinuities, International Journal for Numerical Methods in Engineering
  82~(5) (2010) 564--590.

\bibitem{fries2016}
T.~P. Fries, S.~Omerovi{\'c}, Higher-order accurate integration of implicit
  geometries, International Journal for Numerical Methods in Engineering
  106~(5) (2016) 323--371.

\bibitem{sevilla2008}
R.~Sevilla, S.~Fern{\'a}ndez~M{\'e}ndez, A.~Huerta, {NURBS}-{E}nhanced {F}inite
  {E}lement {M}ethod ({NEFEM}), International Journal for Numerical Methods in
  Engineering 76~(1) (2008) 56--83.

\bibitem{sevilla2011}
R.~Sevilla, S.~Fern{\'a}ndez~M{\'e}ndez, A.~Huerta, Comparison of high-order
  curved finite elements, International Journal for Numerical Methods in
  Engineering 87~(8) (2011) 719--734.

\bibitem{nadal2013}
E.~Nadal, J.~R{\'o}denas, J.~Albelda, M.~Tur, J.~Taranc{\'o}n, F.~Fuenmayor,
  Efficient finite element methodology based on cartesian grids: application to
  structural shape optimization, in: Abstract and applied analysis, Vol. 2013,
  Hindawi, 2013.

\bibitem{kudela2016}
L.~Kudela, N.~Zander, S.~Kollmannsberger, E.~Rank, Smart octrees: Accurately
  integrating discontinuous functions in 3{D}, Computer Methods in Applied
  Mechanics and Engineering 306 (2016) 406--426.

\bibitem{hubrich2017}
S.~Hubrich, P.~Di~Stolfo, L.~Kudela, S.~Kollmannsberger, E.~Rank,
  A.~Schr{\"o}der, A.~D{\"u}ster, Numerical integration of discontinuous
  functions: moment fitting and smart octree, computational Mechanics 60~(5)
  (2017) 863--881.

\bibitem{mousavi2011}
S.~Mousavi, N.~Sukumar, Numerical integration of polynomials and discontinuous
  functions on irregular convex polygons and polyhedrons, computational
  Mechanics 47~(5) (2011) 535--554.

\bibitem{joulaian2016}
M.~Joulaian, S.~Hubrich, A.~D{\"u}ster, Numerical integration of
  discontinuities on arbitrary domains based on moment fitting, computational
  Mechanics 57~(6) (2016) 979--999.

\bibitem{hubrich2019}
S.~Hubrich, A.~D{\"u}ster, Numerical integration for nonlinear problems of the
  finite cell method using an adaptive scheme based on moment fitting,
  Computers \& Mathematics with Applications 77~(7) (2019) 1983--1997.

\bibitem{ventura2006}
G.~Ventura, On the elimination of quadrature subcells for discontinuous
  functions in the extended finite-element method, International Journal for
  Numerical Methods in Engineering 66~(5) (2006) 761--795.

\bibitem{alireza2019}
A.~Abedian, A.~D{\"u}ster, Equivalent legendre polynomials: Numerical
  integration of discontinuous functions in the finite element methods,
  Computer Methods in Applied Mechanics and Engineering 343 (2019) 690--720.

\bibitem{ventura2009}
G.~Ventura, R.~Gracie, T.~Belytschko, Fast integration and weight function
  blending in the extended finite element method, International journal for
  numerical methods in engineering 77~(1) (2009) 1--29.

\bibitem{jonsson2017}
T.~Jonsson, M.~G. Larson, K.~Larsson, Cut finite element methods for elliptic
  problems on multipatch parametric surfaces, Computer Methods in Applied
  Mechanics and Engineering 324 (2017) 366--394.

\bibitem{strang1973}
G.~Strang, G.~J. Fix, An analysis of the finite element method, Vol. 212,
  Prentice-hall Englewood Cliffs, NJ, 1973.

\bibitem{strang1972}
G.~Strang, Variational crimes in the finite element method, in: The
  mathematical foundations of the finite element method with applications to
  partial differential equations, Elsevier, 1972, pp. 689--710.

\bibitem{ern2013}
A.~Ern, J.~L. Guermond, Theory and practice of finite elements, Vol. 159,
  Springer Science \& Business Media, 2013.

\bibitem{burman2016}
E.~Burman, P.~Hansbo, M.~G. Larson, S.~Zahedi, Cut finite element methods for
  coupled bulk--surface problems, Numerische Mathematik 133~(2) (2016)
  203--231.

\bibitem{burman2017}
E.~Burman, P.~Hansbo, M.~G. Larson, A.~Massing, A cut discontinuous galerkin
  method for the laplace--beltrami operator, IMA Journal of Numerical Analysis
  37~(1) (2017) 138--169.

\bibitem{zienkiewicz1971}
O.~Zienkiewicz, R.~Taylor, J.~Too, Reduced integration technique in general
  analysis of plates and shells, International Journal for Numerical Methods in
  Engineering 3~(2) (1971) 275--290.

\bibitem{hughes1978reduced}
T.~J. Hughes, M.~Cohen, M.~Haroun, Reduced and selective integration techniques
  in the finite element analysis of plates, Nuclear Engineering and design
  46~(1) (1978) 203--222.

\bibitem{hughes1978mixed}
D.~S. Malkus, T.~J. Hughes, Mixed finite element methods—reduced and
  selective integration techniques: a unification of concepts, Computer Methods
  in Applied Mechanics and Engineering 15~(1) (1978) 63--81.

\bibitem{taghipour2018}
A.~Taghipour, J.~Parvizian, S.~Heinze, A.~D{\"u}ster, The finite cell method
  for nearly incompressible finite strain plasticity problems with complex
  geometries, Computers \& Mathematics with Applications 75~(9) (2018)
  3298--3316.

\bibitem{piegl2012}
L.~Piegl, W.~Tiller, The {NURBS} book, Springer Science \& Business Media,
  2012.

\bibitem{nitsche1971}
J.~Nitsche, {\"U}ber ein variationsprinzip zur l{\"o}sung von
  dirichlet-problemen bei verwendung von teilr{\"a}umen, die keinen
  randbedingungen unterworfen sind, in: Abhandlungen aus dem mathematischen
  Seminar der Universit{\"a}t Hamburg, Vol.~36, Springer, 1971, pp. 9--15.

\bibitem{nutils}
G.~van Zwieten, J.~van Zwieten, C.~V. Verhoosel, E.~Fonn, T.~van Opstal,
  W.~Hoitinga, \href{http://doi.org/10.5281/zenodo.3243447}{Nutils (version
  5.0)} (2019).
\newline\urlprefix\url{http://doi.org/10.5281/zenodo.3243447}

\bibitem{deprenter2019}
F.~de~Prenter, C.~V. Verhoosel, E.~H. van Brummelen, Preconditioning immersed
  isogeometric finite element methods with application to flow problems,
  Computer Methods in Applied Mechanics and Engineering 348 (2019) 604--631.

\bibitem{hoang2019}
T.~Hoang, C.~V. Verhoosel, C.~Z. Qin, F.~Auricchio, A.~Reali, E.~H. van
  Brummelen, Skeleton-stabilized immersogeometric analysis for incompressible
  viscous flow problems, Computer Methods in Applied Mechanics and Engineering
  344 (2019) 421--450.

\end{thebibliography}
	%
	\appendix
	\section{Midpoint tessellation procedure}\label{sec:tessellation}
As discussed in Section~\ref{sec:octree_partitioning}, at the lowest level of bisectioning we apply a tessellation procedure that provides an order $\mathcal{O}(h^2/2^{2\maxlevel})$ approximation of the interior volume \cite{verhoosel2015}, with $h$ the size of the background element in which the sub-cell resides and with $\maxlevel$ the number of octree bisections. The considered procedure provides an explicit parametrization of both the interior of the trimmed sub-cell and its immersed boundary.

\begin{figure}
	\centering
	\begin{subfigure}{0.32\textwidth}
		\centering
		\begin{tikzpicture}[
    scale=3 ,
    line cap=round,
    line join=round,
    opaque/.style={black,thin,opacity=0.4}
    ]
    
    \coordinate (A) at (0,0,0);
    \coordinate (B) at (1,0,0);
    \coordinate (C) at (1,1,0);
    \coordinate (D) at (0,1,0);
 
    \draw[thick] (A)--(B)--(C)--(D)--cycle;
    
	\coordinate (CD) at ($0.4*(C)+0.6*(D)$);
	\coordinate (AB) at ($0.7*(A)+0.3*(B)$);
	\coordinate (centroid) at (0.5,0.5,0);
	
	\filldraw[draw=green!50!black,very thick,fill=green!20!white] (A) -- (AB) .. controls ($0.3*(AB)+0.3*(CD)$) .. (CD) -- (D) -- cycle;
	\filldraw[draw=red!50!black,very thick,fill=red!20!white] (B) -- (AB) .. controls ($0.3*(AB)+0.3*(CD)$) .. (CD) -- (C) -- cycle;
	\draw[draw=blue!50!black,very thick] (AB) .. controls ($0.3*(AB)+0.3*(CD)$) .. (CD);
    
\end{tikzpicture}
		\caption{}
		\label{fig:cut_cell_2d}
	\end{subfigure}%
	\begin{subfigure}{0.32\textwidth}
		\centering
		\begin{tikzpicture}[
    scale=3 ,
    line cap=round,
    line join=round,
    opaque/.style={black,thin,opacity=0.4}
    ]
    
    \coordinate (A) at (0,0,0);
    \coordinate (B) at (1,0,0);
    \coordinate (C) at (1,1,0);
    \coordinate (D) at (0,1,0);
 
    \draw[thick] (A)--(B)--(C)--(D)--cycle;
    
	\coordinate (CD) at ($0.4*(C)+0.6*(D)$);
	\coordinate (AB) at ($0.7*(A)+0.3*(B)$);
	\coordinate (centroid) at (0.5,0.5,0);
    
    \foreach \c in {(A),(D)}
    \node at \c [circle,draw,fill=white,inner sep=0pt,line width=0.5pt,fill=green!30!white] {$+$};
    
    \foreach \c in {(B),(C)}
    \node at \c [circle,draw,fill=white,inner sep=0pt,line width=0.5pt,fill=red!30!white] {$-$};
    
    \node at (centroid) [circle,draw,fill=white,inner sep=0pt,line width=0.5pt,fill=red!30!white] {$-$};

    \node at (centroid) [circle,draw,fill=white,inner sep=0pt,line width=0.5pt,fill=red!30!white] {$-$};
    
\end{tikzpicture}		
		\caption{}
		\label{fig:levelset_vertex_2d}
	\end{subfigure}%
	\begin{subfigure}{0.32\textwidth}
		\centering
		\begin{tikzpicture}[
    scale=3 ,
    line cap=round,
    line join=round,
    opaque/.style={black,thin,opacity=0.4}
    ]
    
    \coordinate (A) at (0,0,0);
    \coordinate (B) at (1,0,0);
    \coordinate (C) at (1,1,0);
    \coordinate (D) at (0,1,0);
 
    \draw[thick] (A)--(B)--(C)--(D)--cycle;
    
    \coordinate (CD) at ($0.4*(C)+0.6*(D)$);
	\coordinate (AB) at ($0.7*(A)+0.3*(B)$);
	\coordinate (centroid) at (0.5,0.5,0);
	
    \draw[draw=green!50!black,very thick,fill=green!20!white] (A) -- (AB);
    \draw[draw=green!50!black,very thick,fill=green!20!white] (D) -- (CD);
    \draw[draw=green!50!black,very thick] (A)--(D);
    
    \foreach \c in {(A),(D)}
	\node at \c [circle,draw,fill=white,inner sep=0pt,line width=0.5pt,fill=green!30!white] {$+$};

	\foreach \c in {(B),(C)}
	\node at \c [circle,draw,fill=white,inner sep=0pt,line width=0.5pt,fill=red!30!white] {$-$};
    
    \node at (AB) [circle,draw,fill=white,inner sep=1pt,line width=0.5pt,fill=blue!30!white] {$0$};
    \node at (CD) [circle,draw,fill=white,inner sep=1pt,line width=0.5pt,fill=blue!30!white] {$0$};
    \node at (centroid) [circle,draw,fill=white,inner sep=0pt,line width=0.5pt,fill=red!30!white] {$-$};
    
\end{tikzpicture}
		\caption{}
		\label{fig:trim_edges_2d}
	\end{subfigure}\\[12pt]
	\begin{subfigure}{0.32\textwidth}
		\centering
		\begin{tikzpicture}[
scale=3 ,
line cap=round,
line join=round,
opaque/.style={black,thin,opacity=0.4}
]

\coordinate (A) at (0,0,0);
\coordinate (B) at (1,0,0);
\coordinate (C) at (1,1,0);
\coordinate (D) at (0,1,0);

\draw[thick] (A)--(B)--(C)--(D)--cycle;

\coordinate (CD) at ($0.4*(C)+0.6*(D)$);
\coordinate (AB) at ($0.7*(A)+0.3*(B)$);
\coordinate (centroid) at (0.5,0.5,0);

\coordinate (mAc) at ($0.5*(A)+ 0.5*(centroid)$);
\coordinate (mDc) at ($0.4*(D)+ 0.6*(centroid)$);

\draw[draw=green!50!black,very thick,fill=green!20!white] (A) -- (AB);
\draw[draw=green!50!black,very thick,fill=green!20!white] (D) -- (CD);
\draw[draw=green!50!black,very thick] (A)--(D);
\draw[draw=yellow!80!black,very thick,fill=blue!20!white] (A) -- (centroid);
\draw[draw=yellow!80!black,very thick,fill=blue!20!white] (D) -- (centroid);

\foreach \c in {(A),(D)}
\node at \c [circle,draw,fill=white,inner sep=0pt,line width=0.5pt,fill=green!30!white] {$+$};

\foreach \c in {(B),(C)}
\node at \c [circle,draw,fill=white,inner sep=0pt,line width=0.5pt,fill=red!30!white] {$-$};

\node at (AB) [circle,draw,fill=white,inner sep=1pt,line width=0.5pt,fill=blue!30!white] {$0$};
\node at (CD) [circle,draw,fill=white,inner sep=1pt,line width=0.5pt,fill=blue!30!white] {$0$};
\node at (centroid) [circle,draw,fill=white,inner sep=0pt,line width=0.5pt,fill=red!30!white] {$-$};
\node at (centroid) [circle,draw,fill=white,inner sep=0pt,line width=0.5pt,fill=red!30!white] {$-$};
\node at (mAc) [rectangle,draw,fill=white,inner sep=1pt,line width=0.5pt,fill=blue!30!white] {$0$};
\node at (mDc) [rectangle,draw,fill=white,inner sep=1pt,line width=0.5pt,fill=blue!30!white] {$0$};

\end{tikzpicture}
		\caption{}
		\label{fig:levelset_diagonal_2d}
	\end{subfigure}%
	\begin{subfigure}{0.32\textwidth}
		\centering
		\begin{tikzpicture}[
scale=3 ,
line cap=round,
line join=round,
opaque/.style={black,thin,opacity=0.4}
]

\coordinate (A) at (0,0,0);
\coordinate (B) at (1,0,0);
\coordinate (C) at (1,1,0);
\coordinate (D) at (0,1,0);

\draw[thick] (A)--(B)--(C)--(D)--cycle;

\coordinate (CD) at ($0.4*(C)+0.6*(D)$);
\coordinate (AB) at ($0.7*(A)+0.3*(B)$);
\coordinate (centroid) at (0.5,0.5,0);

\coordinate (mAc) at ($0.5*(A)+ 0.5*(centroid)$);
\coordinate (mDc) at ($0.4*(D)+ 0.6*(centroid)$);
\coordinate (midpoint) at ($0.5*(mAc)+0.5*(mDc)$);

\draw[draw=green!50!black,very thick,fill=green!20!white] (A) -- (AB);
\draw[draw=green!50!black,very thick,fill=green!20!white] (D) -- (CD);
\draw[draw=green!50!black,very thick] (A)--(D);
\draw[dotted,draw=yellow!50!black,very thick,fill=blue!20!white] (mAc) -- (mDc); 

\foreach \c in {(A),(D)}
\node at \c [circle,draw,fill=white,inner sep=0pt,line width=0.5pt,fill=green!30!white] {$+$};

\foreach \c in {(B),(C)}
\node at \c [circle,draw,fill=white,inner sep=0pt,line width=0.5pt,fill=red!30!white] {$-$};

\node at (AB) [circle,draw,fill=white,inner sep=1pt,line width=0.5pt,fill=blue!30!white] {$0$};
\node at (CD) [circle,draw,fill=white,inner sep=1pt,line width=0.5pt,fill=blue!30!white] {$0$};
\node at (centroid) [circle,draw,fill=white,inner sep=0pt,line width=0.5pt,fill=red!30!white] {$-$};
\node at (centroid) [circle,draw,fill=white,inner sep=0pt,line width=0.5pt,fill=red!30!white] {$-$};
\node at (mAc) [rectangle,draw,fill=white,inner sep=1pt,line width=0.5pt,fill=blue!30!white] {$0$};
\node at (mDc) [rectangle,draw,fill=white,inner sep=1pt,line width=0.5pt,fill=blue!30!white] {$0$};
\node at (midpoint) [circle,draw,fill=white,inner sep=1pt,line width=0.5pt,fill=blue!30!white] {$0$};
\end{tikzpicture}
		\caption{}
		\label{fig:diagonal_midpoint_2d}
	\end{subfigure}%
	\begin{subfigure}{0.32\textwidth}
		\centering
		\begin{tikzpicture}[
scale=3 ,
line cap=round,
line join=round,
opaque/.style={black,thin,opacity=0.4}
]

\coordinate (A) at (0,0,0);
\coordinate (B) at (1,0,0);
\coordinate (C) at (1,1,0);
\coordinate (D) at (0,1,0);

\draw[thick] (A)--(B)--(C)--(D)--cycle;

\coordinate (CD) at ($0.4*(C)+0.6*(D)$);
\coordinate (AB) at ($0.7*(A)+0.3*(B)$);
\coordinate (centroid) at (0.5,0.5,0);

\coordinate (mAc) at ($0.5*(A)+ 0.5*(centroid)$);
\coordinate (mDc) at ($0.4*(D)+ 0.6*(centroid)$);
\coordinate (midpoint) at ($0.5*(mAc)+0.5*(mDc)$);

\draw[draw=green!50!black,very thick,fill=green!20!white] (A) -- (AB);
\draw[draw=green!50!black,very thick,fill=green!20!white] (D) -- (CD);
\draw[draw=green!50!black,very thick] (A)--(D);
\draw[draw=blue!50!black,very thick,fill=blue!20!white] (AB) -- (midpoint);
\draw[draw=blue!50!black,very thick,fill=blue!20!white] (CD) -- (midpoint);

\foreach \c in {(A),(D)}
\node at \c [circle,draw,fill=white,inner sep=0pt,line width=0.5pt,fill=green!30!white] {$+$};

\foreach \c in {(B),(C)}
\node at \c [circle,draw,fill=white,inner sep=0pt,line width=0.5pt,fill=red!30!white] {$-$};

\node at (AB) [circle,draw,fill=white,inner sep=1pt,line width=0.5pt,fill=blue!30!white] {$0$};
\node at (CD) [circle,draw,fill=white,inner sep=1pt,line width=0.5pt,fill=blue!30!white] {$0$};

\node at (midpoint) [circle,draw,fill=white,inner sep=1pt,line width=0.5pt,fill=blue!30!white] {$0$};
\end{tikzpicture}
		\caption{}
		\label{fig:levelset_midpoint_2d}
	\end{subfigure}\\[12pt]
	\begin{subfigure}{0.32\textwidth}
		\centering
		\begin{tikzpicture}[
scale=3 ,
line cap=round,
line join=round,
opaque/.style={black,thin,opacity=0.4}
]

\coordinate (A) at (0,0,0);
\coordinate (B) at (1,0,0);
\coordinate (C) at (1,1,0);
\coordinate (D) at (0,1,0);

\draw[thick] (A)--(B)--(C)--(D)--cycle;

\coordinate (CD) at ($0.4*(C)+0.6*(D)$);
\coordinate (AB) at ($0.7*(A)+0.3*(B)$);
\coordinate (centroid) at (0.5,0.5,0);

\coordinate (mAc) at ($0.5*(A)+ 0.5*(centroid)$);
\coordinate (mDc) at ($0.4*(D)+ 0.6*(centroid)$);
\coordinate (midpoint) at ($0.5*(mAc)+0.5*(mDc)$);

\draw[draw=green!50!black,very thick,fill=green!20!white] (A) -- (AB);
\draw[draw=green!50!black,very thick,fill=green!20!white] (D) -- (CD);
\draw[draw=green!50!black,very thick] (A)--(D);
\draw[draw=blue!50!black,very thick,fill=blue!20!white] (AB) -- (midpoint);
\draw[draw=blue!50!black,very thick,fill=blue!20!white] (CD) -- (midpoint);
\filldraw[draw=green!50!black,very thick,fill=green!20!white] (A) -- (AB) -- (midpoint) -- cycle;

\foreach \c in {(A),(D)}
\node at \c [circle,draw,fill=white,inner sep=0pt,line width=0.5pt,fill=green!30!white] {$+$};

\foreach \c in {(B),(C)}
\node at \c [circle,draw,fill=white,inner sep=0pt,line width=0.5pt,fill=red!30!white] {$-$};

\node at (AB) [circle,draw,fill=white,inner sep=1pt,line width=0.5pt,fill=blue!30!white] {$0$};
\node at (CD) [circle,draw,fill=white,inner sep=1pt,line width=0.5pt,fill=blue!30!white] {$0$};

\node at (midpoint) [circle,draw,fill=white,inner sep=1pt,line width=0.5pt,fill=blue!30!white] {$0$};
\end{tikzpicture}
		\caption{}
		\label{fig:tessellation1_2d}
	\end{subfigure}%
	\begin{subfigure}{0.32\textwidth}
		\centering
		\begin{tikzpicture}[
scale=3 ,
line cap=round,
line join=round,
opaque/.style={black,thin,opacity=0.4}
]

\coordinate (A) at (0,0,0);
\coordinate (B) at (1,0,0);
\coordinate (C) at (1,1,0);
\coordinate (D) at (0,1,0);

\draw[thick] (A)--(B)--(C)--(D)--cycle;

\coordinate (CD) at ($0.4*(C)+0.6*(D)$);
\coordinate (AB) at ($0.7*(A)+0.3*(B)$);
\coordinate (centroid) at (0.5,0.5,0);

\coordinate (mAc) at ($0.5*(A)+ 0.5*(centroid)$);
\coordinate (mDc) at ($0.4*(D)+ 0.6*(centroid)$);
\coordinate (midpoint) at ($0.5*(mAc)+0.5*(mDc)$);

\draw[draw=green!50!black,very thick,fill=green!20!white] (A) -- (AB);
\draw[draw=green!50!black,very thick,fill=green!20!white] (D) -- (CD);
\draw[draw=green!50!black,very thick] (A)--(D);
\filldraw[draw=green!50!black,very thick,fill=green!20!white] (A) -- (AB) -- (midpoint) -- cycle;
\filldraw[draw=green!50!black,very thick,fill=green!20!white] (A) -- (D) -- (midpoint) -- cycle;
\filldraw[draw=green!50!black,very thick,fill=green!20!white] (D) -- (CD) -- (midpoint) -- cycle;
\draw[draw=blue!50!black,very thick,fill=blue!20!white] (AB) -- (midpoint);
\draw[draw=blue!50!black,very thick,fill=blue!20!white] (CD) -- (midpoint);

\foreach \c in {(A),(D)}
\node at \c [circle,draw,fill=white,inner sep=0pt,line width=0.5pt,fill=green!30!white] {$+$};

\foreach \c in {(B),(C)}
\node at \c [circle,draw,fill=white,inner sep=0pt,line width=0.5pt,fill=red!30!white] {$-$};

\node at (AB) [circle,draw,fill=white,inner sep=1pt,line width=0.5pt,fill=blue!30!white] {$0$};
\node at (CD) [circle,draw,fill=white,inner sep=1pt,line width=0.5pt,fill=blue!30!white] {$0$};

\node at (midpoint) [circle,draw,fill=white,inner sep=1pt,line width=0.5pt,fill=blue!30!white] {$0$};
\end{tikzpicture}
		\caption{}
		\label{fig:tessellation2_2d}
	\end{subfigure}%
	\begin{subfigure}{0.32\textwidth}
		\centering
		\begin{tikzpicture}[
scale=3 ,
line cap=round,
line join=round,
opaque/.style={black,thin,opacity=0.4}
]

\coordinate (A) at (0,0,0);
\coordinate (B) at (1,0,0);
\coordinate (C) at (1,1,0);
\coordinate (D) at (0,1,0);

\draw[thick] (A)--(B)--(C)--(D)--cycle;

\coordinate (CD) at ($0.4*(C)+0.6*(D)$);
\coordinate (AB) at ($0.7*(A)+0.3*(B)$);
\coordinate (centroid) at (0.5,0.5,0);

\coordinate (mAc) at ($0.5*(A)+ 0.5*(centroid)$);
\coordinate (mDc) at ($0.4*(D)+ 0.6*(centroid)$);
\coordinate (midpoint) at ($0.5*(mAc)+0.5*(mDc)$);

\draw[draw=green!50!black,very thick,fill=green!20!white] (A) -- (AB);
\draw[draw=green!50!black,very thick,fill=green!20!white] (D) -- (CD);
\draw[draw=green!50!black,very thick] (A)--(D);
\filldraw[draw=green!50!black,very thick,fill=green!20!white] (A) -- (AB) -- (midpoint) -- cycle;
\filldraw[draw=green!50!black,very thick,fill=green!20!white] (A) -- (D) -- (midpoint) -- cycle;
\filldraw[draw=green!50!black,very thick,fill=green!20!white] (D) -- (CD) -- (midpoint) -- cycle;
\filldraw[draw=red!50!black,very thick,fill=red!20!white] (B) -- (AB) -- (midpoint) -- cycle;
\filldraw[draw=red!50!black,very thick,fill=red!20!white] (B) -- (C) -- (midpoint) -- cycle;
\filldraw[draw=red!50!black,very thick,fill=red!20!white] (C) -- (CD) -- (midpoint) -- cycle;
\draw[draw=blue!50!black,very thick,fill=blue!20!white] (AB) -- (midpoint);
\draw[draw=blue!50!black,very thick,fill=blue!20!white] (CD) -- (midpoint);

\foreach \c in {(A),(D)}
\node at \c [circle,draw,fill=white,inner sep=0pt,line width=0.5pt,fill=green!30!white] {$+$};

\foreach \c in {(B),(C)}
\node at \c [circle,draw,fill=white,inner sep=0pt,line width=0.5pt,fill=red!30!white] {$-$};

\node at (AB) [circle,draw,fill=white,inner sep=1pt,line width=0.5pt,fill=blue!30!white] {$0$};
\node at (CD) [circle,draw,fill=white,inner sep=1pt,line width=0.5pt,fill=blue!30!white] {$0$};

\node at (midpoint) [circle,draw,fill=white,inner sep=1pt,line width=0.5pt,fill=blue!30!white] {$0$};
\end{tikzpicture}
		\caption{}
		\label{fig:tessellation3_2d}
	\end{subfigure}
	\caption{Schematic representation of the mid-point tessellation procedure for a two-dimensional case.}
	\label{fig:midpoint_2d}
\end{figure}
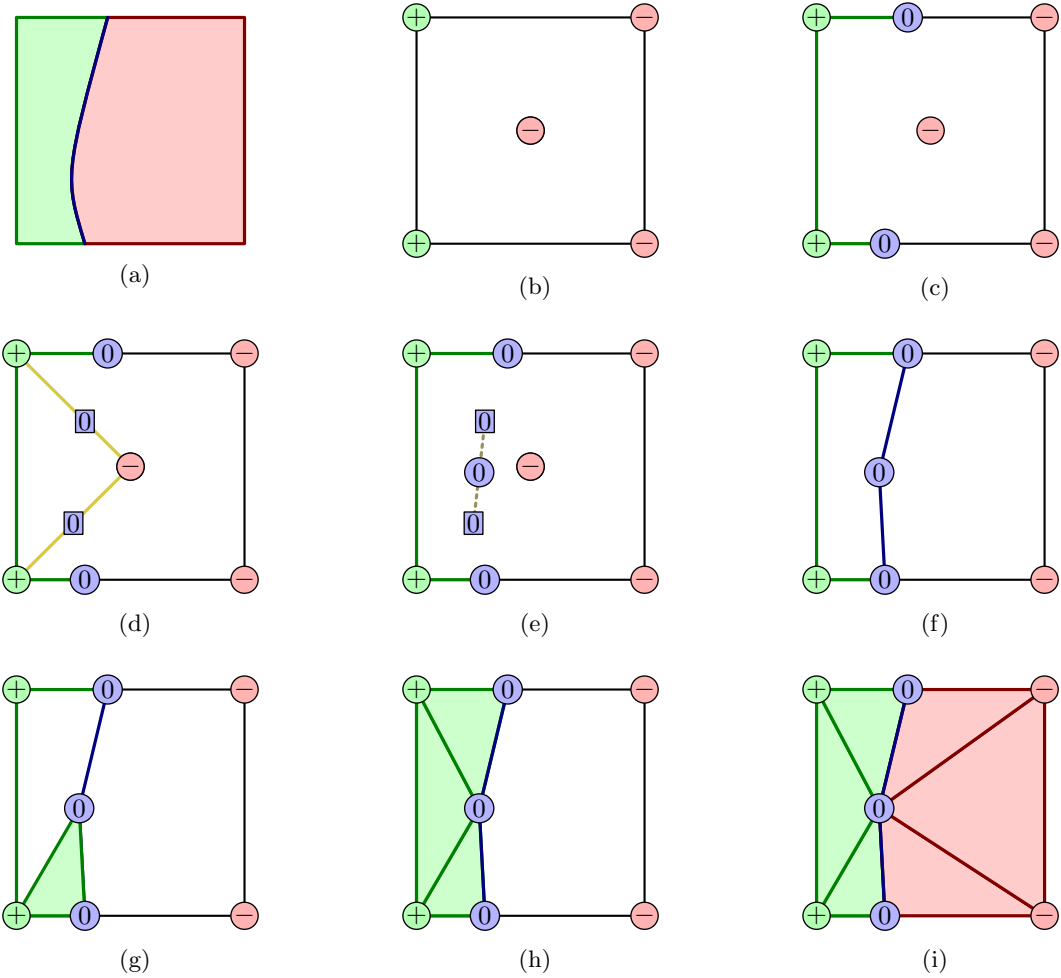

The employed tessellation procedure is illustrated for a two-dimensional sub-cell in Figure~\ref{fig:midpoint_2d}. In order to acquire a tessellation of the immersed boundary and the corresponding trimmed interiors, the following steps are taken:
\begin{itemize}
  \item[(a)] The immersed boundary, \emph{i.e.}, the blue curve in Figure~\ref{fig:cut_cell_2d}, is assumed to corresponds to the zero level set of a continuous function.
  \item[(b)] The level set function is evaluated in the vertices of the quadrilateral cell, and a bi-linear interpolation is used to approximate the level set value in the center point. In Figure~\ref{fig:levelset_vertex_2d} positive level set values are indicated by green circles with a plus sign, and negative ones by a red circle with a minus sign.
  \item[(c)] For all edges of the sub-cell we consider a linear interpolation of the level set function along the edge, and, if the level set values at the edge vertices have opposite signs, we determine the approximate zero level set point based on this linear interpolation. If an edge is intersected by the immersed boundary, it is split up in a positive part and a negative part based on the determined edge intersection. The approximate zero level set points are indicated by the blue circles with a $0$ in Figure~\ref{fig:trim_edges_2d}, and the splitting of the edges is visualized by the edge colors.
  \item[(d)] Along each of the 4 lines connecting the center of the square with its vertices, linearly interpolated zero level set points are determined. For the case considered in Figure~\ref{fig:levelset_diagonal_2d}, the two additionally obtained zero level set points are indicated by the blue squares on the yellow diagonals.
  \item[(e)] An additional zero level set point is formed by taking the arithmetic mean of the coordinates of the zero level set points along the diagonals. This point, which we refer to as the (approximate) midpoint of the trimmed boundary, is illustrated by the blue circle in the interior of the cell in Figure~\ref{fig:diagonal_midpoint_2d}.
  \item[(f)] The trimmed boundary is then constructed by extruding the edge intersections toward the computed approximate midpoint. As can be seen from Figure~\ref{fig:levelset_midpoint_2d}, a piece-wise linear approximation of the trimmed boundary in Figure~\ref{fig:cut_cell_2d} is obtained.
  \item[(g-i)] Interior cells are then constructed by extruding the edges toward the midpoint, thereby creating triangular integration sub-cells. The collection of sub-cells pertaining to the positive side of the trimmed elements are shown in green in Figure~\ref{fig:tessellation2_2d}, whereas the negative side of the element is shown in red in Figure~\ref{fig:tessellation3_2d}.
\end{itemize}

From this tessellation procedure applied in two dimensions it is evident that it traverses through the dimensions of the problem, in the sense that it first trims the edges (a--c), after which the quadrilateral is trimmed through the extrusion of the edges to the computed mid-point (d--i). This approach directly extends to the three-dimensional case, an illustration of which is presented in Figure~\ref{fig:midpoint_3d}. In Figure~\ref{fig:levelset_vertex_3d} the level set function is evaluated in all vertices of the cube. For each of the six faces of the cube the two-dimensional procedure as described above is applied, as illustrated on the unfolded cubes in Figures~\ref{fig:diagonal_3d}--\ref{fig:diagonal_midpoint_3d}. Finally, the trimmed faces are extruded toward the computed mid-point of the cube, as illustrated in Figures~\ref{fig:tessellation1_3d}--\ref{fig:tessellation_3d}. Note that in the three-dimensional case both tetrahedrons and pyramids are created, based on whether or not a face is trimmed. Integration rules are available for both these shapes.

\begin{figure}
	\centering
	\begin{subfigure}{0.33\textwidth}
		\centering
		\begin{tikzpicture}[
  scale=3,
  line cap=round,
  line join=round,
  opaque/.style={black,thin,opacity=0.4}
]

  \coordinate (A) at (0,0,0);
  \coordinate (B) at (1,0,0);
  \coordinate (C) at (1,1,0);
  \coordinate (D) at (0,1,0);
  \coordinate (E) at (0,0,1);
  \coordinate (F) at (1,0,1);
  \coordinate (G) at (1,1,1);
  \coordinate (H) at (0,1,1);
  
  \coordinate (AB) at ($0.6*(A)+0.4*(B)$);
  \coordinate (CD) at ($0.4*(C)+0.6*(D)$);
  \coordinate (GH) at ($0.3*(G)+0.7*(H)$);
  \coordinate (EF) at ($0.7*(E)+0.3*(F)$);

  \draw[dashed,thin,opacity=0.4] (D) -- (A);
  \draw[dashed,thin,opacity=0.4] (B) -- (A);
  \draw[dashed,thin,opacity=0.4] (A) -- (E);
  
  \draw[thick] (B) -- (C);
  \draw[thick] (C) -- (D);
  \draw[thick] (B) -- (F);
  \draw[thick] (C) -- (G);
  \draw[thick] (D) -- (H);
  \draw[thick] (E) -- (F);
  \draw[thick] (F) -- (G);
  \draw[thick] (G) -- (H);
  \draw[thick] (H) -- (E);
  
  \filldraw[draw=green!50!black,very thick,fill=green!20!white] (A) -- (E) -- (H) -- (D) -- cycle;
  \filldraw[draw=green!50!black,very thick,fill=green!20!white] (A) -- (AB).. controls (0.2,0,0.5) .. (EF) -- (E) -- (A);
  \filldraw[draw=green!50!black,very thick,fill=green!20!white] (E) -- (EF) -- (GH) -- (H) -- cycle;
  \filldraw[draw=green!50!black,very thick,fill=green!20!white] (A) -- (AB) .. controls (0.2,0.4,0) .. (CD) -- (D) -- (A);
  \filldraw[draw=blue!20!white,very thick,fill=blue!20!white] (AB) .. controls (0.2,0,0.5) .. (EF) -- (GH) .. controls (0.2,1,0.5) .. (CD) .. controls (0.2,0.4,0) .. (AB) -- cycle;
  \draw[draw=blue!50!black,very thick] (AB) .. controls (0.2,0,0.5) .. (EF);
  \draw[draw=blue!50!black,very thick] (CD) .. controls (0.2,1,0.5) .. (GH);
  \draw[draw=blue!50!black,very thick] (EF) -- (GH);
  \draw[draw=blue!50!black,very thick] (AB) .. controls (0.2,0.4,0) .. (CD);
  
  \draw[thick] (G) -- (GH);
  
\end{tikzpicture}
		\caption{}
		\label{fig:cut_cell_3d}
	\end{subfigure}%
	\begin{subfigure}{0.33\textwidth}
		\centering
		\begin{tikzpicture}[
  scale=3,
  line cap=round,
  line join=round,
  opaque/.style={black,thin,opacity=0.4}
]

  \coordinate (A) at (0,0,0);
  \coordinate (B) at (1,0,0);
  \coordinate (C) at (1,1,0);
  \coordinate (D) at (0,1,0);
  \coordinate (E) at (0,0,1);
  \coordinate (F) at (1,0,1);
  \coordinate (G) at (1,1,1);
  \coordinate (H) at (0,1,1);

  \coordinate (AB) at ($0.6*(A)+0.4*(B)$);
  \coordinate (CD) at ($0.4*(C)+0.6*(D)$);
  \coordinate (GH) at ($0.3*(G)+0.7*(H)$);
  \coordinate (EF) at ($0.7*(E)+0.3*(F)$);
  
  \coordinate (centriod) at (0.5,0.5,0.5);

  \draw[dashed,thin,opacity=0.4] (D) -- (A);
  \draw[dashed,thin,opacity=0.4] (B) -- (A);
  \draw[dashed,thin,opacity=0.4] (A) -- (E);

  \node at (A) [circle,draw,fill=white,inner sep=0pt,opacity=0.4,fill=green!30!white] {+};
  
  \draw[thick] (B) -- (C);
  \draw[thick] (C) -- (D);
  \draw[thick] (B) -- (F);
  \draw[thick] (C) -- (G);
  \draw[thick] (D) -- (H);
  \draw[thick] (E) -- (F);
  \draw[thick] (F) -- (G);
  \draw[thick] (G) -- (H);
  \draw[thick] (H) -- (E);
  
  \foreach \c in {(D),(E),(H)}
  \node at \c [circle,draw,fill=white,inner sep=0pt,line width=0.5pt,fill=green!30!white] {$+$};
  
  \foreach \c in {(F),(G),(C),(B),(centriod)}
  \node at \c [circle,draw,fill=white,inner sep=0pt,line width=0.5pt,fill=red!30!white] {$-$};
  
\end{tikzpicture}
		\caption{}
		\label{fig:levelset_vertex_3d}
	\end{subfigure}%
	\begin{subfigure}{0.33\textwidth}
		\centering
		\begin{tikzpicture}[
  scale=3,
  line cap=round,
  line join=round,
  opaque/.style={black,thin,opacity=0.4}
]

  \coordinate (A) at (0,0,0);
  \coordinate (B) at (1,0,0);
  \coordinate (C) at (1,1,0);
  \coordinate (D) at (0,1,0);
  \coordinate (E) at (0,0,1);
  \coordinate (F) at (1,0,1);
  \coordinate (G) at (1,1,1);
  \coordinate (H) at (0,1,1);

  \coordinate (AB) at ($0.6*(A)+0.4*(B)$);
  \coordinate (CD) at ($0.4*(C)+0.6*(D)$);
  \coordinate (GH) at ($0.3*(G)+0.7*(H)$);
  \coordinate (EF) at ($0.7*(E)+0.3*(F)$);
  
  \coordinate (centriod) at (0.5,0.5,0.5);
  
  \coordinate (m1) at (0.2,0,0.4);
  \coordinate (m2) at (0.3,0.5,1);
  \coordinate (m3) at (0.2,1,0.4);
  \coordinate (m4) at (0.3,0.5,0);

  \draw[dashed,thin,opacity=0.4] (D) -- (A);
  \draw[dashed,thin,opacity=0.4] (B) -- (A);
  \draw[dashed,thin,opacity=0.4] (A) -- (E);

  \node at (A) [circle,draw,fill=white,inner sep=0pt,opacity=0.4,fill=green!30!white] {+};
  
  \draw[thick] (B) -- (C);
  \draw[thick] (C) -- (D);
  \draw[thick] (B) -- (F);
  \draw[thick] (C) -- (G);
  \draw[thick] (D) -- (H);
  \draw[thick] (E) -- (F);
  \draw[thick] (F) -- (G);
  \draw[thick] (G) -- (H);
  \draw[thick] (H) -- (E);
  
  \foreach \c in {(D),(E),(H)}
  \node at \c [circle,draw,fill=white,inner sep=0pt,line width=0.5pt,fill=green!30!white] {$+$};
  
  \foreach \c in {(F),(G),(C),(B),(centriod)}
  \node at \c [circle,draw,fill=white,inner sep=0pt,line width=0.5pt,fill=red!30!white] {$-$};
  
  \node at (AB) [circle,draw,fill=white,inner sep=1pt,line width=0.5pt,fill=blue!30!white] {$0$};
  \node at (CD) [circle,draw,fill=white,inner sep=1pt,line width=0.5pt,fill=blue!30!white] {$0$};
  \node at (EF) [circle,draw,fill=white,inner sep=1pt,line width=0.5pt,fill=blue!30!white] {$0$};
  \node at (GH) [circle,draw,fill=white,inner sep=1pt,line width=0.5pt,fill=blue!30!white] {$0$};
  
\end{tikzpicture}
		\caption{}
		\label{fig:trim_edges_3d}
	\end{subfigure}\\[12pt]
	\begin{subfigure}{0.45\textwidth}
		\centering
		\begin{tikzpicture}[
scale=3 ,
line cap=round,
line join=round,
opaque/.style={black,thin,opacity=0.4}
]

\coordinate (A) at (0,0);
\coordinate (B) at (0.6,0);

\coordinate (C) at (0,0.6);
\coordinate (D) at (0.6,0.6);

\coordinate (E) at (-0.6,0.6);
\coordinate (F) at (0,1.2);
\coordinate (G) at (0.6,1.2);
\coordinate (H) at (1.2,0.6);

\coordinate (I) at (-0.6,1.2);
\coordinate (J) at (0,1.8);
\coordinate (K) at (0.6,1.8);
\coordinate (L) at (1.2,1.2);

\coordinate (M) at (0,2.4);
\coordinate (N) at (0.6,2.4);

\coordinate (c1) at (0.3,0.3);
\coordinate (c2) at (0.3,0.9);
\coordinate (c3) at (0.3,1.5);
\coordinate (c4) at (0.3,2.1);

\coordinate (AB) at ($0.6*(A)+0.4*(B)$);
\coordinate (CD) at ($0.7*(C)+0.3*(D)$);
\coordinate (GF) at ($0.3*(G)+0.7*(F)$);
\coordinate (JK) at ($0.4*(K)+0.6*(J)$);
\coordinate (MN) at ($0.6*(M)+0.4*(N)$);

\coordinate (mAc1) at ($0.5*(A)+ 0.5*(c1)$);
\coordinate (mCc1) at ($0.4*(C)+ 0.6*(c1)$);
\coordinate (mCc2) at ($0.4*(C)+ 0.6*(c2)$);
\coordinate (mFc2) at ($0.4*(F)+ 0.6*(c2)$);
\coordinate (mFc3) at ($0.4*(F)+ 0.6*(c3)$);
\coordinate (mJc3) at ($0.4*(J)+ 0.6*(c3)$);
\coordinate (mJc4) at ($0.4*(J)+ 0.6*(c4)$);
\coordinate (mMc4) at ($0.4*(M)+ 0.6*(c4)$);

\draw[thick] (A)--(B)--(D)--(C)--cycle;
\draw[thick] (C)--(D)--(G)--(F)--cycle;
\draw[thick] (E)--(C)--(F)--(I)--cycle;
\draw[thick] (F)--(G)--(K)--(J)--cycle;
\draw[thick] (D)--(H)--(L)--(G)--cycle;
\draw[thick] (M)--(N)--(K)--(J)--cycle;

\node at (AB) [circle,draw,fill=white,inner sep=1pt,line width=0.5pt,fill=blue!30!white] {$0$};
\node at (CD) [circle,draw,fill=white,inner sep=1pt,line width=0.5pt,fill=blue!30!white] {$0$};
\node at (GF) [circle,draw,fill=white,inner sep=1pt,line width=0.5pt,fill=blue!30!white] {$0$};
\node at (JK) [circle,draw,fill=white,inner sep=1pt,line width=0.5pt,fill=blue!30!white] {$0$};
\node at (MN) [circle,draw,fill=white,inner sep=1pt,line width=0.5pt,fill=blue!30!white] {$0$};

\foreach \c in {(A),(C),(E),(F),(I),(J),(M)}
\node at \c [circle,draw,fill=white,inner sep=0pt,line width=0.5pt,fill=green!30!white] {$+$};

\foreach \c in {(B),(D),(G),(H),(K),(L),(N)}
\node at \c [circle,draw,fill=white,inner sep=0pt,line width=0.5pt,fill=red!30!white] {$-$};

\foreach \c in {(c1),(c2),(c3),(c4)}
\node at \c [circle,draw,fill=white,inner sep=0pt,line width=0.5pt,fill=red!30!white] {$-$};

\end{tikzpicture}
		\caption{}
		\label{fig:diagonal_3d}
	\end{subfigure}%
		\begin{subfigure}{0.45\textwidth}
		\centering
		\begin{tikzpicture}[
scale=3 ,
line cap=round,
line join=round,
opaque/.style={black,thin,opacity=0.4}
]

\coordinate (A) at (0,0);
\coordinate (B) at (0.6,0);

\coordinate (C) at (0,0.6);
\coordinate (D) at (0.6,0.6);

\coordinate (E) at (-0.6,0.6);
\coordinate (F) at (0,1.2);
\coordinate (G) at (0.6,1.2);
\coordinate (H) at (1.2,0.6);

\coordinate (I) at (-0.6,1.2);
\coordinate (J) at (0,1.8);
\coordinate (K) at (0.6,1.8);
\coordinate (L) at (1.2,1.2);

\coordinate (M) at (0,2.4);
\coordinate (N) at (0.6,2.4);

\coordinate (c1) at (0.3,0.3);
\coordinate (c2) at (0.3,0.9);
\coordinate (c3) at (0.3,1.5);
\coordinate (c4) at (0.3,2.1);

\coordinate (AB) at ($0.6*(A)+0.4*(B)$);
\coordinate (CD) at ($0.7*(C)+0.3*(D)$);
\coordinate (GF) at ($0.3*(G)+0.7*(F)$);
\coordinate (JK) at ($0.4*(K)+0.6*(J)$);
\coordinate (MN) at ($0.6*(M)+0.4*(N)$);

\coordinate (mAc1) at ($0.5*(A)+ 0.5*(c1)$);
\coordinate (mCc1) at ($0.4*(C)+ 0.6*(c1)$);
\coordinate (mCc2) at ($0.4*(C)+ 0.6*(c2)$);
\coordinate (mFc2) at ($0.4*(F)+ 0.6*(c2)$);
\coordinate (mFc3) at ($0.4*(F)+ 0.6*(c3)$);
\coordinate (mJc3) at ($0.4*(J)+ 0.6*(c3)$);
\coordinate (mJc4) at ($0.4*(J)+ 0.6*(c4)$);
\coordinate (mMc4) at ($0.4*(M)+ 0.6*(c4)$);

\coordinate (m1) at ($0.5*(mAc1)+0.5*(mCc1)$);
\coordinate (m2) at ($0.5*(mCc2)+0.5*(mFc2)$);
\coordinate (m3) at ($0.5*(mFc3)+0.5*(mJc3)$);
\coordinate (m4) at ($0.5*(mJc4)+0.5*(mMc4)$);

\draw[thick] (A)--(B)--(D)--(C)--cycle;
\draw[thick] (C)--(D)--(G)--(F)--cycle;
\draw[thick] (E)--(C)--(F)--(I)--cycle;
\draw[thick] (F)--(G)--(K)--(J)--cycle;
\draw[thick] (D)--(H)--(L)--(G)--cycle;
\draw[thick] (M)--(N)--(K)--(J)--cycle;

\filldraw[draw=green!50!black,very thick,fill=green!20!white] (A) -- (AB) -- (m1) -- cycle;
\filldraw[draw=green!50!black,very thick,fill=green!20!white] (C) -- (CD) -- (m1) -- cycle;
\filldraw[draw=green!50!black,very thick,fill=green!20!white] (C) -- (A) -- (m1) -- cycle;
\filldraw[draw=green!50!black,very thick,fill=green!20!white] (C) -- (CD) -- (m2) -- cycle;
\filldraw[draw=green!50!black,very thick,fill=green!20!white] (F) -- (GF) -- (m2) -- cycle;
\filldraw[draw=green!50!black,very thick,fill=green!20!white] (F) -- (C) -- (m2) -- cycle;
\filldraw[draw=green!50!black,very thick,fill=green!20!white] (F) -- (GF) -- (m3) -- cycle;
\filldraw[draw=green!50!black,very thick,fill=green!20!white] (J) -- (JK) -- (m3) -- cycle;
\filldraw[draw=green!50!black,very thick,fill=green!20!white] (J) -- (F) -- (m3) -- cycle;
\filldraw[draw=green!50!black,very thick,fill=green!20!white] (J) -- (JK) -- (m4) -- cycle;
\filldraw[draw=green!50!black,very thick,fill=green!20!white] (M) -- (MN) -- (m4) -- cycle;
\filldraw[draw=green!50!black,very thick,fill=green!20!white] (M) -- (J) -- (m4) -- cycle;
\filldraw[draw=green!50!black,very thick,fill=green!20!white] (E) -- (C) -- (F) -- (I) -- cycle;
\draw[draw=blue!50!black,very thick,fill=blue!20!white] (AB) -- (m1);
\draw[draw=blue!50!black,very thick,fill=blue!20!white] (CD) -- (m1);
\draw[draw=blue!50!black,very thick,fill=blue!20!white] (CD) -- (m2);
\draw[draw=blue!50!black,very thick,fill=blue!20!white] (GF) -- (m2);
\draw[draw=blue!50!black,very thick,fill=blue!20!white] (GF) -- (m3);
\draw[draw=blue!50!black,very thick,fill=blue!20!white] (JK) -- (m3);
\draw[draw=blue!50!black,very thick,fill=blue!20!white] (JK) -- (m4);
\draw[draw=blue!50!black,very thick,fill=blue!20!white] (MN) -- (m4);

\node at (AB) [circle,draw,fill=white,inner sep=1pt,line width=0.5pt,fill=blue!30!white] {$0$};
\node at (CD) [circle,draw,fill=white,inner sep=1pt,line width=0.5pt,fill=blue!30!white] {$0$};
\node at (GF) [circle,draw,fill=white,inner sep=1pt,line width=0.5pt,fill=blue!30!white] {$0$};
\node at (JK) [circle,draw,fill=white,inner sep=1pt,line width=0.5pt,fill=blue!30!white] {$0$};
\node at (MN) [circle,draw,fill=white,inner sep=1pt,line width=0.5pt,fill=blue!30!white] {$0$};

\foreach \c in {(A),(C),(E),(F),(I),(J),(M)}
\node at \c [circle,draw,fill=white,inner sep=0pt,line width=0.5pt,fill=green!30!white] {$+$};

\foreach \c in {(B),(D),(G),(H),(K),(L),(N)}
\node at \c [circle,draw,fill=white,inner sep=0pt,line width=0.5pt,fill=red!30!white] {$-$};

\foreach \c in {(m1),(m2),(m3),(m4)}
\node at \c [circle,draw,fill=white,inner sep=1pt,line width=0.5pt,fill=blue!30!white] {$0$};

\end{tikzpicture}
		\caption{}
		\label{fig:diagonal_midpoint_3d}
	\end{subfigure}\\[12pt]
	\begin{subfigure}{0.33\textwidth}
		\centering
		\begin{tikzpicture}[
  scale=3,
  line cap=round,
  line join=round,
  opaque/.style={black,thin,opacity=0.4}
]

  \coordinate (A) at (0,0,0);
  \coordinate (B) at (1,0,0);
  \coordinate (C) at (1,1,0);
  \coordinate (D) at (0,1,0);
  \coordinate (E) at (0,0,1);
  \coordinate (F) at (1,0,1);
  \coordinate (G) at (1,1,1);
  \coordinate (H) at (0,1,1);

  \coordinate (AB) at ($0.6*(A)+0.4*(B)$);
  \coordinate (CD) at ($0.4*(C)+0.6*(D)$);
  \coordinate (GH) at ($0.3*(G)+0.7*(H)$);
  \coordinate (EF) at ($0.7*(E)+0.3*(F)$);

  \coordinate (m1) at (0.2,0,0.4);
  \coordinate (m2) at (0.3,0.5,1);
  \coordinate (m3) at (0.2,1,0.4);
  \coordinate (m4) at (0.3,0.5,0);

  \coordinate (midpoint) at ($0.25*(m1)+0.25*(m2)+0.25*(m3)+0.25*(m4)$);

  \draw[dashed,thin,opacity=0.4] (D) -- (A);
  \draw[dashed,thin,opacity=0.4] (B) -- (A);
  \draw[dashed,thin,opacity=0.4] (A) -- (E);

  \node at (A) [circle,draw,fill=white,inner sep=0pt,opacity=0.4,fill=green!30!white] {+};
  
  \draw[thick] (B) -- (C);
  \draw[thick] (C) -- (D);
  \draw[thick] (B) -- (F);
  \draw[thick] (C) -- (G);
  \draw[thick] (D) -- (H);
  \draw[thick] (E) -- (F);
  \draw[thick] (F) -- (G);
  \draw[thick] (G) -- (H);
  \draw[thick] (H) -- (E);
  
  \filldraw[draw=blue!50!black,very thick,fill=blue!20!white] (AB) -- (m1) -- (EF) -- (m2) -- (GH)-- (m3) -- (CD) -- (m4) -- (AB) -- cycle;
  \draw[draw=blue!50!black,very thick,fill=blue!20!white] (AB) -- (m1);
  \draw[draw=blue!50!black,very thick,fill=blue!20!white] (EF) -- (m1);
  \draw[draw=blue!50!black,very thick,fill=blue!20!white] (EF) -- (m2);
  \draw[draw=blue!50!black,very thick,fill=blue!20!white] (GH) -- (m2);
  \draw[draw=blue!50!black,very thick,fill=blue!20!white] (GH) -- (m3);
  \draw[draw=blue!50!black,very thick,fill=blue!20!white] (CD) -- (m3);  
  \draw[draw=blue!50!black,very thick,fill=blue!20!white] (CD) -- (m4);
  \draw[draw=blue!50!black,very thick,fill=blue!20!white] (AB) -- (m4);
  
  \draw[draw=blue!50!black,very thick,fill=blue!20!white] (AB) -- (midpoint);
  \draw[draw=blue!50!black,very thick,fill=blue!20!white] (EF) -- (midpoint);
  \draw[draw=blue!50!black,very thick,fill=blue!20!white] (GH) -- (midpoint);
  \draw[draw=blue!50!black,very thick,fill=blue!20!white] (CD) -- (midpoint);
  \draw[draw=blue!50!black,very thick,fill=blue!20!white] (m1) -- (midpoint);  
  \draw[draw=blue!50!black,very thick,fill=blue!20!white] (m2) -- (midpoint);    
  \draw[draw=blue!50!black,very thick,fill=blue!20!white] (m3) -- (midpoint);
  \draw[draw=blue!50!black,very thick,fill=blue!20!white] (m4) -- (midpoint);
  
  \draw[thick] (G) -- (H);
  
  \foreach \c in {(D),(E),(H)}
  \node at \c [circle,draw,fill=white,inner sep=0pt,line width=0.5pt,fill=green!30!white] {$+$};
  
  \foreach \c in {(F),(G),(C),(B)}
  \node at \c [circle,draw,fill=white,inner sep=0pt,line width=0.5pt,fill=red!30!white] {$-$};
  
  \node at (AB) [circle,draw,fill=white,inner sep=1pt,line width=0.5pt,fill=blue!30!white] {$0$};
  \node at (CD) [circle,draw,fill=white,inner sep=1pt,line width=0.5pt,fill=blue!30!white] {$0$};
  \node at (EF) [circle,draw,fill=white,inner sep=1pt,line width=0.5pt,fill=blue!30!white] {$0$};
  \node at (GH) [circle,draw,fill=white,inner sep=1pt,line width=0.5pt,fill=blue!30!white] {$0$};
  
  \foreach \c in {(m1),(m2),(m3),(m4)}
  \node at \c [circle,draw,fill=white,inner sep=1pt,line width=0.5pt,fill=blue!30!white] {$0$};
  
  \node at (midpoint) [circle,draw,fill=white,inner sep=1pt,line width=0.5pt,fill=blue!60!white] {$0$}; 
\end{tikzpicture}
		\caption{}
		\label{fig:levelset_midpoint_3d}
	\end{subfigure}%
	\begin{subfigure}{0.33\textwidth}
		\centering
		\begin{tikzpicture}[
  scale=3,
  line cap=round,
  line join=round,
  opaque/.style={black,thin,opacity=0.4}
]

  \coordinate (A) at (0,0,0);
  \coordinate (B) at (1,0,0);
  \coordinate (C) at (1,1,0);
  \coordinate (D) at (0,1,0);
  \coordinate (E) at (0,0,1);
  \coordinate (F) at (1,0,1);
  \coordinate (G) at (1,1,1);
  \coordinate (H) at (0,1,1);

  \coordinate (AB) at ($0.6*(A)+0.4*(B)$);
  \coordinate (CD) at ($0.4*(C)+0.6*(D)$);
  \coordinate (GH) at ($0.3*(G)+0.7*(H)$);
  \coordinate (EF) at ($0.7*(E)+0.3*(F)$);

  \coordinate (m1) at (0.2,0,0.4);
  \coordinate (m2) at (0.3,0.5,1);
  \coordinate (m3) at (0.2,1,0.4);
  \coordinate (m4) at (0.3,0.5,0);

  \coordinate (midpoint) at ($0.25*(m1)+0.25*(m2)+0.25*(m3)+0.25*(m4)$);

  \draw[dashed,thin,opacity=0.4] (D) -- (A);
  \draw[dashed,thin,opacity=0.4] (B) -- (A);
  \draw[dashed,thin,opacity=0.4] (A) -- (E);

  \node at (A) [circle,draw,fill=white,inner sep=0pt,opacity=0.4,fill=green!30!white] {+};
  
  \draw[thick] (B) -- (C);
  \draw[thick] (C) -- (D);
  \draw[thick] (B) -- (F);
  \draw[thick] (C) -- (G);
  \draw[thick] (D) -- (H);
  \draw[thick] (E) -- (F);
  \draw[thick] (F) -- (G);
  \draw[thick] (G) -- (H);
  \draw[thick] (H) -- (E);
  
  \filldraw[draw=green!50!black,very thick,fill=green!20!white] (GH) -- (H) -- (m2) -- cycle;
  \filldraw[draw=green!50!black,very thick,fill=green!20!white] (GH) -- (midpoint) -- (m2) -- cycle;
  \draw[draw=green!50!black, dashed,thin] (H) -- (midpoint);

  \draw[draw=blue!50!black,very thick,fill=blue!20!white] (AB) -- (m1);
  \draw[draw=blue!50!black,very thick,fill=blue!20!white] (EF) -- (m1);
  \draw[draw=blue!50!black,very thick,fill=blue!20!white] (EF) -- (m2);
  \draw[draw=blue!50!black,very thick,fill=blue!20!white] (GH) -- (m2);
  \draw[draw=blue!50!black,very thick,fill=blue!20!white] (GH) -- (m3);
  \draw[draw=blue!50!black,very thick,fill=blue!20!white] (CD) -- (m3);  
  \draw[draw=blue!50!black,very thick,fill=blue!20!white] (CD) -- (m4);
  \draw[draw=blue!50!black,very thick,fill=blue!20!white] (AB) -- (m4);
  
  \foreach \c in {(D),(E),(H)}
  \node at \c [circle,draw,fill=white,inner sep=0pt,line width=0.5pt,fill=green!30!white] {$+$};
  
  \foreach \c in {(F),(G),(C),(B)}
  \node at \c [circle,draw,fill=white,inner sep=0pt,line width=0.5pt,fill=red!30!white] {$-$};
  
  \node at (AB) [circle,draw,fill=white,inner sep=1pt,line width=0.5pt,fill=blue!30!white] {$0$};
  \node at (CD) [circle,draw,fill=white,inner sep=1pt,line width=0.5pt,fill=blue!30!white] {$0$};
  \node at (EF) [circle,draw,fill=white,inner sep=1pt,line width=0.5pt,fill=blue!30!white] {$0$};
  \node at (GH) [circle,draw,fill=white,inner sep=1pt,line width=0.5pt,fill=blue!30!white] {$0$};
  
  \foreach \c in {(m1),(m2),(m3),(m4)}
  \node at \c [circle,draw,fill=white,inner sep=1pt,line width=0.5pt,fill=blue!30!white] {$0$};
  
  \node at (midpoint) [circle,draw,fill=white,inner sep=1pt,line width=0.5pt,fill=blue!60!white] {$0$};
  
\end{tikzpicture}
		\caption{}
		\label{fig:tessellation1_3d}
	\end{subfigure}%
	\begin{subfigure}{0.33\textwidth}
		\centering
		\begin{tikzpicture}[
  scale=3,
  line cap=round,
  line join=round,
  opaque/.style={black,thin,opacity=0.4}
]

  \coordinate (A) at (0,0,0);
  \coordinate (B) at (1,0,0);
  \coordinate (C) at (1,1,0);
  \coordinate (D) at (0,1,0);
  \coordinate (E) at (0,0,1);
  \coordinate (F) at (1,0,1);
  \coordinate (G) at (1,1,1);
  \coordinate (H) at (0,1,1);

  \coordinate (AB) at ($0.6*(A)+0.4*(B)$);
  \coordinate (CD) at ($0.4*(C)+0.6*(D)$);
  \coordinate (GH) at ($0.3*(G)+0.7*(H)$);
  \coordinate (EF) at ($0.7*(E)+0.3*(F)$);

  \coordinate (m1) at (0.2,0,0.4);
  \coordinate (m2) at (0.3,0.5,1);
  \coordinate (m3) at (0.2,1,0.4);
  \coordinate (m4) at (0.3,0.5,0);
  
  \coordinate (midpoint) at ($0.25*(m1)+0.25*(m2)+0.25*(m3)+0.25*(m4)$);

  \draw[dashed,thin,opacity=0.4] (B) -- (A);

  \draw[thick] (B) -- (C);
  \draw[thick] (C) -- (D);
  \draw[thick] (B) -- (F);
  \draw[thick] (C) -- (G);
  \draw[thick] (D) -- (H);
  \draw[thick] (E) -- (F);
  \draw[thick] (F) -- (G);
  \draw[thick] (G) -- (H);
  \draw[thick] (H) -- (E);
  
  \filldraw[draw=green!50!black,very thick,fill=green!20!white] (H) -- (E) -- (A) -- (D) -- cycle;  
  \filldraw[draw=green!50!black,very thick,fill=green!20!white] (AB) -- (A) -- (m1) -- cycle;
  \filldraw[draw=green!50!black,very thick,fill=green!20!white] (A) -- (E) -- (m1) -- cycle;  
  \filldraw[draw=green!50!black,very thick,fill=green!20!white] (EF) -- (E) -- (m1) -- cycle;
  \filldraw[draw=green!50!black,very thick,fill=green!20!white] (EF) -- (E) -- (m2) -- cycle;
  \filldraw[draw=green!50!black,very thick,fill=green!20!white] (E) -- (H) -- (m2) -- cycle;  
  \filldraw[draw=green!50!black,very thick,fill=green!20!white] (GH) -- (H) -- (m2) -- cycle;
  \filldraw[draw=green!50!black,very thick,fill=green!20!white] (GH) -- (H) -- (m3) -- cycle;
  \filldraw[draw=green!50!black,very thick,fill=green!20!white] (D) -- (H) -- (m3) -- cycle;  
  \filldraw[draw=green!50!black,very thick,fill=green!20!white] (CD) -- (D) -- (m3) -- cycle;
  \filldraw[draw=green!50!black,very thick,fill=green!20!white] (CD) -- (D) -- (m4) -- cycle;
  \filldraw[draw=green!50!black,very thick,fill=green!20!white] (A) -- (D) -- (m4) -- cycle;  
  \filldraw[draw=green!50!black,very thick,fill=green!20!white] (AB) -- (A) -- (m4) -- cycle;
  
  \filldraw[draw=blue!20!white,very thick,fill=blue!20!white] (AB) -- (m1) -- (EF) -- (m2) -- (GH)-- (m3) -- (CD) -- (m4) -- (AB) -- cycle;
  
  \draw[draw=blue!50!black,very thick,fill=blue!20!white] (AB) -- (m1);
  \draw[draw=blue!50!black,very thick,fill=blue!20!white] (EF) -- (m1);
  \draw[draw=blue!50!black,very thick,fill=blue!20!white] (EF) -- (m2);
  \draw[draw=blue!50!black,very thick,fill=blue!20!white] (GH) -- (m2);
  \draw[draw=blue!50!black,very thick,fill=blue!20!white] (GH) -- (m3);
  \draw[draw=blue!50!black,very thick,fill=blue!20!white] (CD) -- (m3);  
  \draw[draw=blue!50!black,very thick,fill=blue!20!white] (CD) -- (m4);
  \draw[draw=blue!50!black,very thick,fill=blue!20!white] (AB) -- (m4);

  \draw[draw=blue!50!black,very thick,fill=blue!20!white] (AB) -- (midpoint);
  \draw[draw=blue!50!black,very thick,fill=blue!20!white] (EF) -- (midpoint);
  \draw[draw=blue!50!black,very thick,fill=blue!20!white] (GH) -- (midpoint);
  \draw[draw=blue!50!black,very thick,fill=blue!20!white] (CD) -- (midpoint);
  \draw[draw=blue!50!black,very thick,fill=blue!20!white] (m1) -- (midpoint);  
  \draw[draw=blue!50!black,very thick,fill=blue!20!white] (m2) -- (midpoint);    
  \draw[draw=blue!50!black,very thick,fill=blue!20!white] (m3) -- (midpoint);
  \draw[draw=blue!50!black,very thick,fill=blue!20!white] (m4) -- (midpoint);
  
  \draw[thick] (G) -- (GH);
  
   \foreach \c in {(D),(E),(H)}
  \node at \c [circle,draw,fill=white,inner sep=0pt,line width=0.5pt,fill=green!30!white] {$+$};
  
  \foreach \c in {(F),(G),(C),(B)}
  \node at \c [circle,draw,fill=white,inner sep=0pt,line width=0.5pt,fill=red!30!white] {$-$};
  
  \node at (AB) [circle,draw,fill=white,inner sep=1pt,line width=0.5pt,fill=blue!30!white] {$0$};
  \node at (CD) [circle,draw,fill=white,inner sep=1pt,line width=0.5pt,fill=blue!30!white] {$0$};
  \node at (EF) [circle,draw,fill=white,inner sep=1pt,line width=0.5pt,fill=blue!30!white] {$0$};
  \node at (GH) [circle,draw,fill=white,inner sep=1pt,line width=0.5pt,fill=blue!30!white] {$0$};
  
  \foreach \c in {(m1),(m2),(m3),(m4)}
  \node at \c [circle,draw,fill=white,inner sep=1pt,line width=0.5pt,fill=blue!30!white] {$0$};
  
  \node at (midpoint) [circle,draw,fill=white,inner sep=1pt,line width=0.5pt,fill=blue!60!white] {$0$};  
  
\end{tikzpicture}
		\caption{}
		\label{fig:tessellation_3d}
	\end{subfigure}
	\caption{Schematic representation of the mid-point tessellation procedure for a three-dimensional case.}
	\label{fig:midpoint_3d}
\end{figure}

\end{document}